%% file: MAIN.tex
\documentclass [a4paper,12pt]{book}
\usepackage{geometry}
\usepackage{graphicx}
\usepackage{xcolor}
\usepackage{amssymb}
\usepackage{amsmath}
\usepackage[fleqn]{mathtools}

\usepackage[utf8]{inputenc}

\usepackage{csquotes}

\usepackage{amsthm}
\usepackage{bm}
\usepackage{xcolor}
\usepackage{verbatim}

\newcommand{\bea}{\begin{eqnarray}}
\usepackage{amsthm}
\usepackage{fancyhdr}
\newtheorem{prop}{Proposition}[section]
\newtheorem{thm}{Theorem}[chapter]

\newtheorem*{prf*}{Proof}
\newtheorem{exl}{Example}[chapter]
\newtheorem{defin}{Definition}[chapter]
\newtheorem{rmrk}{Remark}[chapter]
\newtheorem{cor}{Corolary}

\usepackage{MnSymbol}
\usepackage{algorithmic}
\usepackage{algorithm}
\usepackage{cite}
\usepackage{booktabs}
\usepackage{array}

\makeatletter
\newcommand{\doubletilde}[1]{{%
  \mathpalette\double@tilde{#1}%
}}
\newcommand{\double@tilde}[2]{%
  \sbox\z@{$\m@th#1\tilde{#2}$}%
  \ht\z@=.9\ht\z@
  \tilde{\box\z@}%
}
\makeatother

\usepackage{mathdots}
\newcommand{\ciclo}{
\begin{bmatrix}
    c_{0}   & c_{n-1}  & c_{n-2}  & \cdots   & \cdots  & c_{1}   \\
    c_{1}   & \ddots   & \ddots   & \ddots   &         & \vdots  \\
    c_{2}   & \ddots   & \ddots   & \ddots   & \ddots  & \vdots  \\
    \vdots  & \ddots   & \ddots   & \ddots   & \ddots  & c_{n-2} \\ 
    \vdots  &          & \ddots   & \ddots   & \ddots  & c_{n-1} \\
    c_{n-1} & \cdots   & \cdots   & c_{2}    & c_{1}   & c_{0}
\end{bmatrix}
}
\newcommand{\eigDistr}{
  \lim_{n\to \infty} \frac{1}{n} \sum_{j=1}^{n}F(\lambda_j(A_n)) = \frac{1}{\mu(D)}\int_DF(f(x))dx
}
\newcommand{\singDistr}{
  \lim_{n\to \infty} \frac{1}{n} \sum_{j=1}^{n}F(\sigma_j(A_n)) = \frac{1}{\mu(D)}\int_DF(|f(x)|)dx
}
\newcommand{\seqOfSeqs}{ \{\{B_{n,m}\}_n\}_m }
\newcommand{\acsArrow}{\xrightarrow{\text{  a.c.s}}}

\begin{document}
\include{TITLE_PAGE}
\hphantom{ss}
\newpage
\include{ACKNOWLEDGS}
\tableofcontents

\include{INTRODUCTION}
\part{Spectral Analysis and Preconditioned Iterative Solvers for Large Structured Linear Systems}
\include{GLT}

\include{SYMMETRIZED}

\include{FRACTIONAL}

\part{An Iterative Technique for Pricing American Put Options}
\include{AMERICAN}

\newpage
\include{CONCLUSION}

\bibliographystyle{plain}
\bibliography{REFERENCES.bib}
\end{document}

%% file: TITLE_PAGE.tex
\begin{titlepage}
    \begin{center}
        \vspace*{1cm}
            
        \Huge
        \textbf{Spectral Analysis and Preconditioned Iterative Solvers for Large Structured Linear Systems}
            
        \vspace{2.5cm}
        \LARGE
            
        \vspace{4.5cm}
            
        \textbf{Nikos Barakitis}

        \vspace{5.8cm}

        \Large
        ATHENS UNIVERSITY OF ECONOMICS AND BUSINESS \\
        SCHOOL OF INFORMATION SCIENCES AND TECHNOLOGY \\
        DEPARTMENT OF INFORMATICS

    \end{center}
\end{titlepage}

%% file: ACKNOWLEDGS.tex
\section*{Acknowledgment}
\addcontentsline{toc}{section}{Acknowledgment}

I would firstly like to thank my supervisor, Associate Professor Paris Vassalos and my advisor, Professor Stefano Serra-Capizzano, for the unlimited help and support they offered to me during my studies. I would also like to thank Emeritus Professor Evaggelos Mageirou for his throughout invaluable help and guidance, and my advisor Professor Dimitrios Noutsos for being member of the advising commitee and generously sharing ideas with the team.   
\par
I would also like to thank my collaborator researchers, Sven-Erik Erkstrom and Paola Ferrari for the fruitful collaboration we had, and their patience.
\par
Finally, I would like to thank Associate Professor Stavros Toumpis and Professor Panagiotis Katerinis  for their attitude toward me at several moments which helped me to build my self confidence and made me to believe that the completion of this dissertation could be feasible. 
\newline
\newline
Thank you,\\
Nikos Barakitis, Athens 2021.

%% file: INTRODUCTION.tex
\chapter*{Introduction}
\markboth{INTRODUCTION}{INTRODUCTION}

\addcontentsline{toc}{chapter}{Introduction}

The use of iterative methods for solving large structured linear systems has been of interest for more than half a century, and the development of the field has gone together with the improvement of computer systems. In fact, the solution of a large linear system of the form

\begin{align}
Ax=b, \label{linear system}
\end{align}
where the size of $A$ is $n\times n$ and $b$ is a column vector of size $n$, is often the central and most time- and storage-consuming part of the computation.{}

Iterative methods for solving linear systems first appeared in the works of  Gauss, Seidel, and  Jacobi in the 19-th century, with further progress in these methods being made in the first half of the 20-th century. These methods are typically referred to as \emph{stationary}, as opposed to the other  classes of iterative techniques that appeared later and relied on solution searches in Krylov subspaces. For the latter methods, the story began in 1952 with the development of the conjugate gradient (CG) method \cite{Lanczos1952}. This method was proposed for solving symmetric and positive-definite linear systems. Initially, it was considered a direct method, because it was proved analytically to reach the exact solution in at most $n$ steps, or actually in as many steps as the number of distinct eigenvalues of the coefficient matrix. However, in practice, owing to  the limited accuracy of floating-point arithmetic, especially in presence of ill-conditioning  the method requires more iterations than expected for a satisfactory approximation of the solution. Besides this, when considered as a direct method, it required more arithmetic operations than the Gauss elimination.
This method has remained out of interest for two decades. As applications required larger linear systems to be solved, the poor computational escalation of direct methods was overshadowed by the rapid improvement of computers and the evolution of computational methods.{}

This attitude regarding CG changed in the 70s following a publication from J. Reid \cite{reid1971method}. Thereafter, it became clear that for well-conditioned systems, the number of steps that the CG requires to reach the solution with a given accuracy is independent of the size of the system. This work brought Krylov methods back into the focus of the research community. The list of Krylov methods, limited until then, was enriched with methods for non-definite symmetric systems [e.g., the minimum residual method (MINRES) \cite{Paige1975}] and methods for non-symmetric systems [e.g., the generalised minimum residual method (GMRES) \cite{Saad1986}].{}

Published by O. Axelson and G. Lindskog in 1986 \cite{axelsson861}. Since then, the paper has been included in the references of almost every work relating to the solution of linear systems using iterative Krylov-subspace-theory based methods. In this study, it was proved that the efficiency of the preconditioned CG method depends on the clustering of the eigenvalues of the preconditioned matrix, which are \emph{clustered} at $(1)$ and yet certainly far from zero. The notion of eigenvalue clustering will be defined later: Hereafter, preconditioning was officially upgraded to the first research target in the field numerical solution of linear systems.\footnote{It must be mentioned here that for methods such as the Generalized Minimum Residual method, applied to non-symmetric systems, the eigenvalue distribution may not exactly describe the convergence\cite{Greenbaum1996}. However, in every case, a clustered spectrum and a minimal eigenvalue far from zero ensure fast convergence of the method.}. {}

\emph{Preconditioning of a linear system} refers to the replacement of the system (\ref{linear system}) with
\begin{align*}
&M^{-1}Ax=M^{-1}b, \\
&M_1^{-1}AM_2^{-1}y=M_1^{-1}b, \quad x=M_1^{-1}y, \\
&AM^{-1}y=b, \quad x=M^{-1}y, 
\end{align*}
for left, split, and right preconditioning, respectively. In each case, the preconditioned matrix  $M^{-1}A,\:M_1^{-1}AM_2^{-1},\:AM^{-1}$ has a better condition number and superior spectral properties to the original one. {}
For preconditioning to be feasible, the preconditioner must have two somewhat contradictory properties:
\begin{itemize}
	\item The preconditioned system must be easily solvable.
	\item The determination and the application of the preconditioner must be easy. 
\end{itemize}
The first property suggests that the preconditioner must be fairly close to the coefficient matrix of the system; however,  this is generally difficult to solve and contradicts with the second property. To conclude, the next phrase, taken from \cite{Saad2003}, summarises a view widely adopted in the research community:
\begin{displayquote}
"Finding a good preconditioner to solve a given sparse linear system is often viewed as a combination of art and science."
\end{displayquote}

In general, two classes of preconditioning techniques are available. The first includes purely algebraic methods that use only the information contained in the coefficient matrix. Such methods are typically based on some type of \emph{incomplete factorisation } or some type of \emph{sparse approximate inverse} of the coefficient matrix \cite{Benzi2002,Saad2003}. These methods achieve reasonable efficiency for a wide range of problems; however, they might not be the optimal choice for any one particular problem. The other class of methods, primarily applicable to problems arising from PDEs, involves the design of algorithms that are problem-specific. Such methods might be optimal for any specific problem; however, they require complete knowledge of the problem in advance, especially from a spectral point of view. In these methods, the preconditioners are selected by specific classes of matrices; furthermore, for their construction, a detailed spectral analysis of the coefficient matrix is required.{} 

In the first part of this thesis, preconditioning strategies to solve (using Krylov subspace methods) linear systems arising from two specific problems are proposed. In the first problems, the coefficient matrix of the system emerges as  an analytic function of a real Toeplitz matrix. This strategy utilises symmetrisation and preconditioning of the coefficient matrix. Preconditioners are selected from matrix algebras according to the spectral properties of the symmetrised coefficient matrix sequence. The properties of the matrix sequence are extensively analyzed. The second class of problems involves the numerical solution of partial differential equations with a fractional derivative order. These problems have been thoroughly investigated in recent years; however, a new category of preconditioners is here proposed. The new class of preconditioners exhibits optimal behaviour in relation to the proposals given so far in the literature, especially in dimensions of more than one. This behaviour is theoretically confirmed by the numerical results.
 
The first part of this thesis is structured as follows: The first chapter introduces all the necessary definitions and summarises the theory used to analyse the spectral properties of the coefficient matrix sequences. In the second chapter, the problem of symmetrising the large matrices that emerge as analytic functions of real Toeplitz matrices is considered. In the third chapter, is studied the numerical solution of fractional partial differential equations. 
\par

In the second part the numerical solution of a problem arising in finance is considered. In detail a numerical technique based again on an iterative algorithm is used for pricing an American put option. A put option is a financial derivative that gives the right (but not the obligation) to the holder to sell an asset for a pre-specified price $K$, the exercise or strike price. The other party, the writer of the option, must accept the sale for the exercise price regardless of the current price $S$ of the asset at the time of exercise. The right can be exercised either only on the pre-determined expiry or maturity date in case of the European put option, or any day up to the pre-determined expiry date in the American put option. In any case, because the holder does not have an obligation to sell, the right will be exercised only if the current price $S$ is lower than $K$. In this case, the profit of the seller will be $K-S$.

\par
Since they were invented, such rights have been traded on the market as assets; therefore, a fair pricing process is needed. In a viable, where  arbitrage opportunities are not allowed market model, the value of such a right at time $t$ before the maturity date should be equal to the price of a portfolio which has the same expected payoff value as the option at the exercise date at time $t$. For the European put option, it has been shown that the fair price $V(S,t)$ at time $t$ before the maturity date and for asset price $S$ must satisfy the Black--Scholes  equation. At maturity time $T$, the value $V(S,T)$ of the option must be equal to the payoff function. 

\par
Despite the similarities between the two types of put options, the American put, compared with the European put, gives its holder the additional advantage in that it can be exercised any day before the maturity date, thereby offering them additional profit opportunities. It is fair, therefore, that this privilege should be taken into account in the pricing process of the option. Moreover, according to the process in which the price of the option at time $t$ is equal to the price of a portfolio (with the same expected payoff on the exercise day) at time $t$, all the possible exercise times between should be taken into account. Although many of the characteristics of the value function of an American put option have been extensively analysed the pricing of such an options is a problem that has not yet been solved analytically.

Taking advantage of the known characteristics of the optimal value function, the iterative algorithm presented here, which utilises the principles of dynamic programming, iteratively improves exercise \emph{policies}, obtains monotonically increasing value functions and converges quadratically under reasonable assumptions. The exact meaning of a \emph{policy} will be defined in relevant chapter.    

This thesis is based on the published papers \cite{Ferrari2020} and \cite{Magirou2020}, an accepted for publication  \cite{BarErkVas} and an under progress work.

%% file: GLT.tex
\chapter{Generalized Locally Toeplitz Matrix Sequences}

\par

As mentioned in the introduction, the preconditioners for a large class of linear systems are constructed following analysis of the properties of the initial problem, which relate to the effective application of iterative methods. These properties include the condition number, asymptotic distribution of the eigenvalues, and singular values of the system coefficient matrix sequence. The theory used here to study matrix sequences is that of generalised locally Toeplitz (GLT) matrix sequences. This theory unifies, and essentially provides all the tools needed to study the asymptotic behaviour of the eigenvalues and singular values of matrix sequences obtained from the discretisation of differential or integral equations (as well as more besides). 

\par
The theory was introduced in Paollo Tilli’s paper on locally Toeplitz (LT) matrix sequences \cite{Tilli1998}. The idea was further developed by Stefano Serra-Capizzano in \cite{Capizzano2003,glt-2} and also in a series of subsequent papers. A complete presentation of the theory can be found in \cite{MR3674485,MR3674485_vol2}. In the following, we define a matrix sequence as a sequence of the form $\{A_n\}_n$, where $A_n \in \; \mathbb{C}^{n \times n}$ and $n \in \mathbb{N}$. The abbreviations LT and GLT matrix sequences denote locally Toeplitz and generalised locally Toeplitz matrix sequences, respectively.

\section{Singular Value, Eigenvalue Distribution, and Clustering of Matrix Sequences}

\begin{defin}\label{defin:glt eig distr}
Let $f:D\subset \mathbb{R}^n \to \mathbb{C}$ be a function and $\{A_n\}_n$  be a matrix sequence. $\{A_n\}_n$ has an eigenvalue distribution described by $f$, and we write $\{A_n\}_n \sim_{\lambda}f$ if
 
\begin{equation}
\eigDistr, \quad \forall F \in C_c(\mathbb{C}),
\end{equation}
where $\mu(D)$ is the Lebesgue measure of $D$, $\mu(D) \in (0,\infty)$ and $C_c(\mathbb{C})$ is the set of all continuous functions defined on $\mathbb{C}$, whose support\footnote {The support of a function $f$, denoted by $supp f$, is the set $\overline{\{x|f(x)\ne0\}}$} is a closed and bounded subset of $\mathbb{C}$. \newline

We say that $\{A_n\}_n$ has a singular value distribution described by $f$, and we write $\{A_n\}_n \sim_{\sigma}f$ if

\begin{equation}
\singDistr, \quad \forall F \in C_c(\mathbb{R}),
\end{equation} 
where $\mu(D) \in (0,\infty)$ and $C_c(\mathbb{R})$ is the set of all continuous functions defined on $\mathbb{R}$, whose support is a closed and bounded subset of $\mathbb{R}$.
\end{defin}

\par
Intuitively speaking, if the matrix sequence $\{A_n\}_n$ has an eigenvalue distribution described by $f:D\subset \mathbb{R}^n \to \mathbb{C}$, then under the condition that $f$ is continuous $a.e.$\footnote{a.e.: almost anywhere. If $f$ is not continuous over a set of zero measures at most, then $f$ is continuous almost anywhere.}, the properly rearranged eigenvalues are close to a sampling of $f$ on an equispaced grid on $D$. This definition allows eigenvalues to be out of the range of $f$; however, the total number of such eigenvalues is at most $o(n)$. The same is true in the singular value case.

\begin{rmrk}
If $A_n$ is a normal matrix for all $n$, then $\{A_n\}_n \sim_{\lambda}f$ implies that $\{A_n\}_n \sim_{\sigma}f$, because every singular value of each matrix is the absolute value of the corresponding eigenvalue of that matrix.
\end{rmrk}

\par
In the following definition, the notion of the $\epsilon-$expansion of a set $S\subset \mathbb{C}$ is used. The $\epsilon-$expansion of the set $S$ is defined as $D(S,\epsilon)=\cup_{z\in S}D(z,\epsilon)$, where $D(z,\epsilon)$ is the disc centred at $z$ with radius $\epsilon$

\begin{defin}
If $\{A_n\}_n$ is a matrix sequence and $S\subset \mathbb{C}$, we say that the eigenvalues of $\{A_n\}_n$ are strongly clustered at $S$ if, for every  $\epsilon > 0$ and every $n$, the total number of eigenvalues of $A_n$ outside of $D(S,\epsilon)$ is bounded by a constant $C(\epsilon)$ which does not depend on $n$. That is, for every $\epsilon > 0$, we have that 
\begin{equation}
\#\{j \in \{1, \dots, n\}, \lambda_j(A_n) \notin D(S,\epsilon) \} = O(1).
\end{equation}
We say that the eigenvalues of $\{A_n\}_n$ are weakly clustered at $S$ if, for every $\epsilon > 0$ and every $n$, the total number of eigenvalues of $A_n$ outside of $D(S,\epsilon)$ is bounded by a function $g(n,\epsilon)=o(n)$. That is, for every $\epsilon > 0$, we have that 
\begin{equation}
\#\{j \in \{1, \dots, n\}, \lambda_j(A_n) \notin D(S,\epsilon) \} = o(n).
\end{equation}
\end{defin}
We similarly define the notion of the strong and weak clustering of singular values of a matrix sequence on a subset of $\mathbb{R}$.

\begin{thm}\label{distributed_as_f--clustered_at_erf}
If the distribution of eigenvalues of $\{A_n\}_n$ is described by $f$, then the eigenvalues of the matrix sequence are weakly clustered in the essential range of $f$\footnote{The essential range of a function $f:D\subset \mathbb{R}^n \to \mathbb{C}$, denoted by ${\cal E}R(f)$, is the set $z \in \mathbb{C}$, for which $\mu([f\in D(z,\epsilon)])>0 \: \forall \: \epsilon >0$. Therefore, if $f$ takes a value $z_0$ outside ${\cal E}R(f)$, then $\exists \: \epsilon>0 \: \mu([f\in D(z,\epsilon)])=0$.}.
\end{thm}

\section{Approximation in Space of Matrix Sequences}
The basic tool used in the theory of GLT matrix sequences to determine the eigenvalue and singular value distributions and clusters of a matrix sequence is the \emph{closeness} of the sequence with others whose distributions and clusters are already known. The theorems presented in this section indicate the direction in which the concept of \emph{closeness} between two matrix sequences should be defined to obtain identical eigenvalue and singular value distribution and clusters. (See \cite{MR3674485} for more details.)

\par
\begin{defin}
A matrix sequence $\{A_n\}_n$ is said to be sparsely vanishing if for every $M>0$ there exists $n(M)$ such that, for $n>n(M)$,
\[
\frac{\#\{i\:\in\:\{1,\dots,n\}:\sigma_i(A_n)<1/M \}}{n}\leq r(M),
\]
where $\lim_{M\rightarrow \infty}r(M)=0$.
\end{defin}

\begin{thm} \label{thm:GLT1}
Let  $\{A_n\}_n$ and  $\{B_n\}_n$ be two matrix sequences for which, under a sufficiently large $n$, $\|A_n-B_n\|_F^2<c$. Then, the following apply: 
\begin{itemize}
	\item If the singular values of $\{B_n\}_n$ are clustered at $S$, then the singular values of $\{A_n\}_n$ are also clustered at the same set. If the matrices of the two sequences are Hermitian, the same is true for the eigenvalues.
	\item If  $\{B_n\}_n \sim_{\sigma}f$, then $\{A_n\}_n \sim_{\sigma}f$. If the matrices of the sequences are Hermitian,	then $\{B_n\}_n \sim_{\lambda}f$ implies $\{A_n\}_n \sim_{\lambda}f$. The above assertions apply even if $\|A_n-B_n\|_F^2<c(n)=o(n)$.
	\item If the matrices of $\{B_n\}_n$ are invertible, and if furthermore $\|B_n^{-1}\|<M$ for every $n$, the eigenvalues of $\{B_n^{-1}A_n\}_n$ are strongly clustered at $\{1\}$.
	\item If the condition $\|A_n-B_n\|_F^2<c$ is replaced by the relaxed one, 
    \[
	\|A_n-B_n\|_F^2<c(n)=o(n),
	\]
	then all the above apply, with the difference being that the eigenvalues (singular values) of $\{A_n\}_n$ are weakly clustered at sets that are clustered, weakly or strongly, the eigenvalues (singular values) of $\{B_n\}_n$. The eigenvalues (singular values) of $\{B_n^{-1}A_n\}_n$ are weakly clustered at $\{1\}$, if $\{A_n\}_n$ is sparsely vanishing \cite{Taud2}.
\end{itemize}

\begin{prf*}
\par
To prove the first assertion, we consider the case in which the matrices of the sequences are Hermitian. Then, we assume that the eigenvalues of $\{B_n\}_n$ are clustered at a set $S$, and that when $n$ is sufficiently large, $\|A_n-B_n\|_F^2<c$. Thus, we have
 \[
 \sum_{j=1}^{n}|\lambda_j(A_n)-\lambda_j(B_n) |^2 \leq\|A_n-B_n\|_F^2<c, 
 \]
where the first inequality comes from the well-known theorem of Hoffman and Wielandt. 
For any $\epsilon>0$, we define $N_\epsilon = \{j \: \in \{1,\dots,n\}:\:|\lambda_j(A_n)-\lambda_j(B_n)|^2\geq \epsilon\}$. Then, 
\begin{align*}
&\sum_{j=1}^{n}|\lambda_j(A_n)-\lambda_j(B_n) |^2= \\
&\sum_{j \notin N_\epsilon}^{n}|\lambda_j(A_n)-\lambda_j(B_n) |^2 + \sum_{j \in N_\epsilon}^{n}|\lambda_j(A_n)-\lambda_j(B_n) |^2 \leq c \: \Rightarrow \\
&\sum_{j \in N_\epsilon}^{n}|\lambda_j(A_n)-\lambda_j(B_n) |^2 \leq c.
\end{align*}
Let $p$ be the total number of elements of $N_\epsilon$. Then, we have $p \epsilon \leq  c \Rightarrow p \leq c \epsilon^{-1}$.
That is, there are at most $p\leq c \epsilon^{-1}$ pairs of eigenvalues such that
\[
|\lambda_j(A_n)-\lambda_j(B_n)|^2\geq \epsilon \Rightarrow |\lambda_j(A_n)-\lambda_j(B_n)|\geq \epsilon^{1/2}.
\]
If $h(n,\epsilon)$ is the number of eigenvalues of $A_n$ lying outside $D(S, \epsilon)$ and $g(n,\epsilon)$, respectively, for $B_n$, then 
\[
h(n,2\epsilon^{1/2}) \leq p + g(n,\epsilon^{1/2})\leq c\epsilon^{-1} + g(n,\epsilon^{1/2} ).
\]
The bound $c\epsilon^{-1} + g(n,\epsilon^{1/2} )$ does not depend on $n$, because we assume strong clustering at $S$ for the eigenvalues of $\{B_n\}_n$.
\newline
If the matrices of the sequences are non-Hermitian, we define
\begin{align*}
\hat{A}_{2n} = \begin{bmatrix}
0_n  &A_n \\
A_n^*&0_n
\end{bmatrix}
, \quad
\hat{B}_{2n} = \begin{bmatrix}
0_n  &B_n \\
B_n^*&0_n
\end{bmatrix}, \label{A_hat_sequences}
\end{align*}
which are Hermitian, and their eigenvalues are $\pm\sigma_j(A_n)$ and $\pm\sigma_j(B_n)$, respectively. In addition, $\|\hat{A}_{2n}-\hat{B}_{2n}\|_F=2\|A_n-B_n\|_F$\footnote{If $U\Sigma W^*$ is a singular value decomposition of the matrix $A$ and $\hat{A}= \begin{bmatrix}  
   &A \\
A^*&
\end{bmatrix}$, then the unitary matrix that diagonalizes $\hat{A}$ is $\frac{1}{\sqrt{2}}\begin{bmatrix}  
U&-U \\
W&W
\end{bmatrix}$.}. Therefore, because the assumption of the theorem is satisfied, by applying the above arguments to the modified sequences $\{\hat{A}_{2n}\}$ and $\{\hat{B}_{2n}\}$, we can conclude that the singular values of $\{A_n\}_n$ and $\{B_n\}_n$ are clustered in the same sets.

\par
To prove the second inequality  under the condition that $\|A_n-B_n\|_F^2<c$, we first assume that the matrices of the sequences are Hermitian. Again, we define 
\[
N_\epsilon = \{j \: \in \{1,\dots,n\}:\:|\lambda_j(A_n)-\lambda_j(B_n)|^2\geq \epsilon\}.
\]
For $j \in N_\epsilon$, $|\lambda_j(A_n)-\lambda_j(B_n)|\geq \epsilon^{1/2}$, and the total number of its elements is at most  $c\epsilon^{-1}$; for $F \in C_c(\mathbb{C})$, we have
\begin{align*}
&\left| \frac{1}{n}\sum_{j=1}^{n}F(\lambda_j(A_n))-F(\lambda_j(B_n)) \right|= \\
&\left| \frac{1}{n}\sum_{j \notin N_\epsilon}F(\lambda_j(A_n))-F(\lambda_j(B_n)) + \frac{1}{n}\sum_{j \in N_\epsilon}F(\lambda_j(A_n))-F(\lambda_j(B_n)) \right| \leq \\
&\left| \frac{1}{n}\sum_{j \notin N_\epsilon}F(\lambda_j(A_n))-F(\lambda_j(B_n))\right| + \left|\frac{1}{n}\sum_{j \in N_\epsilon}F(\lambda_j(A_n))-F(\lambda_j(B_n))\right| \leq \\
&\omega(\epsilon^{1/2};F) +  \frac{1}{n} 2 \|F\|_\infty c\epsilon^{-1}, 
\end{align*} 
where $\omega(\epsilon;F)$ is the modulus of continuity of $F$, \footnote{The modulus of continuity of a function $F$ is defined as $\omega(\epsilon;F)=sup\{|F(x)-F(y)|: \: |x-y|<\epsilon\}$.}. 
\begin{rmrk}
The existence of $\epsilon^{-1}$ in the representation of $\omega(\epsilon^{1/2};F) +  \frac{1}{n} 2 \|F\|_\infty c\epsilon^{-1}$ 
must be interpreted as follows. If we want $\omega(\epsilon^{1/2};F) +  \frac{1}{n} 2 \|F\|_\infty c\epsilon^{-1}<\delta$, we can take a value of $\epsilon$ that is sufficiently small for $\omega(\epsilon^{1/2};F)<\frac{\delta}{2}$ and a value of $n$ large enough that $n \geq \frac{4 \|F\|_\infty}{\delta \epsilon}$.
\end{rmrk}

If the matrices of the sequences are non-Hermitian, we apply the above conclusion to the modified sequences $\{\hat{A}_{2n}\}_n$ and $\{\hat{B}_{2n}\}_n$, for which---as previously mentioned---it follows that $\|\hat{A}_{2n}-\hat{B}_{2n}\|_F=2\|A_n-B_n\|_F$ and their eigenvalues are $\pm\sigma_j(A_n)$ and $\pm\sigma_j(B_n)$, respectively. Therefore, for every $\epsilon>0$, and for a sufficiently large $n$, we have
\begin{align*}
 &\left| \frac{1}{2n}\sum_{j=1}^{2n}F(\lambda_j(\hat{A}_{2n}))-F(\lambda_j(\hat{B}_{2n}))\right| <\frac{\epsilon}{2} \Rightarrow \\
 &\left| \frac{1}{n}\sum_{j=1}^{n}F(\sigma_j(A_n))-F(\sigma_j(B_n))+\frac{1}{n}\sum_{j=1}^{n}F(-\sigma_j(A_n))-F(-\sigma_j(B_n))\right| <\epsilon, \quad \forall F \: \in \: C_c(\mathbb{C}).
 \end{align*}
 If, in the above equation, we limit the test functions to $C_c(\mathbb{R}^+)$, we obtain the desired result.
 \par
 To prove the third inequality , we need the following propositions.
 \begin{prop}\label{glt:lemma lamda vs sigma}
 For every matrix $A \: \in \: \mathbb{C}^{n\times n}$, $\sum_{i=1}^n |\lambda_i(A)|^2 \leq \sum_{i=1}^n \sigma_i(A)^2$.
 \newline
 If $A=UTU^*$ is the Schur form of the matrix, $T$ is a complex, upper triangular matrix with eigenvalues of $A$ on its diagonal, and  $A=V\Sigma W^*$ is the singular value decomposition, we have 
 \begin{align*}
 &UT^*U^*UTU^*=W\Sigma V^*V\Sigma W^* \: \Rightarrow \: tr(UT^*TU^*)=tr(W\Sigma^2W^*) \: \Rightarrow \:
 tr(T^*T)=tr(\Sigma^2) \: \Rightarrow \\
 &\sum_{j=1}^n \sum_{i=j}^n|T_{i,j}|^2 = \sum_{i=1}^n \sigma_i(A)^2 \: \Rightarrow \sum_{i=1}^n |\lambda_i(A)|^2=\sum_{i=1}^n |T_{i,i}|^2 \leq \sum_{i=1}^n \sigma_i(A)^2.
 \end{align*}
 \end{prop}
 \begin{prop}\label{glt:lemma |AB| frobenious norn}
 If $A,B \: \in \: \mathbb{C}^{n \times n}$, then $\|AB\|_F \leq \|A\|\|B\|_F$.
 \newline
 In this case, if $b_i$ is the $ i $-th column of $B$, we have 
 \[
 \|AB\|_F^2=tr((AB)^*(AB))=\sum_{i=1}^n b_i^*A^*Ab_i = \sum_{i=1}^n \frac{b_i^*A^*Ab_i}{b_i^*b_i}b_i^*b_i \leq \|A\|^2 \sum_{i=1}^n b_i^*b_i = \|A\|^2\|B\|_F^2.
 \]
 \end{prop}
 Thus, we write
 \[
 B_n^{-1}A_n=B_n^{-1}(B_n+(A_n-B_n))=\mathbb{I}_n+B_n^{-1}(A_n-B_n).
 \] 
 Using Proposition \ref{glt:lemma |AB| frobenious norn} and the assumptions of the theorem, we have 
 \[
 \|B_n^{-1}(A_n-B_n)\|_F^2<M^2c;
 \]
 meanwhile, from Proposition \ref{glt:lemma lamda vs sigma}, we have
 \[
 \sum_{i=1}^n |\lambda_i(B_n^{-1}(A_n-B_n))|^2 \leq \sum_{i=1}^n \sigma_i(B_n^{-1}(A_n-B_n))^2 = \|B_n^{-1}(A_n-B_n)\|_F^2<M^2c.
 \]
 Using the same arguments as in the proof of the first part, we deduce that, for every $\epsilon$, a maximum of $Mc\epsilon^{-1}$ (independent of $n$) eigenvalues of $B_n^{-1}(A_n-B_n)$ are greater in absolute value than $\epsilon^{1/2}$. Therefore, the eigenvalues of $\{B_n^{-1}(A_n-B_n)\}_n$ are clustered at $\{0\}$, and the eigenvalues of $B_n^{-1}A_n=\mathbb{I}_n+B_n^{-1}(A_n-B_n)$ are clustered at $\{1\}$.
 \par
 To prove the last part, it only needs to replace the constant $c$ in the three proofs above with a function $c(n)$ which is of $o(n)$. Then, following the same steps as the proof of the first part, we deduce that the eigenvalues (singular values) of $\{A_n\}_n$ are weakly clustered at the same set where the eigenvalues (singular values) of $\{B_n\}_n$ are. Accordingly, we deduce that the eigenvalues of $\{B_n^{-1}(A_n-B_n)\}_n$ are weakly clustered at $\{0\}$ and those of $B_n^{-1}A_n=\mathbb{I}_n+B_n^{-1}(A_n-B_n)$ are clustered at $\{1\}$. The proof of the second part is entirely unaffected, because if  $c$ is constant, and $c(n)=o(n)$ is a function, then
 \[ 
 \lim_{n\rightarrow \infty}\omega(\epsilon^{1/2};F) +  \frac{1}{n} 2 \|F\|_\infty c\epsilon^{-1}=\lim_{n\rightarrow\infty}\omega(\epsilon^{1/2};F) +  \frac{1}{n} 2 \|F\|_\infty c(n)\epsilon^{-1}. 
\]
\end{prf*}
\end{thm}

\begin{thm} \label{thm:GLT2}
Let $\{A_n\}_n$ and $\{B_n\}_n$ be two matrix sequences, for which we have that, for every $n, \: rank(A_n-B_n) \leq r$, where $r$ does not depend on $n$. Then, the following apply:
\begin{itemize}
	\item If the singular values of $\{B_n\}_n$ are clustered at a set $S$, the singular values of $\{A_n\}_n$ are also clustered at the same set. In addition, the singular values of one of the matrix sequences are strongly clustered at a set if and only if the same applies for the singular values of the other sequence. If the matrices of the sequences are Hermitian, the same holds for the eigenvalues of the sequences. 
	\item If  $\{B_n\}_n \sim_{\sigma}f$, then $\{A_n\}_n \sim_{\sigma}f$. If the matrices of the sequences are Hermitian and $\{B_n\}_n \sim_{\lambda}f$, then $\{A_n\}_n \sim_{\lambda}f$.
	\item If every matrix of $\{B_n\}_n$ is invertible, then the eigenvalues of $\{B_n^{-1}A_n\}_n$ are strongly clustered at $\{1\}$.
	\item If the condition $rank(A_n-B_n) \leq r$ is replaced by $rank(A_n-B_n) \leq r(n)=o(n)$, the eigenvalues (singular values) of  $\{A_n\}_n$ are weakly clustered at the sets where the eigenvalues (singular values) of $\{B_n\}_n$ are clustered. The eigenvalues of $\{B_n^{-1}A_n\}_n$ are also weakly clustered at $\{1\}$. Furthermore, 
	\[
	\{B_n\}_n \sim_{\sigma}f \Rightarrow \{A_n\}_n \sim_{\sigma}f;
	\]
	meanwhile, if the matrices of the sequences are Hermitian, 
	\[
	\{B_n\}_n \sim_{\lambda}f \Rightarrow \{A_n\}_n \sim_{\lambda}f.
	\]
\end{itemize}
\begin{prf*}
To prove the first part, we assume that the matrices of the sequences are Hermitian. Then, 
\[
A_n-B_n=V_n^+-V_n^- \Rightarrow A_n=B_n+V_n^+-V_n^-,
\]
where $V_n^+$ and $V_n^-$ are symmetric and positive definite, with
\[
rank(V_n^+)=r_+,\quad rank(V_n^-)=r_-, \quad r_++r_-\leq r.
\]
From Weyl's theorem and conclusions deduced therefrom, we have 
\[
\lambda_{j-r_-}(B_n) \leq \lambda_j(A_n) \leq \lambda_{j+r_+}(B_n) \quad j \: \in \: \{r_-+1,\dots,r_+\}.
\]
Therefore, the number of eigenvalues of $\{A_n\}_n$ lying outside $D(S,\epsilon)$ for some set $S$ can differ from
the number of eigenvalues of $\{B_n\}_n$ lying outside $D(S,\epsilon)$, at a maximum of $r_++r_-\leq r$. If the matrices of the sequences are Hermitian, we use the sequences $\{\hat{A}_n\}_n$ and $\{\hat{B}_n\}_n$, as in the previous theorem. Then, using the same arguments, we deduce that the singular values of the initial sequences are clustered in the same sets.
\newline 
To prove the second assertion, we again assume that the matrices of the sequences are Hermitian. According to Lemma 3.3 in\cite{Tyrtyshnikov19961}, it is enough to prove that the condition 
\[
\lim_{n\rightarrow \infty}\left|\frac{1}{n}\sum_{j=1}^{n}F_{[a,\:b]}(\lambda_j(A_n))-F_{[a,\:b]}(\lambda_j(B_n))\right|=0,
\]
is met for each indicator function $F_{[a,\:b]}$, where $F_{[a,\:b]}(y)=1$ if $y \: \in [a,\:b]$, $F_{[a,\:b]}(y)=0$ else
\footnote{Every continuous function $f$ with bounded suport can be approximated by a simple staircase function of the form $\sum_{j=1}^{m}F_{[a_j,\:b_j]}f(x_j) \:x_j \: \in [a_j,\:b_j] $ }. 
From the first part of the theorem, we deduce that the total number of eigenvalues of the two sequences lying inside an interval $[a,\:b]$ can differ at a maximum of $r$. This means that 
\[
\left|\frac{1}{n}\sum_{j=1}^{n}F_{[a,\:b]}(\lambda_j(A_n))-F_{[a,\:b]}(\lambda_j(B_n))\right| < \frac{r}{n},
\]
and so the conclusion applies to all continuous with bounded support functions. This conclusion is also valid for the case in which $r=r(n)=o(n)$. Analogous to the previous theorems, we can expand the conclusion to the case in which the matrices are non-Hermitian, and we conclude that the sequences have the same singular value distribution.
\newline
For the third assertion, we only have to observe that 
\[
B_n^{-1}A_n=B_n^{-1}(B_n+(A_n-B_n))=\mathbb{I}_n+B_n^{-1}(A_n-B_n),
\]
where $rank(B_n^{-1}(A_n-B_n)) \leq r$.
\newline
If $rank(A_n-B_n) \leq r(n)=o(n)$, the total number of eigenvalues of $\{A_n\}_n$ lying outside of $D(S,\epsilon)$ for a set $S$ can differ from the number of eigenvalues of $\{B_n\}_n$ lying outside of $D(S,\epsilon)$ at a maximum of $r(n)=o(n)$. Thus, we can conclude that $\{A_n\}_n$ has weakly clustered eigenvalues (singular values), even in sets where the eigenvalues (singular values) of $\{B_n\}_n$ are strongly clustered. Similarly, the eigenvalues of $B_n^{-1}A_n$ are weakly clustered at $\{1\}$. The proof of the second part is unaffected, because 
\[
\lim_{n\rightarrow \infty}\frac{r}{n}=\lim_{n\rightarrow \infty}\frac{r(n)}{n}=0.
\]
\end{prf*}
\end{thm}

\begin{thm} \label{thm:GLT3}
Let us suppose that $\{A_n\}_n$ and $\{B_n\}_n$ are two matrix sequences, for which and for a sufficiently large $n$ we have 
\[
A_n-B_n=R_n+N_n, \quad rank(R_n)=r(n)=o(n), \quad \|N_n\|_F^2<c(n) \: \text{, $c(n)=o(n)$}.
\]
Then,
\begin{itemize}
	\item If the singular values of $\{B_n\}_n$ are clustered at a set $S$, then the singular values of $\{A_n\}_n$ are weakly clustered at the same set. If the matrices of the sequences are Hermitian, the same is valid for the eigenvalues.
	\item Furthermore,
	\[
	\{B_n\}_n \sim_{\sigma}f \Rightarrow \{A_n\}_n \sim_{\sigma}f,
	\] 
	whereas for Hermitian sequence matrices, 
	\[
	\{B_n\}_n \sim_{\lambda}f \Rightarrow \{A_n\}_n \sim_{\lambda}f.
	\]
\end{itemize}
\begin{prf*}
Applying Theorem \ref{thm:GLT2} to sequences $\{B_n\}_n$ and $\{B_n+R_n\}_n$, we conclude that the eigenvalues (singular values) of $\{B_n+R_n\}_n$ are clustered in the same sets and have the same distribution as the eigenvalues (singular values) of $\{B_n\}_n$. Then, applying theorem \ref{thm:GLT1} to sequences $\{A_n\}_n$ and $\{B_n+R_n\}_n$, we conclude that the eigenvalues (singular values) of $\{A_n\}_n$ are clustered at the same sets and have the same distribution as the eigenvalues (singular values) of $\{B_n+R_n\}_n$, and consequently with those of $\{B_n\}_n$.
\end{prf*}
\end{thm}

\section{Approximating Classes of Sequences}

\begin{rmrk} \label{rmrk:acs algebric operations on seq}
In the space of matrix sequences, we define the operations of summation, multiplication, and scalar multiplication as an extension of the corresponding matrix operations. Analogously, we define the conjugate transpose of a matrix sequence. In other words, if $\{A_n\}_n, \: \{B_n\}_n$ are matrix sequences, then $A_n, \:B_n \: \in \: \mathbb{C}^{n\times n}$ and $\alpha \: \in \: \mathbb{C}$, then
\begin{align*}
\{A_n\}_n^* &= \{A_n^*\}_n, \\
\{A_n\}_n + \{B_n\}_n &= \{A_n+B_n\}_n, \\
\alpha \{A_n\}_n &= \{\alpha A_n\}_n, \\
\{A_n\}_n\{B_n\}_n &= \{A_nB_n\}_n.
\end{align*}
\end{rmrk}

According to Theorem \ref{thm:GLT3}, two different matrix sequences have the same singular value distribution and clusters (or eigenvalue distribution if their matrices are Hermitian) if their difference can be written as the sum of two matrix sequences, as follows: The rank of the matrices of the first sequence and the Frobenius norm of the matrices of the other sequence are small compared with their size. Thus, we define the concept of convergence of matrix sequences in the context of the theory of approximating classes of sequences. Therefore, the spectral properties of the limit sequence are drawn from the spectral properties of the sequences that constitute the tail of the sequence of matrix sequences. In what follows, the abbreviation a.c.s is used for approximating classes of sequences.

\begin{defin}\label{defin:acs defination}
Let  $\{A_n\}_n$ be a matrix sequence and  $\seqOfSeqs$ be a sequence of matrix sequences. We say that $\seqOfSeqs$ is an  a.c.s.  for $\{A_n\}_n$, and we write $\{B_{n,m}\}_n \xrightarrow{\text{ a.c.s}} \{A_n\}_n$ if, for every $m$, there exists an $n_m$ such that, for $n\geq n_m$,
\begin{align}
A_n=B_{n,m}+R_{n,m}+N_{n,m}, \quad rank(R_{n,m}) \leq c(m)n, \quad \|N_{n,m}\| \leq \omega(m),
\end{align}
where $n_m, \: c(m), \: \omega(m)$ depends only on $m$, and $\lim_{m\rightarrow \infty}c(m)=\lim_{m\rightarrow \infty}\omega(m)=0$. $\| \cdot \|$ denotes the spectral norm.
\end{defin}
\begin{prop}\label{defination_acs_equievalent}
The sequence of matrix sequences $\seqOfSeqs$ is an  a.c.s. for $\{A_n\}_n$ if and only if, for every $\epsilon>0$, there exists an $m(\epsilon)$ such that for every $m>m(\epsilon)$, there exists an $n_m^\epsilon$ such that for every $n>n_m^\epsilon$, we have $A_n=B_{n,m}+R_{n,m}^\epsilon+N_{n,m}^\epsilon \quad rank(R_{n,m}^\epsilon) \leq \epsilon n, \quad \|N_{n,m}^\epsilon\| \leq \epsilon$.
\end{prop}

\par
According to this definition, it is clear that, as $m\rightarrow \infty$, the difference between the sequences $\{A_n\}_n$ and $\{B_{n,m}\}_n$ can be analysed in two sequences, as follows: The rank of the matrices of the first sequence is asymptotically small compared to the matrix size, whereas the size of the matrices of the other sequence is small weighted in the $\|.\|$ norm. That is, as $m\rightarrow \infty$, the conditions of Theorem \ref{thm:GLT3} are exactly met, with the difference that in the theorem, the Frobenius norm is used instead of the spectral norm used in the definition of the a.c.s.. However, this difference is of minor importance, as stated in the following proposition.

\begin{prop}
Let  $\{A_n\}_n$ be a matrix sequence and  $\seqOfSeqs$ be a sequence of matrix sequences. Then,
\[
A_n=B_{n,m}+R_{n,m}+N_{n,m}, \quad rank(R_{n,m}) \leq c(m)n, \quad \|N_{n,m}\| \leq \omega(m),
\]
where $lim_{m\rightarrow \infty}c(m)=lim_{m\rightarrow \infty}\omega(m)=0$, if and only if,
\[
A_n=B_{n,m}+\hat{R}_{n,m}+\hat{N}_{n,m}, \quad rank(\hat{R}_{n,m}) \leq \hat{c}(m,n), \quad \|\hat{N}_{n,m}\|_F^2 \leq \hat{\omega}(m,n),
\]
where $lim_{m\rightarrow \infty}\hat{c}(m,n)=o(n)$ and $lim_{m\rightarrow \infty}\hat{\omega}(m,n)=o(n)$.
\begin{prf*}
To prove the direct, we set
\[
\hat{R}_{n,m}=R_{n,m}, \quad \hat{N}_{n,m}=N_{n,m}, \quad \hat{c}(m,n)=c(m)n, \quad \hat{\omega}(m,n)=\omega(m)^2n.
\]
Applying the equivalent definition of a.c.s. (Proposition \ref{defination_acs_equievalent}) (i.e., for every $\epsilon >0$, there exists an $m(\epsilon)$ such that for $m>m(\epsilon)$, we have $c(n)<\epsilon$ and $\omega(m)<\epsilon$), we find that, for $m>m(\epsilon)$,
\[
\frac{rank(\hat{R}_{n,m})}{n}=\frac{\hat{c}(m,n)}{n}<\epsilon, \quad \frac{\|\hat{N}_{n,m}\|^2}{n}= \frac{\hat{\omega}(m,n)}{n}<\epsilon^2.
\] 
 which completes the proof. To prove the opposite, we observe that if $\|\hat{N}_{n,m}\|_F^2<\hat{\omega}(m,n)$, then 
\[
p=\#\{j \: \in \: \{1,\dots,n\}: \quad \sigma_j^2(A_{n,m})\geq \frac{1}{m}\}< m\hat{\omega}(m,n).
\]
Because p is the total number of singular values of $\hat{N}_{n,m}$ that exceed $\sqrt{\frac{1}{m}}$, we can write $\hat{N}_{n,m}=\hat{N}_{n,m}^R+\hat{N}_{n,m}^N$, where $rank(\hat{N}_{n,m}^R)<m\hat{\omega}(m,n)$ and $\|\hat{N}_{n,m}^N\|<\sqrt{\frac{1}{m}}$. Then, we set
\[
R_{n,m}=\hat{R}_{n,m}+\hat{N}_{n,m}^R, \quad N_{n,m}=\hat{N}_{n,m}^N, \quad nc(m)=\hat{c}(n,m)+m\hat{\omega}(m,n),  \quad \omega(m)=\sqrt{\frac{1}{m}},
\]
and the proof is complete.
\end{prf*}
\end{prop}

The consequences and importance of the above definition are summarised in the two theorems that follow.

\begin{thm} \label{thm:acs sigma}
Let $\{A_n\}_n$ be a matrix sequence and $\seqOfSeqs$ be a sequence of matrix sequences such that $A_n, \: B_{n,m} \: \in \mathbb{C}^{n \times n}$. Let $f,\:f_m:D\subset \mathbb{R}^n \rightarrow \mathbb{C}$ and
\begin{align*}
 &\{B_{n,m}\}_n \sim_{\sigma} f_m, \\
 &\{B_{n,m}\}_n \xrightarrow{\text{ a.c.s}} \{A_n\}_n, \\
 &f_m \rightarrow f  \quad  \text{ in measure}. 
\end{align*}
Thus, $\{A_n\}_n \sim_{\sigma} f $.
\end{thm}

\begin{thm} \label{thm:acs lambda}
Let $\{A_n\}_n$ be a matrix sequence and $\seqOfSeqs$ be a sequence of matrix sequences such that $A_n, \: B_{n,m} \: \in \mathbb{C}^{n \times n}$ Hermitians. Let $f,\:f_m:D\subset \mathbb{R}^n \rightarrow \mathbb{C}$ and
\begin{align*}
&\{B_{n,m}\}_n \sim_{\lambda} f_m, \\
&\{B_{n,m}\}_n \xrightarrow{\text{ a.c.s}} \{A_n\}_n, \\
&f_m \rightarrow f \quad  \text{ in measure }. 
\end{align*}
Then, $\{A_n\}_n \sim_{\lambda} f $.
\end{thm}

The proof of the two theorems is a direct consequence of the definition of the a.c.s sequences and Theorem \ref{thm:GLT3}. As $m\rightarrow \infty$, the requirements of the theorem are satisfied, and the singular value distribution (eigenvalue distribution) of $\{A_n\}_n$ is close to that of $\{B_{n,m}\}_n$. Therefore, for the singular value case and for $F \: \in \: C_c(\mathbb{R})$, we have
\begin{align*}
&\left|\frac{1}{n}\sum_{i=1}^nF(\sigma_i(A_n))- \frac{1}{\mu(D)}\int_D F(|f(\theta)|)d\theta\right| \leq \left|\frac{1}{n}\sum_{i=1}^nF(\sigma_i(A_n))-\frac{1}{n}\sum_{i=1}^nF(\sigma_i(B_{n,m}))\right| \\
&+  \left|\frac{1}{n}\sum_{i=1}^nF(\sigma_i(B_{n,m})) - \frac{1}{\mu(D)}\int_D F(|f_m(\theta)|)d\theta\right| + \left|\frac{1}{\mu(D)}\int_D F(|f_m(\theta)|)d\theta - \frac{1}{\mu(D)}\int_DF(|f(\theta)|)d\theta\right|.
\end{align*}
We have proven that the first term in the above sum has a limit at zero as $m\rightarrow \infty$, whereas the second term has a limit at zero for each $m$ as $n\rightarrow \infty$, by definition. The difference between the two integrals has a limit at zero as $m\rightarrow \infty$, because of the requirement that $f_m \rightarrow f$ in the measure of the two theorems. A similar situation occurs in the eigenvalue case.
\par
It must be mentioned that these two theorems can be proven without changing the norm $\|.\|$ to $\|.\|_F$, though this proof is rather technical and straightforward. The approach used here has been chosen to clarify the two directions and the corresponding tolerance limits of deviation allowed between \emph{nearby} elements in the space of matrix sequences. In the following proposition, are summarised the algebraic properties of the a.c.s. sequences.

\begin{prop}\label{prop:acs algebric properties of acs}
Let $\{A_n\}_n$, $\{\hat{A}_n\}_n$  be matrix sequences, $\alpha,\beta \: \in \: \mathbb{C}$, and 
\begin{itemize}
	\item $\{B_{n,m}\}_n\acsArrow \{A_n\}_n$,
	\item $\{\hat{B}_{n,m}\}_n\acsArrow\{\hat{A}_n\}_n$.
\end{itemize}
Then,

\begin{itemize}
	\item $\{B_{n,m}\}_n^*\acsArrow \{A_n\}_n^*$,
	\item $\{\alpha B_{n,m}+\beta \hat{B}_{n,m}\}_n \acsArrow \{\alpha A_n+\beta \hat{A}_n\}_n$.
	\item Suppose that, for each $M>0$, the total number of singular values of $A_n$ that are further from $M$ as $n\rightarrow \infty$ is bounded by a function $r(M)$, where $\lim_{M\rightarrow \infty}r(M) = 0$; the same applies for $\{\hat{A}_n\}_n$. Then \\
	 $\{B_{n,m} \hat{B}_{n,m}\}_n \acsArrow \{A_n \hat{A}_n\}_n$.
\end{itemize}
\end{prop} 

\subsection{Zero Distributed Sequences}
\par
A central role in the theory of approximation in the space of matrix sequences and in the theory of GLT matrix sequences is played by the class of matrices whose singular values are distributed to zero. The exact definition of the class and the theorem that uses the conclusions of a.c.s are as follows: 

\begin{defin} \label{def:zero distributed seq}
We say that a matrix sequence $\{Z_n\}_n$ is zero distributed if $\{Z_n\}_n \sim_{\sigma} 0$. That is,
\begin{align*}
\lim_{n\rightarrow \infty} \: \frac{1}{n} \sum_{j=1}^{n}F(\sigma_j(Z_n)) = F(0) \: \forall \: F \: \in \: C_c(\mathbb{R}).
\end{align*}
It can be shown that a matrix sequence $\{Z_n\}_n$ is zero-distributed if and only if
\begin{align*}
\lim_{n\rightarrow \infty}\frac{\# \{j\in\{1,\dots n\}, \: \sigma_j(Z_n)>\epsilon\}}{n}=0 \: \forall \: \epsilon>0,
\end{align*} 
or, equivalently, 
\begin{align*}
\forall \epsilon>0, \: \exists \: n(\epsilon), \: \forall \: n>n(\epsilon) \quad \frac{\# \{j\in\{1,\dots n\}, \: \sigma_j(Z_n)>\epsilon\}}{n}<\epsilon . 
\end{align*}
A matrix sequence $Z_n$ is zero distributed if and only if 
\[
\forall \: n, \:\: Z_n= R_n + N_n, \quad \text{where,} \:\:\lim_{n\rightarrow\infty}\frac{rank(R_n)}{n}=\lim_ {n\rightarrow\infty}\|N_n\|=0.
\]
If $\{Z_n\}_n$ is zero distributed, and the matrices of the sequence are Hermitian, then $\{Z_n\}_n \sim_{\lambda} 0$, because in that case, $\sigma_j(Z_n)=|\lambda_j(Z_n)|$. Clearly, in that case, 
\begin{align*}
\lim_{n\rightarrow \infty}\frac{\# \{j\in\{1,\dots n\}, \: |\lambda_j(Z_n)|>\epsilon\}}{n}=0 \: \forall \: \epsilon>0, 
\end{align*}
and
\begin{align*}
\forall \epsilon>0, \: \exists \: n(\epsilon), \: \forall \: n>n(\epsilon) \quad \frac{\# \{j\in\{1,\dots n\}, \: |\lambda_j(Z_n)| > \epsilon \}}{n} < \epsilon.
\end{align*}
\end{defin}

\begin{thm} \label{thm:distr seq plus zero distributed}
Let $\{A_n\}_n$, $\{B_n\}_n$, and $\{Z_n\}_n$ be matrix sequences, where $\{Z_n\}_n$ is zero distributed and also applies that $A_n=B_n+Z_n \: \forall n \in \mathbb{N}$. Then, 
\[
\{B_n\}_n \sim_{\sigma}f \Rightarrow \{A_n\}_n \sim_{\sigma}f.
\]
If the matrices of the sequences are Hermitian,
\[
\{B_n\}_n \sim_{\lambda}f \Rightarrow \{A_n\}_n \sim_{\lambda}f.
\]
\begin{prf*}
Let $Z_n=U_n\Sigma_nV_n^*$ (where $\Sigma_n=diag(\sigma_i) \: \forall \: i=1\dots n$) represent the singular value decomposition of $Z_n$. We define 
\begin{align*}
\Sigma_{n,m}^{R} = diag(x_{[\frac{1}{m},\infty)}(\sigma_i)\sigma_i), \quad \Sigma_{n,m}^{N} = diag(x_{[0, \frac{1}{m})}(\sigma_i)\sigma_i)) \quad i=1\dots n,
\end{align*}
where $x_{[\alpha,\beta]}$ is the indicator function of $[\alpha,\beta]$. In other words, the singular values of $\Sigma_{n,m}^{R}$ are the singular values of $Z_n$ which are greater than or equal to $\frac{1}{m}$, while the remainder are zero. The singular values of $\Sigma_{n,m}^{N}$ are the values of $Z_n$ which are less than $\frac{1}{m}$, whilst the others are zero. We now define 
\begin{align*}
R_{n,m}=U_n\Sigma_{n,m}^{R}V_n^*, \quad N_{n,m}=U_n\Sigma_{n,m}^{N}V_n^*, \quad B_{n,m}=B_n  \quad \forall \: n, m \: \in \mathbb{N}.
\end{align*}
Clearly, $A_n=B_n+Z_n=B_{n,m}+R_{n,m}+N_{n,m}$. Let $\epsilon >0$ and $m(\epsilon)= min\{m \in \mathbb{N}, \: \frac{1}{m}<\epsilon)\}$. By definition, if $m>m(\epsilon)$, then $\|N_{n,m}\|<\epsilon$. In addition, because $\{Z_n\}_n$ is zero distributed, $\exists \: n(\epsilon):=n_m^{\epsilon}$ such that if $n>n_m^{\epsilon}$, $\frac{\# \{j\in\{1,\dots n\}, \: \sigma_j(Z_n)>\epsilon\}}{n}<\epsilon$. Thus, $\frac{rank(R_{n,m})}{n}<\epsilon$.
\newline
From the above, it is clear that $\{B_{n,m}\}_n \acsArrow \{A_n\}_n$. Then, applying Theorem \ref{thm:acs sigma} for $f_m=f$, we deduce that $\{A_n\}_n \sim_{\sigma} f$. If the matrices are Hermitian, we follow the same steps using eigenvalues instead of singular values and apply Theorem \ref{thm:acs lambda}; thus, we deduce that $\{A_n\}_n \sim_{\lambda} f$.
\end{prf*}
\end{thm}

\section{Circulant and Toeplitz Matrices}
Circulant matrices are those of the form
\begin{align*}
C_n = \ciclo. 
\end{align*}
Circulant matrices are diagonalised using the unitary discrete Fourier transform. Their spectrum is known; furthermore, because of their properties, they play a major role in the analysis and design of techniques for solving structured linear systems. Here, in addition to presenting their basic properties, they are used as an example of the use of a.c.s. theory to analyse the distribution of Toeplitz matrix sequences.

\par
Let $F_n$ be a unitary discrete Fourier transform. That is,
\begin{align}
[F_n]_{k,j}=\frac{1}{\sqrt{n}}e^{-i2\pi kj/n}=\frac{1}{\sqrt{n}}w_n^{-kj} \quad k,j=0, \dots n-1 \quad w_n= e^{i2\pi/n}, \label{defin:unitary discrete Fourier }.
\end{align}
It is known that $F_n^{-1}=F_n^*$, where $\phantom{}^*$ is the conjugate transpose operator. 

If $f_{j,n}$ is the $j-th$ column $F_n$, then
\begin{align*}
C_n f_{j,n}=\ciclo \frac{1}{\sqrt{n}} \begin{bmatrix}
w_n^{-0j}     \\
w_n^{-1j}     \\
w_n^{-2j}     \\
\vdots        \\
\vdots        \\
w_n^{-(n-1)j} \\
\end{bmatrix} \\
=\frac{1}{\sqrt{n}}\begin{bmatrix}
&c_0w_n^{-0j}     &+&  c_{n-1}w_n^{-1j}  &+& c_{n-2}w_n^{-2j} &+& \cdots &+& c_1w_n^{-(n-1)j} \\ 
&c_1w_n^{-0j}     &+&  c_{0  }w_n^{-1j}  &+& c_{n-1}w_n^{-2j} &+& \cdots &+& c_2w_n^{-(n-1)j} \\
&c_2w_n^{-0j}     &+&  c_{  1}w_n^{-1j}  &+& c_{2  }w_n^{-2j} &+& \cdots &+& c_3w_n^{-(n-1)j} \\
&\vdots    \\     
&\vdots    \\     
&c_{n-1}w_n^{-0j} &+&  c_{n-2}w_n^{-1j}  &+& c_{n-3}w_n^{-2j} &+& \cdots &+& c_0w_n^{-(n-1)j}
\end{bmatrix} \\
=\frac{1}{\sqrt{n}}\begin{bmatrix}
&c_0w_n^{0j}    &+&  c_{n-1}w_n^{(n-1)j} &+& c_{n-2}w_n^{(n-2)j} &+& \cdots &+& c_1w_n^{1j} \\ 
&c_1w_n^{0j}    &+&  c_{0  }w_n^{(n-1)j} &+& c_{n-1}w_n^{(n-2)j} &+& \cdots &+& c_2w_n^{1j} \\
&c_2w_n^{0j}    &+&  c_{  1}w_n^{(n-1)j} &+& c_{2  }w_n^{(n-2)j} &+& \cdots &+& c_3w_n^{1j} \\
&\vdots    \\  
&\vdots    \\ 
&c_{n-1}w_n^{0j}&+& c_{n-2}w_n^{(n-1)j}  &+& c_{n-3}w_n^{(n-2)j} &+& \cdots &+& c_0w_n^{1j},
\end{bmatrix}
\end{align*}
because $w_n^{-kj}=w_n^{(n-k)j}$. Taking the common factor $w_n^{-kj}$ for $k=0,\dots,n-1$ at the $k$-line, the above becomes

\begin{align*}
\frac{1}{\sqrt{n}}\begin{bmatrix}
&w_n^{-0j}(c_0w_n^{0j}             &+&  c_{n-1}w_n^{(n-1)j} &+& c_{n-2}w_n^{(n-2)j} &+& \cdots  &+& c_1w_n^{1j})  \\ 
&w_n^{-1j}(c_1w_n^{1j}             &+&  c_{0  }w_n^{0j}     &+& c_{n-1}w_n^{(n-1)j} &+& \cdots  &+& c_2w_n^{2j}) \\
&w_n^{-2j}(c_2w_n^{2j}             &+&  c_{  1}w_n^{1j}     &+& c_{2  }w_n^{(2)j}   &+&  \cdots &+& c_3w_n^{3j}) \\
&\vdots    \\  
&\vdots    \\ 
&w_n^{-(n-1)j}(c_{n-1}w_n^{(n-1)j} &+& c_{n-2}w_n^{(n-2)j}  &+& c_{n-3}w_n^{(n-3)j} &+& \cdots  &+& c_0w_n^{0j})
\end{bmatrix}\\
=\frac{1}{\sqrt{n}}\sum_{k=0}^{n-1}c_k(w_n^j)^k \begin{bmatrix}
w_n^{-0j}     \\
w_n^{-1j}     \\
w_n^{-2j}     \\
\vdots        \\
\vdots        \\
w_n^{-(n-1)j} \\
\end{bmatrix} = p(\theta_{j,n})f_{j,n}   \quad j=0,\dots ,n-1,
\end{align*}
where $\theta_{j,n}=\frac{2\pi j}{n}$ and $p(\theta)=\sum_{k=1}^{n-1}c_ke^{ik\theta}$.

\par
Now, we assume a trigonometric polynomial of order $m$, $p_m(\theta)=\sum_{k=-m}^{m}\hat{c}_ke^{ik\theta}$, and a circulant matrix of size $n$ (with $n>2m$), $C_n(p_m)$, whose elements are defined as follows:
\begin{align}
c_k&=\hat{c}_k \quad k=0,\dots m \label{glt:C_n(p_m) defination start} \\
c_k&=0 \quad  m < k < n-m \\
c_k&=\hat{c}_{k-n} \quad k=n-m,\dots n-1. \label{glt:C_n(p_m) defination end}
\end{align}
Then, the eigenvalues of $C_n(p_m)$ according to the above are $\sum_{k=0}^{n-1}c_ke^{i2\pi jk/n}$ for $j=0,\dots n-1$. However, then we have
\begin{align*}
&\sum_{k=0}^{n-1}c_ke^{i2\pi jk/n} = \sum_{k=0}^{m}\hat{c}_ke^{i2\pi jk/n} + \sum_{k=n-m}^{n-1}\hat{c}_{k-n}e^{i2\pi jk/n} = \\
&\sum_{k=0}^{m}\hat{c}_ke^{i2\pi jk/n} + \sum_{k=-m}^{-1}\hat{c}_{k}e^{i2\pi j(k+n)/n} = \sum_{k=-m}^{m}\hat{c}_ke^{i2\pi jk/n}. 
\end{align*}
Finally, the eigenvalues of $C_n(p_m)$ are the values of $p_m(\theta)$ at the $n$ points $\theta_{j,n}$; in other words, the eigenvalues of the matrix  $C_n(p_m)$ constitute a uniform sampling of the function $p_m(\theta)$ at the interval $[0,2\pi]$. It is clear that the eigenvalues of $C_n(p_m)$ are distributed as $p_m(\theta)$. Because the matrix is normal, the same applies to its singular values. To summarise,

\begin{align}
C_n(p_m) \sim_{\sigma,\lambda} p_m.
\end{align}

\subsection{Toeplitz Matrices and Toeplitz Matrix Sequences}
It is known from Fourier analysis that if $f$ is a Lebesgue integral function, defined at $[-\pi,\pi]$, and 
\begin{equation}\label{fourier}
a_{k}=\frac{1}{2\pi} \int_{-\pi} ^{\pi}f(\theta) e^{-i k \theta } \,d\theta,\quad k=0,\pm1,\pm2,\dots,
\end{equation}
then 
\[
f(\theta)=\sum_{k=-\infty}^{\infty}a_{k}e^{ik\theta}.
\]
The above series is called the Fourier series of the function $f$, and it extends periodically across the real line. The coefficients $a_k$ are the  Fourier coefficients of the function.
\begin{defin}\label{glt:defination Toeplitz}
The Toeplitz matrix of size $n$, which is related to the function $f$ via $T_n(f)$, is defined as
\[
T_n(f)=\left[a_{k-j}\right]_{k,j=1}^n,
\]

\[
T_{n}(f)=\begin{bmatrix}{}
a_0 & a_{-1} & \cdots & a_{-n+2} & a_{-n+1} \\
a_1 & a_0 & a_{-1}   &  & a_{-n+2} \\
\vdots & a_1 & a_0 & \ddots & \vdots \\
a_{n-2} &  & \ddots & \ddots & a_{-1} \\
a_{n-1} & a_{n-2} &\cdots & a_1 & a_0
\end{bmatrix}. 
\]
The $f$ is referred to as the generating function of the matrix sequence $\{T_n(f)\}_n$. 
\end{defin}
In the following, several properties of the Toeplitz matrices are given \cite{MR3674485}. Let $f,g \in L^{1}[-\pi,\pi]$, and let $\alpha \in \mathbb{C}$. Then, 
\begin{itemize} 
	\item $T_n(\alpha f)=\alpha T_n(f)$.
	\item $T_n(f+g)=T_n(f)+T_n(g)$.
	\item If $f$ is real, then the $T_n(f)$ is Hermitian and its eigenvalues are contained in the interval $(m_f,M_f)$, where $m_f=ess inf(f)$, $M_f=ess sup(f)$.
	\item $\|T_n(f)\|<M_{|f|}=\|f\|_{L^{\infty}}$ if $|f|$ is not constant almost everywhere. Otherwise $\|T_n(f)\|\leq M_{|f|}$. If $m_f=M_f$ then $f=m_f$ almost everywhere and $T_n(f)=m_f\mathbb{I}_n$.
\end{itemize}
\par
The eigenvalue and singular value distribution of the matrix sequence  $\{T_n(f)\}_n$ have been extensively studied. Initially, Szeg{\H{o}} in \cite{MR890515} proved that if $f\in L^{\infty}([-\pi,\pi])$, then the eigenvalues of the Toeplitz matrix sequence $\{T_n(f)\}_n$ are distributed as $f$. Since then, this conclusion has been extended to sequences with the generating function  $f\in L^{1}([-\pi,\pi])\supset  L^{\infty}([-\pi,\pi])$ complex \cite{MR952991,MR851935,Tyrtyshnikov19961,MR1258226,MR1481397,garoni2015}. The results of this evolutionary research process are summarised in the following theorem, for which a proof borrowed from \cite{MR3674485} is given here as an example of the use of a.c.s. to find the distribution of matrix sequences.

\begin{thm}\label{thm:glt_toeplitz distribution}
Let $f \in L^{1}([-\pi,\pi])$ and $T_n(f)$ be the Toeplitz matrix with generating function $f$. Then,
\[
\{T_n(f)\}_n \sim_{\sigma} f.
\]
If  $f$ is real, then,
\[
\{T_n(f)\}_n \sim_{\lambda} f.
\]
\end{thm}
\begin{prf*}
Let $f_m(\theta)=\sum_{k=-m}^{m}c_ke^{ik\theta}$ be a trigonometric polynomial of order $m$, and let $C_n(f_m)$ be a circulant matrix of size $n$, related to $f_m$ and defined as in (\ref{glt:C_n(p_m) defination start})–(\ref{glt:C_n(p_m) defination end}). Then,
\[
C_n(f_m)-T_n(f_m)=\begin{bmatrix}
    0       & \cdots   & 0        & c_{m}    & \cdots  & c_{1}   \\
    \vdots  & \ddots   & \ddots   & \ddots   & \ddots  & \vdots  \\
    0       & \ddots   & \ddots   & \ddots   & \ddots  & c_{m}   \\
    c_{-m}  & \ddots   & \ddots   & \ddots   & \ddots  & 0       \\ 
    \vdots  & \ddots   & \ddots   & \ddots   & \ddots  & \vdots  \\
    c_{-1}  & \cdots   &  c_{-m}  & 0        & \cdots  & 0
\end{bmatrix},
\]
and $rank(C_n(f_m)-T_n(f_m))\leq 2m$. Then, defining 
\[
T_n(f_m) = C_n(f_m)+(T_n(f_m)-C_n(f_m))=C_n(f_m)+Z_n, \quad Z_n=(T_n(f_m)-C_n(f_m))+0_n.
\]
Clearly, the sequence $\{Z_n \}_n$ is zero distributed, and according to Theorem \ref{thm:distr seq plus zero distributed}, we have that
\[
T _n(f_m) \sim_{\sigma} f_m.
\]
If the polynomial $f_m$ is real, then all the matrices are Hermitian, and using the same theorem again gives 
\[
T _n(f_m) \sim_{\lambda} f_m.
\]
Let now $f \in L^{1}[-\pi, \pi]$, $f(\theta)=\sum_{k=-\infty}^{\infty}c_ke^{ik\theta}$, and $f_m(\theta)=\sum_{k=-m}^{m}c_ke^{ik\theta}$. Then,
\[
\|T_n(f)-T_n(f_m)\|=\|T_n(f-f_m)\|\leq\|f-f_m\|_{L^{\infty}}.
\]

Owing to the uniform convergence of $f_m$ to $f$ at the interval $[-\pi,\pi]$, we have $\lim_{m\rightarrow \infty}\|f-f_m\|_{L^\infty}=0$. In addition, the uniform convergence of $f_m$ to $f$ implies convergence in measure $f_m$ to $f$ in the same interval. So,
\[
T_n(f) = T_n(f_m) + (T_n(f)-T_n(f_m)) = T_n(f_m) + R_{n,m} + N_{n,m} \quad R_{n,m}=0_n, \quad N_{n,m}=(T_n(f)-T_n(f_m)).
\]
Clearly, $T_n(f_m)$ is an a.c.s for $T_n(f)$. Provided that the other conditions of Theorem \ref{thm:acs sigma} also apply, we have that
\[
T_n(f)\sim_{\sigma} f.
\]
If $f$ is real, all the matrices are Hermitian, and all the conditions of Theorem \ref{thm:acs lambda} apply. Consequently,
\[
T_n(f)\sim_{\lambda} f.
\]
\end{prf*}

\section{LT and GLT  Matrix Sequences}

The main source (although certainly not the only one) of problems in the case of large and typically sparse matrices is the discretisation of differential and integral equations. The structures of the matrices appearing in such problems depend on the numerical scheme selected for the discretisation of the specific differential or integral operator. Because this scheme is unchanged in terms of displacement, the matrices produced are Toeplitz. The elements of the coefficient matrix are constant over each diagonal, because the same scheme is chosen for the discretisation of the operator at each point of the unknown function’s domain. 
However, when the unknown function appears in the equation with a non-constant coefficient, all non-zero elements of the Toeplitz matrix are multiplied by the corresponding values of the function. In other words, a sampling of the coefficient function of the differential equation lies along the non-zero diagonals, the coefficient matrix is no longer Toeplitz, and its spectral distribution is not given by the known theorems. In the context of the GLT theory, almost every matrix sequence produced from the discretisation of a differential or integral equation can be approximated in an a.c.s sense by another matrix sequence for which the spectral distribution is known.
\par
The basic definitions and conclusions of the theory are presented in the following subsections.

\subsection{LT  Matrix Sequences}

\begin{defin}
Let $n,m \: \mathbb{N}$, $\alpha:[0,\: 1]\rightarrow \mathbb{C}$, and $f\: \in \: L^1[-\pi, \: \pi]$. Then,
\begin{itemize}
	\item The locally Toeplitz operator is defined as an $n\times n$ matrix,   
    \[
     LT_m^n(\alpha,\: f)=D_m(\alpha)\otimes T_{\lfloor \frac{n}{m}\rfloor}(f)\oplus O_{n\: mod \:m}=diag_{i=1\dots m}[\alpha(\frac{i}{m})T_{\lfloor \frac{n}{m}\rfloor}(f)]\oplus O_{n\: mod \:m},
    \] 
    where $D_m(\alpha)$ is the diagonal matrix of size $m$, and the elements are a uniform sampling of $\alpha$ in $[0,\: 1]$. That is,
    \begin{align} 
    LT_m^n(\alpha,\:f)=
    \begin{bmatrix}
               \alpha(\frac{1}{m})T_{\lfloor \frac{n}{m}\rfloor}(f) &  &  &  &    \\
               &  \alpha(\frac{2}{m})T_{\lfloor \frac{n}{m}\rfloor}(f) &  &  &    \\
               &  &                                              \ddots   &  &    \\
               &  &  &          \alpha(1)T_{\lfloor \frac{n}{m}\rfloor}(f)   &    \\
               &  &  &  &                                        O_{n\:mod\:m}
    \end{bmatrix}.
    \end{align}
    \item Provided that the function $\alpha$ is  Riemann integrable and $\alpha \otimes f=\alpha(x)f(\theta)$, then we can define $\{A_n\}_n$ as a locally Toeplitz sequence, with symbol $\alpha \otimes f$, if 
    \[
    \{LT_m^n(\alpha,f)\}_n \acsArrow \{A_n\}_n.
    \]
     For a locally Toeplitz sequence with symbol $\alpha \otimes f$, we write $\{A_n\} \sim_{LT} \alpha \otimes f$.
\end{itemize}
\end{defin}

\par

\begin{thm}
Let $\alpha_1,\dots\alpha_p:[0,\:1]\rightarrow\mathbb{C}$ and $f_1,\dots f_p \: \in L_1[-\pi,\:\pi]$. Then, for each $m,\:\in \mathbb{N}$ and $F \: \in \: C_c(\mathbb{R})$, we apply
\begin{itemize}
	\item \begin{align*}
	&\lim_{n\rightarrow \infty} \frac{1}{n} \sum_{j=1}^{n}F\left(\sigma_j\left(\sum_{i=1}^pLT_m^n\left(\alpha_i,f_i\right)\right)\right)=\frac{1}{m}\sum_{k=1}^m\frac{1}{2\pi}\int_{-\pi}^{\pi}F\left(\left|\sum_{i=1}^p\alpha_i\left(\frac{k}{m}\right)f_i(\theta)\right|\right)d\theta.
	\end{align*}
	\item \begin{align*}
	&\lim_{n\rightarrow \infty} \frac{1}{n} \sum_{j=1}^{n}F\left(\lambda_j\left(\Re\left(\sum_{i=1}^pLT_m^n\left(\alpha_i,f_i\right)\right)\right)\right)=\frac{1}{m}\sum_{k=1}^m\frac{1}{2\pi}\int_{-\pi}^{\pi}F\left(\Re\left(\sum_{i=1}^p\alpha_i\left(\frac{k}{m}\right)f_i(\theta)\right)\right)d\theta.
	\end{align*}
	\item If $\alpha_1,\dots\alpha_p$ are Riemann integrable\footnote{ 
	The requirement for the function $\alpha$ to be Riemann integrable  is necessary to obtain 
	\[
	\{A_n\} \sim_{LT} \alpha_i \otimes f\: \Rightarrow \:\{A_n\} \sim_\sigma \alpha_i \otimes f.
	\] 
	For example, we define $\alpha:[0,\:1]\rightarrow \mathbb{R}$ with $\alpha(x)=0$ if $x\: \in \: \mathbb{Q}$, and $\alpha(x)=1$ otherwise. In this case, $\alpha$ is not Riemann integrable. Then, ${LT_m^n(\alpha,f)}_n=\mathbb{O}_n$, whilst $\alpha \otimes f = f$ almost everywhere.} and $\{A_n^i\} \sim_{LT} \alpha_i \otimes f_i$; then, 
	\[
	\left\{\sum_{i=1}^pA_n^i\right\}_n\sim_{\sigma}\sum_{i=1}^p\alpha_i \otimes f_i \quad and \quad \left\{\Re\left(\sum_{i=1}^pA_n^i\right)\right\}_n\sim_{\lambda}\Re\left(\sum_{i=1}^p\alpha_i \otimes f_i\right).
	\]
	\item If the matrices of the sequences $\{A_n^i\}_n$ are Hermitian, then $\alpha_i \otimes f_i$ are real almost everywhere, and $\left\{\sum_{i=1}^pA_n^i\right\}_n\sim_{\lambda}\sum_{i=1}^p\alpha_i \otimes f_i$.
\end{itemize}
\end{thm}

The most important consequence of the classification of a matrix sequence as LT is the immediate characterisation of the distribution of its singular values, or its eigenvalues in the Hermitian case.
\newline
It can be proved that the Toeplitz matrix sequences, the sequences of diagonal matrices whose elements are a uniform sampling of a function $\alpha:[0,\:1]\rightarrow \mathbb{C}$, and the zero distributed sequences belong to LT class. More specifically,
\begin{itemize}
	\item $f \: \in \: L^1([-\pi,\:\pi])\:\Rightarrow\:\{T_n(f)\}_n\sim_{LT}1\otimes f=f$,
	\item $\alpha:[0,\:1]\rightarrow \mathbb{C} \: and \: D_n(\alpha)=diag(\alpha(\frac{i}{n}))$ for $i=1,\dots n$, then $\{D_n(\alpha)\}_n\sim_{LT}\alpha \otimes 1 = \alpha$
	\item $\{Z_n\}_n \sim_{\sigma} 0 \: \Rightarrow \: \{Z_n\}_n \sim_{LT} 0$.
\end{itemize}

\subsection{GLT Matrix Sequences}

A matrix sequence is GLT if it is the limit in the a.c.s. sense of a finite sum of LT sequences. That is,

\begin{defin}
Let $\{A_n\}_n$ be a matrix sequence and $\kappa:[0,\:1]\times[-\pi,\:\pi]\rightarrow\mathbb{C}$ be a measurable function. $\{A_n\}_n$ is a GLT sequence, with symbol $\kappa$; then, we write $\{A_n\}_n\sim_{GLT}\kappa$;  if, for each $m \: \in \mathbb{N}$, there exists a finite number of LT sequences  $\{A_n^{i,m}\}_n\sim_{LT}\alpha_{i,m}\otimes f_{i,m}$ such that
\begin{itemize}
	\item $\sum_{i=1}^{N_m}\alpha_{i,m}\otimes f_{i,m}\rightarrow \kappa $ in measure,
	\item $\{\sum_{i=1}^{N_m}A_n^{i,m}\}_n\acsArrow \{A_n\}_n$.
\end{itemize}
\end{defin}

\par
On the one hand, owing to the definition of GLT sequences, it is expected that 

\[
\{A_n\}_n\sim_{GLT}\kappa \:\Rightarrow \:\{A_n\}_n\sim_{\sigma}\kappa.
\]
On the other hand, sequences resulting from basic operations between GLT sequences also belong to the GLT class, with the symbol resulting from the same operations between the symbols of the initial sequences. The most important properties of the GLT sequences are summarised below.

\begin{description}
	\item[{\bf GLT1}] Every GLT sequence is related with a function $k$, $\kappa:[0,\:1]\times[-\pi,\:\pi]\rightarrow \mathbb{C}$, which is the symbol of the sequence. The singular values of the sequence are distributed as the $\kappa$ function. If the matrices of the sequence are Hermitian, the eigenvalues of the sequence are distributed as the $\kappa$.
	\item[{\bf GLT2}] The set of all GLT sequences is an *-algebra. That is, it is closed under linear combinations, multiplications, conjugate transpositions, and inversions, provided that the symbol of the sequence is zero at a set of zero measure. Therefore, a sequence obtained by operations between GLT sequences is GLT with a symbol produced by identical operations between the symbols.
	\item[{\bf GLT3}] Every Toeplitz sequence, with generating function $f \in L^1([-\pi,\pi])$ is GLT, with symbol $\kappa(x,\theta) = f(\theta)$.
	\item[{\bf GLT4}] Every diagonal matrix, whose elements are a uniform sampling of an almost everywhere continuous function $\alpha:[0,\:1]\rightarrow \mathbb{C}$ is GLT with symbol $\kappa(x,\theta) = \alpha(x)$.
	\item[{\bf GLT5}] Every zero-distributed sequence is GLT with symbol $\kappa(x,\theta)=0$.
	\item[{\bf GLT6}] $\{A_n\}_n \sim_{GLT} \kappa$, if and only if there exist GLT sequences $\{B_{n,m}\}_n \sim_{GLT} \kappa_m$ such that $\kappa_m$ converge to $\kappa$ in measure and $\{\{B_{n,m}\}_n\}_m$ is an a.c.s. for $\{A_n\}_n$.
\end{description}

%% file: SYMMETRIZED.tex
\chapter[Symmetrised Matrix Sequences]{Asymptotic spectra of large matrices coming from the symmetrisation of Toeplitz structure functions and applications to preconditioning}

The symmetrisation of a real, non-symmetric Toeplitz system was first proposed by Jennifer Pestana and Andrew Wathen \cite{doi:10.1137/140974213}. The symmetry is obtained by multiplying the system by $Y_n \: \in \: \mathbb{R}^{n \times n}$, where
\begin{align}
Y_n= \begin{bmatrix} \label{anti_identity_matrix}
     &    & 1 \\
     & \iddots   &   \\
1    &    & 
\end{bmatrix}.
\end{align}

The use of Krylov subspace methods in symmetric linear systems offers  significant advantages over the methods used in non-symmetric systems. For the conjugate gradient and minimum residual methods (and other similar methods), which are applied to symmetric positive definite and symmetric non-definite systems, respectively, it is known that the convergence depends on the eigenvalues of the coefficient matrix of the system. Thus, if the eigenvalue distribution of such a system is determined, it is theoretically guaranteed to converge within a number of iterations; however, there is no analogue result for methods applied to non-symmetric systems. In addition, for each iteration of the above methods, the cost is minimal, in the sense that only some table-vector multiplications are required.

\par
The singular value and eigenvalue distribution of symmetrised matrix sequences have been studied in detail in \cite{FFHMS,mazza-pestana,HMS}. The results are presented in the following sections.

\section{Eigenvalue Distribution of Symmetrised Toeplitz Matrix Sequences}

\begin{defin}
Let $g$ be a given function defined in $[0,\: 2\pi]$. We set $\psi_g$ in $[-2\pi,\: 2\pi]$ as follows:
    \begin{equation}\label{def-psi}
    \psi_g(\theta)=\left\{
    \begin{array}{cc}
    g(\theta), & \theta\in [0,2\pi], \\
    -g(\theta+2\pi), & \theta\in [-2\pi,0).
    \end{array}
    \right.
    \end{equation}
\end{defin}

\begin{thm}\label{thm:ferrari symmetrized spectral distribution 1}
Let $f \in L^1([-\pi,\pi])$ with real Fourier coefficients, $Y_n \in \mathbb{R}^{n \times n}$, as in \eqref{anti_identity_matrix}, and let $T_n(f)\in \mathbb{R}^{n \times n}$ be the Toeplitz matrix with generating  $f$. Then,
{}
    \[
    \{Y_nT_n(f)\}_n \sim_{\sigma}  f,
    \]
     
    \[
    \{Y_nT_n(f)\}_n \sim_{\lambda}  \psi_{|f|}.
    \]

\begin{prf*}
The proof of the first result is trivial. Because $Y_n$ is unitary, the singular values of $Y_nT_n(f)$ match those of $T_n(f)$ [Theorem \ref{thm:glt_toeplitz distribution}]. To prove the second result, we first assume that $n=2m$. Then,
\[
Y_nT_n(f)=\begin{bmatrix} 
Y_mH_m(f,+)Y_m   & Y_mT_m(f) \\
   Y_mT_m(f)     & H_m(f,-)  \\
\end{bmatrix}
\]
\[
=\begin{bmatrix} 
                 & Y_mT_m(f) \\
   Y_mT_m(f)     &           \\
\end{bmatrix}+\begin{bmatrix} 
Y_mH_m(f,+)Y_m   &           \\
                 & H_m(f,-)  \\
\end{bmatrix} = B_{2m} + Z_{2m},
\]
where $H_m(f,+)$ is the $m\times m$ Hankel matrix, which contains the Fourier coefficients of the function $f$, starting from $a_1$ at position $(1,1)$ to $a_{2m-1}$ at position $(m,m)$. Analogously, $H_m(f,-)$ is the $m\times m$ Hankel matrix, which contains the Fourier coefficients of the function $f$, starting from $a_{-1}$ at position $(1,1)$ to $a_{-2m+1}$ at position $(m,m)$. $H_m(f,+)$ is exactly the Hankel matrix, with generating $f$, as defined in \cite{Fasino2000}. For this matrix, it was proven that if $f$ is a Lebesgue integrable function, then $\{H_n(f,+)\}_n\sim_{\sigma}0$. Because $H_m(f,-)=H_m(\bar{f},+)$ and $\bar{f}$ is clearly Lebesgue integrable, $\{H_n(f,-)\}_n\sim_{\sigma}0$. In addition , the singular values of $Y_mH_m(f,+)Y_m$ match those of $H_m(f,+)$. The singular values of $Z_{2m}$, because it is a block diagonal, are the singular values of the two blocks.
Let $Y_mT_m(f)=U_m\Sigma_m V_m^*$ be a singular value decomposition of $Y_mT_m(f)$. Then,
\[
\frac{1}{\sqrt{2}}\begin{bmatrix} 
    U_m^*  & V_m^* \\
   -U_m^*  & V_m^* \\
\end{bmatrix}\begin{bmatrix} 
                 & Y_mT_m(f) \\
   Y_mT_m(f)     &           \\
\end{bmatrix}\frac{1}{\sqrt{2}}\begin{bmatrix} 
    U_m    & -U_m \\
    V_m    &  V_m \\
\end{bmatrix}=
\]
\[
\frac{1}{\sqrt{2}}\begin{bmatrix} 
    U_m^*  & V_m^* \\
   -U_m^*  & V_m^* \\
\end{bmatrix}\begin{bmatrix} 
                 & Y_mT_m(f) \\
   (Y_mT_m(f))^*     &           \\
\end{bmatrix}\frac{1}{\sqrt{2}}\begin{bmatrix} 
    U_m    & -U_m \\
    V_m    &  V_m \\
\end{bmatrix}=\begin{bmatrix} 
    \Sigma  &  \\
            & -\Sigma \\
\end{bmatrix}.
\]
Now, let $n=2m+1$. Then,
\[
Y_nT_n(f)=\begin{bmatrix}
Y_mH_m(\hat{f},+)Y_m   & Y_mx_+  & Y_mT_m(f) \\
     x_+^TY_m          & a_0     & x_-^T     \\
   Y_mT_m(f)           & x_-     & H_m(\hat{f},-)  \\
\end{bmatrix}=
\]
\[
\begin{bmatrix}
          &     & Y_mT_m(f) \\
          & 0   &           \\
Y_mT_m(f) &     &           \\
\end{bmatrix}+\begin{bmatrix}
Y_mH_m(\hat{f},+)Y_m   & Y_mx_+  &           \\
     x_+^TY_m          & a_0     & x_-^T     \\
                       & x_-     & H_m(\hat{f},-)  \\
\end{bmatrix}=
\]
\[
B_{2m+1}+Z_{2m+1},
\]
where
\[
x_+=[a_1, a_2 ,\dots,a_m]^T, x_-=[a_{-1},a_{-2},\dots,a_{-m}].
\]
The unitary matrix that diagonalizes $B_{2m+1}$ is
\[
\frac{1}{\sqrt{2}}\begin{bmatrix} 
    U_m    &          & -U_m \\
           & \sqrt{2} &      \\ 
    V_m    &          & V_m  \\
\end{bmatrix},
\]
while its eigenvalues are the same with those of $B_{2m}$ with the addition of 0. Clearly, $\{B_n\}_n\sim_{\lambda}\psi_{|f|}$. Furthermore,
\[
Z_{2m+1}=\begin{bmatrix}
Y_mH_m(\hat{f},+)Y_m   &       &           \\
                       & 0     &     \\
                       &       & H_m(\hat{f},-)  \\
\end{bmatrix}+\begin{bmatrix}
                       & Y_mx_+  &           \\
     x_+^TY_m          & a_0     & x_-^T     \\
                       & x_-     &            \\
\end{bmatrix}.
\]
$Y_mH_m(\hat{f},+)Y_m$ contains the Fourier coefficients of $\hat{f}$, starting from $\hat{a}_{1}$ at position $(1,1)$ to $\hat{a}_{2m-1}$ at position $(m,m)$; thus, we observe that $\hat{a}_i=a_{i+1}$ and $\hat{f}$ is Lebesgue integrable. Therefore, $\{Y_mH_m(\hat{f},+)Y_m\sim_{\sigma}0$ and, analogously, $\{H_m(\hat{f},-)\}_m$. The rank of the second matrix is 2; hence, applying the singular values interlacing theorem for the two matrices, we deduce the sequences $\{Z_{2m}\}_m$ and $\{Z_{2m+1}\}_m$; therefore, $\{Z_n\}_n$ are distributed at zero.
\newline
Then, based on Theorem \ref{thm:distr seq plus zero distributed} for symmetric matrices, $ \{Y_nT_n(f)\}_n \sim_{\lambda}\psi_{|f|}$. 
\end{prf*}
\end{thm}

\section[Symmetrised Functions of Toeplitz Sequences]{Eigenvalue distribution of large matrices produced by the symmetrisation of Toeplitz structure functions}

\subsection{Matrix functions}

Let the $h(z)$ function be analytic at $0$. Then, $h(z)$ is analytic in an open sphere of radius $r$ centred at $0$. Thus,
\[
h(z)=\sum_{k=0}^{\infty}b_kz^k, \quad z \: \in \: \{z \: \in \: \mathbb{C}, \: |z|<r \},
\]
where $b_k=\frac{f^{(k)}(0)}{k!}$. If $A \: \in \: \mathbb{C}^{n\times n}$ with $\|A\|<r$ for a natural matrix norm\footnote{The natural matrix norm is any norm derived from the rule $\|A\|=\sup\limits_{x \neq 0} \frac{\|Ax\|}{\|x\|}$, where $\|.\|$ is a vector norm}, the series $\sum_{k=1}^{\infty}b_kA^k$ converges, and the function $h(A)$ is well defined. If $\rho(A)$ is the spectral radius of $A$, then the condition $\rho(A)<r$ is necessary and sufficient for $\|A\|<r$ for some natural norm $\|.\|$.
\newline
Crucially for symmetrisation, Toeplitz matrices are persymmetric. That is,  
\[
Y_nT_n(f)=T_n(f)^TY_n.
\]
Otherwise, $Y_nT_n(f)$ is symmetric. If $T_n(f)$ is real, then  $Y_nT_n(f)$ is real symmetric and therefore normal. Furthermore, 
\[
Y_nT_n(f)^k = (Y_nT_n(f))T_n(f)^{k-1}=T_n(f)^TY_nT_n(f)^{k-1}=\dots=(T_n(f)^k)^TY_n;
\]
that is, $T_n(f)^k$ is also persymmetric. The above equation applies to all persymmetric matrices.
\begin{prop}
Let $h(z)$ be analytic to $z<r$. If $A \: \in \mathbb{R}^{n\times n}$ is persymmetric and $\rho(A)<r$, then $h(A)=\sum_{k=1}^{\infty}b_kA^k$ is persymmetric.
\end{prop}
In the following, we refer only to series with real coefficients, so that $Y_nh(A)$ is real symmetric.

\subsection{Basic Results}

\begin{prop}\label{lemma:poly_flipped_ftm}
Let $f \in L^\infty([-\pi,\pi])$ with real Fourier coefficients $Y_n \in \mathbb{R}^{n \times n}$, as defined in (\ref{anti_identity_matrix}), and let  $T_n(f)\in \mathbb{R}^{n \times n}$ be the Toeplitz matrix with generating $f$. Then, for every polynomial $p(z)$, we have
        \[
        \{p(T_n(f)\}_n \sim_{\sigma}  p\circ f.
        \]
    
\begin{prf*}
According to Item {\bf GLT3}, every Toeplitz sequence $\{T_n(f)\}_n$ is a GLT with symbol $f$. In addition, according to Item {\bf GLT2}, every sequence obtained with operations between GLT sequences is GLT with a symbol produced via the same operations as between the symbols of the sequences. Finally, according to {\bf GLT1}, the singular values of a GLT sequence with symbol $k$ are distributed as $k$. In the present case,  $\{p(T_n(f))\}_n \sim_{GLT} \tilde f= p\circ f$.
\end{prf*}
\end{prop}

\begin{thm} \label{thm:flipped_ditributions}
Let $f \in L^\infty([-\pi,\pi])$ with real Fourier coefficients $Y_n \in \mathbb{R}^{n \times n}$, as defined in (\ref{anti_identity_matrix}), and let $T_n(f)\in \mathbb{R}^{n \times n}$ be a Toeplitz matrix with generating  $f$. Let  $h(z)$ be an analytic function, with real coefficients and a radius of convergence $r$ such that $\|f\|_\infty<r$.
Then, we have the following distributions:
\begin{equation}\label{eqn:main_flipped_ftm_singular}
    \{h(T_n(f))\}_n \sim_{\sigma} h \circ f,
\end{equation}
and
\begin{equation}\label{eqn:main_flipped_ftm_eig}
\{Y_nh(T_n(f))\}_n \sim_{\lambda}  \psi_{|h\circ f|}.
\end{equation}
\begin{prf*}
The condition $\|f\|_\infty<r$ implies that $\|T_n(f)\|< r$, where $\|\cdot\|$ is the spectral norm; thus, $\rho(T_n(f))<r$. Therefore, matrix $h(T_n(f))$ is well defined.
        
If $|z|<r$, we use the Taylor series at $0$ for $h(z)$. That is, $h(z)=\sum_{k=0}^\infty {b}_k z^k$. For every $m \in \mathbb{N}$, we define the polynomial,
\[
p_m(z) = \sum_{k=0}^m  {b}_k z^k.
\]
Thus, the following properties apply:
\begin{enumerate}
\item $\{p_m(T_n(f))\}_n\sim_{\sigma} p_m \circ f$ for every $m \in \mathbb{N}$,
\item $\{p_m(T_n(f))\}_n \acsArrow \{h(T_n(f))\}_n$,
\item $p_m \circ f \rightarrow h \circ f$ in measure.
\end{enumerate}
The first property is a consequence of Proposition \ref{lemma:poly_flipped_ftm}. To prove the second property, we write  
\[
h(T_n(f) = p_m(T_n(f))+(h(T_n(f))-p_m(T_n(f))). 
\]
Then, we observe that $\|h(T_n(f))-p_m(T_n(f))\|< \epsilon_m$, where $\lim_{m\to\infty} \epsilon_m = 0$. This property follows by setting $R_{n,m}=0_n$ and $N_{n,m}=h(T_n(f))-p_m(T_n(f))$.
For the third property, because we assume that $\|f\|_\infty<r$, $h$ is analytic at the set $f(\theta)$ almost everywhere for $\theta \:\in\:[-\pi,\pi]$. Therefore, $p_m \circ f$ converges almost everywhere at $h \circ f$. In addition, because the domain is bounded, $p_m \circ f$ converges to $h \circ f$. 
Applying Theorem \ref{thm:acs sigma}, it immediately follows that $\{h(T_n(f))\}_n \sim_{\sigma} h \circ f$.
Item {\bf GLT6} implies that the matrix sequence $\{h(T_n(f))\}_n$ is a GLT with symbol $h \circ f$.
\newline
To prove (\ref{eqn:main_flipped_ftm_eig}), we set
\[
\Delta_n(h,f) = h(T_n(f)) - T_n(h\circ f).
\]
Because $h\circ f \: \in \: L^1[-\pi,\pi]$, we have that
\[
\{T_n(h\circ f)\}_n \sim_{\sigma} h\circ f, \quad \{T_n(h\circ f)\}_n \sim_{GLT} h\circ f.
\]
In addition, from the property {\bf GLT2}, we find that the matrix sequence $\{\Delta_n(h,f)\}_n$ is a GLT, since is the difference between the two GLT sequences. Its symbol is the difference between the initial sequences symbols. So, 
\[
\{\Delta_n(h,f)\}_n \sim_{GLT} 0, \quad \{\Delta_n(h,f)\}_n \sim_{\sigma} 0.
\]
Because $Y_n$ is unitary, we have
\[
\{Y_n\Delta_n(h,f)\}_n \sim_{\sigma} 0.
\]
That is, $\{Y_n\Delta_n(h,f)\}_n$ is zero distributed, according to the Definition \ref{def:zero distributed seq}. Then,
\begin{align*} 
\{h(T_n(f))\}_n &= \{T_n(h\circ f)\}_n+\{\Delta_n(h,f)\}_n \Rightarrow \\
\{Y_nh(T_n(f))\}_n &= \{Y_nT_n(h\circ f)\}_n+\{Y_n\Delta_n(h,f)\}_n,
\end{align*}
where $\{Y_nT_n(h\circ f)\}_n \sim_{\lambda} \psi_{|h\circ f|} $, as proved in Theorem \ref{thm:flipped_ditributions}, and $\{Y_n\Delta_n(h,f)\}_n$ is zero distributed. Thus, by applying Theorem \ref{thm:distr seq plus zero distributed}, it immediately follows that
\[
\{Y_nh(T_n(f))\}_n \sim_{\lambda} \psi_{|h\circ f|}.
\]
\end{prf*}
\end{thm}

\section{Numerical Results}

In this section, the eigenvalue and singular value distributions of the sequence $\{Y_nh(T_n(f))\}_n$ are numerically investigated. Considering this distribution, several circulant preconditioners are proposed, and the spectrum of the preconditioned sequence is also investigated.

\subsection{ Spectrum of the $\{Y_nh(T_n(f))\}_n$}
The results presented in this section confirm the claims of Theorem  \ref{thm:flipped_ditributions}. Specifically, for $f$ trigonometric polynomial and $h$, either analytic function either polynomial is confirmed in the following four examples, that the distribution of the eigenvalues of the sequence $ \{Y_nh(T_n(f))\}_n$ is described by $\psi_{|h \circ f|}$. It is also  numerically confirmed that the singular value distribution of $\{h(T_n(f))\}_n$ is described by $|h\circ f|$.

\begin{figure}[!h]
    \centering
    \includegraphics[width=0.8\textwidth]{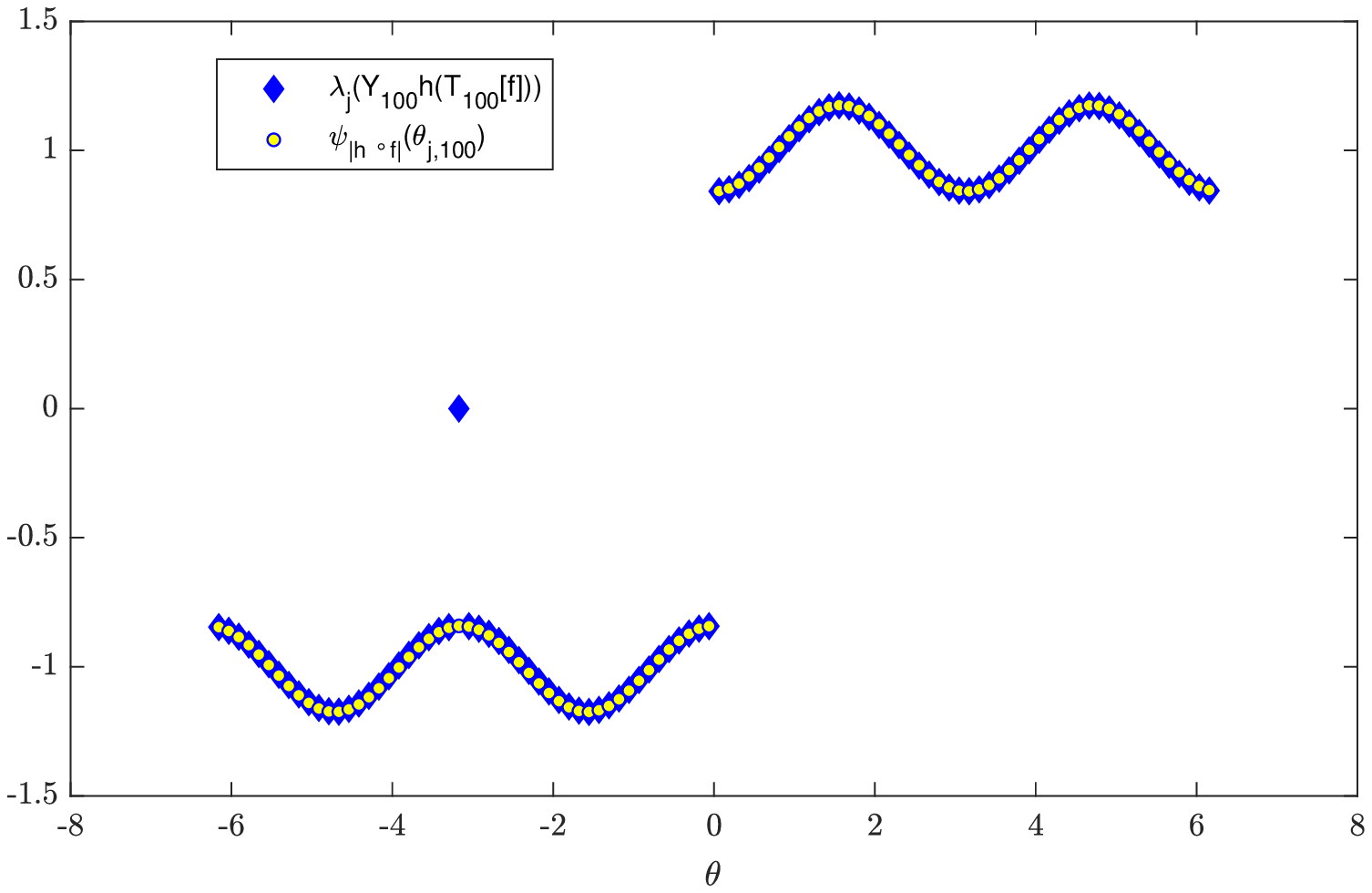}
    \caption{Comparison between the eigenvalues of the symmetrised matrix $Y_{100}h(T_{100}(f))$ and the uniform sampling of $\psi_{|h\circ f|}$ over $[-2\pi,2\pi]$ for $h(z)=\sin(z)$ and  $f(\theta)=e^{{\mathbf{i}}\theta}$.}
    \label{fig:symmetrized h_sin_f_exp}
\end{figure}

\begin{exl}\label{example:symmetrized 1}
In this example the analytic function $h(z)=\sin(z)$ whose Taylor series converges throughout the complex plane, and the trigonometric polynomial $f(\theta)=e^{{\mathbf{i}}\theta}$ are considered. Figure \ref{fig:symmetrized h_sin_f_exp} shows that for $n=100$, the eigenvalues of $Y_nh(T_n(f))$ are well approximated by the uniform sampling of $\psi_{|h\circ f|}$ over $[-2\pi,2\pi]$, except for the presence of an outlier. This indicates that Definition \ref{defin:glt eig distr} does not rule out the existence of such eigenvalues.
\end{exl}

\begin{figure}[!h]
    \centering
    \includegraphics[width=0.8\textwidth]{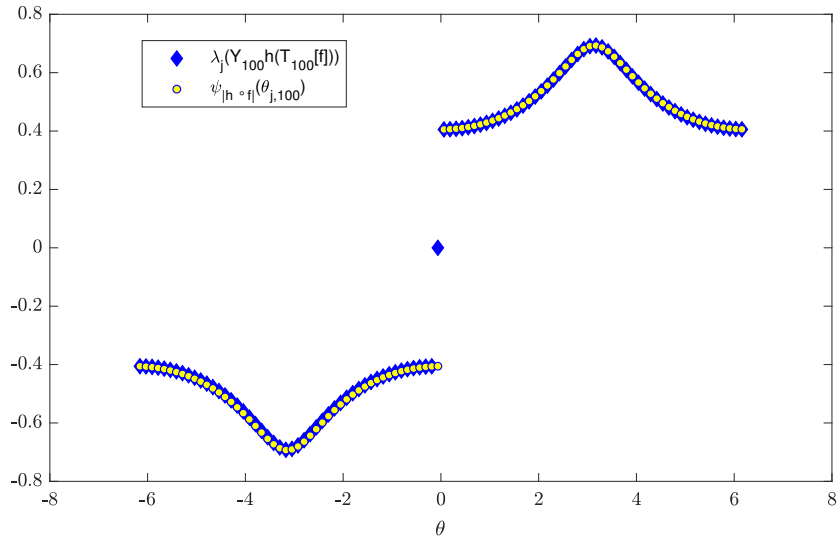}
    \caption{Comparison between the eigenvalues of the symmetrised matrix $Y_{100}h(T_{100}(f))$ and the uniform sampling of $\psi_{|h\circ f|}$ over $[-2\pi,2\pi]$ for $h(z)=\log(1+z)$ and  $f(\theta)=0.5e^{{\mathbf{i}}\theta}$. }
    \label{fig:symmetrized h_log_f_exp}
 \end{figure}

\begin{exl}\label{example:symmetrized 2}
In the second example, for the analytic function $h(z)=\log(1+z)$, whose Taylor series at 0 converges with a radius of convergence equal to 1 is used the trigonometric polynomial $f(\theta)=0.5e^{{\mathbf{i}}\theta}$ with  $\|f\|_\infty<1$, as Theorem \ref{thm:flipped_ditributions} demands.
      Figure \ref{fig:symmetrized h_log_f_exp}, shows that except one outlier, the eigenvalues of $Y_nh(T_n(f))$ for $n=100$ are well approximated by a uniform sampling of $\psi_{|h\circ f|}$ over $[-2\pi,2\pi]$. 
\end{exl}

\begin{figure}[!h]
    \centering
    \includegraphics[width=0.8\textwidth]{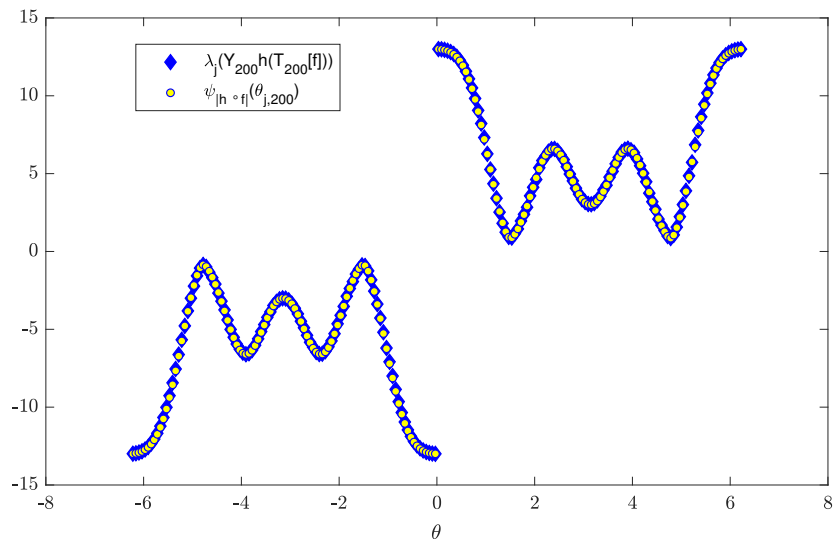}
    \caption{Comparison between the eigenvalues of the symmetrised matrix $Y_{200}h(T_{200}(f))$ and the uniform sampling of $\psi_{|h\circ f|}$ over $[-2\pi,2\pi]$ for $h(z)=1+z+z^2$ and $f(\theta)=-e^{{\mathbf{i}}\theta}+1+e^{-{\mathbf{i}}\theta}+e^{-{\mathbf{i}}2\theta}+e^{-{\mathbf{i}}3\theta}$. }
    \label{fig:symmetrized hon_example_1}
    \end{figure}

\begin{figure}[!h]
    \centering
    \includegraphics[width=0.8\textwidth]{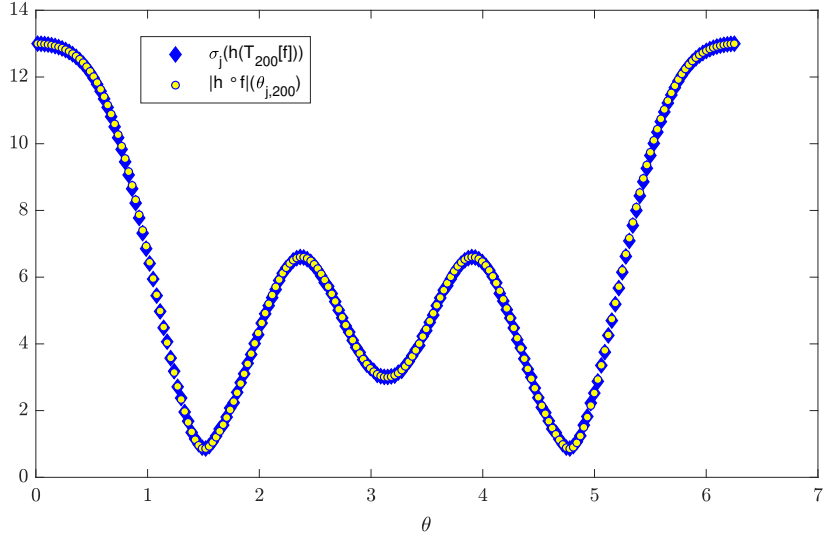}
    \caption{Comparison between the singular values of the matrix $h(T_{200}(f))$ and the uniform sampling of ${|h\circ f|}$ over $[0,2\pi]$ for $h(z)=1+z+z^2$ and $f(\theta)=-e^{{\mathbf{i}}\theta}+1+e^{-{\mathbf{i}}\theta}+e^{-{\mathbf{i}}2\theta}+e^{-{\mathbf{i}}3\theta}$. }
    \label{fig:symmetrized hon_example_1_sigma}
    \end{figure}

 \begin{exl}\label{example:symmetrized 3}
This example was taken from \cite{Hon2018}. Following the same procedure as Examples 1--2 , Figure \ref{fig:symmetrized hon_example_1} shows the spectrum of $Y_nh(T_n(f))$ for $n=200$; the function  $h(z)=1+z+z^2$, whose Taylor series in 0 converges in the whole complex plane; and the trigonometric polynomial $f(\theta)=-e^{{\mathbf{i}}\theta}+1+e^{-{\mathbf{i}}\theta}+e^{-{\mathbf{i}}2\theta}+e^{-{\mathbf{i}}3\theta}$. In the present example, there are no outliers, and the eigenvalues of $Y_{n}h(T_n(f))$ are approximated by the uniform sampling of $\psi_{|h\circ f|}$ over $[-2\pi,2\pi]$. {Moreover, to  numerically confirm the relation \eqref{eqn:main_flipped_ftm_singular} of Theorem \ref{thm:flipped_ditributions}, it is verified that the singular values of the matrix $h(T_n(f))$ can be approximated by a uniform sampling of ${|h\circ f|}$ over $[0,2\pi]$. Indeed, Figure \ref{fig:symmetrized hon_example_1_sigma} shows that the expected approximation holds true already for a moderate size such as $n=200$.}
\end{exl}

\begin{figure}[!h]
    \centering
    \includegraphics[width=0.8\textwidth]{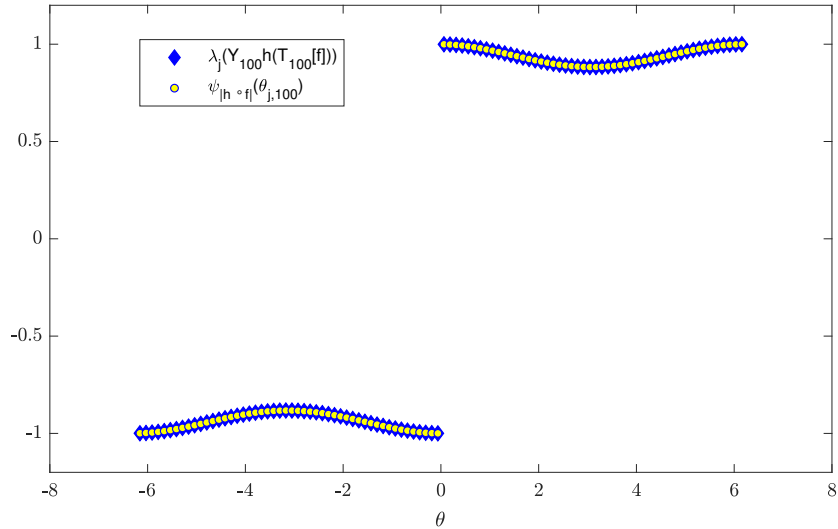}
    \caption{ {Comparison between the eigenvalues of the symmetrised matrix $Y_{100}h(T_{100}(f))$ and the uniform sampling of $\psi_{|h\circ f|}$ over $[-2\pi,2\pi]$ for $h(z)={\rm e}^z$ and $f(\theta)=\sum_{j=-99}^{99}a_{j}e^{\mathbf{i}j\theta}$, with $\lambda = 0.1$, $\mu = -0.9$, $\nu = 0.25$, $\sigma = 0.45$, $r = 0.05$, and $\Delta x=\frac{4}{101}$. }}
    \label{fig:symmetrized hon_example_fin}
\end{figure}

\begin{figure}[!h]
    \centering
    \includegraphics[width=0.8\textwidth]{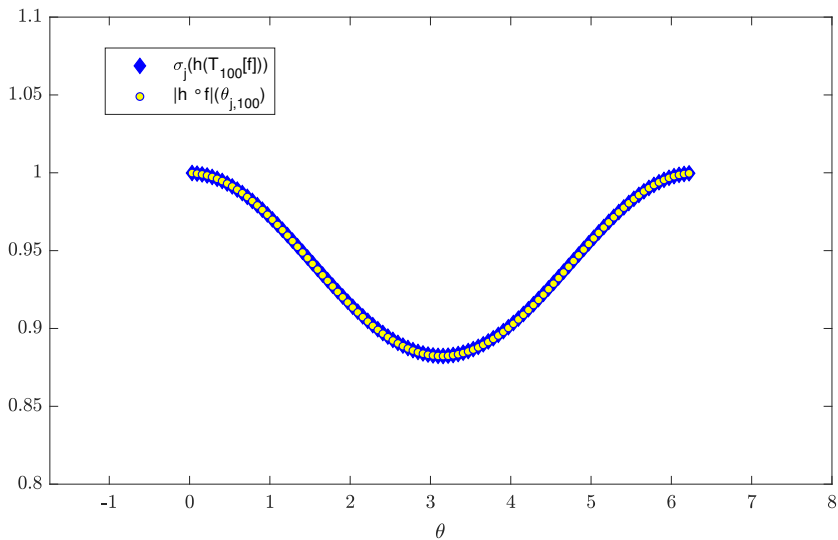}
    \caption{ {Comparison between the singular values of the matrix $h(T_{100}(f))$ and the uniform sampling of $|h\circ f|$ over $[0,2\pi]$  for $h(z)={\rm e}^z$ and $f(\theta)=\sum_{j=-99}^{99}a_{j}e^{\mathbf{i}j\theta}$, with $\lambda = 0.1$, $\mu = -0.9$, $\nu = 0.25$, $\sigma = 0.45$, $r = 0.05$, and $\Delta x=\frac{4}{101}$. }}
    \label{fig:symmetrized hon_example_fin_sigma}
\end{figure}

\begin{exl}\label{example:symmetrized non_prec4}
{ The last example is a practical case taken from \cite{MR2609339, MR2872597}. Here, we have the case of the exponential of a real non-symmetric Toeplitz matrix derived from computational finance (more specifically, from the option pricing framework in jump-diffusion models), where a partial integro-differential equation (PIDE) must be solved. The discretisation of a PIDE can be transformed into a matrix exponential problem which is equivalent to considering the analytic function $h(z)={\rm e}^z$, whose Taylor series centred at $0$ converges in the whole complex plane, as well as a trigonometric polynomial $f(\theta)=\sum_{j=-n+1}^{n-1}a_{j}e^{\mathbf{i}j\theta}$ defined by the following Fourier coefficients:}     

\begin{align}
  a_0&=-{\nu^2}-\Delta x^2(r+\lambda-\lambda w(0)\Delta x);\\\label{eqn:coeff_f0}
  a_1&=\frac{\nu^2}{2}-\Delta x\frac{(2r-2\lambda k-\nu^2)}{4}+\lambda w(-\Delta x)\Delta x^3; \\
  a_{-1}&=\frac{\nu^2}{2}+\Delta x\frac{(2r-2\lambda k-\nu^2)}{4}+\lambda w(\Delta x)\Delta x^3; \\
  a_j&=  \lambda \Delta x^3 w(-j \Delta x),\quad j \in \{-n+1,\dots,-2,\} \cup \{2,\dots, n-1\}. \label{eqn:coeff_fj}
\end{align}

Here, $w(s)=\frac{{\rm e}^{-\frac{(s-\mu)^2}{2\sigma^2}}}{\sqrt{2\pi}\sigma}$ is a normal distribution function with mean $\mu$ and standard deviation $\sigma$, the parameter $k= {\rm e}^{\mu+\frac{\sigma^2}{2}}-1$ is the expectation
of the impulse function, $\Delta x$ is the spatial step size, $\nu$ is the stock return volatility, $r$ is the risk-free interest rate, and $\lambda$ is the arrival intensity of the Poisson process.

{
   Following the same procedure as in Examples 1--3, we plot  in Figure \ref{fig:symmetrized hon_example_fin} the spectrum of $Y_nh(T_n(f))$ for $n=100$. In the present example, we observe that there are no outliers, and the eigenvalues of $Y_{n}h(T_n(f))$ are well  approximated by the uniform sampling of $\psi_{|h\circ f|}$ over $[-2\pi,2\pi]$.}

{   
   In addition, to numerically validate the relation \eqref{eqn:main_flipped_ftm_singular} presented in Figure \ref{fig:symmetrized hon_example_fin_sigma} for $n=100$, we compare the singular values of $h(T_n(f))$ and a uniform sampling of ${|h\circ f|}$ over $[0,2\pi]$. }  

 \end{exl}

\subsection{Circulant Preconditioners for the Symmetrised Toeplitz Sequence}

In the present section, preconditioners for the symmetrised system $Y_{n}h(T_n(f))$ are proposed, and the distribution of the preconditioned sequence is numerically investigated.

{
For the construction of the preconditioners, the approach proposed in \cite{Hon2018} is applied; however, taking into consideration the theoretical results proved here, another circulant preconditioner is also proposed. For the second preconditioner, a theoretical description of the preconditioned sequence distribution is provided. The results are presented in Examples \ref{example:symmetrized prec1}, \ref{example:symmetrized prec2}, and \ref{example:symmetrized prec3}. As mentioned in the introduction, the preferred Krylov method for symmetric, non-definite systems is MINRES. This method has the advantage that the cost per iteration is minimal because only matrix--vector multiplications are required. In addition, if the eigenvalues of the preconditioned coefficient matrix are known, the number of iterations required for convergence with given accuracy is known. The preconditioner must be symmetric and positive-definite. All of the preconditioners proposed here are symmetric and positive definite.
} 

\begin{defin}\cite{doi:10.1137/140974213}\label{circ-modulus}
  For every circulant matrix $C_n \in \mathbb{C}^{n\times n}$, the absolute value circulant matrix $|C_n|$ of $C_n$ is defined as
  \begin{eqnarray}\nonumber
    |C_n|&=& (C_n^* C_n)^{1/2}\\\nonumber
    &=&(C_n C_n^*)^{1/2}\\\nonumber
    &=&F_n|\Lambda_n|F_n^*,
  \end{eqnarray} 

where $F_n$ is defined as in \eqref{defin:unitary discrete Fourier } and $|\Lambda_n|$ is the diagonal matrix of size $n$, whose elements are the absolute values of the eigenvalues of $C_n$.  
\end{defin}

\begin{defin}\label{circ-frob}
The optimal Frobenius preconditioner for a Toeplitz matrix is the circulant $C_n$, defined as

\[ 
c(T_n(f)) = {\rm arg\,min}_{C_n\, }\|T_n(f)-C_n\|_F={\rm arg\,min}_{C_n=F_n\Lambda_n{F_n}^*\,}\|F_n^*T_n(f)F_n-\Lambda_n\|_F,
\]
where $\Lambda_n$ is a diagonal matrix that contains the eigenvalues of  $c(T_n(f))$. It is clear that $\Lambda_n=diag(F_n^*T_n(f)F_n)$. 

\par
Let $\mathbf{c}=[c_0,\:c_1,\:\dots,c_{n-1}]^T$ and $C_{\mathbf{c}}$ be the circulant matrix whose first column is $\mathbf{c}$. To explicitly derive the elements of $c(T_n(f))$, we define $F(\mathbf{c})=\|T_n(f)-C_\mathbf{c}\|_F^2$. We observe that for $k=0,\dots,n-1$, the element $c_k$ of $\mathbf{c}$ appears in two diagonals, in which the $a_k$ and $a_{k-n}$ elements of $T_n(f)$ are located. The first of the two diagonals has $n-k$ elements, and the second has $k$. Thus, we have
\[
F(\mathbf{c})=\|T_n(f)-C_\mathbf{c}\|_F^2= \sum_{k=0}^{n-1}(n-k)(a_k-c_k)^2+k(a_{k-n}-c_k)^2,
\]
where $\mathbf{c}=[c_0,\:c_1,\:\dots,c_{n-1}]^T$ must be determined. To satisfy the necessary conditions for the minimisation of $F$, we require that

\[
\frac{\partial F}{\partial c_k}=0 \Rightarrow -2(n-k)(a_k-c_k)-2k(a_{k-n}-c_k)=0 \Rightarrow c_k =\frac{(n-k)a_k+ka_{k-n}}{n}.
\]

Therefore, the elements $c_k, \; k=0,\dots,n-1$ of the first column of $c(T_n(f))$ are given by $c_k =\frac{(n-k)a_k+ka_{k-n}}{n}$.
\begin{rmrk}
The $\|.\|_F$ norm is produced by a positive inner product which makes $\mathbb{C}^{nxn}$, the space of complex matrices of size n (equipped with $\|.\|_F$ norm) a Hilbert space. The set of circulant matrices of size $n$ constitutes a non-empty, closed, and convex linear subspace of $\mathbb{C}^{nxn}$. Therefore, for every $A\in \mathbb{M}_n(\mathbb{C})$, there exists a unique circulant $c(A)$ [\cite{Rynne2008} Theorem 3.32], such that
\[
\|A-c(A)\|_F = \inf \{\|A-C\|_F : \text{ C is circulant}\}.
\]
\end{rmrk}
\end{defin}
For more properties regarding $c(A)$ see \cite{Skoro}.

As mentioned in Definition \ref{circ-frob}, the diagonal matrix $\Lambda_n$, whose elements are the eigenvalues of the optimal Frobenius preconditioner for $T_n(f)$, is the main diagonal of $F_n^*T_n(f)F_n$. The $j-th$ element in the diagonal is $f_{n,j}^*T_n(f)f_{n,j},\:j=0,\dots,n-1$, where $f_{n,j}$ denotes the $j-th$ column of $F_n$. Therefore, for $j=0,\dots,n-1$, we have that
\begin{align*}
&\lambda_j(c(T_n(f)))= f_{n,j}^*T_n(f)f_{n,j} \\
&= \frac{1}{n} [  w_n^{0j}    ( a_0w_n^{-0j}     + a_{-1}w_n^{-1j}  + a_{-2}w_n^{-2j}  + \dots + a_{-(n-1)}w_n^{-(n-1)j} )   \\
&+               w_n^{1j}    ( a_1w_n^{-0j}     + a_{0}w_n^{-1j}   + a_{-1}w_n^{-2j}  + \dots + a_{-(n-2)}w_n^{-(n-1)j} ) \\
&                \vdots \\
&+               w_n^{(n-1)j}( a_{n-1}w_n^{-0j} + a_{n-2}w_n^{-1j} + a_{n-3}w_n^{-2j}  + \dots + a_{0}w_n^{-(n-1)j} ] \\
&= \frac{1}{n} [ ( a_0w_n^{0j}     + a_{-1}w_n^{-1j}  + a_{-2}w_n^{-2j}  + \dots + a_{-(n-1)}w_n^{-(n-1)j} )   \\
&+               ( a_1w_n^{1j}     + a_{0}w_n^{0j}   + a_{-1}w_n^{-1j}  + \dots + a_{-(n-2)}w_n^{-(n-2)j} ) \\
&                \vdots \\
&+               ( a_{n-1}w_n^{(n-1)j} + a_{n-2}w_n^{(n-2)j} + a_{n-3}w_n^{(n-3)j}  + \dots + a_{0}w_n^{0j} ]. 
\end{align*}
Finally, because $w_n^j=e^{\frac{i2\pi j}{n}}$, the above summation equates to 
\begin{align}
\frac{1}{n}\left[\sum_{k=-n+1}^0a_ke^{\frac{i2\pi jk}{n}}+\sum_{k=-n+2}^1a_ke^{\frac{i2\pi jk}{n}}+\dots+\sum_{k=-n+l+1}^la_ke^{\frac{i2\pi jk}{n}}+\dots+\sum_{k=0}^{n-1}a_ke^{\frac{i2\pi jk}{n}}\right]. \label{lambda_j_of_cf}
\end{align}

Let us suppose that $f(\theta)$ is $l_1$ summable; that is, $\sum_{k=-\infty}^{\infty}|a_k|<\infty$. Then, the partial sum $\sum_{k=-m}^{m}a_ke^{ik\theta}$ uniformly converges to $f(\theta)$. We set $\epsilon_m=\sum_{|k|> m}|a_k|\leq \|f(\theta)-\sum_{k=-m}^m a_ke^{ik\theta}\|_{\infty}$. The partial sum uniformly converges to $f$; thus, for every $\epsilon>0$, there exists an $m$ such that $\epsilon_m < \epsilon$. We observe that for every $m$, the sum \eqref{lambda_j_of_cf} contains $n-2m$ terms, for which the summation starts from a term of order $-k\:\in\:\{-n+m+1,\dots,-m\}$ and ends with a term of order $n-k-1\:\in\:\{m,\dots,n-m-1\}$. All these terms are at most $\epsilon_m$ far from the exact value of $f(\frac{2\pi j}{n})$. We combined the remaining terms in pairs, to take two new terms. One of the orders was higher than $m$, and the other was lower than $m$. For example,
\[
\sum_{k=-1}^{n-2}a_ke^{\frac{i2\pi jk}{n}}+\sum_{k=-n+2}^{1}a_ke^{\frac{i2\pi jk}{n}}=\sum_{k=-n+2}^{n-2}a_ke^{\frac{i2\pi jk}{n}}+\sum_{k=-1}^{1}a_ke^{\frac{i2\pi jk}{n}}.
\]
Hence, we have
\[
|\lambda_j(c(T_n(f)))-f(\frac{2\pi j}{n})|<\frac{\sum_{k=0}^{m-1}\epsilon_k+(n-m)\epsilon_m}{n}.
\]

\begin{prop}\label{prop_optimal_circ_distr} \cite{Estatico2008,Capizzano2000}
Let $\{c(T_n(f))\}_n$ be the sequence of optimal Frobenius preconditioners of the sequence $\{T_n(f)\}_n$. Then, $\{c(T_n(f))\}_n\sim_{\sigma,\lambda,GLT}f$.
\end{prop}

\begin{prop}\label{prop_preconditioners_ditrsi}
Let $f \in L^\infty([-\pi,\pi])$ $l_1$ summable, with real Fourier coefficients. Let $h(z)$ be an analytic function with real coefficients and a radius of convergence $r$ such that $\|f\|_\infty<r$. Then, the circulant matrices $c(T_n(h\circ f))$ and $h(c(T_n(f)))$ are real, and we have that
\begin{align}
\{c(T_n(h\circ f))\}_n\sim_{GLT,\sigma,\lambda}h\circ f, \quad \{h(c(T_n(f)))\}_n\sim_{GLT,\sigma,\lambda}h\circ f. \label{circ_prec_distr}
\end{align}
In addition, the circulant $|c(T_n(h\circ f))|$, $|c(T_n(h\circ f))|^{-1}$, $|h(c(T_n(f)))|$, $|h(c(T_n(f)))|^-{1}$ is real and symmetric, and 
\begin{align}
\{|c(T_n(h\circ f))|\}_n\sim_{GLT,\sigma,\lambda}|h\circ f|,& \quad \{|h(c(T_n(f)))|\}_n\sim_{GLT,\sigma,\lambda}|h\circ f|, \label{abs_circ_prec_distr}\\
\{|c(T_n(h\circ f))|^{-1}\}_n\sim_{GLT,\sigma,\lambda}|h\circ f|^{-1},& \quad \{|h(c(T_n(f)))|^{-1}\}_n\sim_{GLT,\sigma,\lambda}|h\circ f|^{-1}. \label{inv_abs_circ_prec_distr}{}
\end{align}
\begin{prf*}
Under these assumptions, the function $h\circ f$ has real Fourier coefficients and belongs to $L^\infty([-\pi,\pi])\subset L^1([-\pi,\pi])$. The matrix $c(T_n(h\circ f))$ is real for every $n$, because each of its elements represents the weighted average of certain elements of $T_n(h\circ f)$, which is real. The matrix $h(c(T_n(f)))$ is real for every $n$ because $c(T_n(f))$ is real and $h$ has real coefficients. The first distribution at \eqref{circ_prec_distr} is an immediate consequence of the implementation of Proposition \ref{prop_optimal_circ_distr} for $h\circ f$. The second distribution is the consequence of the implementation of the same proposition for $f$, in combination with property {\bf GLT2}. In \cite{HMS}, it was proven that if $C$ is a real circulant, $|C|=(C^*C)^{1/2}$ is real and symmetric. Finally, the distributions at \eqref{abs_circ_prec_distr} and \eqref{inv_abs_circ_prec_distr} arise from the definition of the matrices and the application of the property {\bf GLT2}.
\end{prf*}
\end{prop}

{
The preconditioned matrices $|h(c(T_n(f)))|^{-1}Y_nh(T_n(h\circ f))$ and $|c(T_n(h\circ f))|^{-1}Y_nh(T_n(h\circ f))$ are similar to  real symmetric matrices, and therefore their eigenvalues are real. In fact, because $|h(c(T_n(f)))|^{-1}$ is real, symmetric, and positive definite, it can be written in the form $|h(c(T_n(f)))|^{-1}=Q|\Lambda|Q^T$, where $Q$ is real orthogonal and $|\Lambda|$ is diagonal with no negative elements. Hence,
\[
|h(c(T_n(f)))|^{-1}Y_nh(T_n(h\circ f))=Q|\Lambda|Q^TY_nh(T_n(h\circ f))\sim|\Lambda|^{1/2}QY_nh(T_n(h\circ f))Q^T|\Lambda|^{1/2},
\]
which is symmetric because $Y_nh(T_n(h\circ f))$ is symmetric. Analogously, we apply $|c(T_n(h\circ f))|^{-1}$.
}
The proposition that follows is a restatement (in matrix sequence terms) of Conclusion $3$ at \cite{Hon2018}.
\begin{prop}
\cite{Hon2018} Let $f \in L^\infty([-\pi,\pi])$ and $l_1$ be summable, with real Fourier coefficients. Let $h(z)$ be an analytic function with real coefficients and a radius of convergence  $r$ such that $\|f\|_\infty<r$. If $c(T_n(f))$ is the optimal Frobenius preconditioner for $T_n(f)$, then for the preconditioned sequence $\{|h(c(T_n(f)))|^{-1}Y_nh(T_n(h\circ f))\}_n$, we have that 
\[
\{|h(c(T_n(f)))|^{-1}Y_nh(T_n(h\circ f))\}_n= \{Q_n\}_n + \{Z_n\}_n,
\]
where the matrices of $\{Q_n\}_n$ are real orthogonal, and $\{Z_n\}_n$ is zero-distributed.
\end{prop}
A direct consequence of the proposition above and Theorem \ref{thm:distr seq plus zero distributed} is that the eigenvalues of the    preconditioned sequence are distributed similarly to those of $ \{Q_n \}_n$. The eigenvalues of a real orthogonal matrix can only be $1$ or $-1$. Finally, according to Theorem  \ref{distributed_as_f--clustered_at_erf}, the eigenvalues of the preconditioned sequence are clustered at $\{-1,1\}$.
{
  In the proposition that follows, we prove a similar conclusion for the preconditioned sequence $\{|c(T_n(h\circ f))|^{-1}Y_nh(T_n(h\circ f))\}_n$.
}

\begin{prop}
Let $f \in L^\infty([-\pi,\pi])$ and $l_1$ be summable, with real Fourier coefficients. Let $h(z)$ be an analytic function with real coefficients and a radius of convergence  $r$ such that $\|f\|_\infty<r$. If $c(T_n(h\circ f))$ is the optimal Frobenius preconditioner for $T_n(h\circ f)$, then for the preconditioned matrix sequence $\{|c(T_n(h\circ f))|^{-1}Y_nh(T_n(h\circ f))\}_n$, we have that
\[
\{|c(T_n(h\circ f))|^{-1}Y_nh(T_n(h\circ f))\}_n= \{Q_n\}_n + \{\hat{Z}_n\}_n,
\]
where the matrices of $\{Q_n\}_n$ are real orthogonal, and $\{\hat{Z}_n\}_n$ is zero-distributed.
\begin{prf*}
In Proposition \ref{prop_preconditioners_ditrsi}, it was proved that
\[
\{|c(T_n(h\circ f))|^{-1}\}_n\sim_{GLT}|h\circ f|^{-1}, \quad \{|h(c(T_n(f)))|^{-1}\}_n\sim_{GLT}|h\circ f|^{-1}.
\]
 In addition, as Theorem \ref{thm:flipped_ditributions} states, $h(T_n(f))\sim_{GLT}h\circ f$. By applying the property {\bf GLT2}, we have
\[
\{|c(T_n(h\circ f))|^{-1}h(T_n(h\circ f))\}_n\sim_{GLT}\frac{h\circ f}{|h\circ f|}, \quad \{|h(c(T_n( f)))|^{-1}h(T_n(h\circ f))\}_n\sim_{GLT}\frac{h\circ f}{|h\circ f|}.
\]
So, applying once more the property {\bf GLT2}, we have
\[
\{\Delta_n\}_n= \{|c(T_n(h\circ f))|^{-1}h(T_n(h\circ f))\}_n-\{|h(c(T_n( f)))|^{-1}h(T_n(h\circ f))\}_n \sim_{GLT}0.
\]
Therefore,
\begin{align*}
&\{|c(T_n(h\circ f))|^{-1}h(T_n(h\circ f))\}_n = \{|h(c(T_n( f)))|^{-1}h(T_n(h\circ f))\}_n + \{\Delta_n\}_n \Rightarrow \\
&\{Y_n|c(T_n(h\circ f))|^{-1}h(T_n(h\circ f))\}_n =\{Y_n|h(c(T_n( f)))|^{-1}h(T_n(h\circ f))\}_n + \{Y_n\Delta_n\}_n \Rightarrow \\
&\{|c(T_n(h\circ f))|^{-1}Y_nh(T_n(h\circ f))\}_n = \{|h(c(T_n( f)))|^{-1}Y_nh(T_n(h\circ f))\}_n + \{Y_n\Delta_n\}_n \\
&= \{Q_n\}_n + \{Z_n\}_n + \{Y_n\Delta_n\}_n = \{Q_n\}_n + \{\hat{Z}_n\}_n,
\end{align*}
where $\{\hat{Z}_n\}_n$ is zero-distributed, as the sum of two zero-distributed sequences. The permutation of $Y_n$ with the circulants can be achieved because they are real and Toeplitz.
\end{prf*}{}
\end{prop}

\begin{figure}[!h]
    \centering
    \includegraphics[width=.48\textwidth]{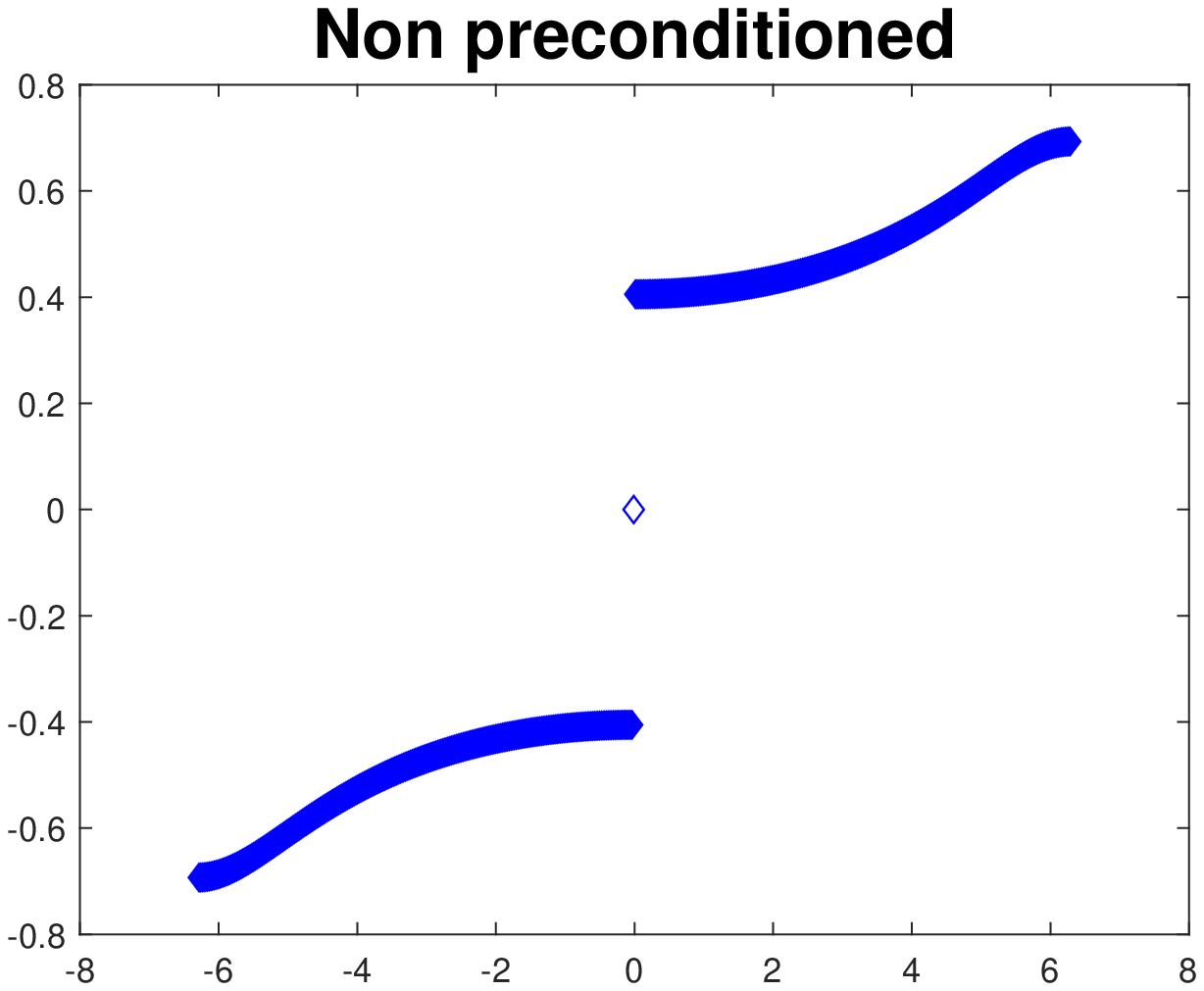}\\
    \includegraphics[width=.48\textwidth]{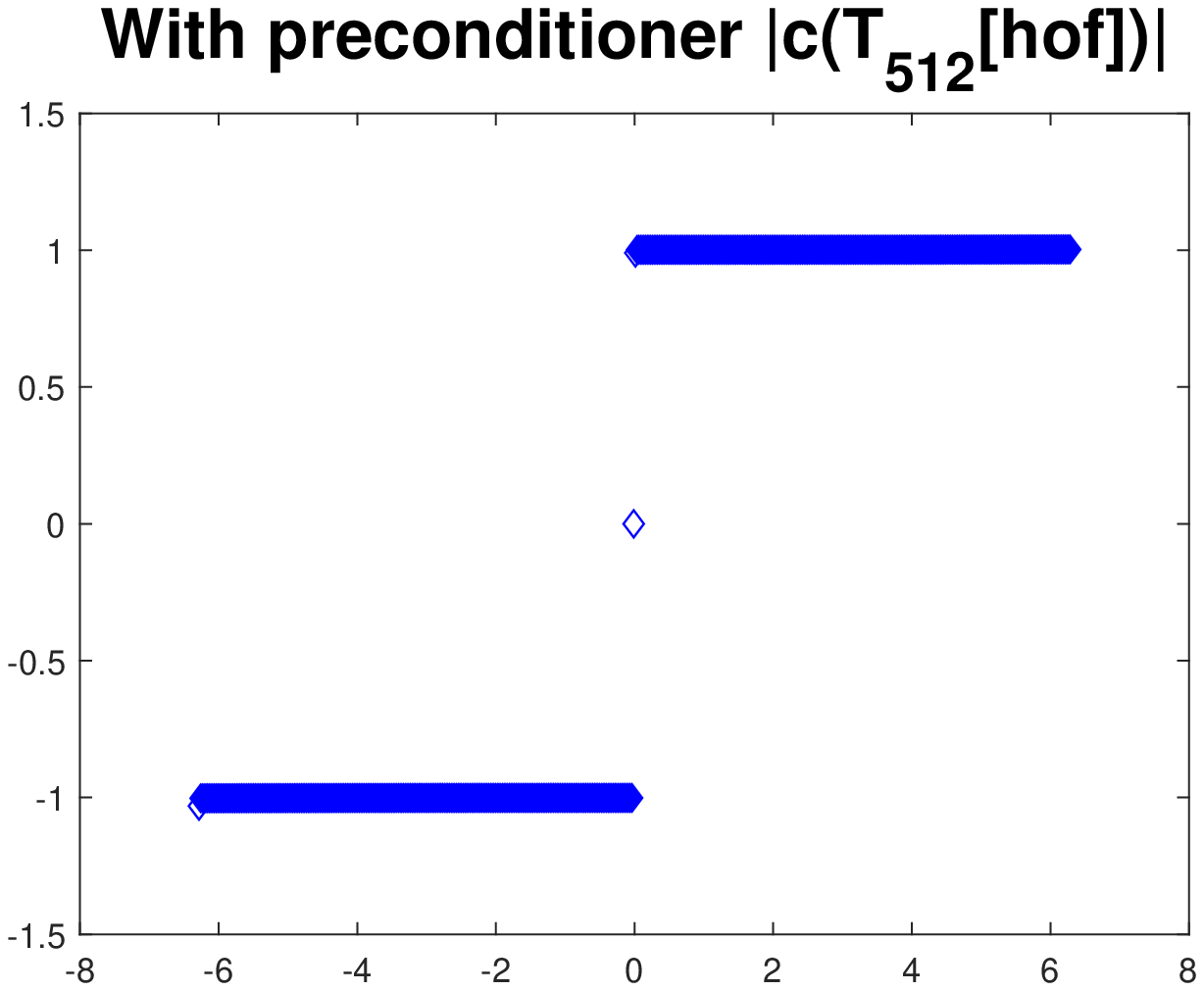}
     \includegraphics[width=.48\textwidth]{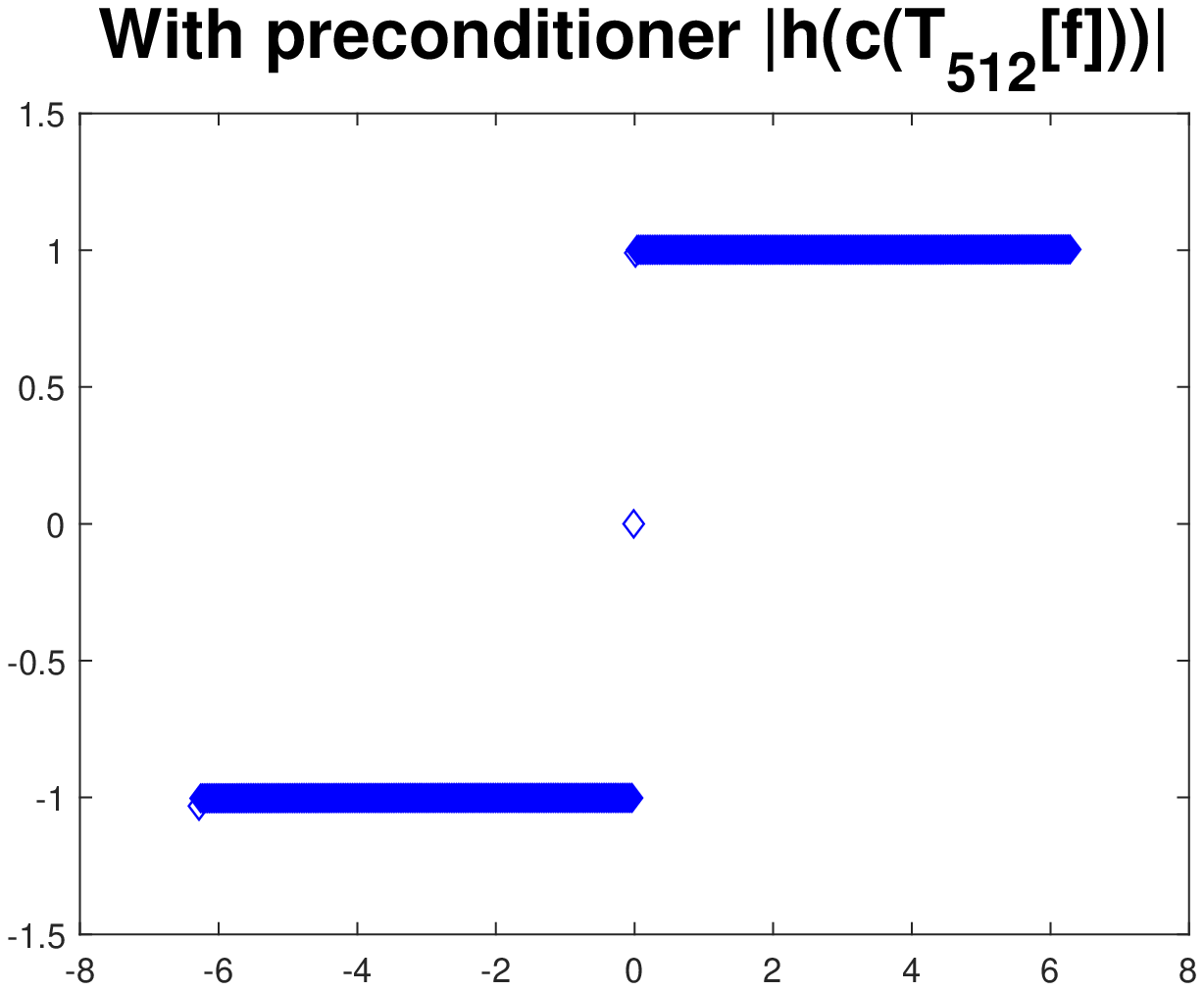}
    \caption{{Spectrum of the symmetrised matrix $Y_{512}h(T_{512}[f])$, for  $h(z)=\log(1+z)$ and $f(\theta)=0.5e^{{\mathbf{i}}\theta}$. Top: without preconditioner; bottom left: preconditioner $P_n=|c(T_n(h\circ f))|$; bottom right: preconditioner $P_n=|h(c(T_n(f)))|$.}}
    \label{fig:prec_h_log_f_exp}
\end{figure}

\begin{exl}\label{example:symmetrized prec1}
In this example, the efficiency of the absolute value circulant matrix $|c(T_n(h\circ f))|$ as a preconditioner  for the symmetrised matrix $Y_nh(T_n(f))$ is tested and compared with $|h(c(T_n(f)))|$; for the functions $h(z)=\log(1+z)$ and $f(\theta)=0.5e^{{\mathbf{i}}\theta}$. In this case, $h\circ f\in L^1([-\pi,\pi])$, and thus, according to Theorem \ref{thm:ferrari symmetrized spectral distribution 1}, it is reasonable to test $P_n=|c(T_n(h\circ f))|$ as a preconditioner for $Y_n(T_n(h\circ f))$ and consequently for $Y_nh(T_n(f))$. 
The efficiency of the two preconditioners is shown in Figure \ref{fig:prec_h_log_f_exp}. In the top panel of the figure, the eigenvalues of the non-preconditioned matrix $Y_nh(T_n(f))$ for $n=512$ are sorted in increasing order. In the two panels that follow, the eigenvalues of the preconditioned matrix are $P_n^{-1}Y_nh(T_n(f))$. For the left graph, the preconditioner is $P_n=|c(T_n(h\circ f))|$, whereas for the right, the preconditioner is $P_n=|h(c(T_n(f)))|$.
\end{exl}

\begin{rmrk}
According to Definition \ref{circ-frob}, for the construction of $|c(T_n(h\circ f))|$, it is necessary to know the Fourier coefficients of the function $h\circ f$. However, these coefficients might not be known and must be calculated. In this case, their calculation should be included in the solution to the problem. These coefficients were not calculated analytically in the examples presented here. For the approximation of the integral (\ref{fourier}) that defines each coefficient, the trapezoidal rule with a uniform partition over $[0,\:2\pi]$ was used. This calculation can be performed using the Fast Fourier Transform. Specifically, the Fourier coefficient $a_k =\frac{1}{2\pi}\int_{0}^{2\pi}h\circ f(\theta)e^{-ik\theta}d\theta $ of the function $h\circ f$ is approximated by
\begin{align*}
\hat{a}_k=\frac{1}{m}\sum_{j=0}^{m-1}h\circ f(\frac{2\pi j}{m})w_{m}^{-kj}, \quad   w_m^k = e^{\frac{i2\pi k}{m}}.
\end{align*}
The vector $[\hat{a}_0,\:\hat{a}_1,\dots,\hat{a}_{m-1}]^T$ is exactly the Fourier transform of the vector
\[
 [h\circ f(0),\:h\circ f(\frac{2\pi}{m}),\dots,h\circ f(\frac{2\pi(m-1)}{m})]^T.
\]
  Given that $\hat{a}_k=a_k+\sum_{|l|\geq 1}a_{k+lm}$, the chosen approximation of $a_{-k}$ is $\hat{a}_{-k+m}$. Setting $m=2n$ and applying the above, we set $2n-1$ necessary coefficients at a total cost of $O(n\log n)$. This indicates that, if the function for which we seek the coefficients is a polynomial of a degree less that $\frac{m}{2}$, this procedure returns the exact coefficients of the function \cite{sauer_T}. 
\end{rmrk}

\begin{exl}\label{example:symmetrized prec2}
In the present example, the functions given in Example \ref{example:symmetrized 3} are considered; that is, $h(z)=1+z+z^2$ and $f(\theta)=-e^{{\mathbf{i}}\theta}+1+e^{-{\mathbf{i}}\theta}+e^{-{\mathbf{i}}2\theta}+e^{-{\mathbf{i}}3\theta}$. In Figure \ref{fig:prec_hon_example1}, {is shown the behaviour of the eigenvalues of the matrix $Y_{512}h(T_{512}(f))$ with and without the use of a preconditioning strategy. In particular, are shown the eigenvalues of the matrix $ Y_{512}h(T_{512}(f))$, sorted in increasing order.} In the bottom-left and bottom-right panels of Figure \ref{fig:prec_hon_example1}, the efficiency of both preconditioning strategies described in the previous example is tested. In both cases, is clear that the eigenvalues of the preconditioned matrix are clustered at -1 and 1, with up to $o(n)$ outliers.

\begin{figure}[!h]
\centering
\includegraphics[width=.48\textwidth]{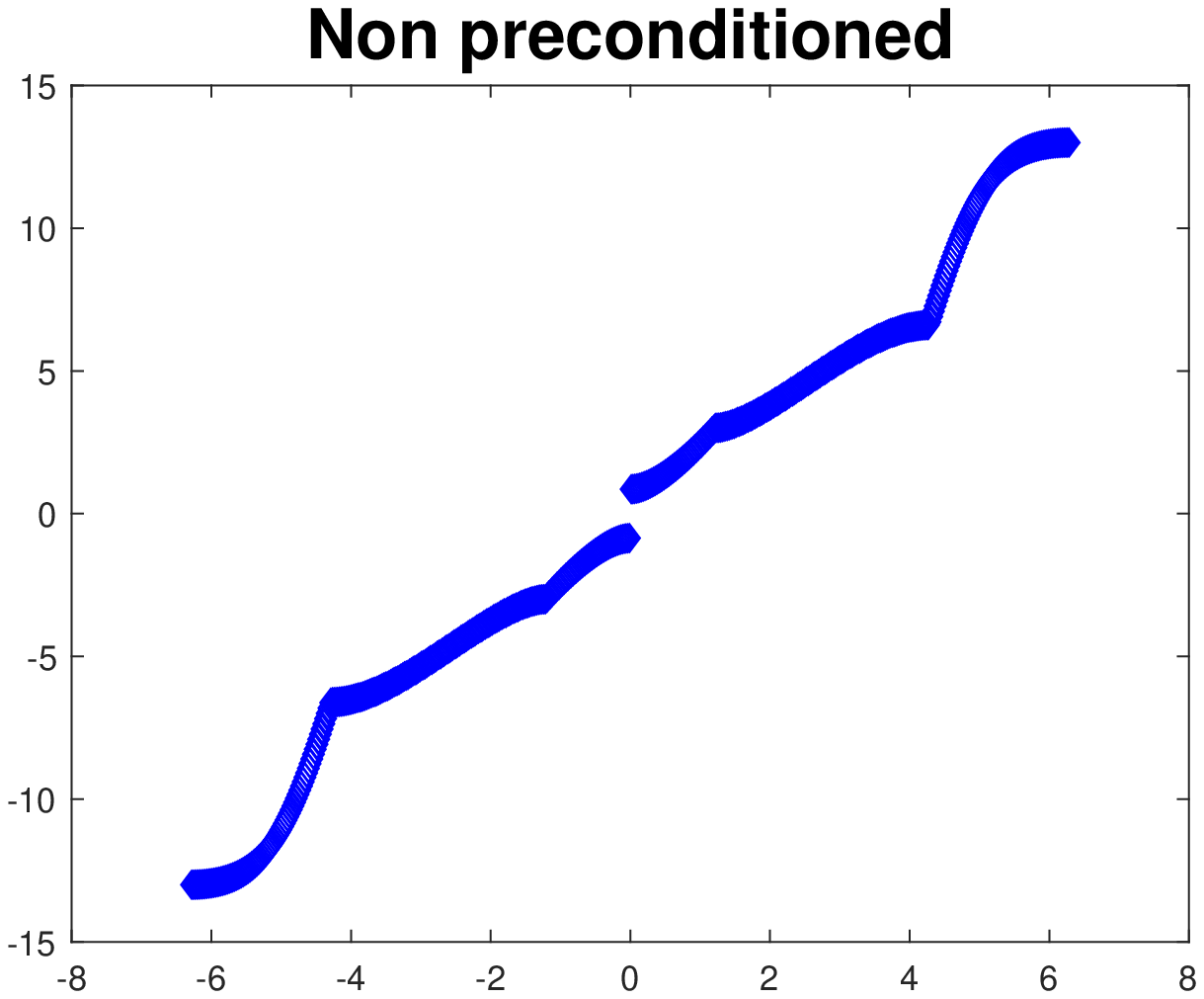} \\
\includegraphics[width=.48\textwidth]{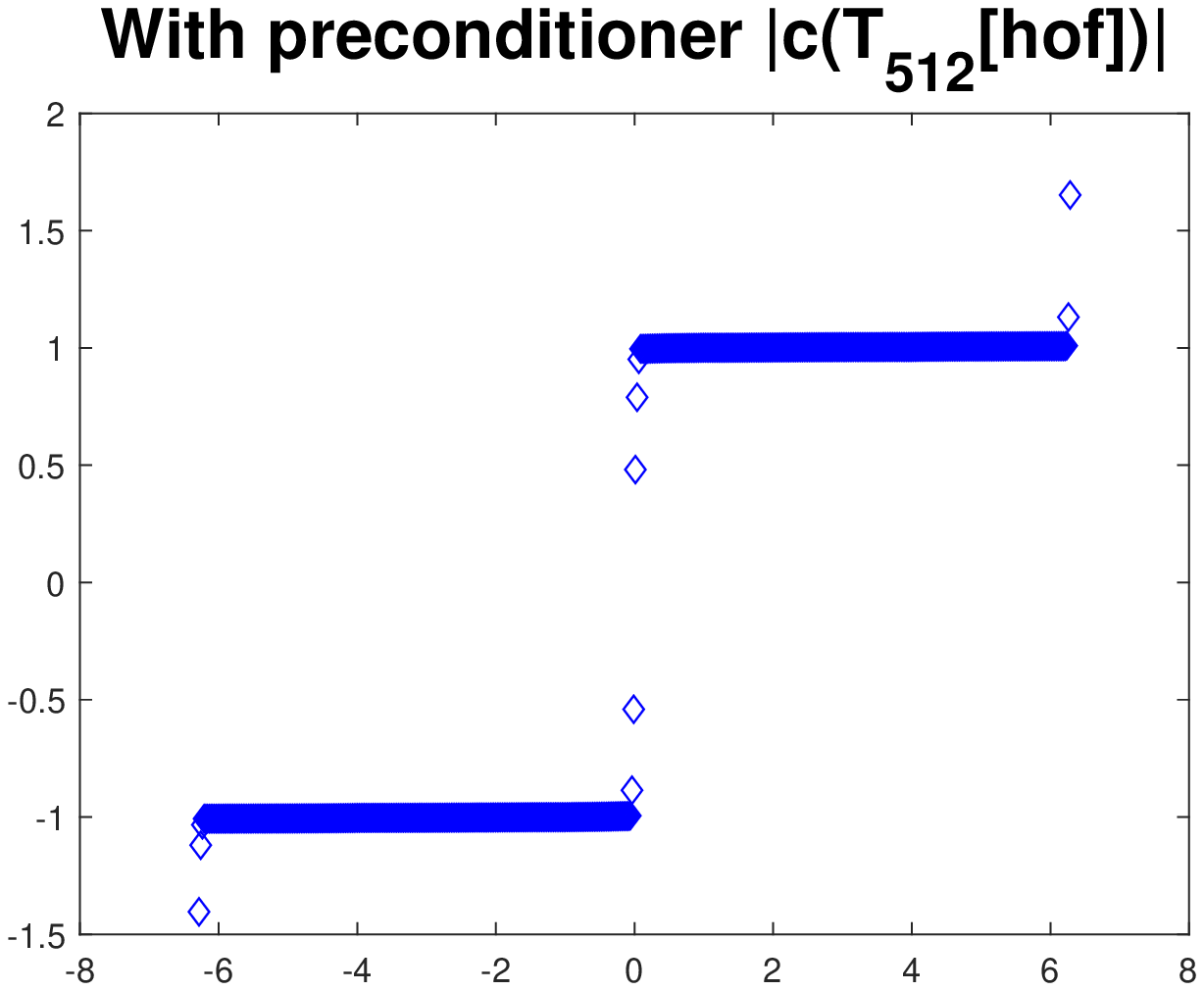}
\includegraphics[width=.48\textwidth]{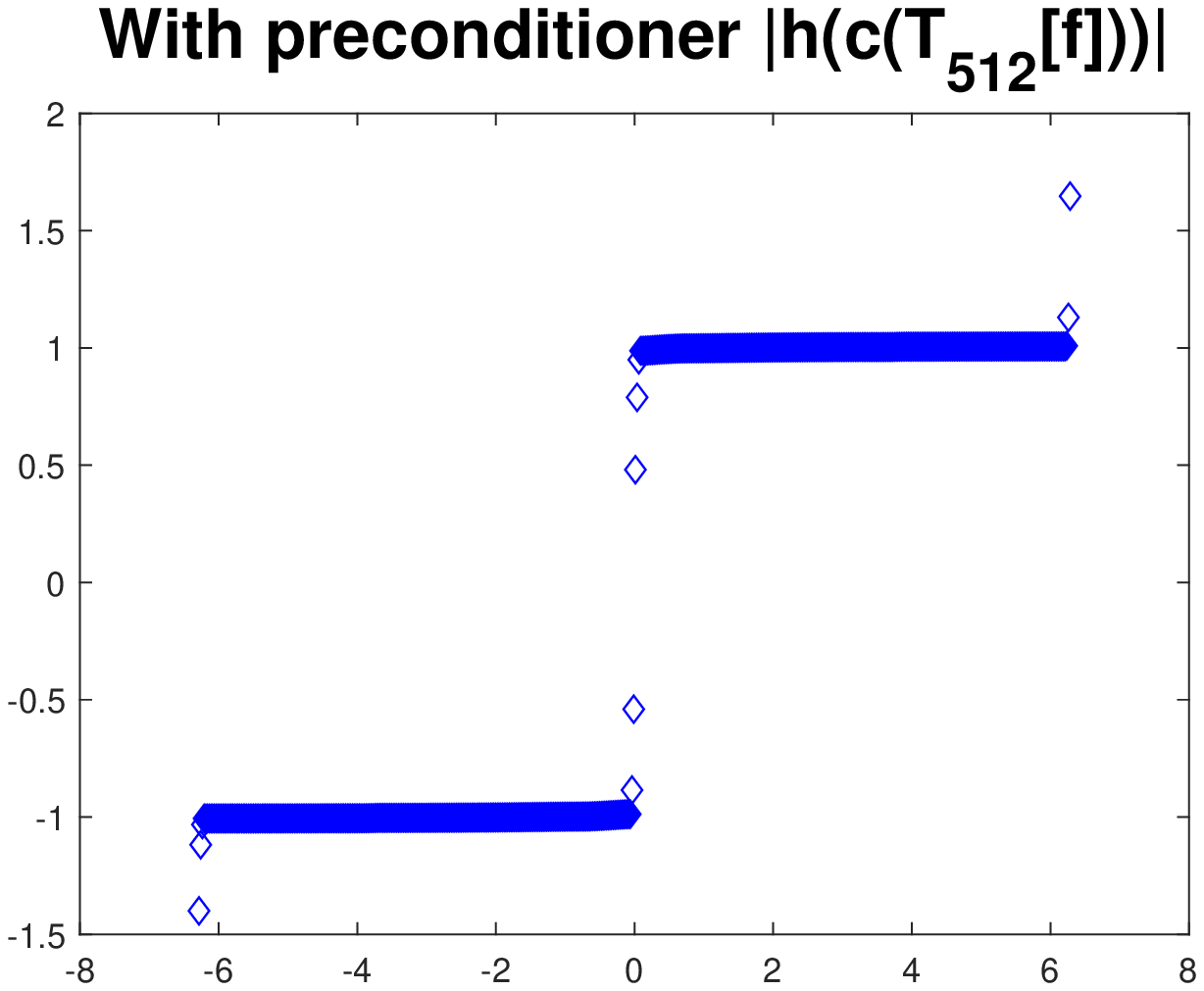}
\caption{{Spectrum of the symmetrised matrix $Y_{512}h(T_{512}(f))$, for $h(z)=1+z+z^2$ and $f(\theta)=-e^{{\mathbf{i}}\theta}+1+e^{-{\mathbf{i}}\theta}+e^{-{\mathbf{i}}2\theta}+e^{-{\mathbf{i}}3\theta}$. Top: without preconditioner; bottom left: preconditioner $P_n=|c(T_n(h\circ f))|$; bottom right: preconditioner $P_n=|h(c(T_n(f)))|$.}}
\label{fig:prec_hon_example1}
\end{figure}
\end{exl}

\begin{figure}[!h]
\centering
\includegraphics[width=.48\textwidth]{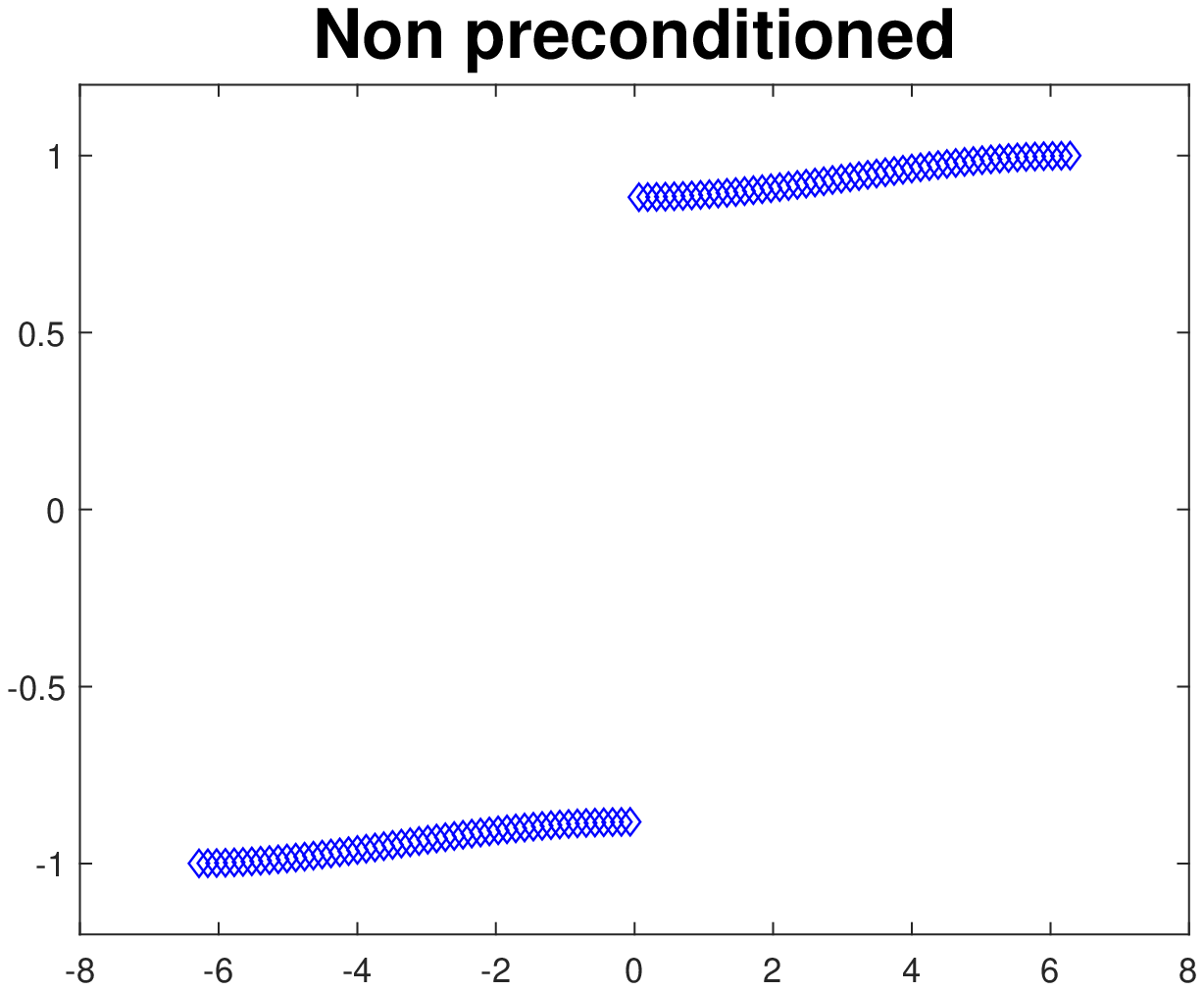} \\
\includegraphics[width=.48\textwidth]{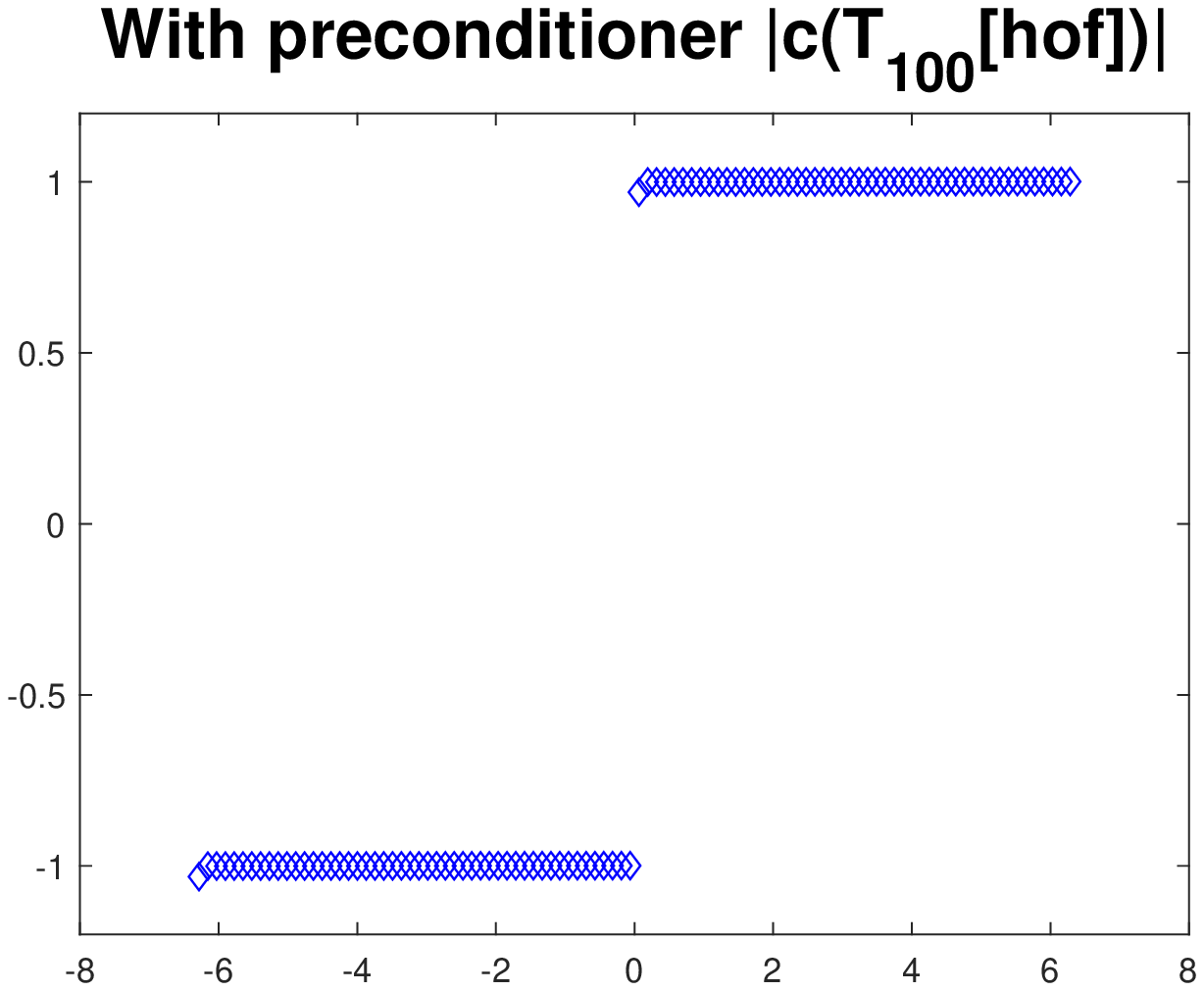}
\includegraphics[width=.48\textwidth]{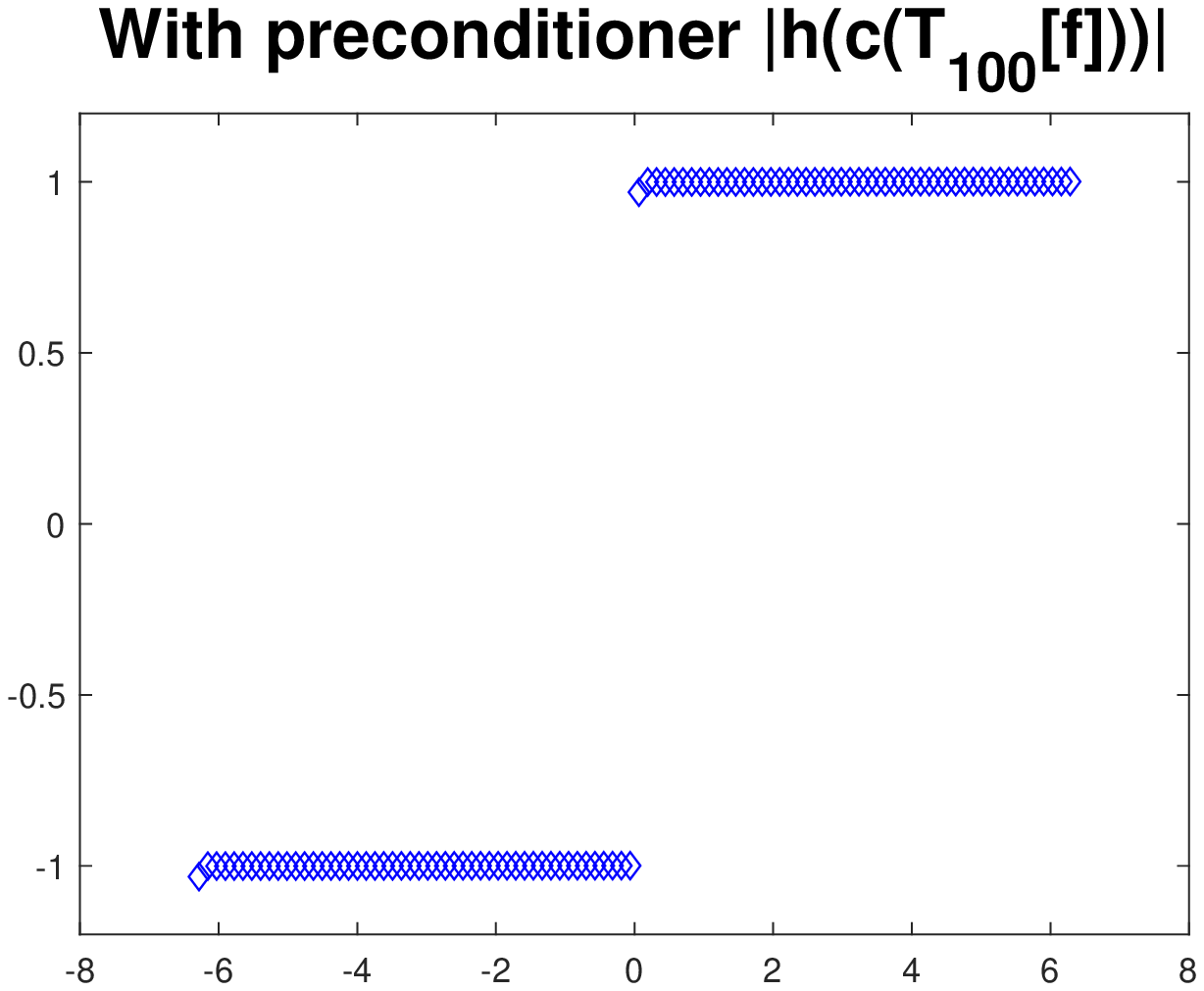}
\caption{{Spectrum of the symmetrised matrix $Y_{100}h(T_{100}(f))$ for $h(z)={\rm e}^z$ and $f(\theta)=\sum_{j=-99}^{99}a_{j}e^{\mathbf{i}j\theta}$, with $\lambda = 0.1$, $\mu = -0.9$, $\nu = 0.25$, $\sigma = 0.45$, $r = 0.05$, and $\Delta x=\frac{4}{101}$. Top: without preconditioner; bottom left: preconditioner $P_n=|c(T_n(h\circ f))|$;  bottom right: preconditioner $P_n=|h(c(T_n(f)))|$.} }
\label{fig:prec_finance}
\end{figure}

\begin{exl}\label{example:symmetrized prec3}
The last preconditioning test is performed on the computational finance case that was studied in Example \ref{example:symmetrized non_prec4}. In other words, we have $h(z)={\rm e}^z$ and $f(\theta)=\sum_{j=-99}^{99}a_{j}e^{\mathbf{i}j\theta}$, with $a_j$ defined as in \eqref{eqn:coeff_f0}-\eqref{eqn:coeff_fj}. First, the preconditioning strategy approach introduced in~\cite{Hon2018} is applied; that is, $P_n=|h(c(T_n(f)))|$. In the right-hand panel of Figure \ref{fig:prec_finance}, the eigenvalues of the preconditioned matrix $P_n^{-1}Y_nh(T_n(f))$ for $n=100$ are shown. The eigenvalues are clustered around -1 and 1, with up to two outliers. 
Analogously, we can study the eigenvalues of the preconditioned matrix $P_{100}^{-1}Y_{100}T_{100}(f)$, where $P_{100}=|c(T_{100}(h\circ f))|$. Indeed, we have $h\circ f\in L^1([-\pi,\pi])$ and, by applying the results in~\cite{FFHMS}, we have that  $P_{100}$ is a valid preconditioner for the matrix $Y_{100}h(T_{100}(f))$. The left-hand panel of Figure \ref{fig:prec_finance} confirms that the eigenvalues of the preconditioned matrix $P_{100}^{-1}Y_{100}h(T_{100}(f))$ are clustered around -1 and 1, with up to two outliers.
\end{exl}

{For each example, the validity of the two different preconditioning strategies was demonstrated. However, we have seen that, for sufficiently large matrices, the spectral results are remarkably similar. Other valid choices of preconditioning are possible; these produce a slightly different effect on the spectrum of the preconditioned matrix. Moreover, is highlighted that the strategy based on the results of \cite[Theorem 5]{FFHMS} provides an entire class of preconditioners suitable for symmetrised Toeplitz structure functions. Indeed, a preconditioner in this class is the absolute value of any circulant matrix $C_n$ such that the following singular value distribution is verified:
\begin{equation}\label{eqn:preconditioner_sv}
    \{C_n^{-1}T_n(h\circ f)\}_n \sim_\sigma 1.
\end{equation}
Concerning the choice of the preconditioning strategy based on this requirement, we used the Frobenius optimal circulant preconditioner because, from the properties of the considered $f$ and $h$, relation (\ref{eqn:preconditioner_sv}) is satisfied. }

{
Finally, we highlight that the choice of the optimal preconditioning strategy between the two approaches analysed in the examples depends on the computational aspects when constructing the matrix $P_n$, which depends on the information available for the specific example. For instance, the computational cost of the construction of the preconditioner $P_n=|c(T_n(h\circ f))|$ decreases if the Fourier coefficients of $h\circ f$ are known.}

%% file: FRACTIONAL.tex
\chapter[Preconditioners for Fractional Diffusion Equations.]{Preconditioners for Fractional Diffusion Equations\\Based on the Spectral Symbol}

\section{Introduction}
\label{fractional:section introduction}

Fractional calculus may be considered as an old and yet novel topic. Old because it dates back to a letter from L'H{\"o}pital to Leibniz in 1695; novel because it has been the object of specialised conferences and treatises for just a little over forty years. In recent years, considerable interest in fractional calculus has been stimulated by its applications in numerical analysis and modelling. Fractional differential equations (FDEs) are used to model anomalous diffusion or dispersion processes. Such  phenomena are ubiquitous in natural and social sciences. Many complex dynamical systems exhibit anomalous diffusion. Fractional kinetic equations are usually an effective method to describe these complex systems, including diffusion, diffusive convection, and Fokker–Planck fractional differential equations. Because analytical solutions are rarely available, these types  of equations are of numerical interest. When the fractional derivative $\alpha= 1$, we obtain the standard diffusion process. With $0 < \alpha < 1$, we obtain a sub-diffusion process or a dispersive, slow diffusion process with an anomalous diffusion index; meanwhile, with $\alpha> 1 $,  an ultra-diffusion process or an increased, fast diffusion process is realised.
\par
Several definitions exist for the fractional derivative, and each definition approaches the ordinary derivative in the integer order limit. In \cite{Meerschaert2004,Meerschaert2006}, the authors proposed two unconditionally stable finite difference schemes of first and second order accuracy based on the shifted Gr{\"u}nwald--Letnikov definition of fractional derivatives. 

In~\cite{Wang2010}, it was shown that once one of these methods is chosen, the coefficient matrix of the generated system can be seen as the sum of two structures, each of which is expressed as a diagonal matrix multiplied by a Toeplitz one.
Because the efficient solution of such systems is of great interest, many iterative solvers have been proposed.
Representative examples include the multigrid method (MGM) scheme proposed by~\cite{Pang2012}, the circulant preconditioner ~\cite{lei131} for the conjugate gradient normal residual (CGNR) method, and two structure-preserving preconditioners proposed in ~\cite{donatelli161}.
In the latter paper, the authors provide a detailed analysis, showing that the sequence of coefficient matrices belongs to the GLT class; furthermore, its spectral symbol, which  describes the asymptotic singular and eigenvalue distributions, is explicitly derived.
In~\cite{moghaderi171}, the analysis was extended to the two-dimensional case, and the authors compared the two-dimensional version of the structure-preserving preconditioner using a decomposition of the Laplacian~\cite{donatelli161} to a preconditioner based on an algebraic MGM.

\par
By studying the simplest (but non-trivial), case of preconditioned Toeplitz systems generated by an even, non-negative function $f$ with zeros of any positive  order, the authors prove \cite{noutsos161} that the essential spectral equivalence between the matrix sequences $\{T_n(f)\}_n$ and $\{\tau_n(f)\}_n $, (where $\{T_n(f)\}_n$ is the sequence of symmetric positive definite (SPD) Toeplitz matrices generated by this function, and $\{\tau_n(f)\}_n$ is the sequence of a specific $\tau$ matrix) is generated as 
\begin{displaymath}
\tau_n(f)=\mathbb{S}_n\mathrm{diag}(f(\boldsymbol{ \theta}))\mathbb{S}_n, \qquad \boldsymbol{ \theta}=\left[\theta_1, \theta_2,\ldots, \theta_n  \right], \qquad  \theta_j= \frac{j\pi}{n+1}=j\pi h, \qquad j=1,\ldots,n, \label{tau_preconditioner}
\end{displaymath}
and
\begin{align}
[\mathbb{S}_n]_{i,j}&=\sqrt{\frac{2}{n+1}}\sin{\left(i\theta_j\right)},  \qquad i,j=1,\ldots, n.  \label{sine_transform}
\end{align}
 We recall here that $\mathbb{S}_n$ is symmetric and orthogonal; therefore, it is the inverse of itself. Furthermore, `essential spectral equivalence' means that all the eigenvalues of $ \{ \tau_n^{-1}(f)T_n(f) \}_n$ belong to an interval $[c, C]$ (except possible $m$ outliers) and do not converge to zero as the matrix size tends to infinity. For generating functions with the order of their zero lying in the interval $[0,3]$, it is worth noting that there are no outliers.

\par
According to the analysis given in the aforementioned works, the coefficient matrix of the system depends on the diffusion coefficients of the fractional DE. In the simplest case (i.e., where they are constant and equal), this is a diagonal times a real SPD Toeplitz matrix with a generating function $\mathcal{F}_\alpha$ that is even, positive, and real, having a zero at zero of real positive order between one and two, plus a positive diagonal with constant entries that asymptotically tend to zero.
The analysis shows that this matrix is present in the more general case, where the diffusion coefficients are neither constant nor equal. In this case, a diagonal times a skew-symmetric real Toeplitz matrix is added to the coefficient matrix.
\par
Taking advantage of this fact, we propose the preconditioner $\mathcal{P}_{\mathcal{F}_\alpha}=D_n\tau_n(\mathcal{F}_\alpha)$, where $D_n$ is a suitable diagonal matrix defined as follows:
We show that this preconditioner can effectively retain the real part of the eigenvalues away from zero, whilst the sine transform maintains the cost per iteration $\mathcal{O}(n\log n)$, using a specific real algorithm or fast Fourier transform (FFT).
It turns out that this preconditioner is very efficient, and although the structure-preserving preconditioners given in \cite{donatelli161} are more efficient in the one-dimensional case, the proposed preconditioner is more efficient in two dimensions than the preconditioners described in \cite{donatelli161} and \cite{moghaderi171}.

\section{Definition of Fractional Derivative}
The fractional derivative has been defined in many ways. Each way has its own physical interpretations and applications. The classic form is given by the Riemann–Liouville integral.
\begin{defin}
As the left Riemann–Liouville integral of order $\alpha>0$, we define the operator $\prescript{}{a}I_x^{\alpha},$ 
\[
\prescript{}{a}I_x^{\alpha}f(x) = \frac{1}{\Gamma(\alpha)} \int_a^x (x-t)^{\alpha-1}f(t)dt,
\]
to be applied on locally integrable functions over $[a,b]$. Analogously, as the right Riemann–Liouville integral of order $\alpha>0$, we define the operator $\prescript{}{x}I_b^{\alpha}$,
\[
\prescript{}{x}I_b^{\alpha}f(x) = \frac{1}{\Gamma(\alpha)} \int_a^x (t-x)^{\alpha-1}f(t)dt,
\]
to be applied in the same class of functions.
\end{defin}
For $\alpha=n \: \in \mathbb{N}$ the left Riemann–Liouville becomes
\[
\prescript{}{a}I_x^{n}f(x)= \frac{1}{\Gamma(n)} \int_a^x (x-t)^{n-1}f(t)dt = \frac{1}{(n-1)!} \int_a^x (x-t)^{n-1}f(t)dt=
\]
\[
 \int_a^{x}\int_a^{s_1}\dots \int_a^{s_{n-1}}f(s_n)ds_n\dots ds_2 ds_1,
 \]
 where the last equality is given by the Cauchy formula for $n-$ times repeated integration. It is immediately apparent that, for $n \in \mathbb{N}$, the operator $\prescript{}{a}I_x^n$ is the $n-$order counter-derivative. That is,
\[
\frac{d^n}{dx^n}\prescript{}{a}I_x^nf=f.
\] 
For the right Riemann–Liouville integral of order $n$, we have that
\[
\frac{d^n}{dx^n}\prescript{}{x}I_b^nf=(-1)^nf.
\]
In general, for $\alpha,\beta >0$, we have that
\[
\frac{d}{dx}\prescript{}{a}I_x^{\alpha+1}f(x)=\frac{d}{dx}\left( \frac{1}{\Gamma{(a+1)}}\int_a^x (x-t)^{\alpha}f(t)dt \right)=\frac{\alpha}{\Gamma{(a+1)}}\int_a^x (x-t)^{\alpha -1 }f(t)dt = \prescript{}{a}I_x^{\alpha}f(x),
\]
and
\[
\prescript{}{a}I_x^{\alpha}\prescript{}{a}I_x^{\beta}=\prescript{}{a}I_x^{\alpha +\beta}.
\]
The same applies for the right Riemann–Liouville integral.
\newline
If we now set $n=\lceil \alpha \rceil, \: \alpha>0$\footnote{The function $\lceil x\rceil$ is the ceiling function of $x$ and is equal to $x$ if $x\in \mathbb{N}$; otherwise,  $\lceil x \rceil=$(the integer part of $x$)+1.
\newline
In addition, we define the floor function of $x$ as $\lfloor x \rfloor=x$ if $x\in \mathbb{N}$; otherwise, $\lfloor x \rfloor=$(the integer part of $x$).}, the operator $\frac{d^n}{dx^n}\prescript{}{a}I_x^{(n-\alpha)}$ is well defined for the locally integrable functions over $[a,b]$. So, if we set
\[
\frac{d^\alpha}{dx^\alpha}f=\frac{d^n}{dx^n}\prescript{}{a}I_x^{(n-\alpha)}f, \quad \alpha>0, \quad n=\lceil \alpha \rceil, 
\]
we have a well-defined operator for which
\[
\lim_{\alpha\rightarrow n }\frac{d^\alpha}{dx^\alpha} = \frac{d^n}{dx^n}, \quad n \: \in \mathbb{N}.
\]
\begin{defin}\label{defin:frac RL fractional derivative defin}
Let $f$ be integrable over $[a,b]$, $\alpha>0$, and $n=\lceil \alpha \rceil$. The left-hand derivative of order $\alpha$, according to Riemann–Liouville, is defined as
\begin{align*}
\frac{d^{\alpha}}{d_{+}x^\alpha}f(x)=\frac{d^n}{dx^n}\prescript{}{a}I_x^{(n-\alpha)}f(x)=\frac{1}{\Gamma(n-a)}\frac{d^n}{dx^n}\int_a^x (x-t)^{n-\alpha-1}f(t)dt.
\end{align*}
The right derivative of order $\alpha$, according to Riemann–Liouville, is defined as
\begin{align*}
\frac{d^{\alpha}}{d_{-}x^\alpha}f(x)=(-1)^n\frac{d^n}{dx^n} \prescript{}{x}I_b^{(n-\alpha)}f(x)=\frac{(-1)^n}{\Gamma(n-a)}\frac{d^n}{dx^n}\int_x^b (t-x)^{n-\alpha-1}f(t)dt.
\end{align*}
\end{defin}

It can be proven that,
\begin{align}
&\frac{d^\alpha}{d_+x^\alpha}\prescript{}{a}I_x^{\alpha}f=\frac{d^\alpha}{d_-x^\alpha}\prescript{}{x}I_b^{\alpha}f=f, \label{inverse operator}\\
&\lim_{\alpha\rightarrow n}\frac{d^\alpha}{d_+x^\alpha}=\lim_{\alpha\rightarrow n}\frac{d^\alpha}{d_-x^\alpha}=\frac{d^n}{dx^n}.
\end{align} 

\par
The definition of the fractional derivative of interest from a numerical point of view is given by Grünwald and is a generalisation of the definition of the derivative of integer order:
\[
\frac{d^n}{dx^n}f(x)=\lim_{h\rightarrow 0}\frac{1}{h^n}\sum_{k=0}^{n}(-1)^k\binom{n}{k}f(x-kh), \quad n \: \in \mathbb{N}.
\]
For $\alpha \: \in \: \mathbb{R}^+$, the left and right derivatives of order $\alpha$ over $[a,b]$ are defined as
\begin{align}\label{defin:fractional derivative}
\frac{d^\alpha}{d_{+}x^\alpha}f(x)=\lim_{h\rightarrow 0}\frac{1}{h^\alpha}\sum_{k=0}^{(x-a)/h}(-1)^k\binom{\alpha}{k}f(x-kh),\\
\frac{d^\alpha}{d_{-}x^\alpha}f(x)=\lim_{h\rightarrow 0}\frac{1}{h^\alpha}\sum_{k=0}^{(b-x)/h}(-1)^k\binom{\alpha}{k}f(x+kh),
\end{align}
respectively. The Grünwald definition of the fractional derivative is equivalent (in the continuous limit) to the Riemann–Liouville definition and immediately provides a method for numerically approximating the fractional derivative of any function.

\section{Fractional Diffusion Equations in One Dimension}

We consider the following initial value problem:
\begin{align}\label{fractional:diff equation two dimensions}
\begin{cases}
\frac{\partial u(x,t)}{\partial t}=
d_+(x,t)\frac{\partial^\alpha u(x,t)}{\partial_+ x^\alpha}+
d_-(x,t)\frac{\partial^\alpha u(x,t)}{\partial_- x^\alpha}+f(x,t), & (x,t) \in (L,R)\times(0,T]\\
u(x,t)=0,                                                         & x\in \{\mathbb{R} \setminus (L,R)\}\times[0,T],\\
u(x,0)=u_0(x)                                                     & x\in[L,R] 
\end{cases}.
\end{align}
Here, $\alpha \: \in \: (1,2)$ is the fractional derivative order, $f(x,t)$ is the source term, and the positive functions $d_{\pm}(x,t)$ are the diffusion coefficients. The left ($\partial_-$) and right ($\partial_+$)  Riemann–Liouville partial fractional derivatives are defined as

\begin{align}
\frac{\partial^\alpha u(x,t)}{\partial_+ x^\alpha}&=
\frac{1}{\Gamma(2-\alpha)}\frac{\partial^2}{\partial x^2}\int_{L}^x (x-\xi)^{1-\alpha}u(\xi,t)d\xi, \\
\frac{\partial^\alpha u(x,t)}{\partial_- x^\alpha}&=
\frac{1}{\Gamma(2-\alpha)}\frac{\partial^2}{\partial x^2}\int_x^{R} (\xi-x)^{1-\alpha}u(\xi,t)d\xi, \nonumber
\end{align}
respectively.
In the present work, to approximate the partial left and right fractional derivatives, two different numerical schemes will be used, and the effectiveness of the method proposed here can be immediately compared with already known methods. These schemes are based on Grünwald’s definition. The scheme is  
adapted in more than one dimension and shifted so that it is consistent and unconditionally stable. More specifically, the left and right partial derivatives (with respect to the spatial variable) of order $\alpha$ are defined as

\begin{align*}
\frac{\partial^\alpha u(x,t)}{\partial_+ x^\alpha}&=\lim_{h\rightarrow 0}\frac{1}{h^\alpha}\sum_{k=0}^{\lfloor(x-L)/h\rfloor}g_k^{(\alpha)}u(x-(k-1)h,t),\\
\frac{\partial^\alpha u(x,t)}{\partial_- x^\alpha}&=\lim_{h\rightarrow 0}\frac{1}{h^\alpha}\sum_{k=0}^{\lfloor(R-x)/h\rfloor}g_k^{(\alpha)}u(x+(k-1)h,t),
\end{align*}
where 
\begin{align}\label{fractional:defination binomial coefficients}
g_k^{(\alpha)}=(-1)^k\binom{\alpha}{k}=\frac{(-1)^k}{k!}\alpha(\alpha-1)\cdots(\alpha-k+1),\quad k= 0,1,\ldots\,
\end{align}
are the fractional binomial coefficients. Then, \cite{Scherer2011}
\begin{align}
\frac{\partial^\alpha u(x,t)}{\partial_+ x^\alpha}&=\frac{1}{h^\alpha}\sum_{k=0}^{\lfloor(x-L)/h\rfloor}g_k^{(\alpha)}u(x-(k-1)h,t)+\mathcal{O}(h), \label{fractional:left deriv discrtization } \\
\frac{\partial^\alpha u(x,t)}{\partial_- x^\alpha}&=\frac{1}{h^\alpha}\sum_{k=0}^{\lfloor(R-x)/h\rfloor}g_k^{(\alpha)}u(x+(k-1)h,t)+\mathcal{O}(h).\label{fractional:right deriv discrtization }
\end{align}

The following method for the discretisation of equation (\ref{fractional:diff equation two dimensions}) was given by Meerschaert and Tadjeran in \cite{Meerschaert2004}. It combines discretisation in time via the implicit Euler method with discretisation of the left and right fractional derivatives (in space) using  formulas (\ref{fractional:left deriv discrtization }) and (\ref{fractional:right deriv discrtization }), respectively. We define

\begin{align*}
x_i=L+ih, \quad     h=\frac{R-L}{n+1},    \quad i=0,\dots,n+1, \\
t_m=m\delta t, \quad \delta t=\frac{T}{M},  \quad m=0,\dots,M.
\end{align*}
We also set 
\[
u_i^{(m)}=u(x_i,t_m), \quad d_{\pm,i}^{(m)}=d_{\pm}(x_i,t_m), \quad f_i^{(m)}=f(x_i,t_m).
\]
Using the implicit Euler method, equation (\ref{fractional:diff equation two dimensions}) becomes
\begin{align*}
\frac{u_i^{(m)}-u_i^{(m-1)}}{\delta t}=d_{+,i}^{(m)}\frac{\partial^{\alpha}u_i^{(m)}}{\partial_+x^{\alpha}}+d_{-,i}^{(m)}\frac{\partial^{\alpha}u_i^{(m)}}{\partial_-x^{\alpha}}+f_i^{(m)}+\mathcal{O}(\delta t).
\end{align*}
Using formulas (\ref{fractional:left deriv discrtization }) and (\ref{fractional:right deriv discrtization }) for the approximation of the left and right fractional derivatives, respectively, we have the following finite difference scheme: 

\begin{align*}
\frac{u_i^{(m)}-u_i^{(m-1)}}{\delta t}=\frac{d_{+,i}^{(m)}}{h^\alpha}\sum_{k=0}^{i}g_k^{\alpha}u_{i-k+1}^{(m)}+\frac{d_{-,i}^{(m)}}{h^\alpha}\sum_{k=0}^{N-i+1}g_k^{\alpha}u_{i+k-1}^{(m)}+f_i^{(m)}\Rightarrow \\
\frac{h^\alpha}{\delta t}u_i^{(m)}-d_{+,i}^{(m)}\sum_{k=0}^{i}g_k^{\alpha}u_{i-k+1}^{(m)}-d_{-,i}^{(m)}\sum_{k=0}^{N-i+1}g_k^{\alpha}u_{i+k-1}^{(m)}=\frac{h^\alpha}{\delta t}u_i^{(m-1)}+h^{\alpha}f_i^m.
\end{align*}
In matrix form, this becomes
\begin{align}
\left(\nu_{M,N}\mathbb{I}_{N}+D_+^{(m)}T_{\alpha,N}+D_-^{(m)}T_{\alpha,N}^{\mathrm{T}}\right)\mathbf{u}^{(m)}&=\nu_{M,N}\mathbf{u}^{(m-1)}+h^\alpha\mathbf{f}^{(m)},
\label{fractional:2dimension analytical matrixform}
\end{align}
where $\mathbb{I}_{N}$ is the identity matrix of size $N$,
\begin{align}
\nu_{M,N}&=\frac{h^\alpha}{\delta_t}, \label{nu_defination}\\
\mathbf{u}^{(m)}&=\left[u_1^{(m)},u_2^{(m)},\ldots,u_{N}^{(m)}\right]^{\mathrm{T}},\nonumber\\
\mathbf{f}^{(m)}&=\left[f_1^{(m)},f_2^{(m)},\ldots,f_{N}^{(m)}\right]^{\mathrm{T}},\nonumber\\
[D_{\pm}^{(m)}]_{i,i}&=d_{\pm,i}^{(m)},\quad i=1,\ldots,N,\nonumber
\end{align}
and
\begin{align}
T_{\alpha,N}=-\left[
\begin{array}{ccccccccc}
g_1^{(\alpha)}&g_0^{(\alpha)}&\phantom{\ddots}\\
g_2^{(\alpha)}&g_1^{(\alpha)}&g_0^{(\alpha)}&\phantom{\ddots}\\
g_3^{(\alpha)}&g_2^{(\alpha)}&g_1^{(\alpha)}&g_0^{(\alpha)}&\phantom{\ddots}\\
\vdots&\ddots&\ddots&\ddots&\ddots\\
\vdots&\ddots&\ddots&\ddots&\ddots&\ddots\\
g_{N-1}^{(\alpha)}&g_{N-2}^{(\alpha)}&\cdots&\ddots&g_2^{(\alpha)}&g_1^{(\alpha)}&g_0^{(\alpha)}\\
g_{N}^{(\alpha)}&g_{N-1}^{(\alpha)}&\cdots\vphantom{\ddots}&\cdots&g_3^{(\alpha)}&g_2^{(\alpha)}&g_1^{(\alpha)}\\
\end{array}
\right].
\label{fractional:defination T_(a,N) matrix}
\end{align}
 If we define,
\begin{align}
\mathcal{M}_{\alpha,N}^{(m)}&=\left(\nu_{M,N}\mathbb{I}_{N}+D_+^{(m)}T_{\alpha,N}+D_-^{(m)}T_{\alpha,N}^{\mathrm{T}}\right), \label{fractional:defination2 dim first order matrix} \\
\mathbf{b}^{(m)}&=\mathbf\nu_{M,N}{u}^{(m-1)}+h^\alpha\mathbf{f}^{(m)}, \nonumber
\end{align}
then the system (\ref{fractional:2dimension analytical matrixform}) becomes
\begin{align}
\mathcal{M}_{\alpha,N}^{(m)}\mathbf{u}^{(m)}&=\mathbf{b}^{(m)}. \label{fractional:equation 2dimension linear system}
\end{align}
For investigating the behaviour of the above system and to design an effective strategy for its solution, it is necessary to investigate the properties of fractional binomial coefficients. 
\begin{prop}\label{fractional:proposition binomial coeff properties}
Let $\alpha \: \in (1,2)$ and $g_k^{(\alpha)}$ be as in (\ref{fractional:defination binomial coefficients}). The following apply:

\[
\left\{
\begin{array}{ccc}
g_0^{(\alpha)}=1,      & g_1^{(\alpha)}=-\alpha,     & g_0^{(\alpha)}>g_2^{(\alpha)}>g_3^{(\alpha)}>\dots>0\\
\sum_{k=0}^{\infty}g_k^{(\alpha)}=0,   & \sum_{k=0}^{n}g_k^{(\alpha)}<0,\: n\geq 1 & \sum_{k=0}^{\infty}|g_k^{(\alpha)}|=2\alpha
\end{array}
\right\}.
\]
That $g_0^{(\alpha)}=1$ and $g_1^{(\alpha)}=-\alpha$ immediately follows from this definition. In addition, because $\alpha  \in  (1,2)$, $1>g_2^{(\alpha)}=(-1)^2\frac{\alpha(\alpha-1)}{2}>0$. If $k>2$, then
\[
g_{k+1}^{(\alpha)}=(-1)^{k+1}\frac{\alpha(\alpha-1)\dots(\alpha-k+1)(\alpha-k)}{(k+1)!}=(-1)\frac{(\alpha-k)}{(k+1)}g_k^{\alpha}.
\]
From the above, it is clear that for $k>2$, the term $g_{k+1}^{(\alpha)}$  retains a positive sign and is less than $g_{k}^{(\alpha)}$. Additionally, using the known identity
\begin{align}
(1-x)^{\alpha}=\sum_{k=0}^{\infty}\binom{\alpha}{k}(-x)^k=\sum_{k=0}^{\infty}g_k^{\alpha}x^k, \label{fractional:binomial identity}
\end{align}
we have,
\[
0=(1-1)^{\alpha}=\sum_{k=0}^{\infty}\binom{\alpha}{k}(-1)^k=\sum_{k=0}^{\infty}g_k^{\alpha}.
\]
Because the only negative term in the above zero sum is $-\alpha$, the sum of the remaining terms must be equal to $\alpha$. Thus, it turns out that, on one hand, any partial sum is less than zero; on the other hand, the sum of the absolute values is equal to $2\alpha$.
\end{prop}

Using Proposition (\ref{fractional:proposition binomial coeff properties}), we find that the matrix $T_{\alpha,N}$ defined in (\ref{fractional:defination T_(a,N) matrix}) is strictly diagonally dominant, and therefore positive and invertible. In \cite{Wang2010}, it was proven that the coefficient matrix of the system (\ref{fractional:equation 2dimension linear system}) $\mathcal{M}_{\alpha,N}^{(m)}$ is also strictly diagonally dominant and invertible. More essential properties of the involved matrices are revealed using the theory of GLT matrix sequences below.
\par

An interesting property of the matrix \eqref{fractional:defination T_(a,N) matrix}, arising from the operator that this matrix discretises, is given in the following proposition.
\begin{prop}
A linear system with a coefficient matrix of $-T_{\alpha,N}$ can be solved by a direct method, with a total cost of $O(n\log n)$ operations. In fact, only one Toeplitz matrix-vector multiplication [of cost $O(n\log n)$] and then a forward substitution [of cost $O(n)$] are needed.
\end{prop}

\begin{prf*}
According to the relationship \eqref{inverse operator}, the left Riemann–Liouville fractional derivative of order $\alpha$ is the left inverse operator of the Riemann–Liouville fractional integral of order $\alpha$ in the space of integrable functions over an interval $[a,b]$. For this integral, in the continuous limit, we apply 
\[
\prescript{}{a}I_x^{\alpha}f(x) = \frac{1}{\Gamma(\alpha)} \int_a^x (x-t)^{\alpha-1}f(t)dt=\lim_{h\rightarrow 0}h^\alpha\sum_{k=0}^{(x-a)/h}(-1)^k\binom{-\alpha}{k}f(x-kh).
\]
 Inspired by this fact, we can pre-multiply any system of the form $-T_{\alpha,N}x=b$ with the matrix that implements the inverse operator according to the above scheme; that is,
 \begin{align}
\hat{T}_{-\alpha,N}=\left[
\begin{array}{ccccccccc}
g_0^{(-\alpha)}    & 0                & \phantom{\ddots}\\
g_1^{(-\alpha)}    & g_0^{(-\alpha)}  & 0               &  \\
\vdots             &\ddots            &\ddots           &\ddots\\
\vdots             &\ddots            &\ddots           &\ddots           & 0 \\
g_{N-1}^{(-\alpha)}&g_{N-1}^{(\alpha)}&\cdots           &g_1^{(-\alpha)}  &g_0^{(-\alpha)}\\
\end{array}
\right].
\label{fractional:defination T_(-a,N) matrix}
\end{align}
It is worth noting that $\hat{T}_{-\alpha,N}$ is the inverse of the lower triangular Toeplitz matrix with the vector $[g_0^{(\alpha)},g_1^{(\alpha)},\dots,g_{N-1}^{(\alpha)}]$ as the first column and it implements the fractional derivative of order $\alpha$ without displacement. Then, $\hat{T}_{-\alpha,N}(-T_{\alpha,N})$ is the following lower Hessenberg matrix:
\begin{align}
\hat{T}_{-\alpha,N}(-T_{\alpha,N})=\left[
\begin{array}{ccccccccc}
-g_1^{(-\alpha)}     & 1               &  \\
-g_2^{(-\alpha)}     & 0               & 1               &   \\
\vdots               &\ddots           &\ddots           &\ddots\\
-g_{N-1}^{(-\alpha)} &0                &\ddots           &\ddots         & 1 \\
-g_{N}^{(-\alpha)}   &0                &\cdots           &0             & 0\\
\end{array}
\right].
\label{fractional:defination T_(-a,N)xT(a,n) matrix}
\end{align}
Because the elements of the first column are all known, it is not necessary to make the multiplication explicitly, and we only need multiply the right-hand vector with $\hat{T}_{-\alpha,N}$. Of course, this multiplication can be performed with a cost of $O(n\log n)$ operations, whilst the system with the coefficient matrix $\hat{T}_{-\alpha,N}(-T_{\alpha,N})$ is clear and can be solved with a forward substitution.

\end{prf*}

\subsection{Second-order Finite Difference Discretisation} \label{section:finite differences 2nd order}
It can be shown \cite{Tian2015} that, 
\begin{align}
\frac{\partial^\alpha u(x,t)}{\partial_+ x^\alpha}&=\frac{1}{h^\alpha}\sum_{k=0}^{\lfloor(x-L)/h\rfloor}w_k^{(\alpha)}u(x-(k-1)h,t)+\mathcal{O}(h^2), \label{fractional:second order left deriv discrtization } \\
\frac{\partial^\alpha u(x,t)}{\partial_- x^\alpha}&=\frac{1}{h^\alpha}\sum_{k=0}^{\lfloor(R-x)/h\rfloor}w_k^{(\alpha)}u(x+(k-1)h,t)+\mathcal{O}(h^2), \label{fractional:second order right deriv discrtization }
\end{align}
where
\begin{align}
    w_0^{(\alpha)}&=\frac{\alpha}{2}g_0^{(\alpha)}, \label{fractional:defination w_0}\\
    w_k^{(\alpha)}&=\frac{\alpha}{2}g_k^{(\alpha)}+\frac{2-\alpha}{2}g_{k-1}^{(\alpha)}, \quad k\geq 1, \label{fractional:defination w_k}
\end{align}
and $g_k^{(\alpha)}$ as defined in (\ref{fractional:defination binomial coefficients}). 
\newline
In this case, the matrix $T_{\alpha,N}$ in the system (\ref{fractional:2dimension analytical matrixform}) must be replaced by the following matrix: 
\begin{align}
S_{\alpha,N}=-\left[
\begin{array}{ccccccccc}
w_1^{(\alpha)}&w_0^{(\alpha)}&\phantom{\ddots}\\
w_2^{(\alpha)}&w_1^{(\alpha)}&w_0^{(\alpha)}&\phantom{\ddots}\\
w_3^{(\alpha)}&w_2^{(\alpha)}&w_1^{(\alpha)}&w_0^{(\alpha)}&\phantom{\ddots}\\
\vdots&\ddots&\ddots&\ddots&\ddots\\
\vdots&\ddots&\ddots&\ddots&\ddots&\ddots\\
w_{n_1-1}^{(\alpha)}&w_{n_1-2}^{(\alpha)}&\cdots&\ddots&w_2^{(\alpha)}&w_1^{(\alpha)}&w_0^{(\alpha)}\\
w_{n_1}^{(\alpha)}&w_{n_1-1}^{(\alpha)}&\cdots\vphantom{\ddots}&\cdots&w_3^{(\alpha)}&w_2^{(\alpha)}&w_1^{(\alpha)}\\
\end{array}
\right].
\label{eq:grunwaldmatrix2}
\end{align}

\begin{prop}\label{fractional:w_k^a properties}
Let $\alpha \in (1,2)$ and $w_k^{(\alpha)}$ be as defined in (\ref{fractional:defination w_0})–(\ref{fractional:defination w_k}). Thus, it is apparent from the definition that
\[
w_0^{(\alpha)}=\frac{\alpha}{2}>0, \quad w_1^{(\alpha)}=\frac{2-\alpha-\alpha^2}{2}<0, \quad w_2^{(\alpha)}=\frac{\alpha(\alpha^2+\alpha-4)}{4}.
\]
Examining the sign of $w_2^{(\alpha)}=\frac{\alpha(\alpha^2+\alpha-4)}{4}$ over $(1,2)$, we find that 
\[
w_2^{(\alpha)}\leq0 \:\: \alpha \in\left(1,\frac{-1+\sqrt{17}}{2}\right], \quad w_2^{(\alpha)}>0 \:\: \alpha \in\left(\frac{-1+\sqrt{17}}{2},2\right).
\]
In addition, if $k>2$, we have that $w_k^{(\alpha)}>0$, because $w_k^{(\alpha)}$ is a weighted average of two positive terms. From the definition and properties of fractional binomial coefficients (\ref{fractional:proposition binomial coeff properties}), we have
\[
1> w_0^{(\alpha)}> w_3^{(\alpha)}> w_4^{(\alpha)}>\dots>0. 
\] 
According to the above,
\[
\sum_{k=0}^{\infty}w_k^{(\alpha)}=\frac{\alpha}{2}g_0^{(\alpha)}+\sum_{k=1}^{\infty}\left(\frac{\alpha}{2}g_k^{\alpha}+\frac{2-\alpha}{2}g_{k-1}^{(\alpha)}\right)=\frac{\alpha}{2}\sum_{k=0}^{\infty}g_k^{(\alpha)}+\frac{2-\alpha}{2}\sum_{k=1}^{\infty}g_{k-1}^{(\alpha)}=0.
\]
If $\alpha \in\left(1,\frac{-1+\sqrt{17}}{2}\right]$, then the two negative terms in the above zero sum are $w_1^{(\alpha)}$ and $w_2^{(\alpha)}$. Therefore, the sum of the remaining terms must be equal to $|w_1^{(\alpha)}|+|w_2^{(\alpha)}|$, and thus,
\[
\sum_{k=0}^{\infty}|w_k^{(\alpha)}|=2(|w_1^{(\alpha)}|+|w_2^{(\alpha)}|).
\]
Analogously, if $\alpha \in\left(\frac{-1+\sqrt{17}}{2},2\right)$, then the only negative term in the zero sum is $w_1^{(\alpha)}$, and we conclude that
\[
\sum_{k=0}^{\infty}|w_k^{(\alpha)}|=2|w_1^{(\alpha)}|.
\]
In any case,
\[
\sum_{k=0}^{\infty}|w_k^{(\alpha)}|<\infty,
\]
and if $n>1$,
\[
\sum_{k=0}^{n}w_k^{(\alpha)}<0.
\]
\end{prop}

\section{Spectral Analysis of the Coefficient Matrices}

In the present section, we present in detail the distribution of eigenvalues and singular values of the matrix sequences $\{T_{\alpha,n}\}_n$, $\{T_{\alpha,n}+T_{\alpha,n}^{\mathrm{T}}\}_n$, and $\{\mathcal{M}_{\alpha,n}^{(m)}\}_n$ appearing on the left-hand side of the system (\ref{fractional:equation 2dimension linear system}). For the first two, which are clearly Toeplitz, their spectral symbol \footnote{In the definition of Toeplitz matrices (Definition \ref{glt:defination Toeplitz}), the term 'generating function' is used instead of 'spectral symbol' for the function whose Fourier coefficients compose the diagonals of the matrices of the sequence. Nevertheless, from the property {\bf GLT3}, every Toeplitz sequence is GLT with a spectral symbol that is the same as the generating function. Hence, 
the term 'spectral symbol' is used instead of 'generating function' for reasons of homogeneity, because the spectral behaviour of all matrix sequences shown here can only be analysed using the GLT theory.} is analyzed. The sequence $\{\mathcal{M}_{\alpha,N}^{(m)}\}_n$ belongs to the GLT class, and its spectral symbol is also considered. In addition, the distribution of the sequence $\{S_{\alpha,n}\}_n$ appears in the second-order finite difference discretisation is considered.

\par
\begin{defin}\label{fractional:defin wiener class}
Let the sequence $\{f_k\}_k$ be such that $\sum_{k=0}^{\infty}|f_k|<\infty$. Then, the series $\sum_{k=0}^{\infty}f_k e^{ik\theta}$ converges uniformly to a continuous $2\pi-$periodic function $f(\theta)$. The set of all these functions is the Wiener class. 
\end{defin}

\begin{prop}\label{fractional:first order symbols}
\cite{donatelli161}
Let $\alpha \in (1,2)$. The matrix sequences $\{T_{\alpha,n}\}_n$, $\{T_{\alpha,n}^{\mathrm{T}}\}_n$, and \phantom{aaaaa} $\{T_{\alpha,n}+T_{\alpha,n}^{\mathrm{T}}\}_n$ are Toeplitz with spectral symbols 
\begin{align}
&g_{\alpha}(\theta)=-e^{-i\theta}(1-e^{i\theta})^{\alpha},  \label{fractional:g_alpha defination1} \\
&g_{\alpha}(-\theta)=\overline{g_{\alpha}(\theta)}, \nonumber \\
&p_{\alpha}(\theta)=g_{\alpha}(\theta)+\overline{g_{\alpha}(\theta)}, \label{fractional:p_alpha defination}
\end{align}
respectively.
\newline
We observe that 
\[
[T_{\alpha,n}]_{k,j}=\begin{cases}
g_{k-j+1}^{(\alpha)} \quad k-j\geq -1, \\{}
\\
0 \quad k-j<-1
\end{cases}.
\]
Based on the properties of the sequence $\{g_n^{\alpha}\}_n$ \eqref{fractional:proposition binomial coeff properties} and the definition of the Wiener class (Definition \ref{fractional:defin wiener class}), $-\sum_{k=-1}^{\infty}g_{k+1}^{(\alpha)}e^{ik\theta}$ is well defined and belongs to that class.
\newline
Then,
\begin{align*}
&-\sum_{k=-1}^{\infty}g_{k+1}^{(\alpha)}e^{ik\theta}=-g_0^{(\alpha)}e^{-i\theta}-g_1^{(\alpha)}-g_2^{(\alpha)}e^{i\theta}- \dots = \\
&-e^{-i\theta}(g_0^{(\alpha)}+g_1^{(\alpha)}e^{i\theta}+g_2^{(\alpha)}e^{2i\theta}+\dots)= \\
&-e^{-i\theta}\sum_{k=0}^{\infty}g_k^{(\alpha)}e^{ik\theta}=-e^{-i\theta}(1-e^{i\theta})^\alpha=g_{\alpha}(\theta).
\end{align*}
Therefore, $T_{\alpha,n}=T_n(g_{\alpha}(\theta))$. It is clear that the spectral symbol of the sequence $T_{\alpha,n}^{\mathrm{T}}$ is the function $g_{\alpha}(-\theta)$, which, because the Fourier coefficients are real, entails that $g_{\alpha}(-\theta)=\overline{g_{\alpha}(\theta)}$. Also, 
\[
T_{\alpha,n}+T_{\alpha,n}^{\mathrm{T}}=T_n(g_{\alpha}(\theta))+T_n(\overline{g_{\alpha}(\theta)})=T_n(g_{\alpha}(\theta)+\overline{g_{\alpha}(\theta)})=T_n(p_{\alpha}(\theta)).
\]
\end{prop}

The spectral distribution of the sequence $\{\mathcal{M}_{\alpha,n}^{(m)}\}_n$ was studied by \cite{donatelli161}. The following propositions summarise the results required to design an effective strategy for solving \ref{fractional:proposition binomial coeff properties}.

\begin{prop}
\cite{donatelli161}
Let $\nu_{M,N}=o(1)$ and $d_+(x,t)=d_+(x)$ and $d_-(x,t)=d_-(t)$ be Riemann integrable over $[L,R]$. For the sequence $\{\mathcal{M}_{\alpha,n}^{(m)}\}_n$, as defined in (\ref{fractional:defination2 dim first order matrix}), we have 
\[
\{\mathcal{M}_{\alpha,n}^{(m)}\}_n \sim_{GLT} \hat{h}_{\alpha}(\hat{x},\theta),
\]
with
\[
\hat{h}_{\alpha}(\hat{x},\theta)=h_{\alpha}(L+(R-L)\hat{x},\theta),\quad h_{\alpha}(x,\theta)=d_+(x)g_{\alpha}(\theta)+d_-(x)g_{\alpha}(-\theta),
\]
where $(\hat{x},\theta)\in[0,1]\times[-\pi,\pi]$ and $(x,\theta)\in[L,R]\times[-\pi,\pi]$. In addition, from property {\bf GLT1}, 
\[
\{\mathcal{M}_{\alpha,n}^{(m)}\}_n \sim_{\sigma}h_{\alpha}(x,\theta).
\]
If $d_+(x)=d_-(x)$, then the function $h_{\alpha}(x,\theta)$ is real, and the matrices $\mathcal{M}_{\alpha,n}^{(m)}$ have real eigenvalues, and
\[
\{\mathcal{M}_{\alpha,n}^{(m)}\}_n \sim_{\lambda}h_{\alpha}(x,\theta).
\]
\end{prop}

\begin{prop}\label{fractional:p_a order of zero}
If $\alpha \in (1,2)$, the function $p_{\alpha}(\theta)$ has a zero of order $\alpha$ at $0$.\footnote{If $f$ is a continuous, non-negative, and real function over $[a,b]$, we say that it has a zero of order $\alpha$ at $\theta_0 \in[a,b]$, if there exist $C_1,C_2>0$ such that $\liminf\limits_{\theta\rightarrow\theta_0}\frac{f(\theta)}{|\theta-\theta_0|^{\alpha}}=C_1$ and $\limsup\limits_{\theta\rightarrow\theta_0}\frac{f(\theta)}{|\theta-\theta_0|^{\alpha}}=C_2$. } For $\alpha=2$, $p_2(\theta)=4-4\cos(\theta)$ and the proposition is true. For $\alpha=1$, $p_1(\theta)=2-2\cos(\theta)$ and the proposition is untrue, because this trigonometric polynomial has a zero of order two.
\end{prop}

\begin{prop}\label{fractional:h_a order of zero}
For the functions $p_{\alpha}(\theta)$ and $h_\alpha(x,\theta)$, as defined above, we apply
\begin{align*}
&\lim_{\theta\rightarrow 0^+}\frac{h_{\alpha}(x,\theta)}{p_{\alpha}(\theta)}=\frac{d_+(x)+d_-(x)}{2}-\mathbf{i}\tan\left(\alpha\frac{\pi}{2}\right)\frac{d_+(x)-d_-(x)}{2} \\
&\lim_{\theta\rightarrow 0^-}\frac{h_{\alpha}(x,\theta)}{p_{\alpha}(\theta)}=\frac{d_+(x)+d_-(x)}{2}+\mathbf{i}\tan\left(\alpha\frac{\pi}{2}\right)\frac{d_+(x)-d_-(x)}{2}.
\end{align*}
\end{prop}

\par
It is evident from the propositions above that the coefficient matrix of the system \eqref{fractional:equation 2dimension linear system}, $\mathcal{M}_{\alpha,n}^{(m)}$, is in a bad condition, because its minimum singular value or eigenvalue if $d_+(x)=d_-(x)$ converges to zero with order $O(n^{-\alpha})$. An effective strategy for preconditioning the system is to keep the singular values or eigenvalues of the system away from zero. It should be noted that, if the preconditioner $\mathcal{C}_n$ is selected from the GLT class (e.g., as a band Toeplitz or circulant with a spectral symbol $f$), then from the property {\bf GLT2},  $\mathcal{C}_n^{-1}\mathcal{M}_{\alpha,n}^{(m)}\sim_{GLT}\frac{h_{\alpha}(x,\theta)}{f(\theta)}$ and $\mathcal{C}_n^{-1}\mathcal{M}_{\alpha,n}^{(m)}\sim_{\sigma}\frac{h_{\alpha}(x,\theta)}{f(\theta)}$. In this case, if $d_+(x)\neq d_-(x)$, the preconditioner is not optimal. This is because the singular values or eigenvalues of the sequence cannot be clustered at $1$, because the spectral symbol $\frac{h_{\alpha}(x,\theta)}{f(\theta)}$ is a nontrivial function of $x$.
\begin{prop}
\cite{moghaderi171}
Let $\alpha \in (1,2)$. The sequences $\{S_{\alpha,n}\}_n$,  $\{S_{\alpha,n}^{\mathrm{T}}\}_n$, and $\{S_{\alpha,n}+S_{\alpha,n}^{\mathrm{T}}\}_n$ are Toeplitz with spectral symbols
\begin{align}
&w_\alpha(\theta)=-\left(\frac{2-\alpha(1-e^{-\mathbf{i}\theta})}{2}\right)\left(1-e^{\mathbf{i}\theta}\right)^\alpha,\label{eq_w_function}  \\
&w_\alpha(-\theta)=\overline{w_\alpha(\theta)}, \nonumber \\
&q_{\alpha}(\theta)=w_\alpha(\theta)+\overline{w_\alpha(\theta)} \label{second_order_symbol},
\end{align}
respectively.
\end{prop}
\begin{prop}
If $\alpha \in (1,2)$, the function $q_{\alpha}(\theta)$ has a zero of order $\alpha$ at 0.
\end{prop}

\section{Fractional Diffusion Equations in Two Dimensions}
\label{sec:introduction:fde2d}
We consider the following initial value problem in two dimensions:

\begin{align}
\begin{cases}
\frac{\partial u(x,y,t)}{\partial t}=
d_+(x,y,t)\frac{\partial^\alpha u(x,y,t)}{\partial_+ x^\alpha}+
d_-(x,y,t)\frac{\partial^\alpha u(x,y,t)}{\partial_- x^\alpha}+\\
\hspace{1.6cm}+ e_+(x,y,t)\frac{\partial^\beta u(x,y,t)}{\partial_+ y^\beta}+
e_-(x,y,t)\frac{\partial^\beta u(x,y,t)}{\partial_- y^\beta}+
f(x,y,t),&(x,y,t)\in\Omega\times (0,T),\\
u(x,y,t)=0,&(x,y,t)\in\mathbb{R}^2\setminus\Omega\times[0,T],\\
u(x,y,0)=u_0(x,y),&x\in\bar{\Omega}.
\end{cases}
\label{eq:fde}
\end{align}
Here, $\Omega=(L_1,R_1)\times (L_2,R_2),$ $\alpha,\beta \in (1,2)$ is the fractional order of the derivative, $f(x,y,t)$ is the source term, and the non-negative functions $d_{\pm}(x,y,t)$ and $e_{\pm}(x,y,t)$ are the diffusion coefficients. 

In this case, the left ($\partial_+$) and right ($\partial_-$)  Riemann–Liouville fractional derivatives are defined as
\begin{align*}
\frac{\partial^\alpha u(x,y,t)}{\partial_+ x^\alpha}&=
\frac{1}{\Gamma(2-\alpha)}\frac{\partial^2}{\partial x^2}\int_{L_1}^x (x-\xi)^{1-\alpha}u(\xi,y,t)d\xi, \\
\frac{\partial^\alpha u(x,y,t)}{\partial_- x^\alpha}&=
\frac{1}{\Gamma(2-\alpha)}\frac{\partial^2}{\partial x^2}\int_x^{R_1} (\xi-x)^{1-\alpha} u(\xi,y,t)d\xi,\\
\frac{\partial^\beta u(x,y,t)}{\partial_+ y^\beta}&=
\frac{1}{\Gamma(2-\beta)}\frac{\partial^2}{\partial y^2}\int_{L_2}^y (y-\eta)^{1-\beta} u(x,\eta,t)d\eta,\\
\frac{\partial^\beta u(x,y,t)}{\partial_- y^\beta}&=
\frac{1}{\Gamma(2-\beta)}\frac{\partial^2}{\partial y^2}\int_y^{R_2} (\eta-y)^{1-\beta}u(x,\eta,t)d\eta.
\end{align*}
For the discretisation of Equation (\ref{eq:fde}), we use a method that combines the Crank--Nicolson method in time with the second-order finite difference in spatial domain scheme (\ref{fractional:second order left deriv discrtization })–(\ref{fractional:second order right deriv discrtization }), adapted for two dimensions. The method was proposed and proven to be consistent and unconditionally stable in \cite{Tian2015}.

We define,
\begin{align}
h_x&=\frac{R_1-L_1}{n_1+1}=(R_1-L_1)h_1\quad\quad x_{i}=L_1+ih_x,\quad i=1,\ldots, n_1,\nonumber\\
h_y&=\frac{R_2-L_2}{n_2+1}=(R_2-L_2)h_2\quad\quad y_{i}=L_2+ih_y,\quad i=1,\ldots, n_2, \nonumber
\end{align}
and $N=n_1n_2$.
For the unknown function $u(x,y,t)$, we set $u_{i,j}^{(m)}=u(x_i,y_j,t^{(m)})$ and 
\begin{align}
\mathbf{u}^{(m)}&=[u_{1,1}^{(m)},\ldots,u_{n_1,1}^{(m)},u_{1,2}^{(m)},\ldots,u_{n_1,2}^{(m)},\ldots, u_{1,n_2}^{(m)},\ldots,u_{n_1,n_2}^{(m)}]^{\mathrm{T}}.\nonumber
\end{align}
For the diffusion coefficients  $d_+(x,y,t)$, $d_-(x,y,t)$,  $e_+(x,y,t)$, and $e_-(x,y,t)$, we set      
\newline $d_{i,j}^{\pm,(m)}=d_\pm(x_i,y_j,t^{(m)})$ and $e_{i,j}^{\pm,(m)}=e_\pm(x_i,y_j,t^{(m)})$. The corresponding discretisation is as follows:.
\begin{align}
\mathbf{d}_{\pm}^{(m)}&=[d_{1,1}^{\pm,(m)},\ldots,d_{n_1,1}^{\pm,(m)},d_{1,2}^{\pm,(m)},\ldots,d_{n_1,2}^{\pm,(m)},\ldots, d_{1,n_2}^{\pm,(m)},\ldots,d_{n_1,n_2}^{\pm,(m)}]^{\mathrm{T}},\nonumber\\
\mathbf{e}_{\pm}^{(m)}&=[e_{1,1}^{\pm,(m)},\ldots,e_{n_1,1}^{\pm,(m)},e_{1,2}^{\pm,(m)},\ldots,e_{n_1,2}^{\pm,(m)},\ldots, e_{1,n_2}^{\pm,(m)},\ldots,e_{n_1,n_2}^{\pm,(m)}]^{\mathrm{T}}. \nonumber
\end{align}
For the discretisation of the source term $f(x,y,t)$, we set $f_{i,j}^{(m)}=f(x_i,y_j,t^{(m)})$ and 
\begin{align}
\mathbf{v}^{(m-1/2)}&=[f_{1,1}^{(m-1/2)},\ldots,f_{n_1,1}^{(m-1/2)},f_{1,2}^{(m-1/2)},\ldots,f_{n_1,2}^{(m-1/2)},\ldots, f_{1,n_2}^{(m-1/2)},\ldots,f_{n_1,n_2}^{(m-1/2)}]^{\mathrm{T}}. \nonumber
\end{align}
We define $D_\pm^{(m)}=\mathrm{diag}(\mathbf{d}_\pm^{(m)})$ and $E_\pm^{(m)}=\mathrm{diag}(\mathbf{e}_\pm^{(m)})$.
In the equation, fractional derivatives of different orders $\alpha$ and $\beta$ appear, and it is also possible to obtain different numbers of discretisation points $n_1,n_2$ in each spatial domain. Thus, we define the matrices $S_{\alpha,n_1}$, $S_{\beta,n_2}$, and the $N\times N$ matrices,
\begin{align}
A_x^{(m)}&=D_+^{(m)}(\mathbb{I}_{n_2}\otimes S_{\alpha,n_1})+D_-^{(m)}(\mathbb{I}_{n_2}\otimes S_{\alpha,n_1}^{\mathrm{T}}),\nonumber\\
A_y^{(m)}&=E_+^{(m)}( S_{\beta,n_2}\otimes \mathbb{I}_{n_1})+E_-^{(m)}( S_{\beta,n_2}^{\mathrm{T}}\otimes \mathbb{I}_{n_1}),\nonumber
\end{align}
where $\mathbb{I}_n$ is the identity matrix of size $n$, and $\otimes$ denotes the Kronecker product. Then, by using the Crank--Nicolson method, we obtain the system

\begin{align}
    \left(\frac{1}{r}\mathbb{I}_{N}+A_{x}^{(m)}+\frac{s}{r}A_{y}^{(m)}\right)\mathbf{u}^{(m)}=\left(\frac{1}{r}\mathbb{I}_{N}-A_{x}^{(m-1)}-\frac{s}{r}A_{y}^{(m-1)}\right)\mathbf{u}^{(m-1)}+2h_x^{\alpha}\mathbf{v}^{(m-1/2)},\nonumber
\end{align}
where $r=\frac{h_t}{2h_x^\alpha}$, $s=\frac{h_t}{2h_y^\beta}$. In compact form, the system is written
\begin{align}
\mathcal{M}_{(\alpha,\beta),N}^{(m)}\mathbf{u}^{(m)}=\mathbf{b}^{(m)}, \nonumber
\end{align}
where
\begin{align}
\mathcal{M}_{(\alpha,\beta),N}^{(m)}&=\frac{1}{r}\mathbb{I}_{N}+A_{x}^{(m)}+\frac{s}{r}A_{y}^{(m)}, \label{defin:2d coefficient matrix}\\
\mathbf{b}^{(m)}&=\left(\frac{1}{r}\mathbb{I}_{N}-A_{x}^{(m-1)}-\frac{s}{r}A_{y}^{(m-1)}\right)\mathbf{u}^{(m-1)}+2h_x^{\alpha}\mathbf{v}^{(m-1/2)}.\nonumber
\end{align}

\begin{prop}
\cite{moghaderi171} We suppose that $\frac{1}{r}=o(1)$ and $\frac{s}{r}=\frac{h_x^{\alpha}}{h_y^{\beta}}=O(1)$. We suppose also that for a given time $t_m$, the functions $d_+(x,y):=d_+(x,y,t_m)$, $d_-(x,y):=d_-(x,y,t_m)$, $e_+(x,y):=e_+(x,y,t_m)$, and $e_-(x,y):=e_-(x,y,t_m)$ are Riemman integrable over $[a_1,b_1]\times[a_2,b_2]$. Then, 
\[
\{ \mathcal{M}_{(\alpha,\beta),N}^{(m)} \}_N \sim_{GLT} \hat{h}_{(\alpha,\beta)}(\hat{\mathbf{x}},\boldsymbol{\theta}), \quad \hat{\mathbf{x}}=(\hat{x},\hat{y}), \: \boldsymbol{\theta}=(\theta_1,\theta_2), 
\]
where
\begin{align*}
&\hat{h}_{(\alpha,\beta)}(\hat{\mathbf{x}},\boldsymbol{\theta})=h_{(\alpha,\beta)}(a_1+(b_1-a_1)\hat{x},a_2+(b_2-a_2)\hat{y},\boldsymbol{\theta}), \\
&h_{(\alpha,\beta)}(x,y,\boldsymbol{\theta})=d_+(x,y)w_{\alpha}(\theta_1)+d_-(x,y)w_{\alpha}(-\theta_1)+\frac{s}{r} \left( e_+(x,y)w_{\beta}(\theta_2)+e_-(x,y)w_{\beta}(-\theta_2) \right), \\
&(\hat{\mathbf{x}},\boldsymbol{\theta})\in[0,1]^2 \times [-\pi,\pi]^2, \quad (x,y,\boldsymbol{\theta}) \in[a_1,b_1]\times[a_2,b_2]\times [-\pi,\pi]^2.
\end{align*}
Therefore 
\[
\{ \mathcal{M}_{(\alpha,\beta),N}^{(m)} \}_N \sim_{\sigma} h_{(\alpha,\beta)}(x,y,\boldsymbol{\theta}).
\]
In addition, if $d_+(x,y)=d_-(x,y)=e_+(x,y)=e_-(x,y)$ we have,
\[
\{ \mathcal{M}_{(\alpha,\beta),N}^{(m)} \}_N \sim_{\lambda} h_{(\alpha,\beta)}(x,y,\boldsymbol{\theta}).
\]
\end{prop}

\section{The $\tau$ Preconditioners}

\subsection{The $\tau$ Preconditioner in one Dimension}
\label{sec:prop1d}
In order that the results are directly combarable with that in  \cite{donatelli161}, in one dimension will be used the first order finite difference scheme. Also, to simplify the notation the time mark will be ommited. Let now,  $T_n=T_{\alpha,n1}$ as in \eqref{fractional:defination T_(a,N) matrix} and  $\mathcal{M}_{n}=\mathcal{M}_{\alpha,n}$ as in \eqref{fractional:equation 2dimension linear system}.
\par
As mentioned in the introduction of the chapter, the proposed preconditioner is a diagonal matrix $D_n$, times a $\tau$ matrix, $\mathcal{P}_{\mathcal{F}_\alpha}=D_n\tau_n(\mathcal{F}_{\alpha}(\boldsymbol{\theta}))$. Such a combination of two matrices as preconditioner is not a new proposal (\cite{NV-2008},\cite{NV2002},\cite{noutsos2005}). 

The form of the coefficient matrix of the system  $\mathcal{M}_n=\nu_{M,n}\mathbb{I}_{n}+D_+T_{n}+D_-T_{n}^{\mathrm{T}}$ suggests for the diagonal $D_n$ the following matrix, 
\begin{align}
D_n&=\frac{1}{2}\left(D_++D_-\right), \nonumber \\
[D_n]_{i,i}&= \frac{d_{+,i}+d_{-,i}}{2}, \label{Diagonal_1D_matrix}
\end{align}
that has been used in other preconditioning strategies also \cite{donatelli161}.
Assuming that the functions $d_{\pm}$ do not have a common zero $x_0 \in [L,R]$ we conclude that the $D_n^{-1}$ is uniformly bounded and
\begin{align}
D_n^{-1}\mathcal{M}_n=\nu_{M,n}D_n^{-1}+D_n^{-1}D_+T_{n}+D_n^{-1}D_-T_{n}^{\mathrm{T}}.\nonumber
\end{align}
If we now define $\delta(x)=\frac{d_+(x)}{d_+(x)+d_-(x)}$, $\delta_i=\delta(x_i)$, $\boldsymbol{\delta}=[\delta_1,\delta_2,\ldots,\delta_n]$,  $G_n=\mathrm{diag}(\boldsymbol{\delta})$ and taking into consideration that the $d_\pm$ are no-negative, we have that $0<\delta(x)<1$ and also,
\begin{align}
D_n^{-1}D_+&=2G_n, \nonumber\\
D_n^{-1}D_-&=2(\mathbb{I}_{n}-G_n).\nonumber
\end{align}

Hence, $D_n^{-1}\mathcal{M}_n$ can be written as
\begin{align}
D_n^{-1}\mathcal{M}_n&=\nu_{M,n}D_n^{-1}+D_n^{-1}D_+T_{n}+D_n^{-1}D_-T_{n}^{\mathrm{T}}\nonumber\\
&=\nu_{M,n}D_n^{-1}+2G_nT_{n}+2(\mathbb{I}_{n}-G_n)T_{n}^{\mathrm{T}} \nonumber \\
&=\nu_{M,n}D_n^{-1}+(T_{n}+T_{n}^{\mathrm{T}})+(2G_n-\mathbb{I}_{n})(T_n-T_{n}^{\mathrm{T}}). \nonumber
\end{align}
Since, from \eqref{fractional:first order symbols}, $T_n\coloneqq T_n(-e^{-\mathbf{i}\theta}\left(1-e^{\mathbf{i}\theta}\right)^\alpha)=T_n(g_{\alpha}(\theta))$ and  $T_n^{\mathrm{T}}\coloneqq T_n(-e^{\mathbf{i}\theta}\left(1-e^{-\mathbf{i}\theta}\right)^\alpha)=T_n(g_{\alpha}(-\theta))$ we have
\begin{align}
D_n^{-1}\mathcal{M}_n&=\nu_{M,n}D_n^{-1}+(T_{n}+T_{n}^{\mathrm{T}})+(2G_n-\mathbb{I}_{n})(T_n-T_{n}^{\mathrm{T}}) \nonumber \\
&=\nu_{M,n}D_n^{-1}+T_n(g_{\alpha}(\theta)+g_{\alpha}(-\theta))+(2G_n-\mathbb{I}_{n})T_n(g_{\alpha}(\theta)-g_{\alpha}(-\theta)) \nonumber  \\
&=\nu_{M,n}D_n^{-1}+T_n(p_{\alpha}(\theta))+(2G_n-\mathbb{I}_{n})T_n(2\mathbf{i}\Im\left\{g_{\alpha}(\theta)\right\}),   \label{DM_analytical_presentation}
\end{align}
where $p_{\alpha}(\theta)$, defined in \eqref{fractional:p_alpha defination}, is real, positive and even.
The above derivation of the $D_n^{-1}\mathcal{M}_n$ matrix is of interest since it makes clear why it is reasonable to use the $\tau$ preconditioner.
The first term of the above matrix, $\nu_{M,n}D_n^{-1}$, is  diagonal with positive and $o(1)$ entries, since we have supposed that the $d_{\pm}$ functions do not have zero at the same point in the domain $[L,R]$ and $\nu_{M,n}=o(1)$.
We mention here that although the entries are $o(1)$, its effect on the eigenvalues of the preconditioned matrix can be  significant. The reason is  explained in the end of this section.
The third term in \eqref{DM_analytical_presentation}  is a diagonal matrix with entries in $[-1,1]$ times a skew-symmetric Toeplitz matrix with generating function $2\mathbf{i}\Im\left\{g_{\alpha}(\theta)\right\},$ and  consequently purely imaginary eigenvalues.
If $d_+=d_-$ this term is vanishing  while if the $d_{\pm}$ are constant but not equal it is a pure skew-symmetric Toeplitz (in that case $(2G_n-\mathbb{I}_{n})=c\mathbb{I}_{n}$ for some constant $c$).

The term in \eqref{DM_analytical_presentation}, which is mainly responsible for the dispersion of the real part of the spectrum, is the second term, that is $T_n(p_{\alpha}(\theta))$. The  $\tau$ preconditioner will effectively  cluster the eigenvalues of this matrix, and consequently the eigenvalues of the whole matrix $D_n^{-1}\mathcal{M}_n$.
Hence, taking advantage of the `essential spectral equivalence' between the matrix sequences $\{\tau_n(f)\}_n$ and $\{T_n(f)\}_n$ proven in \cite{noutsos161}, we propose a preconditioner expressed as
\begin{align}
\mathcal{P}_{\mathcal{F}_\alpha,n}&=D_n\tau_n(p_{\alpha}(\theta))=D_n\mathbb{S}_nF_n\mathbb{S}_n,  \label{eq:prop_1D}
\end{align}
where

\begin{align}
F_n&=\mathrm{diag}(p_{\alpha}(\boldsymbol{\theta})),\qquad  \boldsymbol{ \theta}=\left[\theta_1, \theta_2,\ldots, \theta_n  \right], \qquad  \theta_j= \frac{j\pi}{n+1}=j\pi h, \qquad j=1,\ldots,n, \nonumber
\end{align}
with $D_n$ defined in \eqref{Diagonal_1D_matrix} and $\mathbb{S}_n$ being the sine transform matrix reported in \eqref{sine_transform}. Obviously,  the proposed preconditioner is symmetric and positive definite.


\subsubsection{Case I: $d_{\pm}$ are constants }
  In the case where the diffusion coefficient functions are constants, the \eqref{DM_analytical_presentation} becomes:
  \begin{align*}
  \left(2\frac{\nu_{M,n}}{d_++d_-}\right)\mathbb{I}_n+T_n\left(p_{\alpha}(\theta)\right)+&\left(\frac{d_+-d_-}{d_++d_-}\right)T_n\left(2\mathbf{i}\Im\left\{g_{\alpha}(\theta)\right\}\right)= \\
  &T_n\left(2\frac{\nu_{M,n}}{d_++d_-}+ p_{\alpha}(\theta)\right)+T_n\left(2\left(\frac{d_+-d_-}{d_++d_-}\right)\mathbf{i}\Im\left\{g_{\alpha}(\theta)\right\}\right),
  \end{align*}
  i.e,  is exactly the sum of a symmetric and a skew-symmetric Toeplitz matrix.  
  It is worth noticing that according to the GLT machinery, the term $\frac{2\cdot \nu_{M,n}}{d_++d_-}$ which is added to the symbol  of the first Toeplitz matrix sequence does not change the symbol of the sequence since is of order $o(1)$. However it  affects the speed  in which the minimum eigenvalue of the sequence approaches zero as the dimension of the matrix tends to infinity. Thus, in this special case,  the $\tau$ part of preconditioner  is defined as 
  \begin{align*}
  \tau_{M,n}\left(p_{\alpha}(\theta)+\frac{2\cdot \nu_{M,n}}{d_++d_-}\right)=\mathbb{S}_n \mathrm{diag}\left(p_{\alpha}(\theta)+\frac{2\cdot \nu_{M,n}}{d_++d_-}\right) \mathbb{S}_n=\mathbb{S}_n \hat{F}_n \mathbb{S}_n.
  \end{align*}
  Then,
  \begin{align*}
  \tau_{M,n}^{-1}\left(p_{\alpha}(\theta)+\frac{2\cdot \nu_{M,n}}{d_++d_-}\right) \left[T_n\left(\frac{2\cdot\nu_{M,n}}{d_++d_-}+ p_{\alpha}(\theta)\right)+T_n\left(2\frac{d_+-d_-}{d_++d_-}\mathbf{i}\Im\left\{g_{\alpha}(\theta)\right\}\right)\right] \sim \\
  \hat{F}_n^{-\frac{1}{2}}\mathbb{S}_n\left[ T_n\left(\frac{2\cdot\nu_{M,n}}{d_++d_-}+ p_{\alpha}(\theta)\right)+T_n\left(2\frac{d_+-d_-}{d_++d_-}\mathbf{i}\Im\left\{g_{\alpha}\theta)\right\} \right) \right]\mathbb{S}_n\hat{F}_n^{-\frac{1}{2}}=\\
  \hat{F}_n^{-\frac{1}{2}}\mathbb{S}_nT_n\left(\frac{2\cdot\nu_{M,n}}{d_++d_-}+ p_{\alpha}(\theta)\right)\mathbb{S}_n\hat{F}_n^{-\frac{1}{2}} + \hat{F}_n^{-\frac{1}{2}}\mathbb{S}_nT_n\left(2\frac{d_+-d_-}{d_++d_-}\mathbf{i}\Im\left\{g_{\alpha}(\theta)\right\}\right)\mathbb{S}_n\hat{F}_n^{-\frac{1}{2}}.
  \end{align*}
  The first term in the above sum  is symmetric and its eigenvalues are strongly clustered at 1 since  the conditions of  the main theoretical result of  \cite{noutsos161} are fulfilled concerning   the   spectral equivalence between  a  $\tau$ matrix and  a Toeplitz one. The second term is skew-symmetric and it does not  affect  the real part of the eigenvalues of the whole matrix. Moreover,  it is absent whenever  $d_+=d_-$. Hence, the real parts of the eigenvalues of the preconditioned matrix  are strongly  clustered around  1 and are bounded by constants $c, C$ with $0<c\leq 1\leq C <\infty$.

\subsubsection{Case II.  $d_{-}(x)=d_{+}(x) >0$  }
\label{sec:diffequal}
 In this case,  the term  $2G_n-\mathbb{I}_{n}=\mathbf{0}$ in  \eqref{DM_analytical_presentation} is equal to zero and the preconditioned matrix becomes $\tau_n^{-1}(p_{\alpha}(\theta))(\nu_{M,n}D_n^{-1}+T_n(p_{\alpha}(\theta)))$ which is similar to the SPD
\begin{align}
\tau_n^{-1}(p_{\alpha}(\theta))(\nu_{M,n}D_n^{-1}+T_n(p_{\alpha}(\theta)))&\sim F_n^{-1/2}\mathbb{S}_n(\nu_{M,n}D_n^{-1}+T_n(p_{\alpha}(\theta)))\mathbb{S}_nF_n^{-1/2}\nonumber\\
&=\nu_{M,n}F_n^{-1/2}\mathbb{S}_n(D_n^{-1})\mathbb{S}_nF_n^{-1/2}+F_n^{-1/2}\mathbb{S}_n(T_n(p_{\alpha}(\theta)))\mathbb{S}_nF_n^{-1/2}. \label{eq:main:precond1d}
\end{align}
In the above splitting in positive symmetric terms,   the first  one has ${\it o}(n)$ eigenvalues  tending  to infinity while the second one fulfills the main theoretical result of \cite{noutsos161}  and thus, for every $n$, it has eigenvalues belonging to an interval  [c,C] with $c,C$ constants and  $0<c\leq 1\leq C <\infty$. The claim about the spectrum of the first term can be  proven if we equivalently show that the inverse  of it, i.e. $ F_n (\mathbb{S}_n D_n  \mathbb{S}_n)$ has at most ${\it o}(n)$ eigenvalues tending to 0 as $n\rightarrow \infty$. Since $F_n$ is the diagonal matrix formed by the   values  $p_{\alpha}(j\pi h), \quad j=1,\ldots, n,$ which has a unique zero at zero of order $\alpha$,   there will be an index $\hat{j}$ with $\hat{j}$ of order ${\it o}(n)$ such that $p_{\alpha}(j\pi h)$  being  of order ${\it o}(1)$  for all $j\leq \hat{j}$.  Thus, at most ${\it o}(n)$  eigenvalues of $F_n $  can tend to zero. Using  Rayleigh quotient and taking into account that the matrix $D_n$ is a diagonal matrix with entries bounded from above end below by positive universal constants, the claim is proven.   Consequently, using the Weyl's theorem on (\ref{eq:main:precond1d}) we have that  
\begin{align*}
\lambda_{k}\left(\nu_{M,n}F_n^{-1}\mathbb{S}_n(D_n^{-1})\mathbb{S}_n +F_n^{-1}\mathbb{S}_n(T_n(p_{\alpha}(\theta)))\mathbb{S}_n \right) \leq  \\ 
\nu_{M,n}\lambda_{k}( F_n^{-1}\mathbb{S}_n(D_n^{-1})\mathbb{S}_n) +\lambda_n \left(F_n^{-1}\mathbb{S}_n(T_n(p_{\alpha}(\theta)))\mathbb{S}_n\right). 
\end{align*}

Accordingly, at most $\it o(n)$  eigenvalues of $\tau_n^{-1}(p_{\alpha}(\theta))(\nu_{M,n}D_n^{-1}+T_n(p_{\alpha}(\theta)))$ can tend to infinity. Clearly the term $\nu_{M,n}$ which in general tends to zero as ${\it O}(n^{1-\alpha})$, can further  reduce the number of eigenvalues tending to infinity. 

In the semi elliptic case (see \cite{NSV08} and especially the numerical experiments therein), if the equal functions $d_{\pm}$ have a root then an unpredictable asymptotical behavior of the eigenvalues of  coefficient   matrix $\cal{M_\alpha}$ is expected.

\subsubsection{Case III: General case }
In the case where $d_+\neq d_-$ the term $(2G_n-\mathbb{I}_{n})(T_n-T_{n}^{\mathrm{T}})$ is nonzero and it affects the spectrum of the preconditioned matrix. Specifically,
\begin{align}
&\tau_n^{-1}(p_{\alpha}(\theta))(\nu_{M,n}D_n^{-1}+T_n(p_{\alpha}(\theta))+(2G_n-\mathbb{I}_{n})T_n(2\mathbf{i}\Im\left\{g_{\alpha}(\theta)\right\}))\nonumber\\
&\hspace{1cm}\sim F_n^{-1/2}\mathbb{S}_n(\nu_{M,n}D_n^{-1}+T_n(p_{\alpha}(\theta))+(2G_n-\mathbb{I}_{n})T_n(2\mathbf{i}\Im\left\{g_{\alpha}(\theta)\right\}))\mathbb{S}_nF_n^{-1/2}\nonumber\\
&\hspace{1cm}=F_n^{-1/2}\mathbb{S}_n(\nu_{M,n}D_n^{-1})\mathbb{S}_nF_n^{-1/2}+F_n^{-1/2}\mathbb{S}_n(T_n(p_{\alpha}(\theta)))\mathbb{S}_nF_n^{-1/2}+ \nonumber \\
&\hspace{4cm}F_n^{-1/2}\mathbb{S}_n(2G_n-\mathbb{I}_{n})T_n(2\mathbf{i}\Im\left\{g_{\alpha}(\theta)\right\})\mathbb{S}_nF_n^{-1/2}, \nonumber
\end{align}
where  only the, new,  third term can add imaginary quantity on the eigenvalues.  However, through experimentation it can be shown,   that the effect of this third  term on the real part of the eigenvalues is negligible. In this sense, all  the  numerical experiments given in  Section 4  belong to this case mainly for  showing   the performance of the proposal  there were  the spectral analysis  do not explicitly and  in depth  cover the topic.

\subsection{Proposed Preconditioner: Two Dimensions}
In the two-dimensional case the second order spatial discretization is used, in order to be consistent with \cite{moghaderi171} and be able to readily compare the results. In this case, as reported in Section  \ref{sec:introduction:fde2d}, the coefficient matrix of the system is defined as
\begin{align}
\mathcal{M}_{(\alpha,\beta),N}^{(m)}&=\frac{1}{r}\mathbb{I}_{N}+D_+^{(m)}(\mathbb{I}_{n_2}\otimes S_{\alpha,n_1})+D_-^{(m)}(\mathbb{I}_{n_2}\otimes S_{\alpha,n_1}^{\mathrm{T}})+\frac{s}{r}\left(E_+^{(m)}( S_{\beta,n_2}\otimes \mathbb{I}_{n_1})+E_-^{(m)}( S_{\beta,n_2}^{\mathrm{T}}\otimes \mathbb{I}_{n_1})\right). \label{coef:2d}
\end{align}
 It is reminded that $S_{\alpha,n_1}=T_{n_1}(w_\alpha(\theta))$ and $S_{\beta,n_2}=T_{n_2}(w_\beta(\theta))$.
Again, for simplicity the time dependency in the notation is omitted.

Now let  $\mathcal{F}_{(\alpha,\beta)}(\theta_1,\theta_2)=q_\alpha(\theta_1)+\frac{s}{r}q_\beta(\theta_2)$ where $q$ is the real, nonnegative  and even function defined in \eqref{second_order_symbol}, $\theta_1, \theta_2 \in [-\pi,\pi],$ and $n_1$, $n_2$ the two integers   used for the discretization of  the domain $[L_x, R_x] \times [L_y, R_y]$.
Using the grid in \eqref{tau_preconditioner} we define the diagonal matrices
\begin{align}
F_{n_1,j}=&\mathrm{diag}(\mathcal{F}_{(\alpha,\beta)}(\theta_{i,n_1},\theta_{j,n_2}),  i=1,\ldots,n_1),\quad \nonumber
\end{align}
for each $j=1,\ldots,n_2$. Then, the $N\times N$ diagonal matrix is expressed as
\begin{align}
F_{N}=&
\begin{bmatrix*}
F_{n_1,1}&\\
&F_{n_1,2}&&\\
&&\ddots&\\
&&&\ddots&\\
&&&&F_{n_1,n_2}
\end{bmatrix*}. \label{F_N_matrix}
\end{align}
Let $\mathbb{S}_{n_1}$ and $\mathbb{S}_{n_2}$ be the discrete sine transform matrices of sizes $n_1$ and $n_2$, respectively,  as they defined in \eqref{sine_transform}.
Then, generalizing the idea of  (\ref{eq:prop_1D}),  the proposed preconditioner for this case is
\begin{align}\mathcal{P}_{\mathcal{F}_{(\alpha,\beta)},N}&=D_N\left(\mathbb{S}_{n_2} \otimes \mathbb{S}_{n_1}\right)F_N\left(\mathbb{S}_{n_2} \otimes \mathbb{S}_{n_1}\right), \label{2D_preconditioner}
\end{align}
where
\begin{align}
D_N&=(D_++D_-+E_++E_-)/4. \nonumber
\end{align}

The motivation of the above construction is to create a preconditioner that properly acts on the  different sources  affecting the spectrum of $\mathcal{M}_{(\alpha,\beta),N}.$ Specifically, the  diagonal part  operates   on the spatial space treating  the influence that  the coefficients of the equation have on the matrix, while the $\tau$ matrix  focuses  on the spectral space  and the ill-conditioning generated by the discretization  of the  fractional differential  operator. This observation  is a direct result of the GLT symbol associated to $\mathcal{M}_{(\alpha,\beta),N}$ and has been extensively studied  in \cite{NSV08} and \cite{V18}, for the case of semi elliptic differential equations. Moreover, the spectral  analysis of the preconditioned Conjugate Gradient in 2 dimension is considered in \cite{noutsos2006}.
In the simplest, but not unusual  in applications, case where $d_{\pm}=d$, $e_{\pm}=e$ we can counterbalance  the influence of the term $\frac{1}{r}$ in the spectrum of $\mathcal{M}_{(\alpha,\beta),N}$  incorporate  it  into  the $\tau$ part of the preconditioner.  Particularly,  we define  $\mathcal{\hat{F}}_{(\alpha,\beta)}(\theta_1,\theta_2)=\frac{1}{r}+d\cdot q_\alpha(\theta_1)+\frac{s}{r}e\cdot q_\beta(\theta_2)$ 
 replacing  the sampling of $\mathcal{F}_{(\alpha,\beta)}$ with that of $\mathcal{\hat{F}}_{(\alpha,\beta)}$ for the construction of $\hat{F}_N$ instead of $F_N$ in \eqref{F_N_matrix}.  Accordingly,  the new corresponding preconditioner $\mathcal{P}_{\mathcal{\hat{F}}_{(\alpha,\beta)},N}$  is defined as
 
\begin{align}
\mathcal{P}_{\mathcal{\hat{F}}_{(\alpha,\beta)},N}&=\left(\mathbb{S}_{n_2} \otimes \mathbb{S}_{n_1}\right)\hat{F}_N\left(\mathbb{S}_{n_2} \otimes \mathbb{S}_{n_1}\right). \label{2D_preconditioner_proof_case}
\end{align}

The following theorem shows that in this case,  the spectrum of the preconditioned matrix is bounded by positive constants independent of the size of the matrix.

\begin{thm}\label{prop:main}
Assume that $d_{\pm}=d>0$, $e_{\pm}=e>0$. In this case the coefficient matrix of the system becomes

\begin{align}
A_N=\frac{1}{r}\mathbb{I}_N+(\mathbb{I}_{n_2} \otimes \hat{A}_{n_1}^\alpha)+(A_{n_2}^\beta \otimes \mathbb{I}_{n_1})=\left(\mathbb{I}_{n_2} \otimes (\frac{1}{r}\mathbb{I}_{n_1}+\hat{A}_{n_1}^\alpha)\right)+(A_{n_2}^\beta \otimes \mathbb{I}_{n_1})=\mathbb{I}_{n_2} \otimes A_{n_1}^\alpha+A_{n_2}^\beta \otimes \mathbb{I}_{n_1},
\end{align}

where
\begin{align}
A_{n_1}^\alpha&=\frac{1}{r}\mathbb{I}_N+T_{n_1}\left(d\cdot (w_\alpha(\theta)+w_\alpha(-\theta))\right)=T_{n_1}\left( \frac{1}{r}+d\cdot q_\alpha(\theta)\right),
\nonumber\\
A_{n_2}^\beta&= T_{n_2}\left(e\frac{s}{r}\cdot (w_\beta(\theta)+w_\beta(-\theta))\right)=T_{n_2}\left(e\frac{s}{r}\cdot q_\beta(\theta)\right). \nonumber
\end{align}

Then, the spectrum of the preconditioned matrix sequence $\left\{\mathcal{P}_{\mathcal{\hat{F}}_{(\alpha,\beta)},N}^{-1}A_N\right\}_N$ is bounded by  positive constants $c,C$ independent of  $N$.

\end{thm}

\begin{prf*}
We have that
\begin{align}
h_x=(R_x-L_x)h_1,&\quad
h_y=(R_y-L_y)h_2,\nonumber\\
r=\frac{h_t}{2h_x^\alpha},&\quad s=\frac{h_t}{2h_y^\beta}, \nonumber
\end{align}
and 
\begin{align}
\hat{F}_N&=\mathbb{I}_{n_2}\otimes F_{n_1}^\alpha+F_{n_2}^\beta\otimes \mathbb{I}_{n_1}, \label{D}
\end{align}
where $\mathbb{I}_n$ is the identity matrix of order $n$ and

\begin{align}
F_{n_1}^\alpha&=\mathrm{diag}(d\cdot \mathcal{F}_\alpha(\theta_{i,n_1})+\frac{1}{r}),\quad i=1,\ldots,n_1, \label{h1} \\
F_{n_2}^\beta&=\mathrm{diag}(e\frac{s}{r}\cdot\mathcal{F}_\beta(\theta_{j,n_2})),\quad j=1,\ldots,n_2.\label{h2}
\end{align}
The matrix $A_N$ is  SPD since each of its terms is  a Kronecker product of a diagonal with a SPD Toeplitz matrix.
Hence,
\begin{align}
\mathcal{P}_{N}^{-1}A_{N} = \left(\mathbb{S}_{n_2} \otimes \mathbb{S}_{n_1}\right)\hat{F}_N^{-1}\left(\mathbb{S}_{n_2} \otimes \mathbb{S}_{n_1}\right)A_N,\nonumber
\end{align}
which is similar to the matrix
\begin{align}
 \hat{F}_N^{-1/2}\left(\mathbb{S}_{n_2} \otimes \mathbb{S}_{n_1}\right)A_N \left(\mathbb{S}_{n_2} \otimes \mathbb{S}_{n_1}\right)\hat{F}_N^{-1/2}.\nonumber
\end{align}
Thus,

\begin{align}
&\hat{F}_N^{-1/2}(S_{n_2} \otimes \mathbb{S}_{n_1})\left((\mathbb{I}_{n_2} \otimes A_{n_1}^\alpha) +(A_{n_2}^\beta \otimes \mathbb{I}_{n_1})\right)(\mathbb{S}_{n_2} \otimes \mathbb{S}_{n_1})\hat{F}_N^{-1/2}\nonumber\\
&\hspace{1cm}=\hat{F}_N^{-1/2}\left((\mathbb{S}_{n_2} \otimes \mathbb{S}_{n_1})(\mathbb{I}_{n_2} \otimes A_{n_1}^\alpha)(\mathbb{S}_{n_2} \otimes \mathbb{S}_{n_1})+(\mathbb{S}_{n_2} \otimes \mathbb{S}_{n_1})(A_{n2}^\beta \otimes \mathbb{I}_{n_1})(\mathbb{S}_{n_2} \otimes \mathbb{S}_{n_1})\right)\hat{F}_N^{-1/2}\nonumber\\
&\hspace{1cm}=\hat{F}_N^{-1/2}\left(\mathbb{I}_{n_2}\otimes \mathbb{S}_{n_1}A_{n_1}^\alpha S_{n_1} + \mathbb{S}_{n_2}A_{n_2}^\beta \mathbb{S}_{n_2} \otimes \mathbb{I}_{n_1}\right)\hat{F}_N^{-1/2}\nonumber\\
&\hspace{1cm}=\hat{F}_N^{-1/2}\left(\mathbb{I}_{n_2}\otimes (F_{n_1}^\alpha)^{1/2} (F_{n_1}^\alpha)^{-1/2} \mathbb{S}_{n_1}A_{n_1}^\alpha \mathbb{S}_{n_1}(F_{n_1}^\alpha)^{-1/2} (F_{n_1}^\alpha)^{1/2}\right. +\nonumber\\
&\hspace{2.13cm} +\left.(F_{n_2}^\beta)^{1/2} (F_{n_2}^\beta)^{-1/2}\mathbb{S}_{n_2}A_{n_2}^\beta \mathbb{S}_{n_2}(F_{n_2}^\beta)^{-1/2} (F_{n_2}^\beta)^{1/2} \otimes \mathbb{I}_{n_1}\vphantom{\frac{1}{r}}\right)\hat{F}_N^{-1/2}\nonumber\\
&\hspace{1cm}=\hat{F}_N^{-1/2}\left((\mathbb{I}_{n_2}\otimes (F_{n_1}^\alpha)^{1/2})\underbrace{(\mathbb{I}_{n_2}\otimes  (F_{n_1}^\alpha)^{-1/2} \mathbb{S}_{n_1}A_{n_1}^\alpha \mathbb{S}_{n_1}(F_{n_1}^\alpha)^{-1/2}}_{=L})(\mathbb{I}_{n_2}\otimes (F_{n_1}^\alpha)^{1/2})\right. + \nonumber\\
&\hspace{2.13cm} +\left.((F_{n_2}^\beta)^{1/2} \otimes \mathbb{I}_{n_1})\underbrace{((F_{n_2}^\beta)^{-1/2}\mathbb{S}_{n_2}A_{n_2}^\beta \mathbb{S}_{n_2}(F_{n_2}^\beta)^{-1/2}) \otimes \mathbb{I}_{n_1})}_{=R}((F_{n_2}^\beta)^{1/2} \otimes \mathbb{I}_{n_1})\vphantom{\frac{1}{r}}\right)\hat{F}_N^{-1/2}\nonumber\\
&\hspace{1cm}=\underbrace{\hat{F}_N^{-1/2}(\mathbb{I}_{n_2}\otimes (F_{n_1}^\alpha)^{1/2}) L (\mathbb{I}_{n_2}\otimes (F_{n_1}^\alpha)^{1/2})\hat{F}_N^{-1/2}}_{=A_L}+\underbrace{\hat{F}_N^{-1/2}((F_{n_2}^\beta)^{1/2} \otimes \mathbb{I}_{n_1})R((F_{n_2}^\beta)^{1/2} \otimes \mathbb{I}_{n_1})\hat{F}_N^{-1/2}}_{=A_R}. \label{precMatrix}
\end{align}

Let
\begin{align}
P_{n_1}^{\alpha}&=\mathbb{S}_{n_1}F_{n_1}^\alpha \mathbb{S}_{n_1},\nonumber\\
P_{n_2}^{\beta} &=\mathbb{S}_{n_2}F_{n_2}^\beta \mathbb{S}_{n_2}. \nonumber
\end{align}
Then, (see \cite{noutsos161}),  there exist positive  constants $c$ and $C$ independent of $n_1,n_2$, such that 
\begin{displaymath}
c<\sigma\left(\left(P_{n_1}^{\alpha}\right)^{-1}A_{n_1}^\alpha\right)<C \Rightarrow c<\sigma\left((F_{n_1}^\alpha)^{-1/2} \mathbb{S}_{n_1}A_{n_1}^\alpha \mathbb{S}_{n_1}(F_{n_1}^\alpha)^{-1/2}\right)<C
\end{displaymath}
 and 
 \begin{displaymath}
 c<\sigma\left(\left(P_{n_2}^{\beta}\right)^{-1}A_{n_2}^\beta\right)<C\Rightarrow c<(F_{n_2}^\beta)^{-1/2}\mathbb{S}_{n_2}A_{n_2}^\beta \mathbb{S}_{n_2}(F_{n_2}^\beta)^{-1/2}<C.
\end{displaymath} 
  Consequently, for every normalized vector $x \in \mathbb{R}^N$ we find that: 
\begin{align}
c<x^{\mathrm{T}} L x<C, \qquad  c< x^{\mathrm{T}} R x<C. \nonumber
\end{align}
Since the  matrices $A_L,$ $A_R$ that form \eqref{precMatrix}  are SPD, some properties concerning such kind of matrices are used here.  Specifically, the inequality $A>B$ for $A,B$ SPD matrices if $A-B >0$  is positive definite is used and in addition  if $A$, $B$, $C$, $D$, and $E$ are SPD, then, by the Sylvester inertia law and by the definition of Rayleigh quotient
\begin{align}
A>B &\Leftrightarrow EAE>EBE, \label{eq:spd:1}\\
A>B \text{ and } C>D &\Leftrightarrow A+C>B+D.\label{eq:spd:2}
\end{align}

Therefore, we infer
\begin{align}
&\begin{cases}
c\mathbb{I}_N<L<C\mathbb{I}_N, \\
c\mathbb{I}_N<R<C\mathbb{I}_N,
\end{cases}\nonumber
\end{align}
and, using \eqref{eq:spd:1} and \eqref{eq:spd:2}, we deduce
\begin{align}
&\begin{cases}
 c\hat{F}_N^{-1}(\mathbb{I}_{n_2}\otimes F_{n_1}^\alpha)<A_L<C\hat{F}_N^{-1}(\mathbb{I}_{n_2}\otimes F_{n_1}^\alpha),\\
 c\hat{F}_N^{-1}(F_{n_2}^\beta \otimes \mathbb{I}_{n_1})<A_R<C\hat{F}_N^{-1}(F_{n_2}^\beta \otimes \mathbb{I}_{n_1}).
 \end{cases}\label{threeinequalities}
 \end{align}
Using again \eqref{eq:spd:1} and \eqref{eq:spd:2}, taking into account the two inequalities of \eqref{threeinequalities}, and \eqref{D}, we have
 \begin{align}
 c\hat{F}_N^{-1}(\mathbb{I}_{n_2}\otimes F_{n_1}^\alpha)+c\hat{F}_N^{-1}(F_{n_2}^\beta \otimes \mathbb{I}_{n_1}) &= c \hat{F}_N^{-1}\hat{F}_N = c\mathbb{I}_N,\nonumber\\
C\hat{F}_N^{-1}(\mathbb{I}_{n_2}\otimes F_{n_1}^\alpha)+C\hat{F}_N^{-1}(F_{n_2}^\beta \otimes \mathbb{I}_{n_1}) &= C \hat{F}_N^{-1}\hat{F}_N=C\mathbb{I}_N. \nonumber
\end{align}
Consequently, we conclude that 
\begin{displaymath} c\mathbb{I}_N\leq F_N^{-1/2}(\mathbb{S}_{n_1} \otimes \mathbb{S}_{n_2})A_N(\mathbb{S}_{n_1} \otimes \mathbb{S}_{n_2})F_N^{-1/2}\leq C \mathbb{I}_N. \end{displaymath}
 Therefore, the spectrum of the preconditioned matrix, which is similar to the $F_N^{-1/2}(\mathbb{S}_{n_1} \otimes \mathbb{S}_{n_2})A_N(\mathbb{S}_{n_1} \otimes \mathbb{S}_{n_2})F_N^{-1/2}$, lies in $[c,C]$. Moreover, from \cite{noutsos161} we expect all the eigenvalues to be clustered  around 1, something that is  numerically confirmed in the next section. 

 \end{prf*}
 

 \begin{cor}
Let the functions $d_{+}(x,y,t),d_{-}(x,y,t),e_{+}(x,y,t),e_{-}(x,y,t)$ being  strictly positive functions on $\Omega$, with $d_{+}(x,y,t)=d_{-}(x,y,t)=e_{+}(x,y,t)=e_{-}(x,y,t).$  Then,  the preconditioned matrix  sequence $\left\{\mathcal{P}_{\mathcal{\hat{F}}_{(\alpha,\beta)},N}^{-1} \mathcal{M}_{(\alpha,\beta),N}^{(m)} \right\}_N$ is bounded by  positive constants $c,C$ independent of  $N$.
\end{cor}

\begin{prf*}

The  proof  can be easily obtained  from the results of Theorem \ref{prop:main}  and the observation that the coefficient matrix in (\ref{coef:2d}) can be bounded by  

\begin{displaymath}
A^{c}_{N}\leq \mathcal{M}_{(\alpha,\beta),N}^{(m)} \leq A^{C}_N,
\end{displaymath}

where 

 \begin{align}
A_{N}^c&=\frac{1}{r}\mathbb{I}_{N}+c(\mathbb{I}_{n_2}\otimes S_{\alpha,n_1})+c(\mathbb{I}_{n_2}\otimes S_{\alpha,n_1}^{\mathrm{T}})+\frac{s\cdot c}{r}\left( (S_{\beta,n_2}\otimes \mathbb{I}_{n_1})+( S_{\beta,n_2}^{\mathrm{T}}\otimes \mathbb{I}_{n_1})\right), \nonumber
\end{align}
 \begin{align}
A_{N}^{C} &=\frac{1}{r}\mathbb{I}_{N}+C(\mathbb{I}_{n_2}\otimes S_{\alpha,n_1})+C (\mathbb{I}_{n_2}\otimes S_{\alpha,n_1}^{\mathrm{T}})+\frac{s\cdot C}{r}\left ( (S_{\beta,n_2}\otimes \mathbb{I}_{n_1})+( S_{\beta,n_2}^{\mathrm{T}}\otimes \mathbb{I}_{n_1})\right), \nonumber
\end{align}

and 
 \begin{displaymath}
c =\min_{(x,y,t)\in \Omega}\{d_{+}(x,y,t),d_{-}(x,y,t),e_{+}(x,y,t),e_{-}(x,y,t)\}, 
\end{displaymath}
\begin{displaymath}
 C=\max_{(x,y,t)\in \Omega}\{d_{+}(x,y,t),d_{-}(x,y,t), e_{+}(x,y,t),e_{-}(x,y,t)\}.
  \end{displaymath} 

Then, using Rayleigh quotient we obtain

\begin{displaymath}
\mathcal{P}_{\mathcal{\hat{F}}_N}^{-1} A^{c}_{N}  \leq \mathcal{P}_{\mathcal{\hat{F}}_N}^{-1} \mathcal{M}_{(\alpha,\beta),N}^{(m)} \leq \mathcal{P}_{\mathcal{\hat{F}}_N}^{-1} A^{C}_N
\end{displaymath}


\begin{displaymath}
\lambda_1(\mathcal{P}_{\mathcal{\hat{F}}_N}^{-1} A^{c}_{N})  \leq \lambda_1(\mathcal{P}_{\mathcal{\hat{F}}_N}^{-1} \mathcal{M}_{(\alpha,\beta),N}^{(m)})\leq \lambda_N(\mathcal{P}_{\mathcal{\hat{F}}_N}^{-1} \mathcal{M}_{(\alpha,\beta),N}^{(m)}) \leq \lambda_N(\mathcal{P}_{\mathcal{\hat{F}}_N}^{-1} A^{C}_N),
\end{displaymath}

and the proof is completed.
\end{prf*}

\section{Numerical Examples}
\label{sec:numerical}
In this section we present three numerical examples to show the efficiency of the proposed preconditioners, compared with preconditioners discussed in~\cite{donatelli161} (one dimension) and \cite{moghaderi171} (two dimensions).

\begin{itemize}
\item Example 1 is a one-dimensional problem, taken from~\cite[Example 1.]{donatelli161}, and we compare and discuss the preconditioners therein with the proposed $\mathcal{P}_{\mathcal{F}_\alpha,n}$, and a few variations based on the spectral symbol. The fractional derivatives are of order $\alpha\in \{1.2,1.5,1.8\}$.

\item Example 2 is a two-dimensional problem, taken from~\cite[Example 1.]{moghaderi171}, and we compare and discuss the preconditioners therein with the proposed $\mathcal{P}_{\mathcal{F}_{(\alpha,\beta)},N}$. The fractional derivatives are $\alpha=1.8$ and $\beta=1.6$.

\item Example 3 is the same experiment as Example 2, but with the fractional derivatives $\alpha=1.8$ and $\beta=1.2$.
\end{itemize}
The numerical experiments presented in Tables~\ref{tbl:table1}--\ref{tbl:table4} were implemented in \textsc{Julia} v1.1.0, using GMRES from the package \textsc{IterativeSolvers.jl} (GMRES tolerance is set to $10^{-7}$) and the \textsc{FFTW.jl} package. Benchmarking is done with \textsc{BenchmarkTools.jl} with 100 samplings and minimum time is presented in milliseconds. Experiments were run, in serial, on a computer with dual Intel Xeon E5 2630 v4 2.20 GHz (10 cores each) cpus, and with 128 GB of RAM.

The Figures~\ref{fig:1a}--\ref{fig:1c} (and Figures~\ref{fig:2} and \ref{fig:3}) show the scaled spectra of the preconditioned coefficient matrix $\mathcal{P}^{-1}\mathcal{M}_{\alpha,n_1}$ (and $\mathcal{P}^{-1}\mathcal{M}_{(\alpha,\beta),N}$) for different preconditioners $\mathcal{P}$, fractional derivatives $\alpha$, and matrix orders $n_1$ (and $\beta$, $N=n_1,n_2$). The scaling by a constant $c_0$ is performed the following way: find the smallest enclosing circle over all the eigenvalues of the matrix of interest $A$. The center is denoted $c_0$ and the radius is $r$. Then, the spectrum is scaled as $\lambda_j(A)/c_0$ and the circle scaled and centered in $(1,0)$. The Julia package \textsc{BoundingSphere.jl} was used to compute $c_0$ and $r$ for all figures. The current scaling of the eigenvalues of preconditioned coefficient matrices is a visualization of the important effect  for the convergence rate of GMRES  of both the clustering and of the shape of the clustering.

In Tables~\ref{tbl:table1}--\ref{tbl:table4}, for each preconditioner, we present the number of iterations [it], minimal timing [ms], and the condition number of the preconditioned  matrix $\kappa$.  Best results are highlighted in bold.

\subsection{Example 1}
We compare the proposed preconditioner $\mathcal{P}_{\mathcal{F}_\alpha,n}$ with the ones presented in Example 1 from~\cite{donatelli161} (and two alternative symbol based preconditioners).
We consider the one-dimensional form of~\eqref{eq:fde} in  the domain $[L_1,R_1]\times[t_0,T]= [0,2]\times [0,1]$,  where  the diffusion coefficients 
\begin{align}
d_+(x)&=\Gamma(3-\alpha)x^\alpha,\nonumber\\
d_-(x)&=\Gamma(3-\alpha)(2-x)^\alpha,\nonumber
\end{align}
are non-constant in space. Furthermore, the source term is
\begin{align}
f(x,t)=-32e^{-t}\left(x^2+\frac{(2-x)^2(8+x^2)}{8}-\frac{3(x^3+(2-x)^3)}{3-\alpha}+\frac{3(x^4+(2-x)^4)}{(4-\alpha)(3-\alpha)}\right), \nonumber
\end{align}
and the initial condition is
\begin{align}
u(x,0)=4x^2(2-x)^2,\nonumber
\end{align}
which yield an analytical solution $u(x,t)=4e^{-t}x^2(2-x)^2$. We assume $h_x=h_t=2/(n_1+1)$, that is, $\nu_{M,n_1}=h_x^{\alpha-1}$ and number of time steps $M=(n_1+1)T/(R_1-L_1)=(n_1+1)/2$.
The set of fractional derivatives $\alpha$, for which a solution is computed for, is $\{1.2,1.5,1.8\}$ and in addition we consider the following set of partial dimensions for $n_1$, that is $\{2^6-1,2^7-1,2^8-1,2^9-1\}$.

In Table~\ref{tbl:table1} we present the results for the following preconditioners
\begin{itemize}
\item Identity ($\mathbb{I}_{n_1}$): GMRES without any preconditioner.
\item Circulant ($\mathcal{P}_{C,n_1}$): Described in~\cite{lei131} and implemented using FFT.
\item ``Full'' symbol ($\mathcal{P}_{\textsc{full},n_1}$): Defined as 
\[
\mathbb{S}_{n_1}\mathrm{diag}\left(\nu_{M,n_1}+d_{+,i}g_\alpha(\theta_{j,n_1})+d_{-,i}g_\alpha(-\theta_{j,n_1}), ~j=1,2,\ldots n_1\right)\mathbb{S}_{n_1}
\] and implemented using FFT.
\item Symbol ($\mathcal{P}_{\mathcal{F}_\alpha,n_1}$): Proposed in Section~\ref{sec:prop1d}, $D_{n_1}\mathbb{S}_{n_1}\mathrm{diag}\left( p_\alpha(\theta_{j,n_1}), j=1,2,\ldots n_1\right)\mathbb{S}_{n_1}$, and implemented using FFT.
\end{itemize}

\begin{figure}[!ht]
\centering
\includegraphics[width=0.32\textwidth]{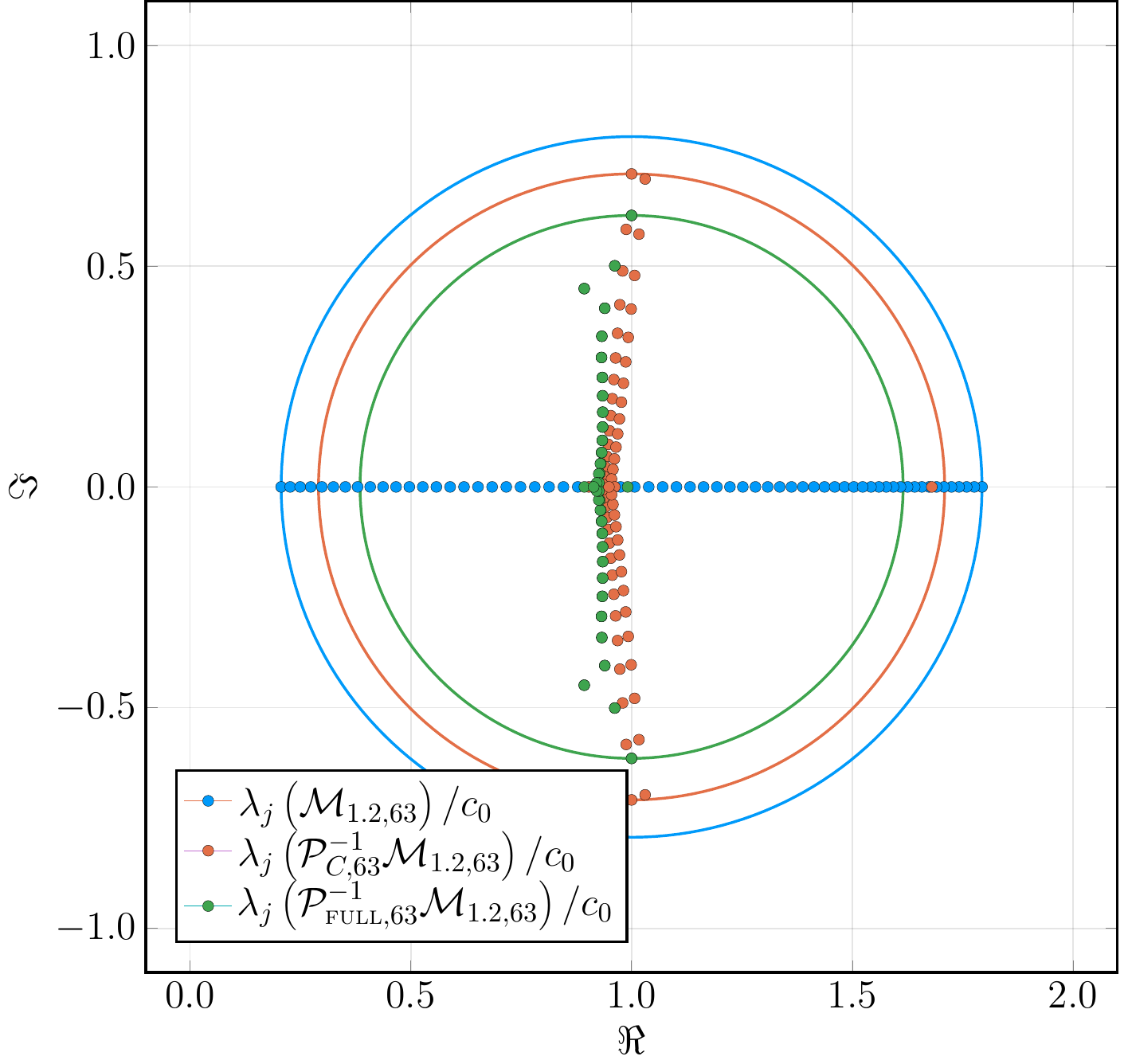}
\includegraphics[width=0.32\textwidth]{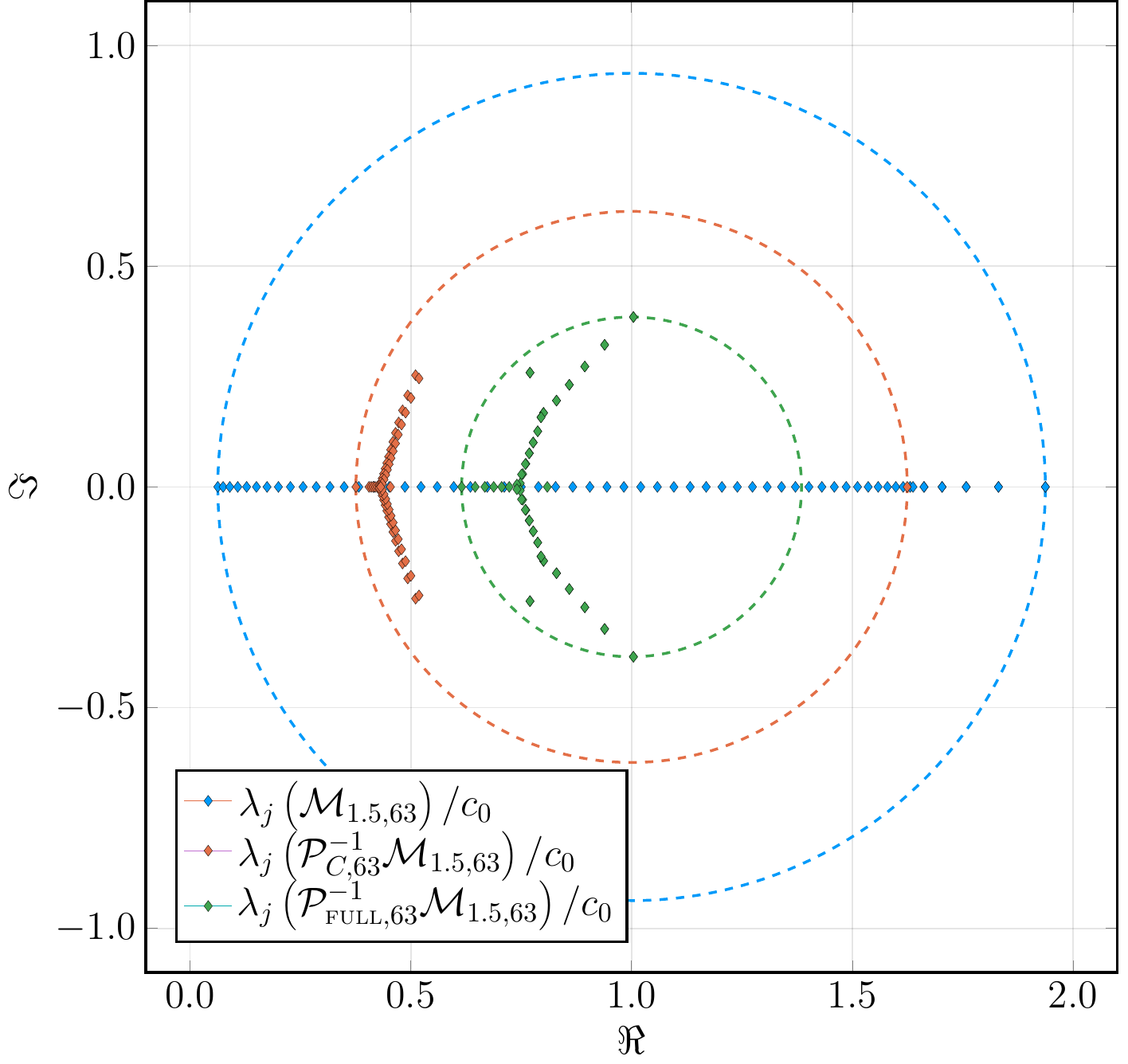}
\includegraphics[width=0.32\textwidth]{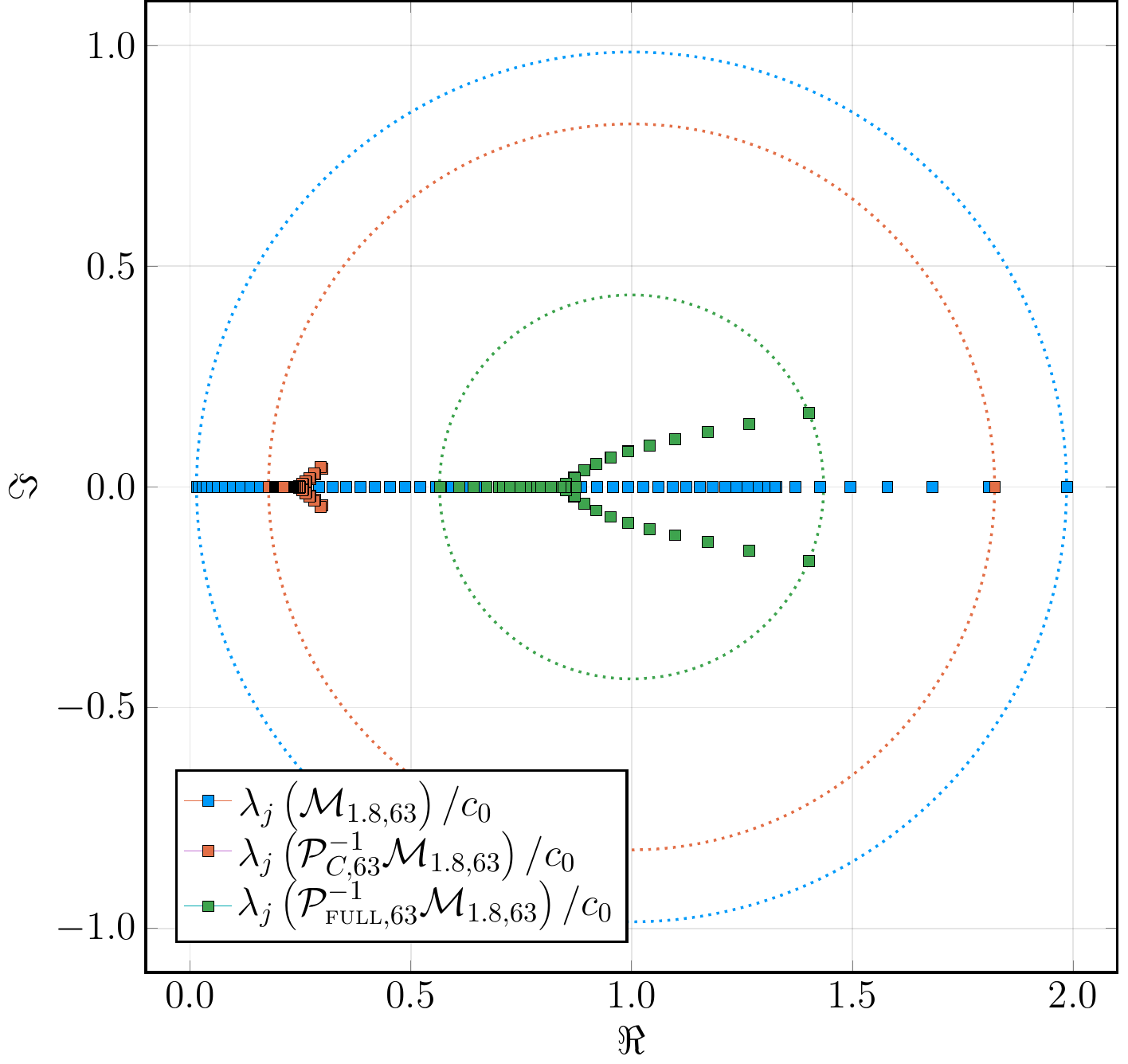}
\caption{[Example 1: 1D, $\alpha=\{1.2,1.5,1.8\}$] Scaled spectra of the resulting matrices when the preconditioners $\mathbb{I}_{n_1}$, $\mathcal{P}_{C,n_1}$, and
$\mathcal{P}_{\textsc{full},n_1}$
are applied to the coefficient matrices $\mathcal{M}_{\alpha,n_1}$ and $n_1=2^6-1$.
\textbf{Left:} $\alpha=1.2$.  \textbf{Middle:} $\alpha=1.5$. \textbf{Right:} $\alpha=1.8$.}
\label{fig:1a}
\end{figure}
In Figure~\ref{fig:1a} we present the scaled spectra of the resulting matrices, when the preconditioners $\mathbb{I}_{n_1}$, $\mathcal{P}_{C,n_1}$, and
$\mathcal{P}_{\textsc{full},n_1}$
are applied to the coefficient matrices $\mathcal{M}_{\alpha,n_1}$ when $n_1=2^6-1$ and $\alpha=1.2$ (left), $\alpha=1.5$ (middle), and $\alpha=1.8$ (right). We conclude that the spectral behavior resulting from the circulant and ``full'' symbol preconditioner resemble each other, but the condition number is lower for the ``full'' symbol preconditioner, as seen in Table~\ref{tbl:table1}.
\begin{figure}[!ht]
\centering
\includegraphics[width=0.48\textwidth]{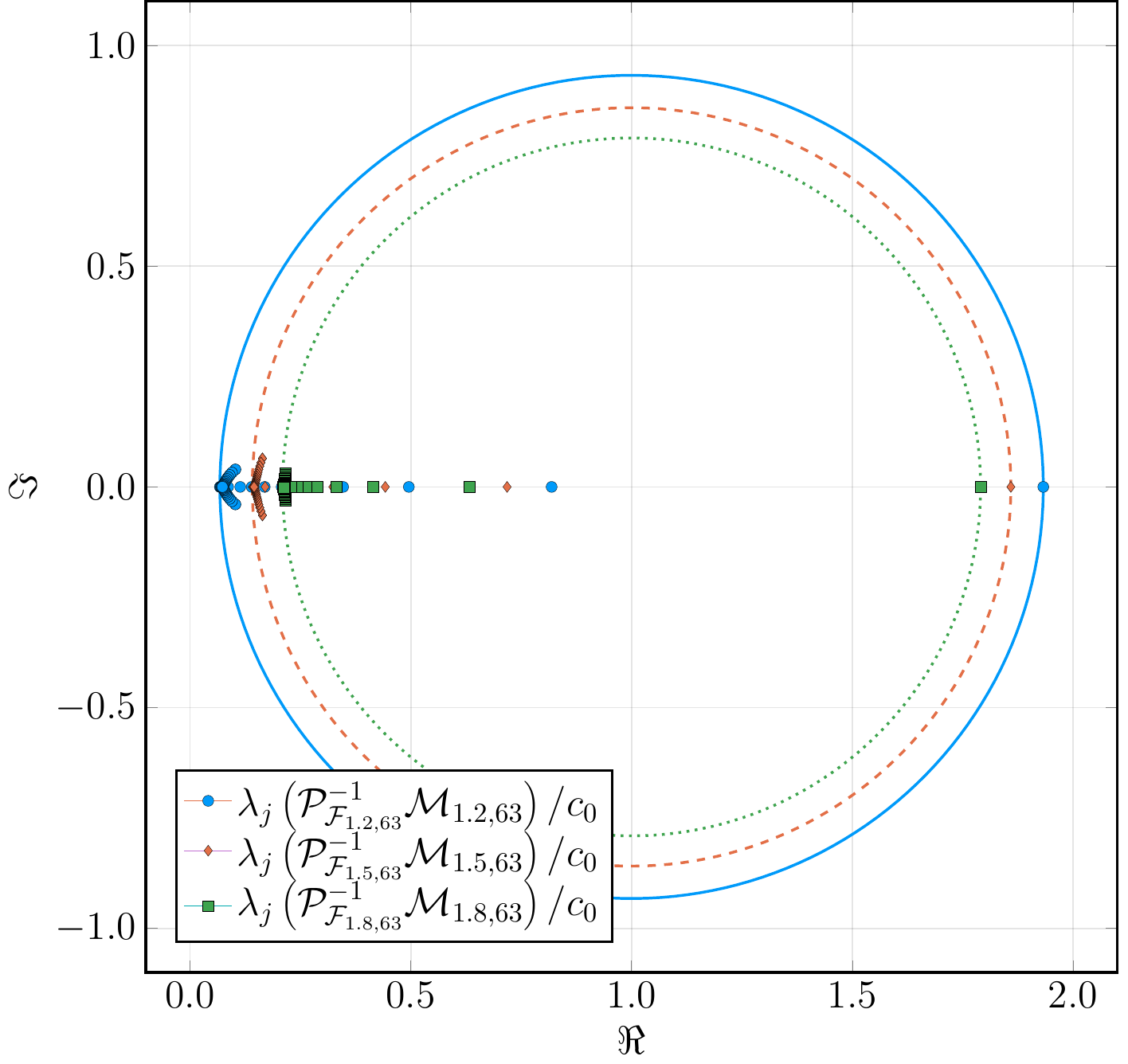}
\caption{[Example 1: 1D, $\alpha=\{1.2,1.5,1.8\}$] Scaled spectra of the resulting matrices when the preconditioners $\mathcal{P}_{\mathcal{F}_\alpha,n_1}$ are applied to the coefficient matrices $\mathcal{M}_{\alpha,n_1}$ for $n_1=2^6-1$.}
\label{fig:1b}
\end{figure}
In Figure~\ref{fig:1b} we show the scaled spectra of the resulting matrices when the preconditioners $\mathcal{P}_{\mathcal{F}_\alpha,n_1}$ are applied to the coefficient matrices $\mathcal{M}_{\alpha,n_1}$ with $n_1=2^6-1$ and $\alpha=\{1.2,1.5,1.8\}$.
We note that the clustering of the eigenvalues of the preconditioned matrices is very good except for a few large eigenvalues, especially one for any given $\alpha$. The condition number is higher for the symbol preconditioner, compared to the ``full'' symbol preconditioner, however, as seen in Table~\ref{tbl:table1} both the number of iterations and execution time is lower for the symbol preconditioner.
This confirms numerically that the term $\nu_{M,n}\mathbb{I}_n$ in the ``full'' preconditioner, has a negative impact on the performance of the preconditioner, as stated in Section~\ref{sec:diffequal}.
This is due to the fact the GMRES convergence rate largely depends on the clustering of the spectrum, and a few large eigenvalues, which might give higher condition numbers, do not degrade the convergence rate, see \cite{axelsson861}.
In Table~\ref{tbl:table2} we present the results for the following preconditioners
\begin{itemize}
\item First derivative ($\mathcal{P}_{1,n_1}$): Tridiagonal preconditioner based on the finite difference discretization of the first derivative, proposed in~\cite{donatelli161} and implemented using the Thomas algorithm.
\item Second derivative ($\mathcal{P}_{2,n_1}$): Tridiagonal preconditioner based on the finite difference discretization of the second derivative, proposed in~\cite{donatelli161} and implemented using the Thomas algorithm.
\item Tridiagonal ($\mathcal{P}_{\textsc{tri},n_1}$): Tridiagonal preconditioner based on the three main diagonals of the coefficient matrix and implemented using the Thomas algorithm.
\item Alternative symbol based ($\mathcal{P}_{\tilde{\mathcal{F}}_{\alpha},n_1}$): Constructed by $\mathbb{S}_{n_1}D_{n_1}\mathrm{diag}(p_\alpha(\theta_{j,n_1}))\mathbb{S}_{n_1}$ and implemented using FFT.
\end{itemize}
\begin{table}[!ht]
\centering
\caption{[Example 1: 1D, $\alpha=\{1.2,1.5,1.8\}$] Numerical experiments with GMRES and different preconditioners. For each preconditioner we present: average number of iterations for one time step [it], total timing in milliseconds [ms] to attain the approximate solution at time $T$, and the condition number $\kappa$ of the preconditioned  matrix, $\mathcal{P}^{-1}\mathcal{M}_{\alpha,n_1}$. Best results are highlighted in bold.}
\begin{tabular}{cc|ccc|ccc|ccc|ccccc}
\toprule
$\alpha$ & $n_1+1$&&$\mathbb{I}_{n_1}$&&&$\mathcal{P}_{\textsc{C},n_1}$&&&$\mathcal{P}_{\textsc{full},n_1}$&&&$\mathcal{P}_{\mathcal{F}_\alpha,n_1}$\\
&&\footnotesize [it]&\footnotesize [ms]&\footnotesize $\kappa$&\footnotesize [it]&\footnotesize [ms]&\footnotesize $\kappa$&\footnotesize [it]&\footnotesize [ms]&\footnotesize $\kappa$&\footnotesize [it]&\footnotesize [ms]&\footnotesize $\kappa$\\
\midrule
\footnotesize 1.2&\footnotesize $2^6$&\footnotesize \hfill 28.0&\footnotesize \hfill \textbf{1.7}&\footnotesize \hfill 9.6&\footnotesize \hfill 13.0&\footnotesize \hfill 9.6&\footnotesize \hfill 3.3&\footnotesize 14.0&\footnotesize \hfill 3.8&\footnotesize \textbf{1.6}&\footnotesize \hfill \textbf{7.2}&\footnotesize \hfill 2.3&\footnotesize \hfill 30.8\\
&\footnotesize $2^7$&\footnotesize \hfill 39.0&\footnotesize \hfill 24.3&\footnotesize \hfill 11.5&\footnotesize \hfill 14.0&\footnotesize \hfill 53.5&\footnotesize \hfill 3.6&\footnotesize 14.0&\footnotesize \hfill 17.6&\footnotesize \textbf{1.8}&\footnotesize \hfill \textbf{8.6}&\footnotesize \hfill \textbf{13.3}&\footnotesize \hfill 63.7\\
&\footnotesize $2^8$&\footnotesize \hfill 46.0&\footnotesize \hfill 114.9&\footnotesize \hfill 13.4&\footnotesize \hfill 13.0&\footnotesize \hfill 119.8&\footnotesize \hfill 3.8&\footnotesize 14.0&\footnotesize \hfill 68.8&\footnotesize \textbf{2.0}&\footnotesize \hfill \textbf{9.9}&\footnotesize \hfill \textbf{58.2}&\footnotesize \hfill 132.2\\
&\footnotesize $2^9$&\footnotesize \hfill 51.0&\footnotesize \hfill 594.5&\footnotesize \hfill 15.5&\footnotesize \hfill 12.0&\footnotesize \hfill 574.0&\footnotesize \hfill 4.2&\footnotesize 13.0&\footnotesize \hfill 312.7&\footnotesize \textbf{2.2}&\footnotesize \hfill \textbf{9.9}&\footnotesize \hfill \textbf{285.2}&\footnotesize \hfill 274.7\\
&\footnotesize $2^{10}$&\footnotesize \hfill 54.0&\footnotesize \hfill 2882.0 &\footnotesize \hfill 17.9&\footnotesize \hfill 11.0&\footnotesize \hfill 1927.0&\footnotesize \hfill 4.5&\footnotesize 12.0&\footnotesize \hfill \textbf{1415.0}&\footnotesize \hfill \textbf{2.4} &\footnotesize \hfill \textbf{10.9}&\footnotesize \hfill 1450.0&\footnotesize \hfill 571.4\\
&\footnotesize $2^{11}$&\footnotesize \hfill 56.0&\footnotesize \hfill  18569.0&\footnotesize \hfill 20.5&\footnotesize \hfill \textbf{10.0}&\footnotesize \hfill 11749.0&\footnotesize \hfill 4.9&\footnotesize 11.0&\footnotesize \hfill \textbf{8840.0}&\footnotesize \hfill \textbf{2.5}&\footnotesize \hfill 12.8&\footnotesize \hfill 9773.0&\footnotesize \hfill 1189.7 \\
\midrule
\footnotesize 1.5&\footnotesize $2^6$&\footnotesize \hfill 32.0&\footnotesize \hfill 2.0&\footnotesize \hfill 33.4&\footnotesize \hfill 12.0&\footnotesize \hfill 8.8&\footnotesize \hfill 7.1&\footnotesize 13.0&\footnotesize \hfill 3.2&\footnotesize \textbf{1.8}&\footnotesize \hfill \textbf{6.7}&\footnotesize \hfill \textbf{2.2}&\footnotesize \hfill 16.1\\
&\footnotesize $2^7$&\footnotesize \hfill 60.0&\footnotesize \hfill 37.2 &\footnotesize \hfill 51.2&\footnotesize \hfill 12.0&\footnotesize \hfill 46.7&\footnotesize \hfill 9.2&\footnotesize 13.0&\footnotesize \hfill 16.4&\footnotesize \textbf{2.1}&\footnotesize \hfill \textbf{8.0}&\footnotesize \hfill \textbf{12.5}&\footnotesize \hfill 33.3\\
&\footnotesize $2^8$&\footnotesize \hfill 89.0&\footnotesize \hfill 213.1&\footnotesize \hfill 75.8&\footnotesize \hfill 12.0&\footnotesize \hfill 111.3&\footnotesize \hfill 12.0&\footnotesize 13.0&\footnotesize \hfill 64.5&\footnotesize \textbf{2.3}&\footnotesize \hfill \textbf{8.5}&\footnotesize \hfill \textbf{52.6}&\footnotesize \hfill 70.9\\
&\footnotesize $2^9$&\footnotesize \hfill 122.0&\footnotesize \hfill 1389.0&\footnotesize \hfill 109.9&\footnotesize \hfill 12.0&\footnotesize \hfill 544.2&\footnotesize \hfill 15.8&\footnotesize 12.0&\footnotesize \hfill 288.9&\footnotesize \textbf{2.6}&\footnotesize \hfill \textbf{10.0}&\footnotesize \hfill \textbf{280.2}&\footnotesize \hfill 152.7\\
&\footnotesize $2^{10}$&\footnotesize \hfill 158.0&\footnotesize \hfill 8007.0&\footnotesize \hfill 157.7 &\footnotesize \hfill 11.0&\footnotesize \hfill 1779.0&\footnotesize \hfill 21.2&\footnotesize \hfill 11.0&\footnotesize \hfill \textbf{1366.0}&\footnotesize \hfill \textbf{2.9}&\footnotesize \hfill \textbf{10.0}&\footnotesize \hfill 1386.0&\footnotesize \hfill 331.8\\
&\footnotesize $2^{11}$&\footnotesize \hfill 195.0&\footnotesize \hfill 56266.0&\footnotesize \hfill 224.7&\footnotesize \hfill \textbf{10.0}&\footnotesize \hfill 11538.0&\footnotesize \hfill 28.6&\footnotesize \hfill \textbf{10.0} &\footnotesize \hfill \textbf{8551.0}&\footnotesize \hfill \textbf{3.2}&\footnotesize \hfill 11.0&\footnotesize \hfill 9142.0&\footnotesize \hfill 724.3\\
\midrule
\footnotesize 1.8&\footnotesize $2^6$&\footnotesize \hfill 32.0&\footnotesize \hfill 2.1&\footnotesize \hfill 136.5&\footnotesize \hfill 9.0&\footnotesize \hfill 6.6&\footnotesize \hfill 23.0&\footnotesize 10.0&\footnotesize \hfill 2.6&\footnotesize \textbf{2.6}&\footnotesize \hfill \textbf{6.1}&\footnotesize \hfill \textbf{2.2}&\footnotesize \hfill 9.7\\
&\footnotesize $2^7$&\footnotesize \hfill 67.0&\footnotesize \hfill 42.2&\footnotesize \hfill 266.3&\footnotesize \hfill 9.0&\footnotesize \hfill 36.1&\footnotesize \hfill 37.8&\footnotesize 11.0&\footnotesize \hfill 14.5&\footnotesize \textbf{2.8}&\footnotesize \hfill \textbf{6.8}&\footnotesize \hfill \textbf{11.2}&\footnotesize \hfill 19.5\\
&\footnotesize $2^8$&\footnotesize \hfill 131.0&\footnotesize \hfill 332.3&\footnotesize \hfill 494.8&\footnotesize \hfill 9.0&\footnotesize \hfill 89.8&\footnotesize \hfill 63.0&\footnotesize 10.0&\footnotesize \hfill 53.6&\footnotesize \textbf{2.9}&\footnotesize \hfill \textbf{7.0}&\footnotesize \hfill \textbf{47.2}&\footnotesize \hfill 40.8\\
&\footnotesize $2^9$&\footnotesize \hfill 231.2&\footnotesize \hfill 3085.0&\footnotesize \hfill 893.8&\footnotesize \hfill 9.0&\footnotesize \hfill 446.8&\footnotesize \hfill 106.3&\footnotesize \hfill 9.0&\footnotesize \hfill \textbf{257.9}&\footnotesize \textbf{2.9}&\footnotesize \hfill \textbf{8.6}&\footnotesize \hfill 262.8&\footnotesize \hfill 86.9\\
&\footnotesize $2^{10}$&\footnotesize \hfill 341.0&\footnotesize \hfill 20620.0&\footnotesize \hfill 1589.3&\footnotesize \hfill \textbf{8.0}&\footnotesize \hfill 1503.0&\footnotesize \hfill 180.5&\footnotesize \hfill \textbf{8.0}&\footnotesize \hfill \textbf{1191.0}&\footnotesize \hfill \textbf{3.0}&\footnotesize \hfill 10.0&\footnotesize \hfill 1370.0&\footnotesize \hfill 187.5\\
&\footnotesize $2^{11}$&\footnotesize \hfill 470.0&\footnotesize \hfill 163700.0&\footnotesize \hfill 2800.9&\footnotesize \hfill 8.0&\footnotesize \hfill 10197.0&\footnotesize \hfill 308.3&\footnotesize \hfill \textbf{7.0} &\footnotesize \hfill \textbf{7759.0}&\footnotesize \hfill \textbf{3.0}&\footnotesize \hfill 11.0&\footnotesize \hfill 9125.0&\footnotesize \hfill 408.1\\
\bottomrule
\end{tabular}
\label{tbl:table1}
\end{table}

\begin{table}[!ht]
\centering
\caption{[Example 1: 1D, $\alpha=\{1.2,1.5,1.8\}$] Numerical experiments with GMRES and different preconditioners. For each preconditioner we present: average number of iterations for one time step [it], total timing in milliseconds [ms] to attain the approximate solution at time $T$, and the condition number $\kappa$ of the preconditioned mass matrix, $\mathcal{P}^{-1}\mathcal{M}_{\alpha,n_1}$. Best results are highlighted in bold.}
\begin{tabular}{cc|ccc|ccc|ccc|cccccc}
\toprule
$\alpha$ & $n_1+1$&&$\mathcal{P}_{1,n_1}$&&&$\mathcal{P}_{2,n_1}$&&&$\mathcal{P}_{\textsc{tri},n_1}$&&&$\mathcal{P}_{\tilde{\mathcal{F}}_{\alpha},n_1}$\\
&&\footnotesize [it]&\footnotesize [ms]&\footnotesize $\kappa$&\footnotesize [it]&\footnotesize [ms]&\footnotesize $\kappa$&\footnotesize [it]&\footnotesize [ms]&\footnotesize $\kappa$&\footnotesize [it]&\footnotesize [ms]&\footnotesize $\kappa$\\
\midrule
\footnotesize 1.2&\footnotesize $2^6$&\footnotesize \hfill 8.0&\footnotesize \hfill 1.1&\footnotesize \hfill \textbf{1.2}&\footnotesize \hfill 9.0&\footnotesize \hfill 1.0&\footnotesize \hfill 2.1&\footnotesize \hfill \textbf{5.0}&\footnotesize \hfill \textbf{0.7}&\footnotesize \hfill 1.3&\footnotesize \hfill 7.5&\footnotesize \hfill 2.1&\footnotesize \hfill 29.2\\
&\footnotesize $2^7$&\footnotesize \hfill 8.0&\footnotesize \hfill 7.5&\footnotesize \hfill \textbf{1.3}&\footnotesize \hfill 10.0&\footnotesize \hfill 8.6&\footnotesize \hfill 2.2&\footnotesize \hfill \textbf{5.0}&\footnotesize \hfill \textbf{5.9}&\footnotesize \hfill 1.4&\footnotesize \hfill 8.5&\footnotesize \hfill 12.2&\footnotesize \hfill 58.7\\
&\footnotesize $2^8$&\footnotesize \hfill 7.0&\footnotesize \hfill 32.0&\footnotesize \hfill \textbf{1.3}&\footnotesize \hfill 10.0&\footnotesize \hfill 37.4&\footnotesize \hfill 2.4&\footnotesize \hfill \textbf{5.0}&\footnotesize \hfill \textbf{32.0}&\footnotesize \hfill 1.5&\footnotesize \hfill 9.9&\footnotesize \hfill 52.0&\footnotesize \hfill 118.6\\
&\footnotesize $2^9$&\footnotesize \hfill 7.0&\footnotesize \hfill 180.9&\footnotesize \hfill \textbf{1.4}&\footnotesize \hfill 10.0&\footnotesize \hfill 191.2&\footnotesize \hfill 2.6&\footnotesize \hfill \textbf{5.0}&\footnotesize \hfill \textbf{171.0}&\footnotesize \hfill 1.5&\footnotesize \hfill 9.9&\footnotesize \hfill 254.3&\footnotesize \hfill 239.7\\
&\footnotesize $2^{10}$&\footnotesize \hfill 6.0&\footnotesize \hfill 959.7 &\footnotesize \hfill \textbf{1.4}&\footnotesize \hfill 9.0&\footnotesize \hfill 1066.0&\footnotesize \hfill 2.8&\footnotesize \hfill \textbf{5.0}&\footnotesize \hfill \textbf{928.7}&\footnotesize \hfill 1.6&\footnotesize \hfill 11.0&\footnotesize \hfill 1363.0&\footnotesize \hfill 484.0\\
&\footnotesize $2^{11}$&\footnotesize \hfill 6.0&\footnotesize \hfill 7026.0&\footnotesize \hfill \textbf{1.5}&\footnotesize \hfill 9.0&\footnotesize \hfill 7675.0&\footnotesize \hfill 3.0&\footnotesize \hfill \textbf{5.0}&\footnotesize \hfill \textbf{6914.0}&\footnotesize \hfill 1.7&\footnotesize \hfill 12.0&\footnotesize \hfill 10787.0&\footnotesize \hfill 976.3\\
\midrule
\footnotesize 1.5&\footnotesize $2^6$&\footnotesize \hfill 16.0&\footnotesize \hfill 1.5&\footnotesize \hfill 2.5&\footnotesize \hfill 8.0&\footnotesize \hfill 1.0&\footnotesize \hfill \textbf{2.1}&\footnotesize \hfill \textbf{7.0}&\footnotesize \hfill \textbf{1.0}&\footnotesize \hfill 2.4&\footnotesize \hfill 8.7&\footnotesize \hfill 2.7&\footnotesize \hfill 13.6\\
&\footnotesize $2^7$&\footnotesize \hfill 20.0&\footnotesize \hfill 14.4&\footnotesize \hfill 3.1&\footnotesize \hfill 9.0&\footnotesize \hfill 8.1&\footnotesize \hfill \textbf{2.3}&\footnotesize \hfill \textbf{8.0}&\footnotesize \hfill \textbf{7.5}&\footnotesize \hfill 3.0&\footnotesize \hfill 8.0&\footnotesize \hfill 12.1&\footnotesize \hfill 26.3\\
&\footnotesize $2^8$&\footnotesize \hfill 24.0&\footnotesize \hfill 67.9&\footnotesize \hfill 4.0&\footnotesize \hfill 9.0&\footnotesize \hfill \textbf{35.3}&\footnotesize \hfill \textbf{2.7}&\footnotesize \hfill 11.0&\footnotesize \hfill 40.2&\footnotesize \hfill 4.0&\footnotesize \hfill \textbf{8.4}&\footnotesize \hfill 47.7&\footnotesize \hfill 51.8\\
&\footnotesize $2^9$&\footnotesize \hfill 26.0&\footnotesize \hfill 366.7&\footnotesize \hfill 5.2&\footnotesize \hfill 10.0&\footnotesize \hfill \textbf{197.5}&\footnotesize \hfill \textbf{3.0}&\footnotesize \hfill 13.0&\footnotesize \hfill 227.3&\footnotesize \hfill 5.4&\footnotesize \hfill \textbf{9.9}&\footnotesize \hfill 248.1&\footnotesize \hfill 103.0\\
&\footnotesize $2^{10}$&\footnotesize \hfill 27.0&\footnotesize \hfill 1810.0&\footnotesize \hfill 6.9&\footnotesize \hfill \textbf{10.0}&\footnotesize \hfill \textbf{1105.0}&\footnotesize \hfill \textbf{3.5}&\footnotesize \hfill 15.0&\footnotesize \hfill 1331.0&\footnotesize \hfill 7.4&\footnotesize \hfill \textbf{10.0}&\footnotesize \hfill 1636.0&\footnotesize \hfill 205.9\\
&\footnotesize $2^{11}$&\footnotesize \hfill 25.4&\footnotesize \hfill 11212.0&\footnotesize \hfill 9.0&\footnotesize \hfill \textbf{11.0}&\footnotesize \hfill \textbf{8179.0}&\footnotesize \hfill \textbf{4.0}&\footnotesize \hfill 18.0&\footnotesize \hfill 9684.0&\footnotesize \hfill 10.4&\footnotesize \hfill \textbf{11.0}&\footnotesize \hfill 10563.0&\footnotesize \hfill 424.5\\
\midrule
\footnotesize 1.8&\footnotesize $2^6$&\footnotesize \hfill 25.0&\footnotesize \hfill 2.5&\footnotesize \hfill 8.4&\footnotesize \hfill \textbf{6.0}&\footnotesize \hfill \textbf{0.8}&\footnotesize \hfill \textbf{1.6}&\footnotesize \hfill 7.0&\footnotesize \hfill 1.0&\footnotesize \hfill 3.5&\footnotesize \hfill 8.0&\footnotesize \hfill 2.3&\footnotesize \hfill 9.0\\
&\footnotesize $2^7$&\footnotesize \hfill 40.0&\footnotesize \hfill 27.3&\footnotesize \hfill 14.3&\footnotesize \hfill \textbf{6.0}&\footnotesize \hfill \textbf{6.3}&\footnotesize \hfill \textbf{1.7}&\footnotesize \hfill 10.0&\footnotesize \hfill 8.7&\footnotesize \hfill 5.6&\footnotesize \hfill 7.8&\footnotesize \hfill 11.3&\footnotesize \hfill 17.0\\
&\footnotesize $2^8$&\footnotesize \hfill 61.0&\footnotesize \hfill 159.8&\footnotesize \hfill 25.3&\footnotesize \hfill 7.0&\footnotesize \hfill \textbf{31.0}&\footnotesize \hfill \textbf{1.8}&\footnotesize \hfill 15.0&\footnotesize \hfill 48.3&\footnotesize \hfill 9.4&\footnotesize \hfill \textbf{6.9}&\footnotesize \hfill 43.3&\footnotesize \hfill 33.1\\
&\footnotesize $2^9$&\footnotesize \hfill 88.0&\footnotesize \hfill 1083.0&\footnotesize \hfill 44.7&\footnotesize \hfill \textbf{7.0}&\footnotesize \hfill \textbf{170.1}&\footnotesize \hfill \textbf{2.0}&\footnotesize \hfill 22.0&\footnotesize \hfill 325.4&\footnotesize \hfill 16.6&\footnotesize \hfill \textbf{7.0}&\footnotesize \hfill 222.4&\footnotesize \hfill 65.4\\
&\footnotesize $2^{10}$&\footnotesize \hfill 120.0&\footnotesize \hfill 6277.0&\footnotesize \hfill 78.8&\footnotesize \hfill \textbf{7.0}&\footnotesize \hfill \textbf{999.3}&\footnotesize \hfill \textbf{2.3}&\footnotesize \hfill 31.0&\footnotesize \hfill 1983.0&\footnotesize \hfill 30.0&\footnotesize \hfill 8.9&\footnotesize \hfill 1569.0&\footnotesize \hfill 130.1\\
&\footnotesize $2^{11}$&\footnotesize \hfill 158.0&\footnotesize \hfill 46716.0&\footnotesize \hfill 138.2&\footnotesize \hfill \textbf{7.0}&\footnotesize \hfill \textbf{7309.0}&\footnotesize \hfill \textbf{2.6}&\footnotesize \hfill 44.7&\footnotesize \hfill 15756.0&\footnotesize \hfill 54.6&\footnotesize \hfill 10.0&\footnotesize \hfill 10249.0&\footnotesize \hfill 259.8\\
\bottomrule
\end{tabular}
\label{tbl:table2}
\end{table}

\begin{figure}[!ht]
\centering
\includegraphics[width=0.32\textwidth]{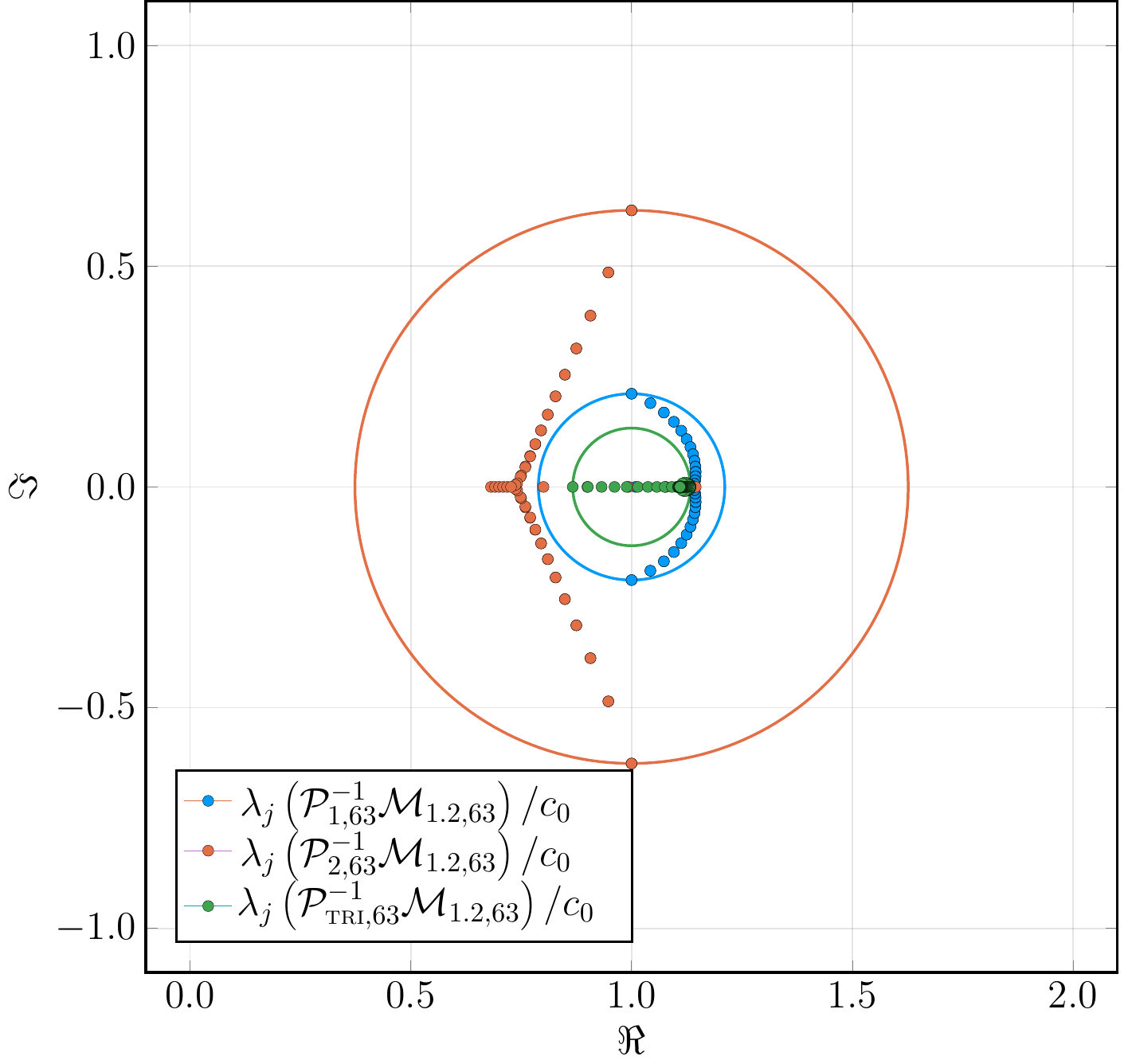}
\includegraphics[width=0.32\textwidth]{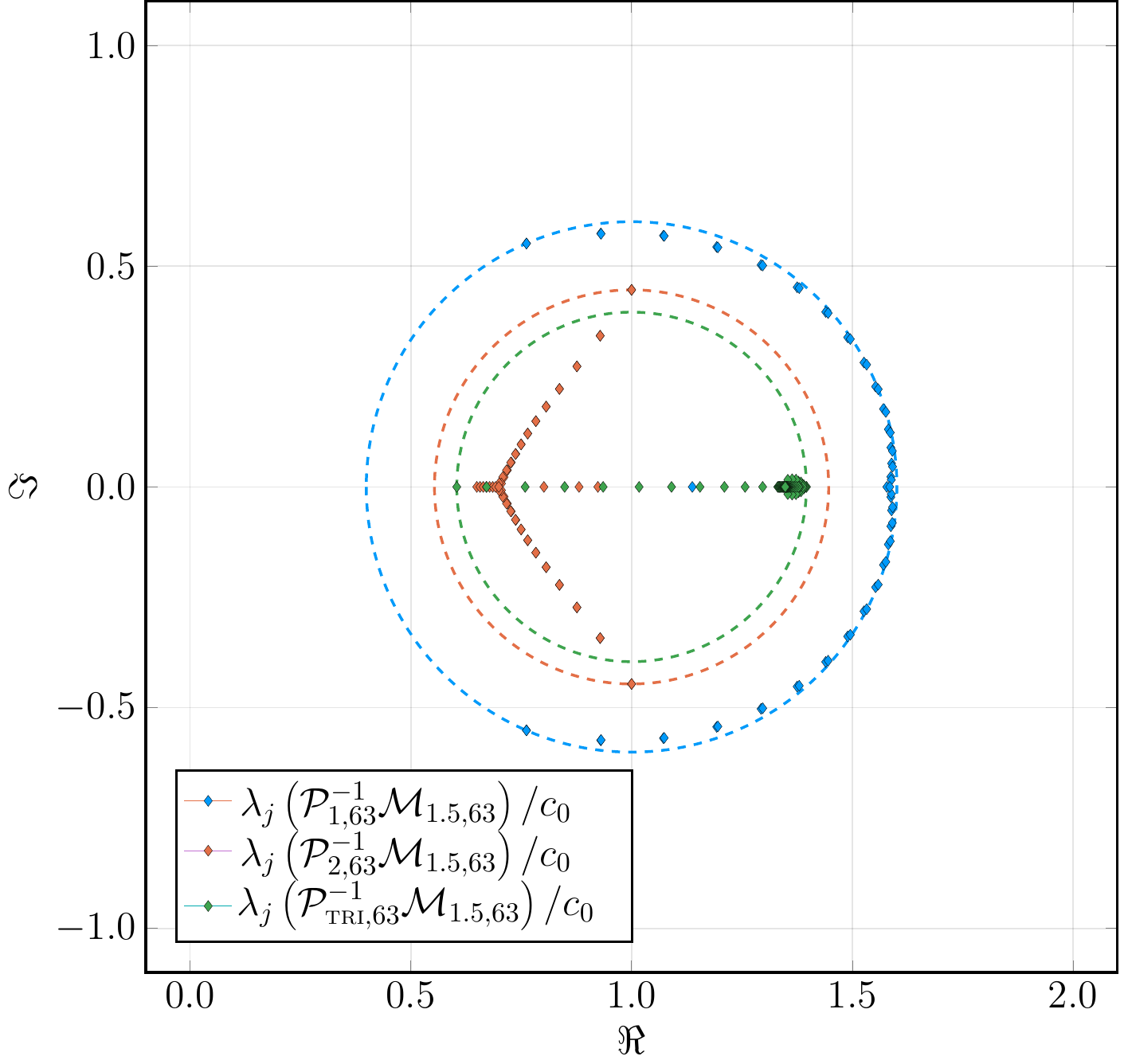}
\includegraphics[width=0.32\textwidth]{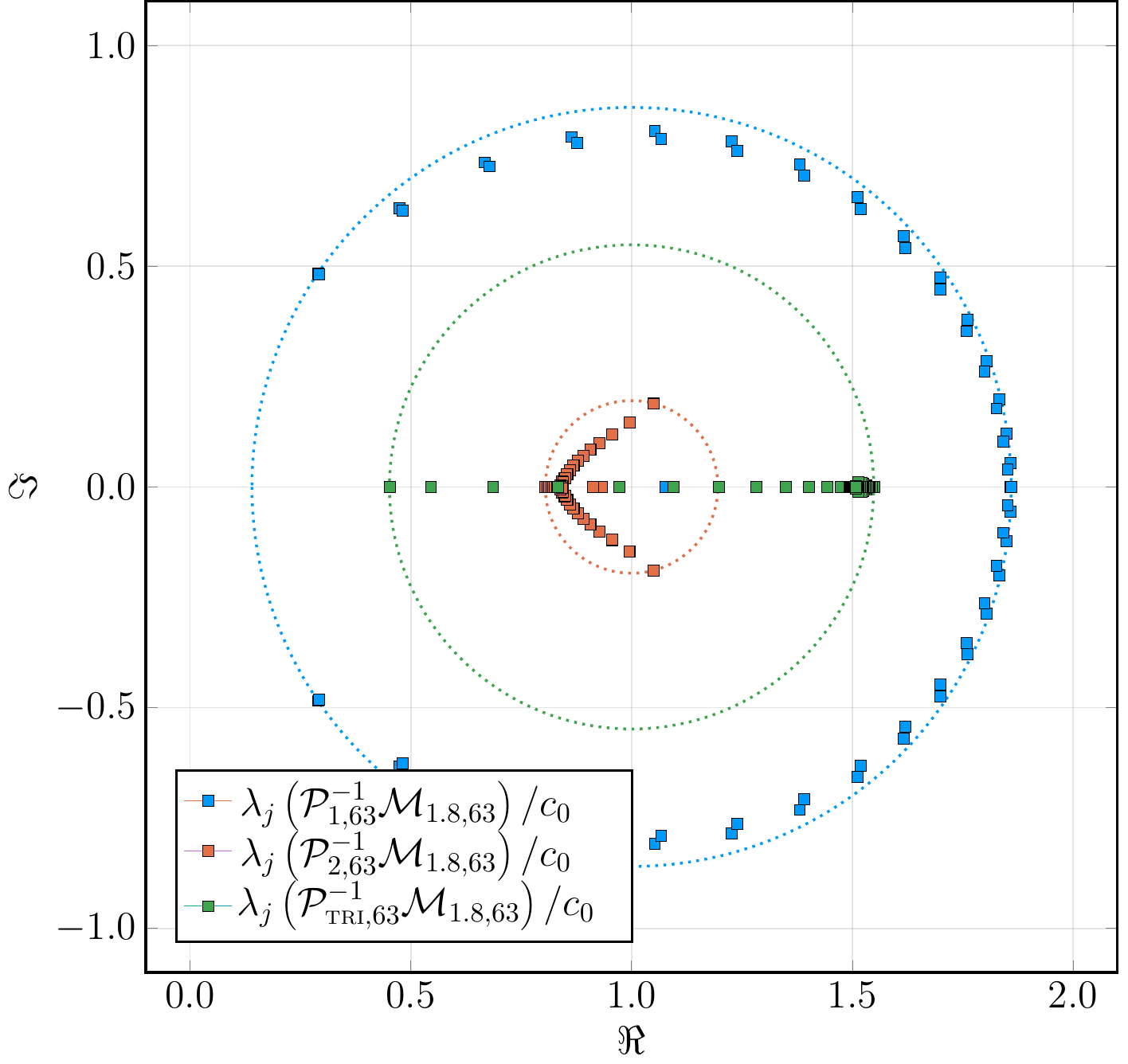}
\caption{[Example 1: 1D, $\alpha=\{1.2,1.5,1.8\}$] Scaled spectra of the resulting matrices when the preconditioners $\mathcal{P}_{1,n_1}$, $\mathcal{P}_{2,n_1}$, and $\mathcal{P}_{\textsc{tri},n_1}$ are applied to the  matrices $\mathcal{M}_{\alpha,n_1}$ and $n_1=2^6-1$.
\textbf{Left:} $\alpha=1.2$.  \textbf{Middle:} $\alpha=1.5$. \textbf{Right:} $\alpha=1.8$.}
\label{fig:1c}
\end{figure}

\begin{figure}[!ht]
\centering
\includegraphics[width=0.48\textwidth]{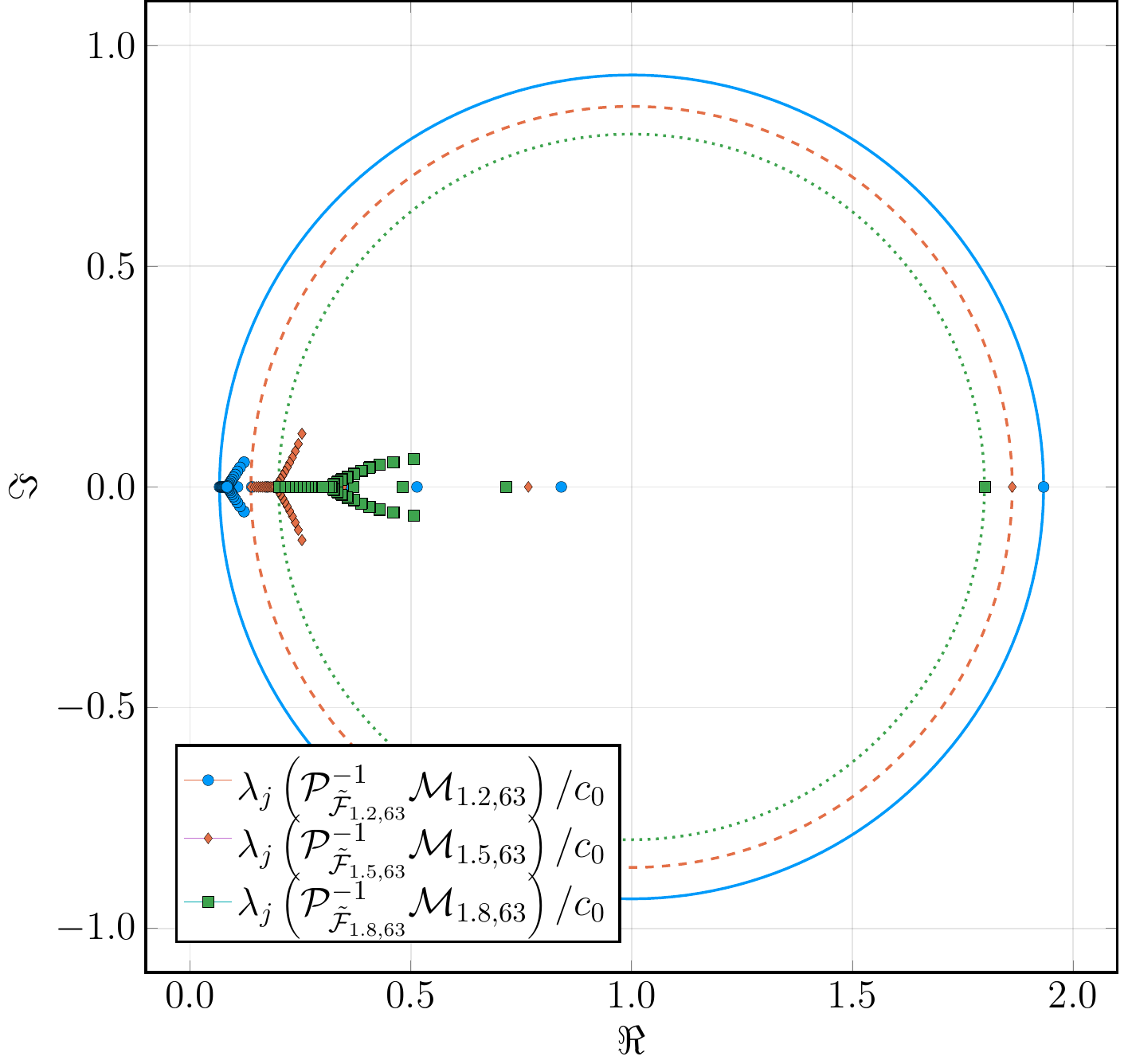}
\caption{[Example 1: 1D, $\alpha=\{1.2,1.5,1.8\}$] Scaled spectra of the resulting matrices when the preconditioners $\mathcal{P}_{\tilde{\mathcal{F}}_\alpha,n_1}$ are applied to the matrices $\mathcal{M}_{\alpha,n_1}$ for $n_1=2^6-1$.}
\label{fig:1d}
\end{figure}

Like in Figure~\ref{fig:1a}, in Figure~\ref{fig:1c} we present the scaled spectra of the preconditioned matrix. The spectral behavior of the three preconditioners (first and second derivative and the tridiagonal) for different $\alpha$ correlate well with the results presented in Table~\ref{tbl:table2}. In the left panel of Figure~\ref{fig:1c} the best clustering is obtained using the tridiagonal preconditioner, followed by the first derivative, and then by the second derivative. Since $\alpha=1.2$, a value close to one, this behavior is expected. When $\alpha=1.5$, as presented in the middle panel of Figure~\ref{fig:1c}, the results are similar for the three preconditioners, but the second derivative preconditioner performs best as $n_1$ increases. In the right panel of Figure~\ref{fig:1c} we see that the best clustering is observed for the second derivative preconditioner, and it also performs best for all $n_1$ and all reported quantities (iterations, timings, and condition numbers).
The better performance of the preconditioners reported in Table~\ref{tbl:table2} as opposed the ones in Table~\ref{tbl:table1} is expected: this is due to the computational complexity of $\mathcal{O}(n)$ for the Thomas algorithm, as opposed to $\mathcal{O}(n\log n)$ for the DFT. However,  due to the inherit parallel nature of FFT opposed to serial  Thomas algorithm, this disadvantage will turn to be a significant  benefit for our proposal if  a parallel environment  is  used.        
In Figure~\ref{fig:1d} we present the scaled spectrum of an alternative symbol based preconditioner, $\mathcal{P}_{\tilde{\mathcal{F}}_\alpha,n_1}$, which performs slightly better than the proposed preconditioner $\mathcal{P}_{\mathcal{F}_\alpha,n_1}$ in Section~\ref{sec:prop1d} (compare Tables~\ref{tbl:table1} and \ref{tbl:table2}). This is mainly due to the avoided multiplication with the inverse of $D_n$ for $\mathcal{P}_{\tilde{\mathcal{F}}_\alpha,n_1}$, since the spectrum  of the resulted preconditioned matrices using  $\mathcal{P}_{\tilde{\mathcal{F}}_\alpha,n_1}$ and $\mathcal{P}_{\mathcal{F}_\alpha,n_1}$ are comparable. Furthermore, in this case it seems that the  most efficient choice of preconditioner is problem specific, depending on $d_{\pm}$.

\subsection{Example 2}
The considered two-dimensional example is originally from~\cite[Example 4.]{pang151} and also discussed in~\cite[Example 1.]{moghaderi171}. In \eqref{eq:fde}, define $\alpha=1.8$, $\beta=1.6$, and
\begin{align}
    d_+(x,y)=\Gamma(3-\alpha)(1+x)^\alpha(1+y)^2, \qquad d_-(x,y)=\Gamma(3-\alpha)(3-x)^\alpha(3-y)^2,\nonumber \\
    e_+(x,y)=\Gamma(3-\beta)(1+x)^2(1+y)^\beta,  \qquad e_-(x,y)=\Gamma(3-\beta)(3-x)^2(3-y)^\beta.\nonumber
\end{align}
The spatial domain is $\Omega=[0,2]\times[0,2]$ and the time interval is $[t_0,T]=[0,1]$.
The initial condition is
\begin{align}
    u(x,y,0)=u_0(x,y)=x^2y^2(2-x)^2(2-y)^2,\nonumber
\end{align}
and the source term is
\begin{align}
f(x,y,t)&=-16e^{-t}\left(x^2(2-x)^2y^2(2-y)^2
+f_\alpha(x,y)
+f_\alpha(2-x,2-y)
+f_\beta(y,x)
+f_\beta(2-y,2-x)\right),\nonumber\\
f_\gamma(x,y)&=\left(8x^{2-\gamma}-\frac{24x^{3-\gamma}}{3-\gamma}+\frac{24x^{4-\gamma}}{(4-\gamma)(3-\gamma)}\right)(1+x)^\gamma(1+y)^2y^2(2-y)^2,\nonumber
\end{align}
such that the solution to the FDE is given by $u(x,y,t)=16e^{-t}x^2(2-x)^2y^2(2-y)^2$.
Let $h=h_x=h_y=2/(n+1)$, with $n=n_1=n_2=M$, and $h_t=1/(M+1)$. Then,
\begin{align}
    \frac{1}{r}=\frac{2h^\alpha}{h_t}=\frac{2^{\alpha+1}M}{(n+1)^\alpha}=\frac{2^{\alpha+1}n}{(n+1)^\alpha}, \qquad \frac{s}{r}=\frac{h^\alpha}{h^\beta}=2^{\alpha-\beta}(n+1)^{\beta-\alpha}.\nonumber
\end{align}
In Table~\ref{tbl:table3} (and also Table~\ref{tbl:table4}) we present the results for the following preconditioners:
\begin{itemize}
\item Second derivative ($\mathcal{P}_{2,N}$): Preconditioner based on the finite difference discretization of the second derivative, proposed in~\cite{moghaderi171} and implemented using one Galerkin projection multigrid V-cycle.
\item Algebraic multigrid ($\mathcal{P}_{\textsc{MGM},N}$): Preconditioner based on algebraic multigrid, proposed in~\cite{moghaderi171} and implemented using one algebraic multigrid V-cycle.
\item Symbol ($\mathcal{P}_{\mathcal{F}_{(\alpha,\beta)},N}$): Proposed preconditioner and implemented using FFT.
\end{itemize}

We mention that in  multi dimensional setting,   holds  a negative results concerning the optimality of circulant algebra  when it is used for preconditioning  Toeplitz matrices  generated by  function with  zeros of order greater than 1 (see \cite{NSV_tcs}, \cite{NSV_con}).  Thus, we consider unnecessary  a comparison with such kind of preconditioners. 
\begin{table}[!ht]
\centering
\caption{[Example 2: 2D, $\alpha=1.8, \beta=1.6$] Numerical experiments with GMRES and different preconditioners. For each preconditioner we present: average number of iterations for one time step [it], total timing in milliseconds [ms] to attain the approximate solution at time $T$, and the condition number $\kappa$ of the preconditioned  matrix, $\mathcal{P}^{-1}\mathcal{M}_{(\alpha,\beta),N}$. Best results are highlighted in bold.}
\begin{tabular}{c|ccc|ccc|ccc|cccc}
 \toprule
  $n_1=n_2$&&$\mathbb{I}_N$&&&$\mathcal{P}_{2,N}$&&&$\mathcal{P}_{\textsc{mgm},N}$&&&$\mathcal{P}_{\mathcal{F}_{(\alpha,\beta)},N}$ \\
  &\footnotesize [it]&\footnotesize [ms]&\footnotesize $\kappa$&\footnotesize [it]&\footnotesize [ms]&\footnotesize $\kappa$&\footnotesize [it]&\footnotesize [ms]&\footnotesize $\kappa$&\footnotesize [it]&\footnotesize [ms]&\footnotesize $\kappa$\\
 \midrule
 \footnotesize $2^4$&\footnotesize \hfill 37.0&\footnotesize \hfill 32.2&\footnotesize \hfill 57.4&\footnotesize \hfill 21.0&\footnotesize \hfill 64.8&\footnotesize \hfill 48.6&\footnotesize \hfill 10.0&\footnotesize \hfill 40.8&\footnotesize \hfill 3.7&\footnotesize \hfill \textbf{8.0}&\footnotesize \hfill \textbf{35.1}&\footnotesize \hfill \textbf{1.9}\\
 \footnotesize $2^5$&\footnotesize \hfill 73.0 &\footnotesize \hfill 331.4&\footnotesize \hfill 167.4&\footnotesize \hfill 17.6&\footnotesize \hfill 551.1&\footnotesize \hfill 31.7&\footnotesize \hfill 11.0&\footnotesize \hfill 383.1&\footnotesize \hfill 5.4&\footnotesize \textbf{8.0}&\footnotesize \hfill \textbf{296.8}&\footnotesize \hfill \textbf{2.7}\\
 \footnotesize $2^6$&\footnotesize \hfill 137.0 &\footnotesize \hfill 35440.0&\footnotesize \hfill 429.4&\footnotesize \hfill 17.0&\footnotesize \hfill 10465.0&\footnotesize \hfill 310.7&\footnotesize \hfill 11.0&\footnotesize \hfill 16146.0&\footnotesize \hfill 8.2&\footnotesize \hfill \textbf{9.0}&\footnotesize \hfill \textbf{6569.0}&\footnotesize \hfill \textbf{4.3}\\
 \footnotesize $2^7$&\footnotesize \hfill 251.0 &\footnotesize \hfill 1644134.0&\footnotesize \hfill 966.8&\footnotesize \hfill 17.0&\footnotesize \hfill 213713.0&\footnotesize \hfill 678.4&\footnotesize \hfill 10.0&\footnotesize 352471.0&\footnotesize \hfill 12.2&\footnotesize \hfill \textbf{9.0}&\footnotesize \hfill \textbf{135535.0}&\footnotesize \hfill \textbf{7.7}\\
 \bottomrule
\end{tabular}
\label{tbl:table3}
\end{table}

\begin{figure}[!ht]
\centering
\includegraphics[width=0.48\textwidth]{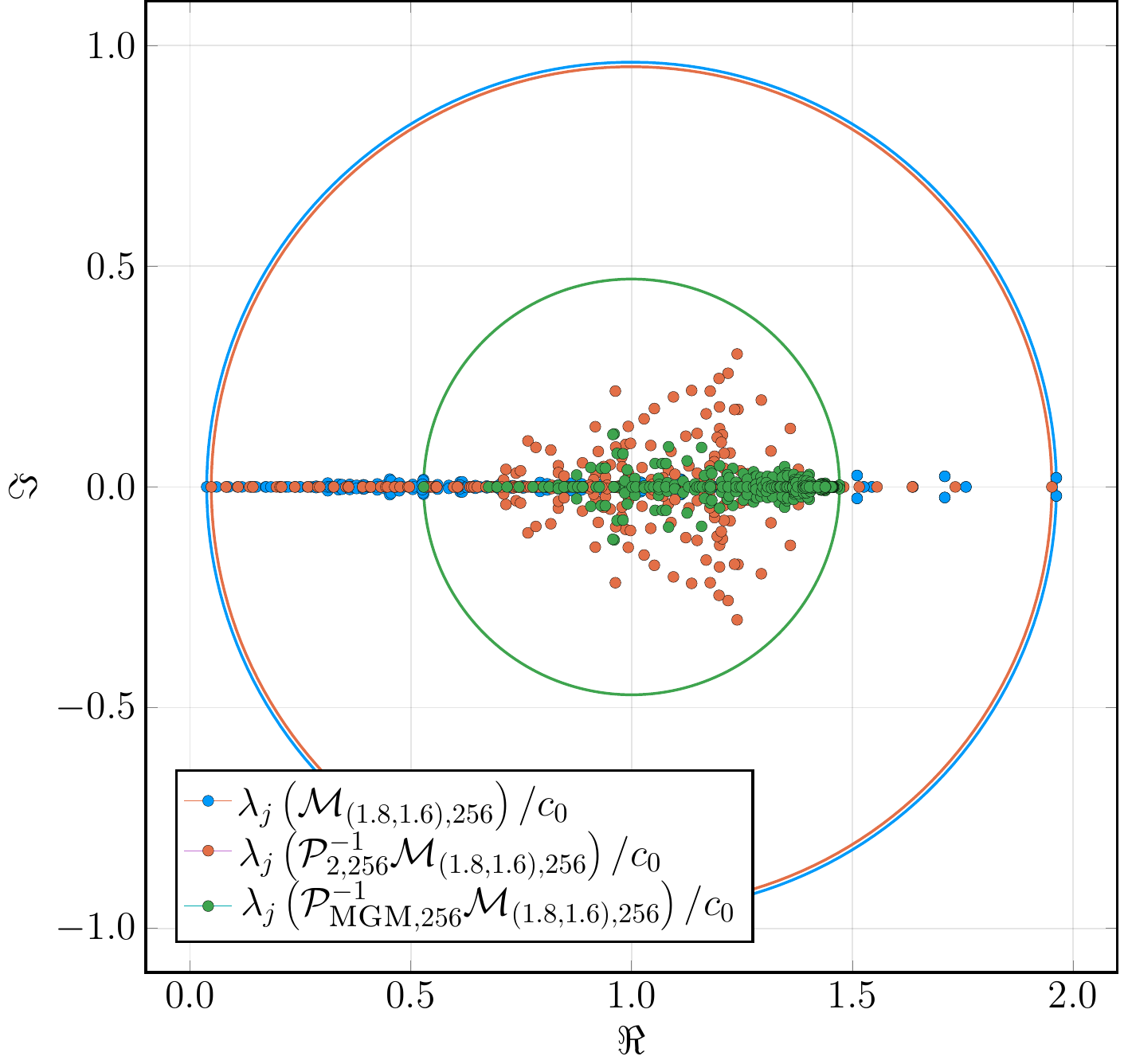}
\includegraphics[width=0.48\textwidth]{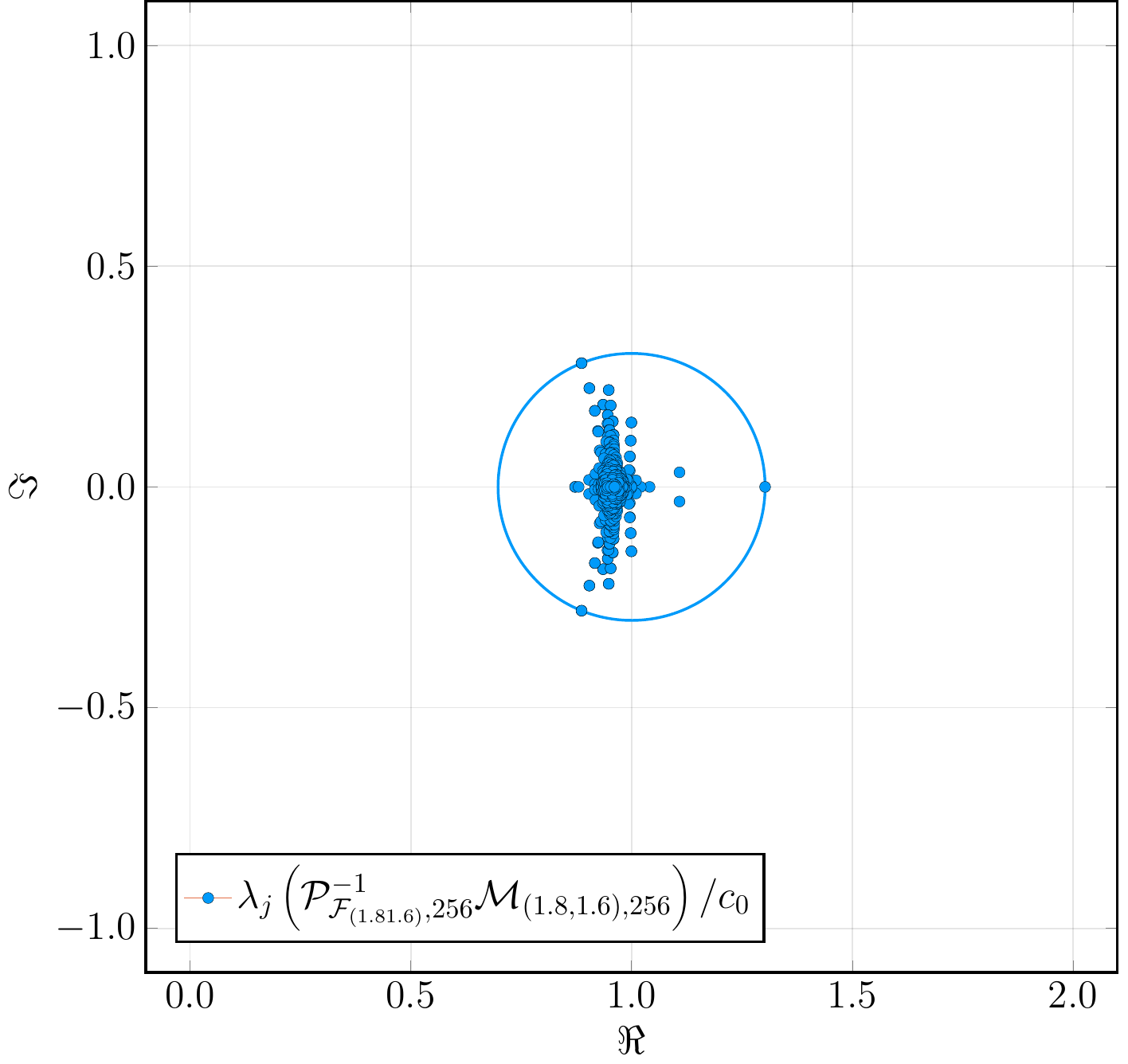}
\caption{[Example 2: 2D, $\alpha=1.8, \beta=1.6$] Scaled spectra of the resulting matrices when the preconditioners are applied to the  matrices $\mathcal{M}_{(\alpha,\beta),n_1^2}$ and $n_1=2^4$. \textbf{Left:} Preconditioners $\mathbb{I}_N$, $\mathcal{P}_{2,N}$, and $\mathcal{P}_{\textsc{mgm},N}$ \textbf{Right:} Preconditioner $\mathcal{P}_{\mathcal{F}_{(\alpha,\beta)},N}$.}
\label{fig:2}
\end{figure}

For details on the multigrid based preconditioners, $\mathcal{P}_{2,N}$ (Galerkin projection multigrid) and $\mathcal{P}_{\textsc{mgm},N}$ (algebraic multigrid), see~\cite{moghaderi171}.
The proposed symbol-based preconditioner, $\mathcal{P}_{\mathcal{F}_{(\alpha,\beta)},N}$, performs better than the multigrid-based preconditioners, as seen in Table~\ref{tbl:table3}. In Figure~\ref{fig:2} we present the scaled spectra of the preconditioned  matrices for $N=n_1n_2=2^8$. The clustering is better for the proposed symbol-based preconditioners than the other three, as seen comparing the left and right panels. We note in Table~\ref{tbl:table3} that the number of iterations are essentially constant both for the algebraic multigrid and the symbol-based preconditioners.

By fine tuning parameters for the multigrid-based preconditioners, such as number of smoothing steps, W-cycles etc, these results might be improved. However, the simplicity of the proposed preconditioner, where no fine-tunings are required, is advantageous.

\subsection{Example 3}
By modifying the coefficients $\alpha=1.8$ and $\beta=1.6$ in Example 2, to $\alpha=1.8$ and $\beta=1.2$ we obtain Example 3. In Table~\ref{tbl:table4} we present the same type of computations as in Table~\ref{tbl:table3}.
As discussed in~\cite{moghaderi171}, the performance of the proposed multigrid-based preconditioners depend on the fractional derivatives $\alpha$ and $\beta$. Since, in this example, $\alpha$ and $\beta$ differ more than in Example 2, and $\beta$ is far away from two, we clearly see in Table~\ref{tbl:table4} that the multigrid-based preconditioners perform worse than in Example 2. Especially note the worse behavior of the condition number for the algebraic multigrid-based preconditioner $\mathcal{P}_{\textsc{mgm},N}$. The condition numbers are essentially the same for the symbol-based preconditioner $\mathcal{P}_{\mathcal{F}_{(\alpha,\beta)},N}$ in Examples 2 and 3.

\begin{table}[!ht]
\centering
\caption{[Example 3: 2D, $\alpha=1.8, \beta=1.2$] Numerical experiments with GMRES and different preconditioners. For each preconditioner we present: average number of iterations for one time step [it], total timing in milliseconds [ms] to attain the approximate solution at time $T$, and the condition number $\kappa$ of the preconditioned  matrix, $\mathcal{P}^{-1}\mathcal{M}_{(\alpha,\beta),N}$. Best results are highlighted in bold.}
\begin{tabular}{c|ccc|ccc|ccc|cccc}
 \toprule
  $n_1=n_2$&&$\mathbb{I}_N$&&&$\mathcal{P}_{2,N}$&&&$\mathcal{P}_{\textsc{mgm},N}$&&&$\mathcal{P}_{\mathcal{F}_{(\alpha,\beta)},N}$ \\
  &\footnotesize [it]&\footnotesize [ms]&\footnotesize $\kappa$&\footnotesize [it]&\footnotesize [ms]&\footnotesize $\kappa$&\footnotesize [it]&\footnotesize [ms]&\footnotesize $\kappa$&\footnotesize [it]&\footnotesize [ms]&\footnotesize $\kappa$\\
 \midrule
 \footnotesize $2^4$ &\footnotesize \hfill 49.0 &\footnotesize \hfill 37.1 &\footnotesize \hfill 57.8&\footnotesize \hfill 26.5 &\footnotesize \hfill 79.7 &\footnotesize \hfill 42.8&\footnotesize \hfill 18.0 &\footnotesize \hfill 39.0 &\footnotesize \hfill 8.2&\footnotesize \hfill \textbf{10.0} &\footnotesize \hfill \textbf{37.0}&\footnotesize \hfill \textbf{1.9} \\
 \footnotesize $2^5$ &\footnotesize \hfill 92.0 &\footnotesize \hfill 394.0 &\footnotesize \hfill 162.9 &\footnotesize \hfill 32.0 &\footnotesize \hfill 713.8  &\footnotesize \hfill 104.0&\footnotesize \hfill 26.0 &\footnotesize \hfill 450.7 &\footnotesize \hfill 16.7&\footnotesize \textbf{12.0} &\footnotesize \hfill \textbf{329.0}&\footnotesize \hfill \textbf{2.7}    \\
 \footnotesize $2^6$ &\footnotesize \hfill  173.0 &\footnotesize \hfill 44532.0&\footnotesize \hfill  401.7 &\footnotesize \hfill 41.0 &\footnotesize \hfill 17197.0   &\footnotesize \hfill 231.6&\footnotesize \hfill 33.0 &\footnotesize \hfill 35021.0 &\footnotesize 32.8&\footnotesize \hfill \textbf{13.0} &\footnotesize \hfill \textbf{7493.0}&\footnotesize \hfill \textbf{4.4}   \\
 \footnotesize $2^7$ &\footnotesize  \hfill 316.0&\footnotesize \hfill 2070478.0&\footnotesize \hfill 876.4 &\footnotesize \hfill 51.0 &\footnotesize \hfill 438344.0 &\footnotesize \hfill 515.8&\footnotesize \hfill 41.0 &\footnotesize \hfill 1107711.0&\footnotesize 62.9&\footnotesize \hfill \textbf{14.5} &\footnotesize \hfill \textbf{171500.0}& \footnotesize \hfill \textbf{7.9}\\
 \bottomrule
\end{tabular}
\label{tbl:table4}
\end{table}

\begin{figure}[!ht]
\centering
\includegraphics[width=0.48\textwidth]{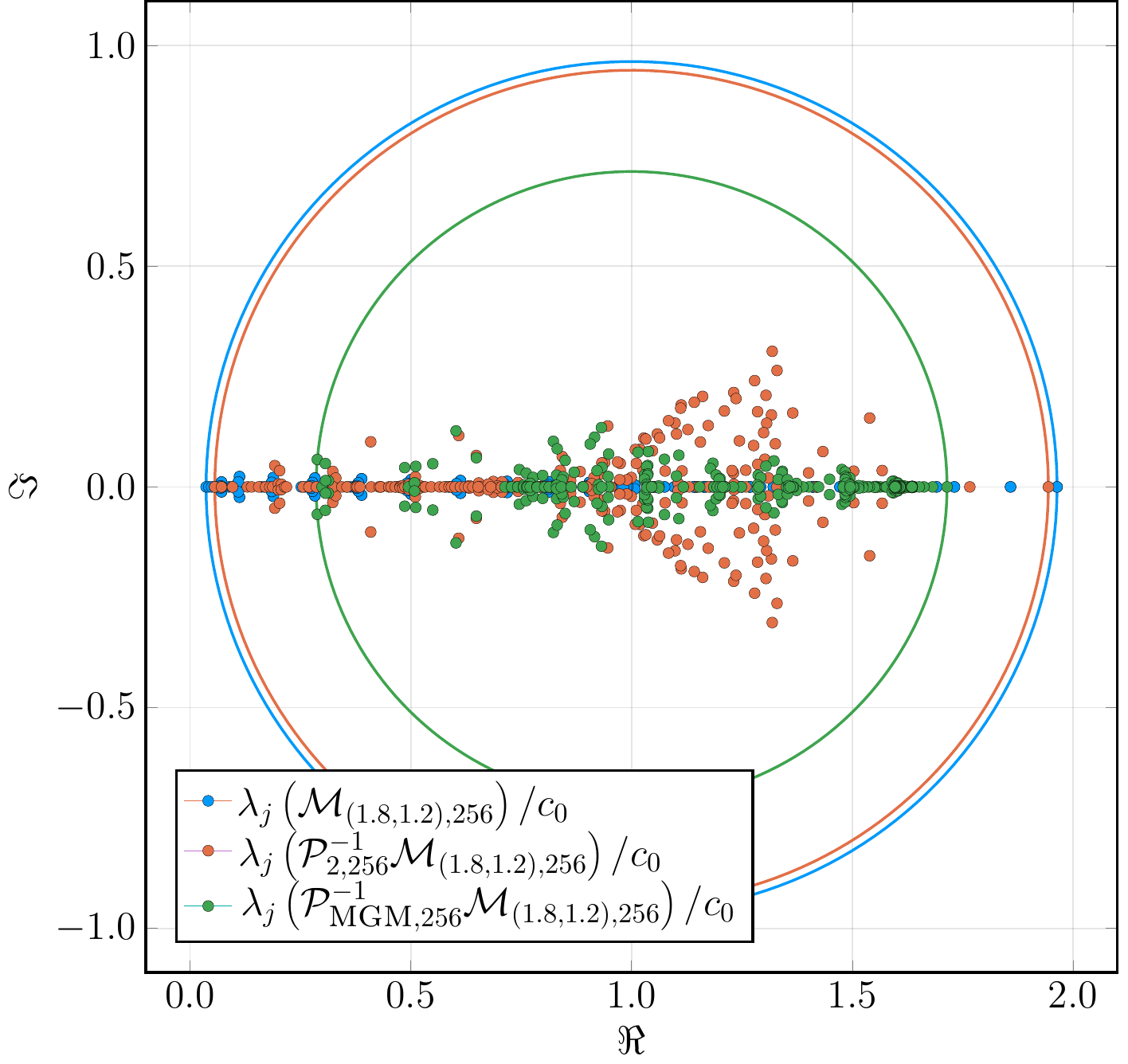}
\includegraphics[width=0.48\textwidth]{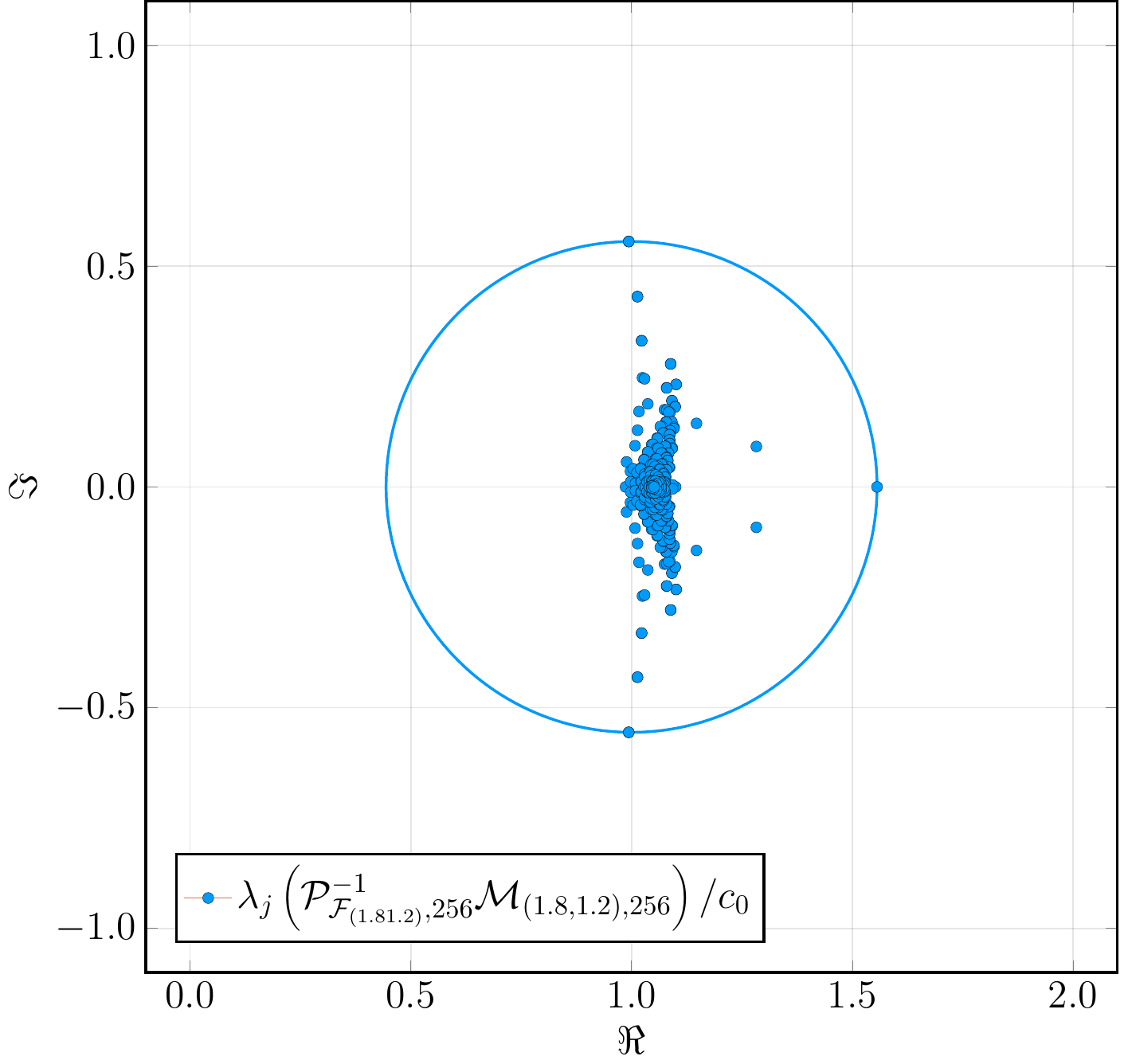}
\caption{[Example 3: 2D, $\alpha=1.8, \beta=1.2$] Scaled spectra of the resulting matrices when the preconditioners are applied to the coefficient  matrices $\mathcal{M}_{(\alpha,\beta),n_1^2}$, and $n_1=2^4$. \textbf{Left:} Preconditioners $\mathbb{I}_N$, $\mathcal{P}_{2,N}$, and $\mathcal{P}_{\textsc{mgm},N}$ \textbf{Right:} Preconditioner $\mathcal{P}_{\mathcal{F}_{(\alpha,\beta)},N}$.}
\label{fig:3}
\end{figure}

In Figure~\ref{fig:3} we present the same scaled spectra as in Figure~\ref{fig:2}, but regarding Example 3. Again we note the advantageous clustering properties of the proposed symbol-based preconditioner in the right panel.

%% file: AMERICAN.tex
\chapter{Policy Iteration Algorithm}

\section{Introduction}

\par
As mentioned in the introduction the algorithm developed here iteratively improves exercise \emph{policies} untill no further improvement is feasible. A policy for the American put option is specified by an exercise boundary $b=b(t)$, the underlying price for which the holder has pre-decided that they will exercise the option at moment $t$. If starting above the boundary, the option is exercised at the moment the underlying asset first reaches it; otherwise, it is exercised immediately. The optimal policy is one such exercise policy \cite{H.McK:65} and can therefore be obtained by appropriately selecting the boundary; meanwhile, the option value above the boundary is determined by solving the Black--Scholes equation with the value specified on the boundary. Intuitively, the optimal boundary that corresponds to the optimal policy maximises the option value. A closed-form expression of the optimal boundary for the American put is not available; however, it can be characterised in terms of nonlinear integral equations. The first such equation was given in \cite{H.McK:65} but presented computational difficulties because it involved the derivative of a boundary which becomes infinite towards the end of the exercise interval. Other integral equations do not include the derivative, (see, e.g., \cite{Jacka}, \cite{CarrJarrowMyn}, \cite{Kim}, and \cite{Zhu18}). In \cite{J.Goo:02}, an alternative nonlinear integral equation (involving boundary derivatives) of the second type was developed to facilitate iterative improvement; this equation was used to refine the boundary estimates at the critical horizon end of the procedure proposed in \cite{Barles}.

\par
The policy iteration algorithm operates directly by appropriately modifying the boundary. The arbitrage value of a policy specified by a given boundary in the continuation region (above the boundary) is the solution of the Black--Scholes equation with the appropriate boundary condition. Even though a policy modification in which the profit from an immediate exercise exceeds the continuation would be beneficial, the scope of such modifications is limited to the continuation region. A more efficient approach would be to examine the Black--Scholes solution that satisfies the same boundary condition for the entire $x,t$ region, assuming its existence, and to then update the boundary to one in which immediate exercise is advantageous either in the continuation or in the stopping region. Via the maximum principle property, the new policy will represent an improvement in the continuation region, and by a careful selection of the modification, the new policy can be shown to be an improvement for all $x,t$. Such modifications are possible until the smooth pasting condition is satisfied. By selecting the new boundary greedily to maximise the benefit from the change, fast convergence to the optimal boundary can be achieved.

\section{The Black--Scholes Market Model}

\begin{defin}\label{american:black sholes model market}
The Black--Scholes market model consists of two assets $B$ and $x$, whose dynamics are given by
\begin{align} 
dB(t) = rB(t)dt, \\
dx(t) = \mu x(t)dt+\sigma x(t)dw(t), \label{american:asset price dif equatio}
\end{align}
where $r$, $\mu$, and $\sigma$ are known constants. $B$ is a bond with interest rate $r$, $x$ is the price of an asset, and $w$ describes the standard Brownian motion. This model does not permit arbitrage opportunities. 
\end{defin}
\begin{thm}
Suppose that in the market model above, we want to price a contingent claim of the form  $\xi = \Phi(x(T))$. Then, the only price function $V(x,t)$ for this that is consistent with the absence of arbitrage is the solution of the following boundary value problem:
\begin{align}
rV=V_t+rxV_x+\frac{1}{2}\sigma^2x^2V_{xx}, \\
V(x,T)=\Phi(x).
\end{align}
The solution of this system in the domain  $[0 \: T]\times \mathbb{R}$ gives for the European-type option, 
\begin{align*}
V(s,t)=e^{-r(T-t)}E^{Q}[\Phi(x(T))|x(t)=s],
\end{align*}
where $Q$ is a probability measure defined such that the price of the underlying asset has dynamics of the form
\begin{align*}
dx(t) = r x(t)dt+\sigma x(t)dW(t).
\end{align*}
$W(t)$ is used for Brownian motion specified by the probability measure $Q$.
\end{thm}
That is, the price of the European option, at the moment $t$ and for the underlying price $s$, is the discounted mean value of its payoff function under the probability measure $Q$, given that the asset price at time $t$ is $s$. Because the American-type option can be exercised at any time $\tau$ up to the expiration date $T$, the holder can choose a time of exercise that maximises their profit. Therefore, a fair price in that case should satisfy the following condition:
\begin{align*}
V(s,t)=\underset{\tau}{sup}\:e^{-r(\tau -t)}E^{Q}[\Phi(x(\tau))|x(t)=s], \: \tau \in [t\:T].
\end{align*}

\section{Policy Iteration Algorithm}

As mentioned before, an American put option with exercise price $K$ and expiration date $T$ gives its holder the right to sell a specified asset, any time before the expiration date at the price of $K$. Assuming that the holder of the option does not own the asset, but buys it at the current market price $x_t$, their profit will be $K-x_t$. Because the right will be exercised only if the current market price is lower than $K$, the profit can be expressed as $[K-x_t]^+$, where $[z]^+$ is zero if $z$ is negative and $z$ otherwise. The asset price $x_t$ is given by Equation \eqref{american:asset price dif equatio}. 

\par
A policy $\pi$ is defined by a boundary; that is, a function $b(t)$ over $[0,T]$ with $b(t)\leq K$. The policy implies that if the price of the asset at time $t$ is greater than the boundary $b(t)$, the holder does not exercise his right until the stopping time
\[ \tau^{\pi} = \min \{t : x_t \leq b(t) \}, \]
whereas they exercise their right immediately if the asset price is lower than the boundary $b(t)$.
\par
Now, we define the {\it continuation region} and {\it stopping region} as $\{(x,t):b(t)\leq x,\: t\in [0,\:T] \}$ and $\{(x,t):x<b(t),\:t\in [0,\:T] \}$, respectively.
In the Black--Scholes market model (Definition \ref{american:black sholes model market}), such a right, with specified stopping time $\tau^{\pi}$, can be replicated at time $t$ by a portfolio whose initial value is $V^{\pi}(x,t)$ ($x$ is the asset price at time $t$) where $V^{\pi}$ satisfies the Black--Scholes equation in the continuation region:

\begin{equation}\label{BS}
rV^{\pi}=V^{\pi}_t+rxV^{\pi}_x+\frac{1}{2}\sigma^2x^2V^{\pi}_{xx},
\end{equation}
\[
V^{\pi}(b(t),t)=K-b(t), \quad V^{\pi}(\infty,t)=0, \quad V(x,T)=0, \quad \forall \{T,x: x\geq b(T)\}.
\]

Thus, we define the replication value for policy $\pi$ \textbf{  $U^{\pi}(x,t)$} 
\[
U^\pi (x,t)=\left\{
\begin{array}{lrr}
V^{\pi}(x,t),&x > b(t) & \text{Continuation region}\\
 \left[ K-x\right]^+, &x\leq b(t)& \text{Stopping region.}
\end{array}\right\}
\]
The {\it optimal policy} is defined by the {\it optimal boundary } $b^*$, which must satisfy a {\it smooth pasting condition} \cite{H.McK:65}. $V^*$ is the solution of \eqref{BS} with boundary $b^*$; thus,  $V^*_x(b^*(t),t)=-1$ for every $t$. Therefore, $U^{\pi^*}$ has a continuous $x$ derivative.
\par
Using the maximum principle property for parabolic equations, a better policy can be obtained by modifying the boundary $b$ to a new one $\gamma$ for which immediate exercise is better than continuation; that is, $V^\pi (\gamma(t),t)\leq [K-\gamma(t)]^+$ for all $t$. This stopping time policy will be referred as $q$. In this formulation, the new boundary must be within the continuation region where $V^\pi$ is defined. However, assuming that the solution of \eqref{BS} can be smoothly extended to all $x\ge 0$ and $t\le T$, modifications can be considered in either the continuation or stopping regions. This extension assumption means that a function $\tilde{V}$ satisfying \eqref{BS} for all $x\ge0\quad t\le T$ such that $\tilde{V}=V^\pi$ for $x\ge b(t)$ can be identified. 

It is expected that even for the boundary modification obtained using this extended solution of \eqref{BS}, the continuation value $V^q$ dominates that of $V^\pi$ in the continuation region of $q$, as shown in the following theorem, which shows that by a proper choice of $\gamma$, the improvement is global in the replication values; that is, $U^q(x,t) \geq U^\pi(x,t),\;\forall x,t.$

\begin{thm}\label{Monot}
Let policy $\pi$ be specified by a boundary $b$, and assume that the extended Black--Scholes problem, 
\begin{equation}\label{BSBY}
\begin{split}
rV =V_t+rxV_x+\frac{1}{2}\sigma^2x^2V_{xx}, &\quad \mbox{with}  \\ V^{\pi}(b(t),t)=K-b(t), \quad V^{\pi}(\infty,t)=0, \quad &V(x,T)=0, \quad \forall \{t,x: x\geq b(T)\} ,
\end{split}
\end{equation}	
has a solution $V^{\pi}$ defined for all \textbf{$x \ge 0, \; 0 \leq t \leq T$} .

Consider a policy $q$ with boundary $\gamma$ such that  $V^\pi(\gamma(t),t)\leq [K-\gamma(t)]^+,\;\forall t$. Then, assuming the existence of the solution $V^q$ of the extended Black--Scholes problem for $\gamma$, the following holds:
 \begin{itemize}
\item [\rm{i)}]  $V^\pi (x,t) \leq V^q(x,t),$  for $x\geq \gamma(t),$ 
\item [\rm{ii)}] For a proper selection of $\gamma$:\; $U^\pi (x,t)\leq U^q(x,t),\;\;\forall (x,t)$.
\end{itemize}
\end{thm}

\begin{prf*}
From the definition of $\gamma$, we have that $K-\gamma(t)\geq V^\pi (\gamma(t),t)$ for every $t$.  Let us consider $V^q$ and $V^\pi$ at the boundary $ \gamma$. From the previous inequality, we have
\begin{align*} 
V^q(\gamma(t),t) = K-\gamma(t)- V^\pi(\gamma(t),t)+ V^\pi(\gamma(t),t)  \geq  V^\pi(\gamma(t),t) \Rightarrow \\
V^q(\gamma(t),t) -V^\pi(\gamma(t),t)\geq0.
\end{align*}
The functions $V^q,\;V^\pi$ are solutions of the Black--Scholes equation that can be transformed to the heat equation by a monotonic transformation \cite{Koh:14}. 
For the transformed equations (and, subsequently, for $V^q(\gamma(t),t) -V^\pi(\gamma(t),t)$) we apply the maximum and minimum principles in the region 
\[
\Omega_{\gamma}=\{(x,y):x\ge \gamma(t) , 0\le t \le T\}.
\] 
In the rest of the boundary, applies $V^q(x,t) -V^\pi(x,t)=0$ because  
\[
V^q(\infty ,T)=V^\pi(\infty ,t)=0,\;V^q(x,T)=V^\pi(x,T)=0 \quad x\ge K.
\] 
Therefore, according to the minimum principle, we have $V^q(\gamma(t),t) -V^\pi(\gamma(t),t)\geq 0$ everywhere in the interior of $\Omega_{\gamma}$\footnote{The $\Omega_{\gamma}$ is not bounded or closed; however, we can limit the bounded subset because $\lim_{x\rightarrow \infty}V(x,t)=0$}.

\rm{(ii)}  Supposing that $b$ is not optimal, we have that $V_x^\pi(b(t),t) \ne -1$ for some $t$. Then, we obtain the intervals  $(\bar{x},b(t))$ or $(b(t),\breve{x}),$ where $V^\pi(x,t) \le \left[K-x\right]^+$. Specifically, if $V_x^\pi(b(t_0),t_0 > -1$ applies for some $t_0$; then, the continuity of $V_x$ in an open set around $t_0$ means that $V_x^\pi(b(t),t) > -1$. In this set, we have the extended function $V^\pi(x,t)<K-x$ for $x<b(t)$; thus, for this set, we choose $\gamma(t)<b(t)$. Then, for $x\leq \gamma(t)$,
\[
U^q(x,t)=K-x=U^\pi(x,t),
\]
whereas for $\gamma(t)< x \leq b(t)$  
 \[
 U^q(x,t)>K-x=U^\pi(x,t).
 \]
 The inequality $U^q(x,t)>K-x$ in this specific set applies because $U^q(\gamma(t),t)=V^q(\gamma(t),t)=K-\gamma(t)$, and because $\gamma$ is optimal, we have that $V^q_x(\gamma(t),t)=-1$. Considering the last equation and the fact that $V^q(x,t)$ is convex \cite{BrenSchw}, we deduce that 
 \[
 U^q(x,t)=V^q(x,t)>K-x=U^\pi(x,t)
 \]
 for $\gamma(t)< x \leq b(t)$.

If applies that $V_x^\pi(b(t_0),t_0) < -1$ for some $t$, the continuity of $V_x$ in an open set around $t_0$ means that $V_x^\pi(b(t),t) < -1$ in all set. In this set, we have for the extended function $V^\pi(x,t)<K-x$ for $b(t)<x\leq \gamma(t)$, and we choose $b(t)<\gamma(t)$. Then, for $x\leq b(t)$, we have that 
\[
U^q(x,t)=K-x=U^\pi(x,t);
\] 
Meanwhile, when $b(t)<x\leq \gamma(t)$ applies, 
\[
U^q(x,t)=K-x>V^\pi(x,t)=U^\pi(x,t)
\].

In the regions where $V_x^\pi(b(t),t) = -1$, we set $\gamma(t)=b(t)$. Finally, to verify the inequality $U^q(x,t)\geq U^\pi(x,t)$ in the region $[\max(\gamma(t),b(t)),\infty]\times [0,\:T]$, we apply the minimum principle for the function $V^q(x,t)- V^q(x,t)$ in that region, as in \rm{(i)}.

\end{prf*}

A greedy improving strategy would be to exercise upon first reaching an asset value that locally maximises the improvement; namely, $\arg \max_{x,locally} \{ [K-x]^+ -V^{\pi}(x,t)\}$. It is not mandatory to take into account the properties stated in part ii) of Theorem \ref{Monot} because a strict maximisation $V$ will lead to an optimal boundary. Sequentially applying this boundary update, we obtain the following algorithm, which is in the spirit of  policy iteration \cite{Bellman} because the updates rely on the solution of the extended Black--Scholes problem \eqref{BSBY}, which is closely related to the replication value of the current policy.  The algorithm is described in the following steps: \\

\textbf{Policy Iteration Algorithm - PIA}
\begin{enumerate}
        \item Select an arbitrary stopping time policy $\pi_0$ for a boundary $b_0.$
\item  Compute $V_0$ for all $(x,t)$ via the extended Black--Scholes problem \eqref{BSBY} on $b_0$.
\item $i \gets 0 $ 
\item Repeat
\begin{enumerate}
\item $b_{i+1}(t)\gets \arg\max_x \{ [K-x]^{+}\!\!-V_i(x,t) \} $. \textit{Comment}: A local maximum is chosen, preferably but not necessarily satisfying the conditions in Theorem \ref{Monot}(b).
\item Compute $V_{i+1}$ via the extended Black--Scholes problem \eqref{BSBY} on $b_{i+1}$.
\item $i\gets i+1$
\end{enumerate}  
\item {Until $\|V_{i+1}-V_i\|_{\infty} \leq \epsilon$ \textbf{or} $\| b_{i+1}-b_i\|_{\infty} \leq  \epsilon $ } with $\epsilon$ being a desired tolerance.
\end{enumerate}

The $V_i$ values obtained when using the above algorithm monotonically improve in the respective continuation regions as a consequence of Theorem \ref{Monot}, and so do the corresponding replication values $U_i$, provided that the boundaries are properly chosen. Because the $V_i$ functions are bounded, it is expected that they converge to a function satisfying the smooth pasting condition; however, no formal arguments to that effect are presented here. Instead, it is shown that under reasonable assumptions, the sequence of the generated boundaries converges fast to the optimal one. 

\section{Convergence Properties of the {\bf PIA}}\label{ConvergSect}
In the context of Markovian decision processes and control theory, it has been highlighted that the policy iteration algorithm is related to Newton's root-finding method (see, for instance, \cite{Put:04}) and \cite{E.Mag:77}). Given the quadratic convergence of Newton's method, it is reasonable to ask whether {\bf PIA} shows similar characteristics. We will  present such a property, whose justification hinges on the assumption we state in the following paragraphs. The results of the quadratic convergence of a policy iteration algorithm are also given in \cite{JackaVolat}; however, our results differ because they concern the convergence of boundaries and not of values.

Consider a smooth boundary $b(t)$ and two solutions of the Black--Scholes equation which vanish at infinity and on their boundaries take values given by two functions $f(t),g(t)$: $V^f(b(t),t)=f(t)$ and $V^g(b(t),t)=g(t)$. Consider the partial derivatives $\frac{\partial{V^f(b(t),t)}}{\partial{x}}$ and $\frac{\partial{V^g(b(t),t)}}{\partial{x}}$ on the boundary. For functions $ f,g $ mutually close in (say) the maximum norm, we consider the difference in their partial derivatives on the common boundary. We assume that for functions $f,g$, 
\begin{equation}\label{partials}
\left\| \frac{\partial{V^f(b(t),t)}}{\partial{x}}-\frac{\partial{V^g(b(t),t)}}{\partial{x}}\right\|_{\infty}=O(\|f-g\|_{\infty}),
\end{equation}
where we recall that the Landau symbol $O$ indicates that there is a positive constant $c$ such that for a sufficiently small $\|f-g\|_{\infty}$, the derivative difference $\| \frac{\partial{V^f(b(t),t)}}{\partial{x}}-\frac{\partial{V^g(b(t),t)}}{\partial{x}}\|_{\infty}$ is less than or equal  to $c\|f-g\|_{\infty}$. For the remainder this assumption will be referred to as \emph{The Assumption of the Partials} ({\bf AP}). The proof of the quadratic convergence of the algorithm relies on {\bf AP}, which is not valid in general. However, the assumption is valid when the functions $f,g$ at the boundary $b$ are multiples of some function $h$; that is, $f(t)-g(t)=\lambda h(t)$, as in the numerical results shown below. 

\begin{thm} \label{Converg}
We considered the optimal exercise boundary $b^*$. Under {\bf AP}, boundaries $b_i\!$ s generated by {\bf PIA} satisfy 
\[
\|b_{i+1}(t)-b^*(t)\|_{\infty}=O(\|b_i(t)-b^*(t)\|_{\infty}^2);
\]
that is, the convergence of the algorithm is quadratic.
\end{thm}

\begin{prf*}
        \begin{figure}[!ht]
                \centering
                \includegraphics[width=10cm, height=8cm]{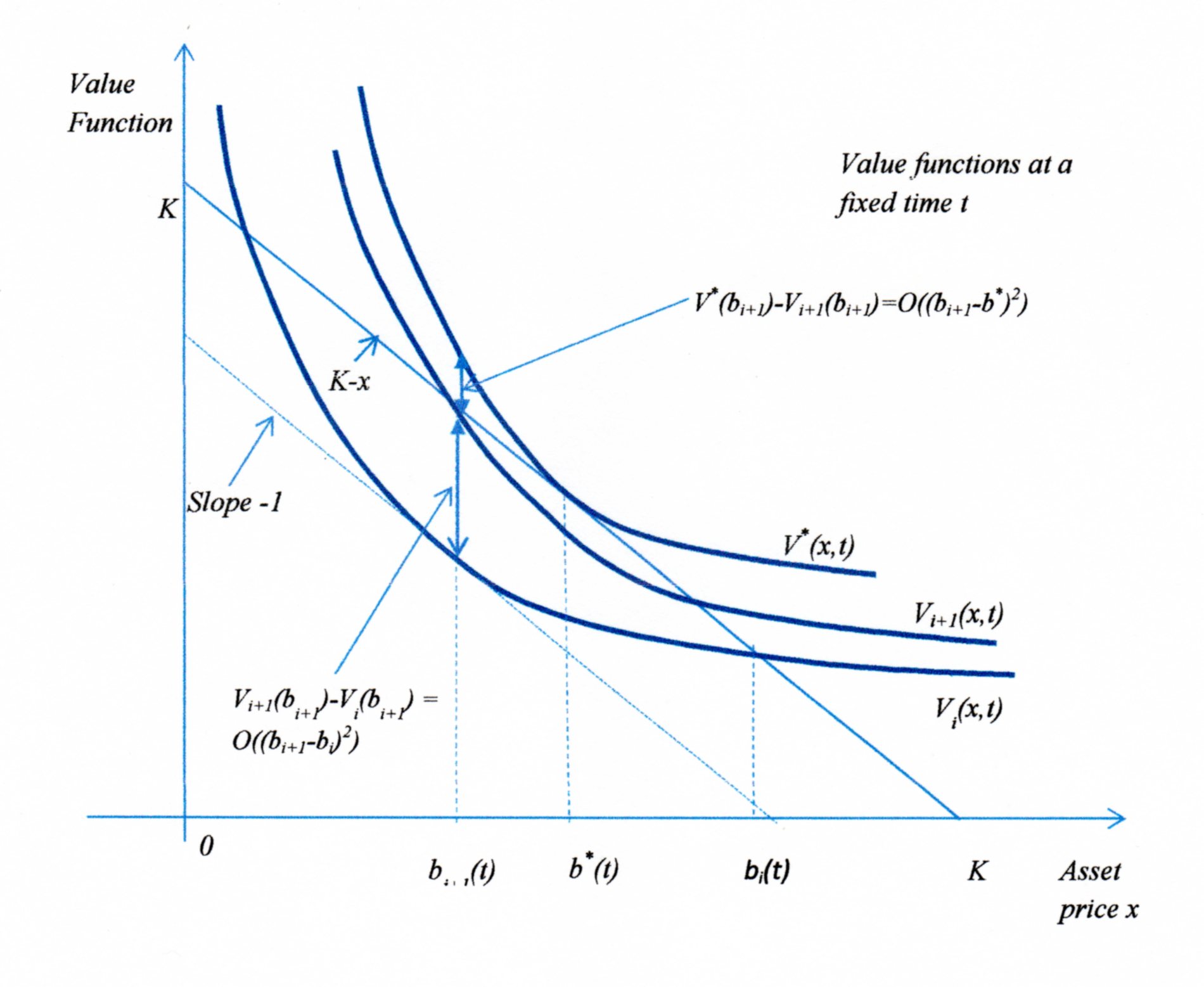}
                \caption{The steps of the \textbf{PIA} algorithm}
                \label{fig:1}
        \end{figure}
The solution of the Black--Scholes problem $V^*$ with boundary condition $V^*(b^*(t),t)=K-b^*(t)$ satisfies the smooth pasting $\frac{\partial{V^*(b^*(t),t)}}{\partial{x}}=-1$, and vanishes as $x$ tends to infinity. (See Figure {\ref{fig:1}}). Let $b_i$ and $b_{i+1}$ be two successive boundaries, as specified by the \textbf{PIA} algorithm.  Then, from the construction of $b_{i+1}$, we apply $V_{i,x}(b_{i+1}(t),t)=-1$. From the Taylor theorem for some $\hat{b}$, we have that
\[
V_i(b_i(t),t)=V_i(b_{i+1}(t),t)+ V_{i,x}(b_{i+1}(t),t)(b_i(t)-b_{i+1}(t))+\frac{1}{2}V_{i,xx}(\hat{b},t)(b_i(t)-b_{i+1}(t))^2.
\]
Thus, 
\begin{align*}
V_i(b_{i+1}(t),t)&=K-b_i(t)+(b_i(t)-b_{i+1}(t))-\frac{1}{2}V_{i,xx}(\hat{b},t)(b_i(t)-b_{i+1}(t))^2\\
                 &= K-b_{i+1}(t)-\frac{1}{2}V_{i,xx}(\hat{b},t)(b_i(t)-b_{i+1}(t))^2,
\end{align*}
and finally
\begin{equation}\label{con1}
V_{i+1}(b_{i+1}(t),t)-V_i(b_{i+1}(t),t)=\frac{1}{2}V_{i,xx}(\hat{b},t)(b_i(t)-b_{i+1}(t))^2.
\end{equation}

Furthermore, for some $\dot{b}$, we have from Taylor's theorem that
\[
V^*(b_{i+1}(t),t)=V^*(b^*(t),t)+V^*_x(b^*(t),t)(b_{i+1}(t)-b^*(t))+\frac{1}{2}V^*_{xx}(\dot{b},t)(b_{i+1}(t)-b^*(t))^2.
\]
Because $V_x^*(b^*,t)=-1$ and $V^*(b^*,t)=K-b^*$ from the previous equation, we take
\begin{equation}\label{con2}
V^*(b_{i+1}(t),t)-V_{i+1}(b_{i+1},t)=\frac{1}{2}V^*_{xx}(\dot{b},t)(b_{i+1}(t)-b^*(t))^2.
\end{equation}
Applying Taylor's theory to the partial $x$ derivatives, we have that for some $\bar{b}$
\[
 V^*_x (b_{i+1}(t),t)=V^*_x (b^*(t),t)+V_{xx}^* (\bar{b},t)(b_{i+1}(t)-b^*(t)).
\]   

Because $V^*_x (b^*(t),t)$ = $V_{i,x}(b_{i+1}(t),t)$ = $-1$, subtracting the last equation from $V_{i+1}$ at $b_{i+1}$, we have
\begin{align*}
V_{i+1,x}(b_{i+1}(t),t)&=V_{i+1,x}(b_{i+1}(t),t)-V^*_x (b_{i+1}(t),t)+V_{i,x} (b_{i+1}(t),t) \\
&+V_{xx}^* (\bar{b},t)(b_{i+1}(t)-b^*(t)),
\end{align*}
and so 
\begin{align*}
 V_{xx}^* (\bar{b},t)(b_{i+1}(t)-b^*(t))&=V_{i+1,x}(b_{i+1}(t),t)-V_{i,x} (b_{i+1}(t),t)\\
 &+V^*_x (b_{i+1}(t),t)-V_{i+1,x}(b_{i+1}(t),t).
 \end{align*}
Taking absolute values
\begin{align*}
\left | V_{xx}^* (\bar{b},t)(b_{i+1}(t)-b^*(t))\right | &\leq \left | V_{i+1,x}(b_{i+1}(t),t)-V_{i,x} (b_{i+1}(t),t)\right |\\
                                                        &+ \left | V_{i+1,x}(b_{i+1}(t),t)-V^*_x (b_{i+1}(t),t)\right |.
\end{align*}   

Equations \eqref{con1} and \eqref{con2} show that on the boundary $b_{i+1}$, the functions $V_{i+1}$ and $V^*$ differ by an order of  $[b_{i+1}(t)-b^*(t)]^2$, whilst $V_{i}$ and $V_{i+1}$ differ by an order of $[b_{i+1}(t)-b_i(t)]^2$ . Note (see \cite{J.Goo:02}) that at $b^*$ the second derivative of the optimal value function equals $V^*_{xx}(b^*(t),t)=\frac{2rK}{\sigma^2 b^{*2}}$ and is therefore bounded in the vicinity of $b^*$ for the interval $[0,T]$. Thus, for sufficiently small differences in the boundaries, and by virtue of {\bf AP} \eqref{partials}, there exist constants $\alpha$, $c_1$, and $c_2$ such that

\begin{equation}\label{con3}
\alpha\lvert b_{i+1}(t) -b^*(t) \rvert \leq c_1 (b_{i+1}(t) -b_i(t) )^2+c_2(b_{i+1}(t) -b^*(t) )^2.
\end{equation}
By writing $\left(b_{i+1}(t)-b_i(t)\right)^2$ as $\left(b_{i+1}(t) -b^*(t)+b^*(t)-b_i(t)\right)^2$ and expanding the previous inequality, Equation \eqref{con3} becomes
\[
\alpha\lvert b_{i+1}(t) -b^*(t) \rvert - (c_1+c_2)(b_{i+1}(t) -b^*(t) )^2 +2c_1 (b_{i+1}(t) -b^*(t))(b_i(t) -b^*(t)) \leq c_1 (b_i(t) -b^*(t) )^2.
\]
Because $|b_{i+1}(t) -b^*(t)|$ and $|b_i(t) -b^*(t)|$ are sufficiently small, the above inequality becomes
\begin{equation}\label{con4}
\alpha\lvert b_{i+1}(t) -b^*(t) \rvert  \leq c_1 (b_{i}(t) -b^*(t) )^2.
\end{equation}
Supposing that $ \alpha$, $c_1$, and $c_2$ are bounded away from $0$ and $\infty$, the above inequality can be written as
\begin{equation}
 \left\| b_{i+1}(t) -b^*(t) \right\|_{\infty} \leq c \left\| b_{i}(t) -b^*(t) \right\|^2_{\infty},
\end{equation}
where $c=\frac{c_1}{\alpha}$. 

\end{prf*}

As stated earlier, the application of \textbf{AP} in the previous proof is justified if the values $V_i(x,t)$ generated are eventually of the form $V_i(x,t)=V^*(x,t)+\lambda_i h(x,t)$ for some function $h$. The constants involved in the \textbf{AP} depend on the boundary $b_i$ but are independent of $i$ when the boundaries are close to $b^*$. The computations reported in Section \ref{SectComput} corroborate the applicability of this assumption.

\section[Policy Iteration Algorithm for Stochastic Control Problems]{A Policy Iteration Algorithm for Free Boundary Stochastic Control Problems} \label{FBProblem}

The policy improvement used in the previous section relied on examining the extension of the value function inside the stopping region, which makes sense if we want to assess delaying the exercise.  The same principle can be applied  to general free  boundary control problems and such a procedure is presented here, first for deterministic systems and then for stochastic ones.  Several simplifying assumptions are made and only a monotonicity result is given while value or policy convergence is not examined.

Consider the control problem to determine
\begin{equation*}
\sup_{u,T}\int_{t_0}^{T}f(x,\tau,u)d\tau+F(x_T,T),
\end{equation*}
with $x\in \mathbb{R}^n$  satisfying the DE 
\begin{equation}\label{detsys}
\frac{dx}{dt}=g(x,t,u)\;\; \; x(t_0)=x_0.
\end{equation}

The end point $T$ is freely chosen. It is assumed that $u=u(x,t)$ is an acceptable control in that the resulting differential equation \eqref{detsys} has a solution.  

A \textit{stopping policy} $\pi^k = \left\lbrace u_k(x,t),\Delta^k_t\right\rbrace $ consists of an acceptable  control $u_k$  and a collection of stopping regions $\Delta^k_t\subset \mathbb{R}^n$ for $t\ge t_0$ as well as continuation regions $C^k_t=\mathbb{R}^n-\Delta^k_t$.  Applying $\pi^k$ for a $x_0$ in the continuation region for $t_0$ consists of using the specified control to obtain a trajectory $x^k_t$ and stopping when it first enters a stopping region, i.e. at $T^k=\min_t\left\lbrace t \ge t_0\;|\;x^k_t\in \Delta^k_t\right\rbrace$, assuming   $T^k$ to be finite. The value of  $\pi  _k,\;V^k(x,t)$, is defined in the continuation region as 
\begin{equation}\label{PIValuefcn}
V^k(x_0,t_0)=\int_{t_0}^{T^k}f(x^k,\tau,u^k) d\tau+F(x_{T^k},T^k).
\end{equation}
For a $x_0$ in the stopping region of $t_0$ the value is by definition $F(x_0,t_0)$ and thus the value $U^k(x,t)$ achieved by the stopping policy $\pi^k$ is 
\[
U^k (x,t)=\left\{
\begin{array}{lrr}
V^k(x,t),&x \in C^k_t, & Continuation\;Region,\\
F(x,t),&x \in \Delta^k_t, &Stopping\;Region.
\end{array}\right.
\] 
To extend the continuation value $V^k$ in the stopping region we assume that for $(x_0,t_0)$ in the stopping region there is a prior time $T^k$ and a corresponding state $x_{T^k}$ on the boundary of $\Delta^k_{T^k}$ such that the process moves using the specidied control from $x_{T^k}$ to $x_0$ at time $t_0>T^k$ while staying inside the stopping region.  The $V^k$ is again given by \eqref{PIValuefcn}.  
 
A new policy $\pi^{k+1}= \left\lbrace u_{k+1}(x,t),\Delta^{k+1}_t\right\rbrace $ is an {\it improvement} provided it is both a control and a stopping region improvement, namely the following conditions hold:
 \begin{description}
 \item [a) ] {\it Control improvement.}
 \begin{equation}
 \!\!\!\!\!\! f(x,t,u^{k+1})+g(x,t,u^{k+1})V^k_x(x,t) \ge  f(x,t,u^{k})+g(x,t,u^{k})V^k_x(x,t) \; \; \forall x,t. 
 \end{equation}
 A greedy choice would be to select the control that maximizes the right hand side.
 \item [b)]  {\it Stopping region improvement}.
 In the new termination region the immediate exercise value must be greater than the continuation using $\pi_k$, namely 
  $$\Delta^{k+1}_t \subseteq \left\lbrace x \mid F(x,t) \ge V^k(x,t)\right\rbrace.$$
A greedy choice would be   $$\Delta^{k+1}_t= \left\lbrace x\mid x=\arg\max_x [F(x,t)-V^k(x,t)]\right\rbrace.$$ 
This particular choice makes stopping difficult and must be proven consistent with the previously stated requirement of finite stopping times $T^k < \infty$.

 \item [c) ] {\it Termination.}
 The algorithm stops if it is not possible to improve on either criterion, namely the control $u^k$ maximizes the Hamiltonian $f+gV^k_x$ and $ F(x,t)\le V^k(x,t) \; \forall x,t$. Approximate termination criteria could involve $\|u_{k+1}-u_k\|$ and/or $\|V^{k+1}-V^k\|$.
\end{description}
The value of the updated policy $V^{k+1}$  is an improvement, i.e. $V^{k+1}(x,t) \ge V^{k}(x,t)$ for all $(x,t)$.  To show this,  it is noted that for acceptable $u^{k+1},u^k$ the functions $V^{k+1},V^k$ satisfy the value PDE's:

\begin{equation} \label{valuepde}
 f(x,t,u^l)+g(x,t,u^l)V^l_x(x,t)+V^l_t(x,t)=0\quad l=k,k+1.
\end{equation}
Given the choice of $u^{k+1}$  we have the following relations:

\begin{equation}
\begin{split}
f(x,t,u^{k+1})+&g(x,t,u^{k+1})V^k_x(x,t)+V^k_t(x,t) \ge  \\
&f(x,t,u^{k})+g(x,t,u^{k})V^k_x(x,t)+V^k_t(x,t) =0
\end{split}
\end{equation}

Substituting $f$ from \eqref{valuepde} with $l=k+1,$ we have
\begin{equation} \label{contrinequ}
\left[V^k_t(x,t)-V^{k+1}_t(x,t)\right] + g(x,t,u^{k+1})\left[V^k_x(x,t)-V^{k+1}_x (x,t)  \right] \ge 0.
\end{equation}
Applying the control $u^{k+1}$ starting at $(t_0,x_0)$ in the continuation region, we obtain a trajectory $x^{k+1}$ on which the difference $\left[V^k_t(x,t^{k+1})-V^{k+1}_t(x,t^{k+1})\right]$ equals  the left hand side of \eqref{contrinequ}.  Recalling  the assumption  that the termination time $T^{k+1}$ is finite and  integrating  \eqref{contrinequ} on $x^{k+1}$ from $t_0$ to  $T^{k+1}$ we obtain
\begin{equation}
V^k (x^{k+1}_{T^{k+1}},T^{k+1})-V^{k+1}(x^{k+1}_{T^{k+1}},T^{k+1})-V^k(x_0,t_0)+V^{k+1}(x_0,t_0)\ge 0.
\end{equation}

By the choice of $T^{k+1}$ the difference of the first two terms is non-positive and thus  
$V^{k+1}(x_0,t_0)\ge V^k(x_0,t_0)$.  To show the inequality for an  initial point $(x_0,t_0)$ inside the stopping region the  calculation can be repeated starting from the prior point $(x^{k+1}_{T^{k+1}},T^{k+1})$ from which the control $u^{k+1}$ drives the system to $(t_0,x_0)$. 

As in Theorem \ref{Monot} ii), we can choose the stopping region improvement such that the value of the stopping policies is everywhere nondecreasing i.e. $U^{k+1}(x,t)\ge U^k(x,t)$.  This can be done if the stopping region modifications $\left[\Delta^{k+1}_t-\Delta^{k}_t\right] \cup \left[ \Delta^{k}_t-\Delta^{k+1}_t \right]$ are a subset of the connected region in which $V^k\le F$ and includes the boundary of $\Delta^{k}_t$, and then, if necessary, perform the second modification in Theorem \ref{Monot} ii.

The stochastic control case is similar, but the assumptions required are stricter.  Consider a stochastic control system
\begin{equation} \label{sde}
dx=g(x,t,u)dt+\sigma(x,t,u)dz\;\;\;x(t_0)=x_0.
\end{equation}
We want to determine 
\begin{equation}
\sup_{u,T}E\left[\int_{t_o}^{T}f(x,\tau,u)d\tau+F(x_T,T)\right].
\end{equation}
We use the same definition of a policy $\pi^k = \left\lbrace u_k(x,t),\Delta^k_t\right\rbrace$  as in the deterministic case and consider the value function 
\begin{equation}\label{PIStochValfcn}
V^k(x_0,t_0)=E\left[ \int_{t_0}^{T^k}f(x^k,\tau,u^k)d\tau+F(x_{T^k},T^k) \mid x(t_0)=x_0 \right].
\end{equation}
Starting outside the termination region the stopping time is given by $$T^k=\min_t\left\lbrace t \ge t_0 \mid x^k_t\in \Delta^k_t \right\rbrace .$$ 
On the other hand if $(x_0,t_0)$ is inside the termination region, we consider the set of points $(x_{T^k},T^k)$ on the termination region boundaries from which the specified control leads $x^{k}$ to  $(x_0,t_0)$  while staying inside the termination region. Expectation is taken over the set consisting of these points $(x_{T^k},T^k)$ conditional on the process  reaching$(x_0,t_0)$.  

As in the deterministic case  a policy $\pi  _{k+1}= \left\lbrace u_{k+1}(x,t),\Delta^{k+1}_t\right\rbrace $ is an \textit{improvement}  provided it is both a control and a stopping region improvement, namely: 
\begin{description}
 \item [a) ] {\it Control improvement.}
\begin{align}
f(x,t,u^{k+1})+g(x,t,u^{k+1})V^k_x(x,t)+\frac{1}{2}\sigma^2(x,t,u^{k+1})V^{k}_{xx}(x,t)  &\ge& \nonumber \\
 f(x,t,u^{k})+g(x,t,u^{k})V^k_x(x,t) +\frac{1}{2}\sigma^2(x,t,u^{k})V^{k}_{xx}(x,t) . &&
 \end{align}
\item [b) ] Stopping region improvement. 
The immediate termination value must be greater than the continuation using $\pi_k$, namely 
\[
\Delta^{k+1}_t \subseteq \left\lbrace x\mid F(x,t) \ge V^k(x,t)\right\rbrace .
\]
If $\||F-V^k\|_{\infty}=M,$ then one could select  as stopping region the $x$'s for which the difference $F-V^k$ is greater than $M-\epsilon$. A  greedy choice would be   $$\Delta^{k+1}_t= \left\lbrace x \mid x=\arg \max_x [F(x,t)\ge V^k(x,t)]\right\rbrace .$$ 

This choice would make stopping difficult and must be shown consistent with the requirement of finite expected stopping times $E(T^k)<\infty$ which later in this section will be shown necessary for the algorithm to be improving.

\item [c) ]Termination: As in the deterministic system.

\end{description} 

We assume that the controls $u^k$ lead to strong solutions of the stochastic system, a complication dealt in detail in \cite{JackaPI}.  Then, the values $V^{k+1},V^k$ corresponding to the policies $\pi_{k+1},\pi_k,$ satisfy the PDE's 
\begin{align} 
f(x,t,u^k)+g(x,t,u^k)V^k_x(x,t)+\frac{1}{2}\sigma^2(x,t,u^k)V^l_{xx}(x,t)+V^k_t(x,t)=0, \label{svaluepde_k} \\
f(x,t,u^{k+1})+g(x,t,u^{k+1})V^{k+1}_x(x,t)+\frac{1}{2}\sigma^2(x,t,u^{k+1})V^l_{xx}(x,t)+V^{k+1}_t(x,t)=0. \label{svaluepde_k_1}
\end{align}
The conditions for the boundary of the corresponding stopping region for the above equations are $F(x,t)=V^k(x,t)$ and $F(x,t)=V^{k+1}(x,t)$ respectively. We assume there exist such solutions. From the control boundary \eqref{contrinequ}  applies that, 
\begin{align*}
&f(x,t,u^{k+1})+g(x,t,u^{k+1})V^k_x(x,t)+\frac{1}{2}\sigma^2(x,t,u^{k+1})V^{k}_{xx}(x,t) +V^k_t(x,t) \geq \\
&f(x,t,u^{k})+g(x,t,u^{k})V^k_x(x,t) +\frac{1}{2}\sigma^2(x,t,u^{k})V^{k}_{xx}(x,t)+V^k_t(x,t)=0.
\end{align*}
If we replace the term $f(x,t,u^{k+1})$ according to \eqref{svaluepde_k_1} in the left of the above inequality we have
\begin{align}
\left[V^k_t(x,t)-V^{k+1}_t(x,t)\right] + g(x,t,u^{k+1})\left[V^k_x(x,t)-V^{k+1}_x(x,t) \right] &+ \nonumber \\
 \sigma^2(x,t,u^{k+1})/2\left[V^k_{xx}(x,t)-V^{k+1}_{xx}(x,t)\right]&\ge0. \label{scontrinequ}
\end{align}
If $x^{k+1}$ is the solution of the stochastic system \eqref{sde}, when the control is  $u^{k+1}$, we then consider the difference to the successive values  $W^k(x^{k+1},t)=V^k(x^{k+1},t)-V^{k+1}(x^{k+1},t)$. Then, from Ito's Lemma and assuming a strong solution for the differential equation
\begin{equation}\label{sdp}
dW^k=(W^k_t+gW^k_x+\frac{\sigma^2}{2}W^k_{xx})dt+\sigma W^k_xdz.
\end{equation}
By virtue of \eqref{scontrinequ}, the coefficient appearing in $dt$ is nonnegative. Therefore,  integrating \eqref{sdp} from $t_0$ to $T^{k+1}$ we obtain the inequality 
\begin{equation}
W^k(x_{T^{k+1}},T^{k+1})-W^k(x_0,t_0)\ge \int_{t_0}^{T^{k+1}}\sigma W^k_xdz.
\end{equation}
Again by the construction of the stopping time $T^{k+1},$ the first term is non positive. Furthermore, taking expectations and assuming $E(T^{k+1}) <\infty$ we have by Dynkin's lemma \cite{Wil91}  that the expectation of the right hand integral vanishes. Consequently,  $W^k(x_0,t_0) \le 0$ leading to the desired inequality $$V^{k+1}(x_0,t_0)\ge V^{k}(x_0,t_0),$$ everywhere, inside or outside the stopping region.  As in the deterministic case we can choose the modifications so that the stopping policy values $U^k$ are everywhere nondecreasing.

We collect the assumptions made in the exposition of the {\bf PIA} for the free boundary stochastic system:
\begin{description}
\item [a)] The controls $u^k$ lead to strong solutions.
\item [b)] The PDE's \eqref{valuepde} with boundary values at the stopping region boundaries have smooth solutions allowing the use of Ito's Lemma.
\item [c)] The corresponding stopping times have finite expectations.
\end{description}
Subject to all these assumptions the {\bf PIA} leads to monotonically improving policies.  In its application to the American Put problem the control improvement is superfluous, since it is only the end point value that determines the payoff.  Note also that the algorithm provides yet another proof of the necessity of smooth pasting.

\section{Implementation, Accuracy and Computational Results}\label{SectComput}
The crucial step in the implementation of the proposed algorithm is the  solution of Equation \eqref{BS} for a given boundary specified at discrete points, but not necessarily coinciding with the standard, uniform discretisation. In order to solve the equation above the boundary the popular implicit Euler method is used (see for example  \cite{Hull}, \cite{BrenSchw}), modified in the vicinity of the boundary where the derivative estimates are adjusted to take into account the grid non uniformity.  A piecewise linear boundary is implemented although a smooth interpolation would have been more accurate; however this simple approach suffices for the fine discretisation used.   As in \cite{BrenSchw} the derivatives of the value function (and not the function itself) are set to zero at a large asset price and verified that increasing it does not affect the results.  Using standard arguments as in \cite{IsaacKell}, the error in the value function is $O(\Delta t+\Delta x^2)$ in the region above the boundary. The reason that the Crank-Nicolson method wasn't used is twofold: the Crank Nicolson method  is unconditional stable in $l_2$ norm.  This, together with consistency, ensures convergence in the $l_2$ norm for initial data which lies also in $l_2$. Moreover  the order of convergence may be less than the second order achieved for smooth initial data, see for example \cite{Tho70}.  Even though there are modifications of  Crank-Nicolson method like the  Rannacher one \cite{Ran84} that retain the second order accuracy also for $l_2$ smooth initial data,  they would require a more complicated adjustment near the boundary. Our simpler implementation is sufficient  to clarify the algorithm's features. In the single asset case these refinements were not necessary for a successful implementation but may prove crucial in a multidimensional treatment. 

The algorithm requires the solution of \eqref{BS} below as well as above the boundary.  As stated earlier, extending the solution smoothly below the boundary is an initial value problem, whose error estimate increases exponentially and the implementations that tried to extend the solution below the boundary diverged fast.  In the cases where the boundary had to be decreased the local updating that was used in Equation \eqref{locupd} proved to be satisfactory.

The computations show that the  algorithm converges at a boundary $\tilde{b}$ on which the asset derivative equals $-1$ to several decimals.  Consequently, the boundary error  $\parallel\tilde{b}-b^*\parallel_\infty$ is of order $O(\Delta x)$ provided that the discretisation satisfies  $\Delta t\le \Delta x^2$.  To show this,  note that  the partial derivative estimate under this condition is accurate to $O(\Delta x)$.  Moreover,  let $V,b$ be the current value and boundary, and $V^*,b^*$ the corresponding optimal ones.  Then, by the arguments in  Section \ref{ConvergSect}, we have 
\begin{align*}
V(b(t),t)-V^*(b(t),t)&=V^*_{xx}(b^*,t)(b-b^*)^2+O(\Delta b)^3,\\
V^*_x(b,t)&=-1+V^*_{xx}(b^*,t)(b-b^*)+O(\Delta b)^2,
\end{align*}
and $V_x(b,t)=-1+O(\Delta x)$.  By  {\bf AP} (\ref{partials}) applied on the boundary \textit{b} and noting that $V,V^*$ differ by a quadratic term, so will their derivatives, and thus 
$$V^*_x(b,t)=-1+V^*_{xx}(b^*,t)(b-b^*)+O(\Delta b)^2=V_x(b,t)=-1+O(\Delta x).$$ 
Consequently,  $\|b-b^*\|_{\infty}=O(\Delta x)$.  It is mentioned that for a Crank-Nicolson type scheme and smooth enough initial data a discretisation of $\Delta x= \Delta t$ would have sufficed for an $O(\Delta x)$ derivative error, and the computational burden would have improved by an order of magnitude.

The boundary updating stipulated by the algorithm is implemented in a simplified fashion as follows: At the $k$-th iteration with a boundary $b_k$ we  calculate the value  $V^k(x,t)$ for $x \ge b_k(t)$.  We then compute $\max\{ [K-x]^+-V^k(x,t)\}$ and if it is positive we set the new boundary $b_{k+1}(t)$ at the value attaining the maximum.  If it occurs at an $x$ below the boundary is updated by maximizing the local quadratic approximation and obtain the updating formula  
\begin{equation}\label{locupd}
b_{k+1}(t)=b_k(t)-(1+V^k_{x})/V^k_{xx}.
\end{equation}
The derivatives $V^k_{x},V^k_{xx}$ are to be evaluated at the boundary and standard one sided formulas were applied.  This simplified updating did not affect the claimed speed of convergence  since these refer only to the vicinity of the optimal boundary where the approximate updating is accurate.  It is noted that the second derivative expression has the same accuracy as the first derivative one since it depends on it through Equation \eqref{BS}.  

The implementation displayed the value monotonicity property claimed, and all values generated on the $x,t$ grid were  monotonically increasing for the chosen  discretisation.  This property would probably not hold for a coarser discretisation, and would be of interest to prove   for the corresponding methods, as in \cite{Forsyth1},  the necessary conditions under which  the  monotonicity holds.

The algorithm is applied to a recent example appearing in Zhu et. al. \cite{Zhu18} that uses a novel integral equation for the boundary.  The example has  a risk free interest rate $r=10\%$, volatility $\sigma=30\%$, exercise price $K=100$ and horizon $T=1$.  We applied the {\bf PIA} algorithm with $\Delta x= 0.05\;\;\Delta t=0.0025$ and the computations are in agreement with those in \cite{Zhu18} as shown in Table \ref{Results}.  A linear interpolation was used to obtain Zhu's times to expiry.  We also extended the horizon calculation to $T=10$ units with $\Delta x = 0.1$, and we obtained an exercise boundary of  69.2371, slightly above the perpetual put value 68.9655, obtained from the expression \cite{Koh:14} $K_{exer}\frac{\phi}{1-\phi},$  with $\phi=\sigma^{-2}\left(\sigma^2/2-r-\sqrt{(r-\sigma^2/2)^2+2r\sigma^2}\right).$

\begin{table}
\centering
\caption{Computational Results}     
\label{Results}
\begin{small}
\begin{tabular}{||c|c|c|c||}
\hline\hline
Time    &       PIA& Zhu &      Adjusted \\
to Expiry       &       & (average)& Brennan Schwartz\\         
\hline
0.0868  &       87.3735 &       87.3548 &       87.3842 \\
0.1515  &       85.0142 &       84.9176 &       85.0140 \\
0.2321  &       83.0725 &       82.9635 &       83.0649 \\
0.3039  &       81.8029 &       81.7002 &       81.7972 \\
0.3697  &       80.8666 &       80.7677 &       80.8589 \\
0.4480  &       79.9438 &       79.8531 &       79.9364 \\
0.5083  &       79.3373 &       79.2523 &       79.3312 \\
0.5761  &       78.7375 &       78.6575 &       78.7328 \\
0.6521  &       78.1472 &       78.0710 &       78.1428 \\
0.7376  &       77.5655 &       77.4928 &       77.5623 \\
0.8335  &       76.9949 &       76.9246 &       76.9919 \\
0.9413  &       76.4356 &       76.3635 &       76.4336 \\
\hline\hline
\end{tabular}
\end{small} 
\end{table}

The computational cost  per boundary calculation is $O(mn),$ with $m,n$ being the number of grid points in time and asset price, respectively.  Since $\Delta x^2=\Delta t$ the aforementioned cost is   $O(n^3)$ flops. The number of iterations required to achieve smooth pasting i.e. the derivative estimate being -1 to several significant figures, was of the order of 10.  However, if we stop when the derivative is $O(\Delta x^2)$ close to -1, we have reached the expected error level and might as well stop. This is achieved in about 5 iterations even when starting from an inaccurate initial policy.   

\begin{figure}[h]
\centering
\includegraphics[width=8cm, height=6cm]{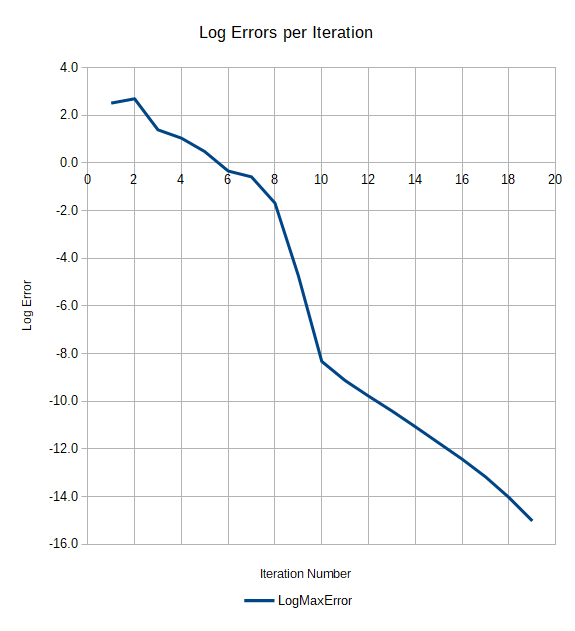}
\caption{ \textbf{Log Errors, all iterations}  }
\label{AllErrors}
\end{figure}

\begin{figure}[h]
\centering
\includegraphics[width=8cm, height=6cm]{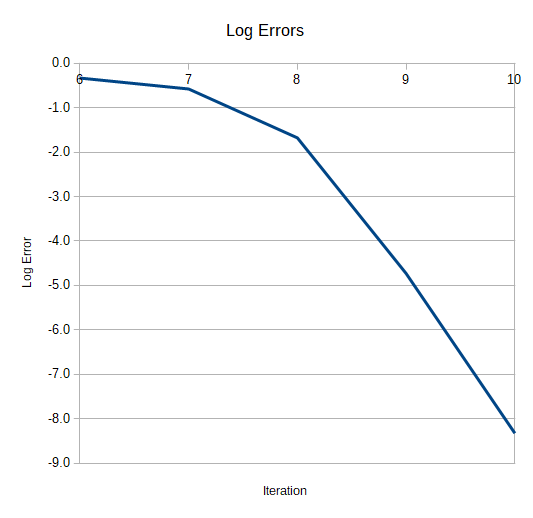}
\caption{ \textbf{Log Errors selected iterations}}
\label{ LogErrors}
\end{figure}

\begin{figure}[h]
\includegraphics[width=8cm, height=6cm]{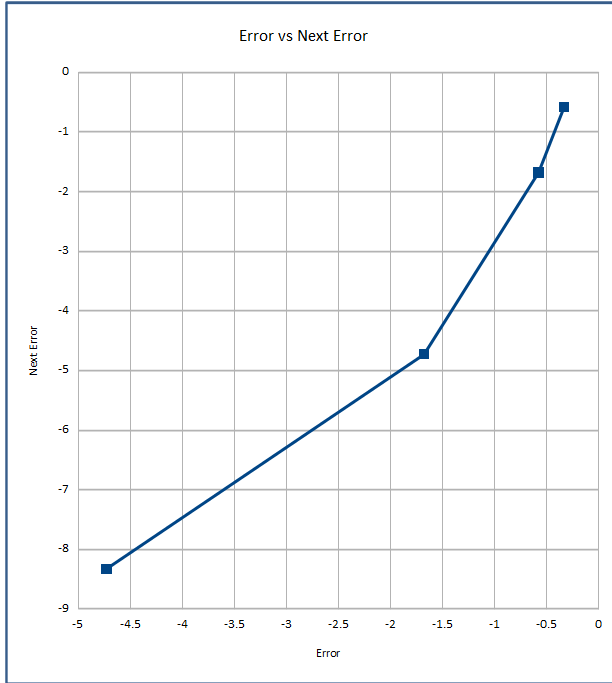}
\centering
\caption{ \textbf{Log Error vs. Next  Error}  }
\label{CorelErrors}
\end{figure}

\begin{figure}[h]
\includegraphics[width=8cm, height=6cm]{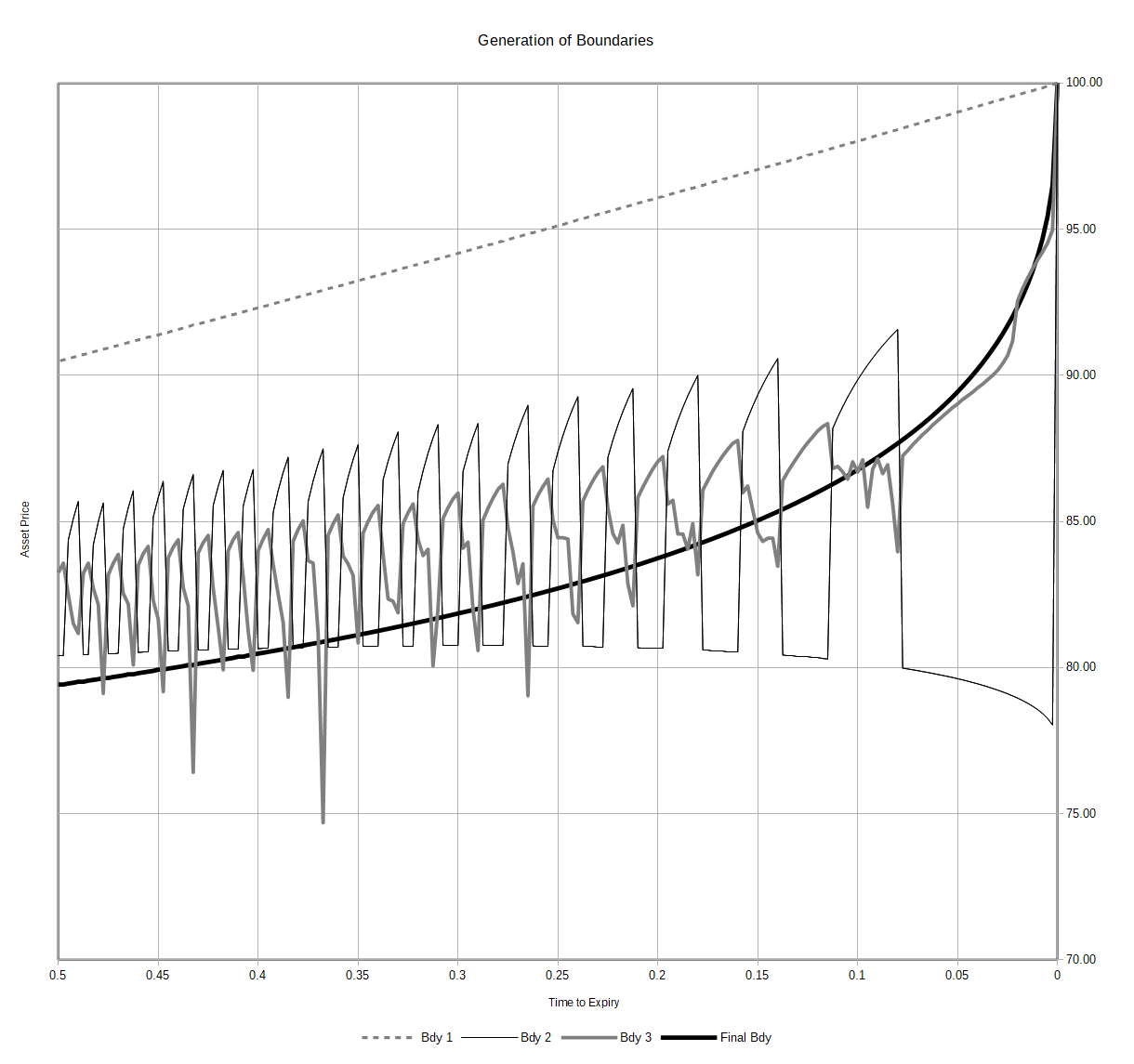}
\centering
\caption{ \textbf{Sequence of Boundaries Generated by the PIA}  }
\label{BoundGen}
\end{figure}

The claimed quadratic convergence of the algorithm was analyzed in the above computations.  We examine whether the error $e_n= \|b_n(t)-b^*(t)\|_{\infty}$ satisfies $e_n\le \alpha e^2_{n-1}$ or in logarithmic terms  $\eta_n=a+2\eta_{n-1}$  with $\eta_n=ln(e_n)\quad a=ln(\alpha)$.   We estimated $e_n$ by sampling over equidistant $t_i$'s, and present the results in  Figures \ref{AllErrors}, \ref{ LogErrors}, and  \ref{CorelErrors}. We expect the quadratic relation to hold whenever the errors are sufficiently small, but not smaller than $O(\Delta x)$.  Specifically,  Figure \ref{AllErrors} shows that  the errors satisfying these conditions are those in iterations 6-10. The results of these iterations  are shown in Figures \ref{ LogErrors} and  \ref{CorelErrors}.  In Figure \ref{ LogErrors} we present the logarithmic error behavior  and  note  the expected concave shape resulting from the solution of the  difference equation $\eta_n=2^n(\eta_0 + a)-a$, the coefficient $\eta_o+a$ of the power term assumed negative, a reasonable assumption.  In Figure  \ref{CorelErrors} we observe that the plot of $\eta_n$ vs $\eta_{n-1}$ has a slope close to 2 for the above iterations. For smaller error values we see in Figure \ref{AllErrors} that they decrease linearly and not quadratically, but this is of little interest since the  convergence to the boundary  is only $O(\Delta x) $ accurate.
In Figure \ref{BoundGen} we illustrate  the  generated boundary sequence.  Starting with an initial boundary far from the optimal, the algorithm generates wildly fluctuating ones for a few iterations, and then settles on the limiting one. This behavior is not inconsistent with the monotonicity property which refers to the values and not the boundaries generated.

A  fast method like Brennan and Schwartz \cite{BrenSchw},  can be easily modified by the Policy Iteration principle to produce more accurate results.  Let the value obtained by that method at  $x_j,t_i$ be $\tilde{V}(x_j,t_i)$. A point $x_j$ is considered as a boundary one if $K-x_j \le\tilde{V}(x_j,t_i)$ but at the immediately lower point $x_{j-1}$ the opposite holds.  It is clear that the derivative at these points is greater than $-1$, and hence the actual boundary should be at a lower point. Having solved \eqref{BS} above the boundary just once the local update \eqref{locupd} is then applied.  These calculations are presented in the last column of Table \ref{Results}; the method was applied with $\Delta x=0.05 ,\; \Delta t=\Delta x^2$ and the results are in close agreement with those obtained by the other methods.  The non adjusted Brennan and Schwartz boundaries  differ by about $1\%$ for $\Delta x=0.05$.   Finally, for a horizon of $T=10$ the boundary is  at 70.04, relatively close to the perpetual put value 68.97. 

The above computational results are further verified by a simple simulation.  The asset process $z_t$ was generated by the exponential formula $z_{t+\delta t}=z_t \exp\left((r-\frac{\sigma^2}{2})\delta t+\sigma\epsilon_t\sqrt{\delta t}\right)$ with $\epsilon_t$ independent samples from a standard normal distribution.  The  time step $\delta t$ used was much smaller than  the boundary discretisation $\Delta t$, since a crossing of the linearized boundary by $z_t$ is sufficient but not necessary for a crossing by the asset process.  Thus, it is expected the value obtained in the simulation to be smaller than the one calculated by PIA. This is indeed observed, see Table \ref{Simulation}.

\begin{table}
        \centering
        \caption{Simulation Results}     
        \label{Simulation}
        \begin{small}
        \begin{tabular}{||c|c|c|c||}
        \hline\hline
\multicolumn{4}{||c||}{Time to Expiry =  0.5}           \\      
\hline                          
Asset Price     &       Option  &       Simulation      &       Standard        \\
&       Value   &       result  &       Deviation       \\
\hline
86.56   &       15.236  &       14.860  &       0.215   \\
96.94   &       9.558   &       9.422   &       0.214   \\
107.32  &       5.883   &       5.965   &       0.184   \\
117.71  &       3.563   &       3.510   &       0.149   \\
128.09  &       2.142   &       2.216   &       0.119   \\
138.47  &       1.271   &       1.280   &       0.091   \\
148.85  &       0.751   &       0.487   &       0.055   \\
159.24  &       0.445   &       0.321   &       0.043   \\
169.62  &       0.267   &       0.223   &       0.035   \\
180.00  &       0.168   &       0.146   &       0.030   \\
\hline\hline
\end{tabular}
\end{small} 
\end{table}

%% file: CONCLUSION.tex
\chapter*{Conclusions}
\markboth{CONCLUSIONS}{CONCLUSIONS}
\addcontentsline{toc}{chapter}{Conclusions}

In this thesis, the numerical solution of three different classes of problems  have been studied. Specifically, we have proposed new techniques and their theoretical analysis has been performed, accompanied by a wide set of numerical experiments, for investigating further the effectiveness and performance of our approach. The first two belong to the research area of numerical linear algebra and concern the spectral analysis and preconditioning of the coefficient matrix of large structured linear systems. The third concerns a problem from economics namely the pricing of an American put option.
\par
In the first set of problems the singular value distribution of the matrix sequence $\{h(T_n(f))\}_n$ and the eigenvalue distribution of the symmetrised matrix sequence $\{Y_nh(T_n(f))\}_n$ was provided. The spectral assymptotic behavior of this matrix sequences was studied under the assumption that $h$ is an analytic function at 0, with radious of convergance $r$, and $f\in L^\infty([-\pi,\pi])$ with $\|f\|_\infty<r$   so that $h\circ f \in L^\infty([-\pi,\pi])\subset  L^1([-\pi,\pi])$. Taking advantage of this analysis circulant preconditioners were proposed, and the eigenvalue distribution of the preconditioned matrix sequences was given. All theoretical results were numerically confirmed. A desirable future development in this class of problems is the investigation on the possibility of relaxation of the given conditions:
\begin{description} 
\item[A)] If the function $h$ is analytic in a given disk centered at $z_0\neq 0$, then the arguments working in the case $z_0=0$ can be repeated verbatim also in the new setting. 
\item[B)] When $f\in L^1([-\pi,\pi])$ (but $f$ does not belong to $ L^\infty([-\pi,\pi])$), then the situation is more complicated and a further step of analysis is required. It could be used the cut--off argument as in \cite{MR1481397,TyL1} and the versatility of the a.c.s. notion. An alternative to the cut-off idea is the use of polynomials such as the Cesaro sum $f$ converging to $f$ in the $ L^1([-\pi,\pi])$ metric plus the trace-norm estimates of $T_n(f)$ derived in \cite{SeTi}.
\end{description}

\par
The second problem that we studied concerned the theoretical and numerical exploration of proper preconditioners based on the spectral symbols of the coefficient matrix for FDE problems. Beside the theoretical study, a comparison between the new and old preconditioners was conducted, especially those presented in~\cite{donatelli161,moghaderi171}.
As expected and  numerically shown in Example 1 which concerned  the one dimensional case, the proposed preconditioners performed slightly worse, at least in a sequential model of computation,  than the tridiagonal preconditions proposed in \cite{donatelli161}, because of the computational complexity. However, in the two dimensional case as discussed in Examples 2 and 3, the proposed preconditioners did indeed perform better than the previously proposed multigrid-based or band preconditioners proposed and studied in \cite{moghaderi171}.
Future directions of research may include more complex problems, further analysis, and more extensive numerical experimentation. Also, problems where the fractional derivatives are greater than 2 may be considered, since then it is expected the symbol-based preconditioners to be even more advantageous, maybe even in the one dimensional case.

\par
For the pricing of an American put option an iterative algorithm was used. Taking advantage of the already known characteristics of the optimal value function it was shown theoretically and confirmed numerically that the proposed algorithm obtains monotonically increasing value functions and converges to the optimal one. Issues regarding the efficiency in multiple asset cases or how can the method be applied in a finite element context are to be concerned in future works.
\par
In the same spirit of the iterative algorithm used here for pricing the American put option, we are working on a further work on the solution of Stefan problem. It consists of a boundary value problem which describes the evolution of the boundary between two faces--that is currently being developed. A preprint of this work is available \cite{magstef}.

%% file: MAIN.bbl
\begin{thebibliography}{10}

\bibitem{MR952991}
F.~Avram.
\newblock On bilinear forms in {G}aussian random variables and {T}oeplitz
  matrices.
\newblock {\em Probab. Theory Related Fields}, 79(1):37--45, 1988.

\bibitem{axelsson861}
O.~Axelsson and G.~Lindskog.
\newblock On the rate of convergence of the preconditioned conjugate gradient
  method.
\newblock {\em Numer. Math.}, 48(5):499--523, 1986.

\bibitem{BarErkVas}
N.~Barakitis, S.-E. Ekstr{\"o}m, and P.~Vassalos.
\newblock Preconditioners for fractional diffusion equations based on the
  spectral symbol.
\newblock {\em Numer. Linear Algebra Appl.}, (accepted).

\bibitem{Barles}
G.~Barles, J.~Burdeau, K.~Romano, and N.~Samsoen.
\newblock Critical stock price near expiration.
\newblock {\em Appl. Math. Finance}, 5:77--95, 1995.

\bibitem{Bellman}
R.~Bellman.
\newblock {\em Dynamic Programming}.
\newblock Princeton University Press, 1957.

\bibitem{Benzi2002}
M.~Benzi.
\newblock Preconditioning techniques for large linear systems: A survey.
\newblock {\em J. Comput. Phys.}, 182(2):418--477, 2002.

\bibitem{BrenSchw}
M.~Brennan and E.~Schwartz.
\newblock The valuation of \uppercase{A}merican put options.
\newblock {\em J. Finance}, 32:449--462, 1977.

\bibitem{CarrJarrowMyn}
P.~Carr, R.~Jarrow, and R.~Myneni.
\newblock Alternative characterizations of american put options.
\newblock {\em Math. Financ.}, 2:87--106, 1992.

\bibitem{donatelli161}
M.~Donatelli, M.~Mazza, and S.~Serra-Capizzano.
\newblock Spectral analysis and structure preserving preconditioners for
  fractional diffusion equations.
\newblock {\em J. Comput. Phys.}, 307:262--279, 2016.

\bibitem{Estatico2008}
C.~Estatico and S.~Serra-Capizzano.
\newblock Superoptimal approximation for unbounded symbols.
\newblock {\em Linear Algebra Appl.}, 428(2-3):564--585, 2008.

\bibitem{Fasino2000}
D.~Fasino and P.~Tilli.
\newblock Spectral clustering properties of block multilevel {H}ankel matrices.
\newblock {\em Linear Algebra Appl.}, 306(1-3):155--163, 2000.

\bibitem{Ferrari2020}
P.~Ferrari, N.~Barakitis, and S.~Serra-Capizzano.
\newblock Asymptotic spectra of large matrices coming from the symmetrization
  of {T}oeplitz structure functions and applications to preconditioning.
\newblock {\em Numer. Linear Algebra Appl.}, 28(1), 2020.

\bibitem{FFHMS}
P.~Ferrari, I.~Furci, S.~Hon, M.~A. Mursaleen, and S.~Serra-Capizzano.
\newblock The eigenvalue distribution of special 2-by-2 block matrix-sequences
  with applications to the case of symmetrized {T}oeplitz structures.
\newblock {\em SIAM J. Matrix Anal. Appl.}, 40(3):1066--1086, 2019.

\bibitem{Forsyth1}
{P. A.} Forsyth and G.~Labahn.
\newblock Numerical methods for controlled hamilton-jacobi bellman pdes in
  finance.
\newblock {\em J. Comput. Finance}, 11:1--44, 2007.

\bibitem{MR3674485}
C.~Garoni and S.~Serra-Capizzano.
\newblock {\em Generalized locally {T}oeplitz sequences: theory and
  applications. {V}ol. {I}}.
\newblock Springer, Cham, 2017.

\bibitem{MR3674485_vol2}
C.~Garoni and S.~Serra-Capizzano.
\newblock {\em Generalized locally {T}oeplitz sequences: theory and
  applications. {V}ol. {II}}.
\newblock Springer, Cham, 2018.

\bibitem{garoni2015}
C.~Garoni, S.~{Serra-Capizzano}, and P.~Vassalos.
\newblock A general tool for determining the asymptotic spectral distribution
  of {H}ermitian matrix-sequences.
\newblock {\em Operators and Matrices}, 9:549--561, 2015.

\bibitem{J.Goo:02}
J.~Goodman and D.~Ostrov.
\newblock On the early exercise boundary of the \uppercase{A}merican put
  option.
\newblock {\em SIAM J. Appl. Math.}, 62:1823--1835, 2002.

\bibitem{Greenbaum1996}
A.~Greenbaum, V.~Pt{\'{a}}k, and Z.~Strako{\v{s}}.
\newblock Any nonincreasing convergence curve is possible for {GMRES}.
\newblock {\em {SIAM} J Matrix Anal. Appl.}, 17(3):465--469, 1996.

\bibitem{MR890515}
U~Grenander and G~Szeg\H{o}.
\newblock {\em {T}oeplitz forms and their applications}.
\newblock Chelsea Publishing Co., New York, second edition, 1984.

\bibitem{HMS}
S.~Hon, M.~A. Mursaleen, and S.~Serra-Capizzano.
\newblock A note on the spectral distribution of symmetrized {T}oeplitz
  sequences.
\newblock {\em Linear Algebra Appl.}, 579:32--50, 2019.

\bibitem{Hon2018}
S.~Hon and A.~Wathen.
\newblock Circulant preconditioners for analytic functions of {T}oeplitz
  matrices.
\newblock {\em Numer. Algorithms}, 79(4):1211--1230, 2018.

\bibitem{Hull}
J.~Hull.
\newblock {\em Options, Futures and Other Derivatives}.
\newblock Pearson, 9th edition, 2017.

\bibitem{IsaacKell}
E.~Isaacson and H.~Keller.
\newblock {\em Analysis of Numerical Methods}.
\newblock Dover, 1994.

\bibitem{Jacka}
S.~D. Jacka.
\newblock Optimal stopping and the \uppercase{A}merican put.
\newblock {\em Math. Financ.}, 1:1--14, 1991.

\bibitem{JackaPI}
{S. D.} Jacka and A.~Mijatovic.
\newblock On the policy improvement algorithm in continuous time.
\newblock {\em Stochastics}, 89(1):348--359, 2017.

\bibitem{Kim}
{I.J.} Kim.
\newblock The analytic valuation of american options.
\newblock {\em Rev. Financial Stud.}, 3:547--572, 1990.

\bibitem{Koh:14}
R.~Kohn.
\newblock Lecture notes in \uppercase{PDE'}s for \uppercase{F}inance, 2014.

\bibitem{Tho70}
{H. O.} Kreiss, V.~Thom{\'e}e, and O.~Widlund.
\newblock Smoothing of initial data and rates of convergence for parabolic
  difference equations.
\newblock {\em Comm. Pure Appl. Math.}, 43:241--259, 1970.

\bibitem{Lanczos1952}
C.~Lanczos.
\newblock Solution of systems of linear equations by minimized iterations.
\newblock {\em J. Res. Natl. Bur. Stand.}, 49(1):33, 1952.

\bibitem{MR2872597}
S.~T. Lee, X.~Liu, and H.~Sun.
\newblock Fast exponential time integration scheme for option pricing with
  jumps.
\newblock {\em Numer. Linear Algebra Appl.}, 19(1):87--101, 2012.

\bibitem{MR2609339}
S.~T. Lee, H.~Pang, and H.~Sun.
\newblock Shift-invert {A}rnoldi approximation to the {T}oeplitz matrix
  exponential.
\newblock {\em SIAM J. Sci. Comput.}, 32(2):774--792, 2010.

\bibitem{lei131}
S.~Lei and H.~Sun.
\newblock A circulant preconditioner for fractional diffusion equations.
\newblock {\em J. Comput. Phys.}, 242:715--725, 2013.

\bibitem{JackaVolat}
J.~Maeda and S.~D. Jacka.
\newblock Market driver volatility model via policy improvement algorithm.
\newblock {\em arXiv:1612.0078v1}, 2016.

\bibitem{E.Mag:77}
E.~Mageirou.
\newblock Iterative techniques for \uppercase{R}icatti game equations.
\newblock {\em J. Optim. Theory Appl.}, 22:51--61, 1977.

\bibitem{Magirou2020}
E.~Magirou, P.~Vassalos, and N.~Barakitis.
\newblock A policy iteration algorithm for the {A}merican put option and free
  boundary control problems.
\newblock {\em J. Comput. Appl. Math}, 373:112544, 2020.

\bibitem{magstef}
E.~Magirou, P.~Vassalos, and N.~Barakitis.
\newblock On a boundary updating method for the scalar {S}tefan problem.
\newblock {\em arXiv preprint arXiv:2202.06418}, 2022.

\bibitem{mazza-pestana}
M.~Mazza and J.~Pestana.
\newblock Spectral properties of flipped {T}oeplitz matrices and related
  preconditioning.
\newblock {\em BIT}, 59(2):463--482, 2019.

\bibitem{H.McK:65}
H.~McKean.
\newblock A free boundary problem for the heat equation arising from a problem
  in mathematical economics.
\newblock {\em Ind. Manag. Rev.}, 6:32--39, 1965.

\bibitem{Meerschaert2004}
M.~Meerschaert and C.~Tadjeran.
\newblock Finite difference approximations for fractional
  advection{\textendash}dispersion flow equations.
\newblock {\em J. Comput. Appl. Math.}, 172(1):65--77, 2004.

\bibitem{Meerschaert2006}
M.~Meerschaert and C.~Tadjeran.
\newblock Finite difference approximations for two-sided space-fractional
  partial differential equations.
\newblock {\em Appl. Numer. Math.}, 56(1):80--90, 2006.

\bibitem{moghaderi171}
H.~Moghaderi, M.~Dehghan, M.~Donatelli, and M.~Mazza.
\newblock Spectral analysis and multigrid preconditioners for two-dimensional
  space-fractional diffusion equations.
\newblock {\em J. Comput. Phys.}, 350:992--1011, 2017.

\bibitem{NSV_con}
D.~Noutsos, S.~Serra-Capizzano, and P.~Vassalos.
\newblock {\em Spectral Equivalence and Matrix Algebra Preconditioners for
  Multilevel {T}oeplitz Systems: A Negative Result.}, pages 313--322.
\newblock AMS, USA, 2001.

\bibitem{NSV_tcs}
D.~Noutsos, S.~Serra-Capizzano, and P.~Vassalos.
\newblock Matrix algebra preconditioners for multilevel {T}oeplitz systems do
  not insure optimal convergence rate.
\newblock {\em Theor. Comput. Sci}, 315(2):557--579, 2004.

\bibitem{noutsos2005}
D.~Noutsos, S.~{Serra-Capizzano}, and P.~Vassalos.
\newblock A preconditioning proposal for ill-conditioned {H}ermitian two-level
  {T}oeplitz systems.
\newblock {\em Numer. Linear Algebra Appl.}, 12(2-3):231--239, 2005.

\bibitem{noutsos2006}
D.~Noutsos, S.~Serra-Capizzano, and P.~Vassalos.
\newblock Block band {T}oeplitz preconditioners derived from generating
  function approximations: analysis and applications.
\newblock {\em Numer. Math.}, 104(3):339--376, 2006.

\bibitem{NSV08}
D.~Noutsos, S.~{Serra-Capizzano}, and P.~Vassalos.
\newblock The conditioning of {FD} matrix sequences coming from semi-elliptic
  differential equations.
\newblock {\em Linear Algebra Appl.}, 428(2-3):600--624, 2008.

\bibitem{noutsos161}
D.~Noutsos, S.~Serra-Capizzano, and P.~Vassalos.
\newblock Essential spectral equivalence via multiple step preconditioning and
  applications to ill conditioned {T}oeplitz matrices.
\newblock {\em Linear Algebra Appl.}, 491:276--291, 2016.

\bibitem{NV2002}
D.~Noutsos and P.~Vassalos.
\newblock New band {T}oeplitz preconditioners for ill-conditioned symmetric
  positive definite {T}oeplitz systems.
\newblock {\em SIMAX}, 23(3):728--743, 2002.

\bibitem{NV-2008}
D.~Noutsos and P.~Vassalos.
\newblock Superlinear convergence for {PCG} using band plus algebra
  preconditioners for {T}oeplitz systems.
\newblock {\em Comput. Math. with Appl.}, 56(5):1255 -- 1270, 2008.

\bibitem{Paige1975}
C.~C. Paige and M.~A. Saunders.
\newblock Solution of sparse indefinite systems of linear equations.
\newblock {\em {SIAM} J. Numer. Anal.}, 12(4):617--629, 1975.

\bibitem{pang151}
H.~Pang and H.~Hai-Wei~Sun.
\newblock Fast {N}umerical {C}ontour {I}ntegral {M}ethod for {F}ractional
  {D}iffusion {E}quations.
\newblock {\em J. Sci. Comput.}, 66(1):41--66, 2015.

\bibitem{Pang2012}
H.~Pang and H.~Sun.
\newblock Multigrid method for fractional diffusion equations.
\newblock {\em J. Comput. Phys.}, 231(2):693--703, 2012.

\bibitem{MR851935}
S.~V. Parter.
\newblock On the distribution of the singular values of {T}oeplitz matrices.
\newblock {\em Linear Algebra Appl.}, 80:115--130, 1986.

\bibitem{doi:10.1137/140974213}
J.~Pestana and A.~Wathen.
\newblock A preconditioned {MINRES} method for nonsymmetric {T}oeplitz
  matrices.
\newblock {\em SIAM J. Matrix Anal. Appl.}, 36(1):273--288, 2015.

\bibitem{Put:04}
M.~Puterman.
\newblock {\em Markov \uppercase{D}ecision \uppercase{P}rocesses:
  \uppercase{D}iscrete \uppercase{S}tochastic \uppercase{D}ynamic
  \uppercase{P}rogramming}.
\newblock John Wiley \& Sons, Inc., New Jersey, 2005.

\bibitem{Ran84}
R.~Rannacher.
\newblock Finite element solution of diffusion problems with irregular data.
\newblock {\em Numer. Math}, 43:309--327, 1984.

\bibitem{reid1971method}
J.~K. Reid.
\newblock On the method of conjugate gradients for the solution of large sparse
  systems of linear equations.
\newblock In {\em Pro. the Oxford conference of institute of mathmatics and its
  applications}, pages 231--254, 1971.

\bibitem{Rynne2008}
B.~Rynne and M.~Martin~Youngson.
\newblock {\em Linear Functional Analysis}.
\newblock Springer London, 2008.

\bibitem{Saad2003}
Y.~Saad.
\newblock {\em Iterative Methods for Sparse Linear Systems}.
\newblock Society for Industrial and Applied Mathematics, 2003.

\bibitem{Saad1986}
Y.~Saad and M.~H. Schultz.
\newblock {GMRES}: A generalized minimal residual algorithm for solving
  nonsymmetric linear systems.
\newblock {\em {SIAM} J. Sci. Comput.}, 7(3):856--869, 1986.

\bibitem{sauer_T}
T.~Sauer.
\newblock {\em Numerical Analysis}.
\newblock Addison-Wesley Publishing Company, USA, 2nd edition, 2011.

\bibitem{Scherer2011}
R.~Scherer, S.~L. Kalla, Y.~Tang, and J.~Huang.
\newblock The {G}r\"{u}nwald{\textendash}{L}etnikov method for fractional
  differential equations.
\newblock {\em Comput. Math. with Appl.}, 62(3):902--917, 2011.

\bibitem{Skoro}
S.~Serra-Capizzano.
\newblock A korovkin-type theory for finite {T}oeplitz operators via matrix
  algebras.
\newblock {\em Numerische Mathematik}, 82(1):117--142, 1999.

\bibitem{Capizzano2000}
S.~Serra-Capizzano.
\newblock Korovkin tests, approximation, and ergodic theory.
\newblock {\em Mathematics of Computation}, 69(232):1533--1559, 2000.

\bibitem{Taud2}
S.~Serra-Capizzano.
\newblock Spectral behavior of matrix sequences and discretized boundary value
  problems.
\newblock {\em Linear Algebra Appl.}, 337:37--78, 2001.

\bibitem{Capizzano2003}
S.~Serra-Capizzano.
\newblock Generalized locally {T}oeplitz sequences: spectral analysis and
  applications to discretized partial differential equations.
\newblock {\em Linear Algebra Appl.}, 366:371--402, 2003.

\bibitem{glt-2}
S.~Serra-Capizzano.
\newblock The {GLT} class as a generalized {F}ourier analysis and applications.
\newblock {\em Linear Algebra Appl.}, 419(1):180--233, 2006.

\bibitem{SeTi}
S.~Serra-Capizzano and P.~Tilli.
\newblock On unitarily invariant norms of matrix-valued linear positive
  operators.
\newblock {\em J. Inequal. Appl.}, 7(3):309--330, 2002.

\bibitem{Tian2015}
W.Y. Tian, H.~Zhou, and W.~Deng.
\newblock A class of second order difference approximations for solving space
  fractional diffusion equations.
\newblock {\em Math. Comput.}, 84(294):1703--1727, 2015.

\bibitem{Tilli1998}
P.~Tilli.
\newblock Locally {T}oeplitz sequences: spectral properties and applications.
\newblock {\em Linear Algebra Appl.}, 278(1-3):91--120, 1998.

\bibitem{MR1258226}
E.~E. Tyrtyshnikov.
\newblock New theorems on the distribution of eigenvalues and singular values
  of multilevel {T}oeplitz matrices.
\newblock {\em Dokl. Akad. Nauk}, 333(3):300--303, 1993.

\bibitem{Tyrtyshnikov19961}
E.~E. Tyrtyshnikov.
\newblock A unifying approach to some old and new theorems on distribution and
  clustering.
\newblock {\em Linear Algebra Appl.}, 232:1--43, 1996.

\bibitem{TyL1}
E.~E. Tyrtyshnikov and N.~L. Zamarashkin.
\newblock Spectra of multilevel {T}oeplitz matrices: advanced theory via simple
  matrix relationships.
\newblock {\em Linear Algebra Appl.}, 270:15--27, 1998.

\bibitem{V18}
P.~Vassalos.
\newblock Asymptotic results on the condition number of {FD} matrices
  approximating semi-elliptic {PDEs}.
\newblock {\em Electron. J. Linear Algebra}, 34:566--581, 2018.

\bibitem{Wang2010}
H.~Wang, K.~Wang, and T.~Sircar.
\newblock A direct $\mathcal{O}({N}\log^2{N})$ finite difference method for
  fractional diffusion equations.
\newblock {\em J. Comput. Phys.}, 229(21):8095--8104, 2010.

\bibitem{Wil91}
D.~Williams.
\newblock {\em Probability with martingales}.
\newblock Cambridge University Press, 1991.

\bibitem{MR1481397}
N.~L. Zamarashkin and E.~E. Tyrtyshnikov.
\newblock Distribution of the eigenvalues and singular numbers of {T}oeplitz
  matrices under weakened requirements on the generating function.
\newblock {\em Mat. Sb.}, 188(8):83--92, 1997.

\bibitem{Zhu18}
{S. P.} Zhu, {X. J.} He, and {X. P.} Lu.
\newblock A new integral equation formulation for \uppercase{A}merican put
  options.
\newblock {\em Quant. Finance}, 18:483--490, 2018.

\end{thebibliography}
